\documentclass[preprint,12pt]{elsarticle}

\usepackage[english]{babel}
\usepackage{amsmath}
\usepackage{amsthm}
\usepackage{amssymb}
\usepackage[scr=boondoxo]{mathalfa}
\usepackage{graphicx}
\usepackage{epstopdf}
\usepackage{subfigure}
\usepackage[dvipsnames]{xcolor}
\usepackage{geometry,fullpage}
\newcommand{\bpr}{\begin{trivlist} \item[]{\bf Proof. }}
	\newcommand{\epr}{\hspace*{\fill} \qed \end{trivlist}}

\newcommand{\tr}{\textrm{tr}\,}
\newcommand{\be}{\begin{eqnarray}}
\newcommand{\ee}{\end{eqnarray}}
\newcommand{\ba}{\begin{align}}
\newcommand{\ea}{\end{align}}
\newcommand{\bi}{\begin{itemize}}
	\newcommand{\ei}{\end{itemize}}

\newcommand{\secref}[1]{Section~\ref{sec:#1}}

\newcommand{\seclab}[1]{\label{sec:#1}}
\newcommand{\eqlab}[1]{\label{eq:#1}}
\renewcommand{\eqref}[1]{(\ref{eq:#1})}

\newcommand{\figref}[1]{Fig.~\ref{fig:#1}}
\newcommand{\figlab}[1]{\label{fig:#1}}
\newcommand{\propref}[1]{Proposition~\ref{proposition:#1}}
\newcommand{\proplab}[1]{\label{proposition:#1}}
\newcommand{\lemmaref}[1]{Lemma~\ref{lemma:#1}}
\newcommand{\lemmalab}[1]{\label{lemma:#1}}
\newcommand{\assumptionref}[1]{Assumption~\ref{assumption:#1}}
\newcommand{\assumptionlab}[1]{\label{assumption:#1}}
\newcommand{\remref}[1]{Remark~\ref{remark:#1}}
\newcommand{\remlab}[1]{\label{remark:#1}}
\newcommand{\thmref}[1]{Theorem~\ref{theorem:#1}}
\newcommand{\thmlab}[1]{\label{theorem:#1}}

\newcommand{\defnlab}[1]{\label{defn:#1}}
\newcommand{\defnref}[1]{Definition~\ref{defn:#1}}
\newcommand{\appref}[1]{Appendix~\ref{app:#1}}

\newcommand{\applab}[1]{\label{app:#1}}

\newcommand{\sgn}{{\mathrm sgn}}

\newtheorem{theorem}{Theorem}[section]
\newtheorem{proposition}[theorem]{Proposition}
\newtheorem{definition}[theorem]{Definition}
\newtheorem{lemma}[theorem]{Lemma}

\newtheorem{remark}[theorem]{Remark}

\newtheorem{assumption}{Assumption}
\numberwithin{equation}{section}
\newtheorem*{classification}{Classification}

\usepackage{todonotes}
\newcommand\SJ[1]{{\color{black}{#1}}} 
\newcommand\note[1]{{\color{black}{#1}}} 
\newcommand\SJnew[1]{{\color{black}{#1}}} 
\definecolor{mygreen}{RGB}{0,175,0}
\newcommand\W[1]{{\color{black}{#1}}} 

\journal{Journal of Differential Equations}

\begin{document}

\begin{frontmatter}




\title{Singularly Perturbed Boundary-Focus Bifurcations}


\cortext[]{Corresponding author}

\author{Samuel Jelbart${}^{a,\ast}$\corref{Corresponding author}}
\ead{sjel7292@uni.sydney.edu.au}
\author{Kristian Uldall Kristiansen${}^{b}$}
\ead{krkri@dtu.dk}
\author{Martin Wechselberger${}^{a}$}
\ead{martin.wechselberger@sydney.edu.au}

\address{
	${}^{a}$School of Mathematics \& Statistics, University of Sydney, Camperdown, NSW 2006, Australia
	${}^{b}$Department of Applied Mathematics and Computer Science, Technical University of Denmark, Lyngby, Kgs. 2800, Denmark
}

\begin{abstract}
We consider smooth systems limiting as $\epsilon \to 0$ to piecewise-smooth (PWS) systems with a boundary-focus (BF) bifurcation. After deriving a suitable local normal form, we study the dynamics for the \textit{smooth} system \note{with sufficiently small but non-zero $\epsilon$,} using a combination of \textit{geometric singular perturbation theory} and \textit{blow-up}. We show that the type of BF bifurcation in the PWS system determines the bifurcation structure for the smooth system within an $\epsilon-$dependent domain which shrinks to zero as $\epsilon \to 0$, identifying a supercritical Andronov-Hopf bifurcation in one case, and a supercritical Bogdanov-Takens bifurcation in two other cases. We also show that PWS cycles associated with BF bifurcations persist as relaxation oscillations in the smooth system, and prove existence of a family of stable limit cycles which connects the relaxation oscillations to regular cycles within the $\epsilon-$dependent domain described above. Our results are applied to models for Gause predator-prey interaction and mechanical oscillation subject to friction.
\end{abstract}




\end{frontmatter}


\section{Introduction}
\seclab{introduction}

Boundary equilibrium (BE) bifurcations are piecewise-smooth (PWS) bifurcations which occur when an equilibrium collides with a discontinuity manifold $\Sigma$ under variation of a system parameter $\mu$. There is a growing literature devoted to the study of such problems, including (but not limited to) the work contained in \cite{Carvalho2018,Carvalho2011,Hogan2016,Kuznetsov2003}, owing much to the work of Filippov \cite{Filippov1960}. In particular, a complete (PWS) topological classification
exists for codimension-1 BE bifurcations in planar PWS systems \cite{Hogan2016,Kuznetsov2003}.
Since PWS systems often constitute approximations for smooth systems with sharp transitions, many authors have chosen to adopt an approach in which PWS systems are viewed as occurring in the singular limit $\epsilon \to 0$ of a smooth, \textit{regularised} system; see \cite{Bonet2016,Buzzi2006,Guglielmi2015,Jeffrey2018,Kristiansen2015b,Kristiansen2015a,Kristiansen2019,Kristiansen2017,Kristiansen2019d,Llibre2007,Sotomayor1996,Teixeira2012} 
and in particular \cite{Carvalho2018,Carvalho2011} for applications in the context of BE bifurcations. On this approach, the `transition' over $\Sigma$ is governed by singularly perturbed dynamics occurring within a narrow domain which scales with $\epsilon$, known as the \textit{switching layer}. For a complete description of the dynamics, one must be able to `connect', or `match' the dynamics at the boundary to the switching layer. This matching is particularly important for understanding BE phenomena, since these are inherently linked to topological changes in the dynamics occurring as the equilibrium enters the switching layer. 

A particularly powerful approach to the study of such problems has been the combination of \textit{geometric singular perturbation theory} (GSPT) \cite{Fenichel1979,Jones1995,Kuehn2015,Wechselberger2019} with a method for geometric desingularisation known as \textit{blow-up} \cite{Dum1996,Krupa2001b}; see, e.g. \cite{Buzzi2006,Carvalho2018,Carvalho2011,Kristiansen2015a,Kristiansen2015b,Kristiansen2019,Kristiansen2019c,Llibre2008,Llibre2007}. 
In the context of smooth systems having a PWS system in their singular limit $\epsilon \to 0$, the {\em blow-up} method can resolve the degeneracy associated with a loss of smoothness at the switching manifold $\Sigma$ by `blowing it up' to a higher dimensional manifold which has improved \W{smoothness and} hyperbolicity properties when viewed within an extended, or `blown-up' phase space.
On this approach, 
the dynamics outside and within the switching layer can be represented within a single (extended) 
phase space, thus providing a framework for analysing the connection between the two regimes.


\

In this work we apply this combination of GSPT and blow-up to the study of smooth systems which limit to PWS systems with a codimension-1 \textit{boundary-focus} (BF) bifurcation, i.e.~a BE bifurcation in which the incident equilibrium is a hyperbolic focus. \note{The remaining cases are treated in the sequel \cite{Jelbart2021}.} Previous work in this direction (e.g. \cite{Carvalho2018,Carvalho2011}) has typically involved the use of `Sotomayor-Teixeira' (ST) regularisation functions, originally introduced in \cite{Sotomayor1996} for the study of certain singularities and as part of Peixoto's program for the study of structural stability in PWS systems \cite{Peixoto1959}. Since ST regularisations involve the somewhat artificial cutoff at the boundary to the switching layer, however, they do not occur naturally in applications. We consider an alternative class of regularisation functions without abrupt cessation at the switching layer boundary, with an emphasis on detailed analysis of the \textit{smooth} dynamics for $0 < \epsilon \ll 1$, as is relevant in applications; see also \cite{Jelbart2019c,Kristiansen2019e,Kristiansen2019c,Kristiansen2019d}. 
Following this approach, we derive a local normal form for smooth systems which limit to PWS systems having a codimension-1 BF bifurcation, allowing for a rigorous study of its unfolding in the corresponding smooth system with $0 < \epsilon \ll 1$.

In particular cases, BF bifurcation is known to correlate with the termination of \textit{PWS cycles} -- closed PWS orbits consisting (but not entirely so) of orbit segments contained within $\Sigma$ -- at the bifurcation value $\mu = \mu_{bf}$. We show that these PWS cycles persist as `relaxation oscillations' in the smooth system with $0 < \epsilon \ll 1$ for $\mu-$values bounded away from $\mu_{bf}$. These oscillations are \textit{not slow-fast}, but rather of relaxation type in the broader sense, i.e. components of the oscillation exhibit qualitatively distinct dynamics as $\epsilon \to 0$; see e.g. \cite{Jelbart2019c,Jelbart2020a,Kosiuk2016,Kristiansen2019b,Kristiansen2019d} for examples of \W{such `relaxation' oscillations in applications}. In order to understand the dynamics near $\epsilon = 0$, $\mu = \mu_{bf}$, multiple blow-up transformations are applied in order to obtain a {\em desingularised system} in which the unfolding of the BF point occurs. Here we identify different mechanisms, depending on the type of BF bifurcation exhibited by the original PWS system obtained in the limit $\epsilon \to 0$. In those cases known to correlate with the birth/termination of PWS cycles, we identify a supercritical Andronov-Hopf bifurcation as a mechanism for the onset of stable oscillations. In a subset of cases, saddle-node and homoclinic bifurcations are also identified, and shown to be organised by the presence of a codimension-two Bogdanov-Takens point. In each case, the identified bifurcations are `singular' in the sense that they occur on a domain which shrinks to zero in the singular limit $\epsilon \to 0$. Importantly for applications, it is shown that the size of this domain (and hence, e.g. any oscillations contained within it) is quantitatively linked both to $\epsilon$ and a `transition coefficient' $k$ associated with the gradient of the transition between the switching and outer layers. It is important to emphasise that the unfolding of the BF point occurs within a parameter regime which lies at the interface between outer and switching layer dynamics, and that detailed analysis of the matching problems obtained at the boundary to both the outer and switching layer is necessary for a complete description of the dynamics. We present such an analysis for a restricted set of cases, proving (in particular) a connection between regular oscillations identified in the desingularised problem and the relaxation oscillations. We show and discuss the application of our results in the context of models for predator-prey interactions and mechanical oscillators subject to friction.
In addition to the direct implications of our results for the normal form for the qualitative dynamics of these systems, we also show how \textit{quantitative} information (concerning e.g.~bifurcations) can be obtained by the application of a suitable (explicit) coordinate transformation, without the need to transform the system into local normal form.

\

The manuscript is structured as follows: In \secref{normal_form} we introduce the problem, provide the relevant background on regularisation and PWS theory, and present the local normal form. 
\SJnew{Main results on bifurcations, persistence of PWS cylces as relaxation oscillations, and a connection for the dynamics associated with the unfolding of the BF point with the dynamics identified in both outer and switching layers, are presented in \secref{main_results}.} 
\SJnew{In \secref{bf_bifurcation_in_application} we apply our results to models for Gause predator-prey interaction and a mechanical oscillator subject to friction.}
We summarise and conclude in \secref{summary_and_conclusion}, deferring the remaining details and proofs to the appendix. 

\section{\SJnew{Setup} for the smooth BF bifurcation}
\seclab{normal_form}

We consider planar systems of the general form
\begin{equation}
\eqlab{main}
\dot z = Z \left(z, \phi\left( y \epsilon^{-1} \right), \alpha \right) ,
\end{equation}
where $z=(x,y) \in \mathbb R^2$, $\phi : \mathbb R \to \mathbb{R}$, $\epsilon \in (0,\epsilon_0]$, and $\alpha \in I \subset \mathbb{R}$. We assume that $Z : \mathbb R^2 \times \mathbb R \times I \to \mathbb R^2$ is smooth in all arguments. Following \cite{Kristiansen2019c}, we impose the following assumptions on \eqref{main}:

\begin{assumption}
	\assumptionlab{ass1}
	The map $p \mapsto Z(z, p, \alpha)$ is affine, i.e.
	\[
	Z(z, p, \alpha) = p Z^+(z, \alpha) + (1 - p) Z^-(z, \alpha) ,
	\]
	where the vector fields $Z^\pm : \mathbb R^2 \times I \to \mathbb R^2$ are smooth.
\end{assumption}
\begin{assumption} \assumptionlab{ass1b}
	The smooth `regularisation function' $\phi : \mathbb R \to \mathbb R$ satisfies the monotonicity condition
	\[
	\frac{\partial \phi(s)}{\partial s} > 0 ,
	\]
	for all $s \in \mathbb R$ and, moreover,
	\begin{equation}
	\eqlab{mono}
	\phi(s) \to
	\begin{cases}
	1  & \text{for } s \to \infty , \\
	0  & \text{for } s \to -\infty .
	\end{cases}
	\end{equation}
\end{assumption}

As a direct consequence of these assumptions, the pointwise limit as $\epsilon \to 0^+$ in \eqref{main} yields the following PWS system on $\mathbb R^2 \setminus \{y=0\}$:
\begin{equation}
\eqlab{PWS}
\dot z =
\begin{cases}
Z^+(z, \alpha) , & \text{for } y>0, \\
Z^-(z, \alpha) , & \text{for } y<0, \\
\end{cases}
\end{equation}
where the two open half-planes $\{y>0\}$ and $\{y<0\}$ are separated by the line 
\begin{equation}
\eqlab{switching_manifold}
\Sigma = \left\{(x,y) \in \mathbb R^2 : y = 0 \right\} ,
\end{equation}
often referred to in the PWS literature as the \textit{switching manifold} \cite{Bernardo2008}.
\begin{figure}[t!]
	\centering
	\includegraphics[scale=0.75]{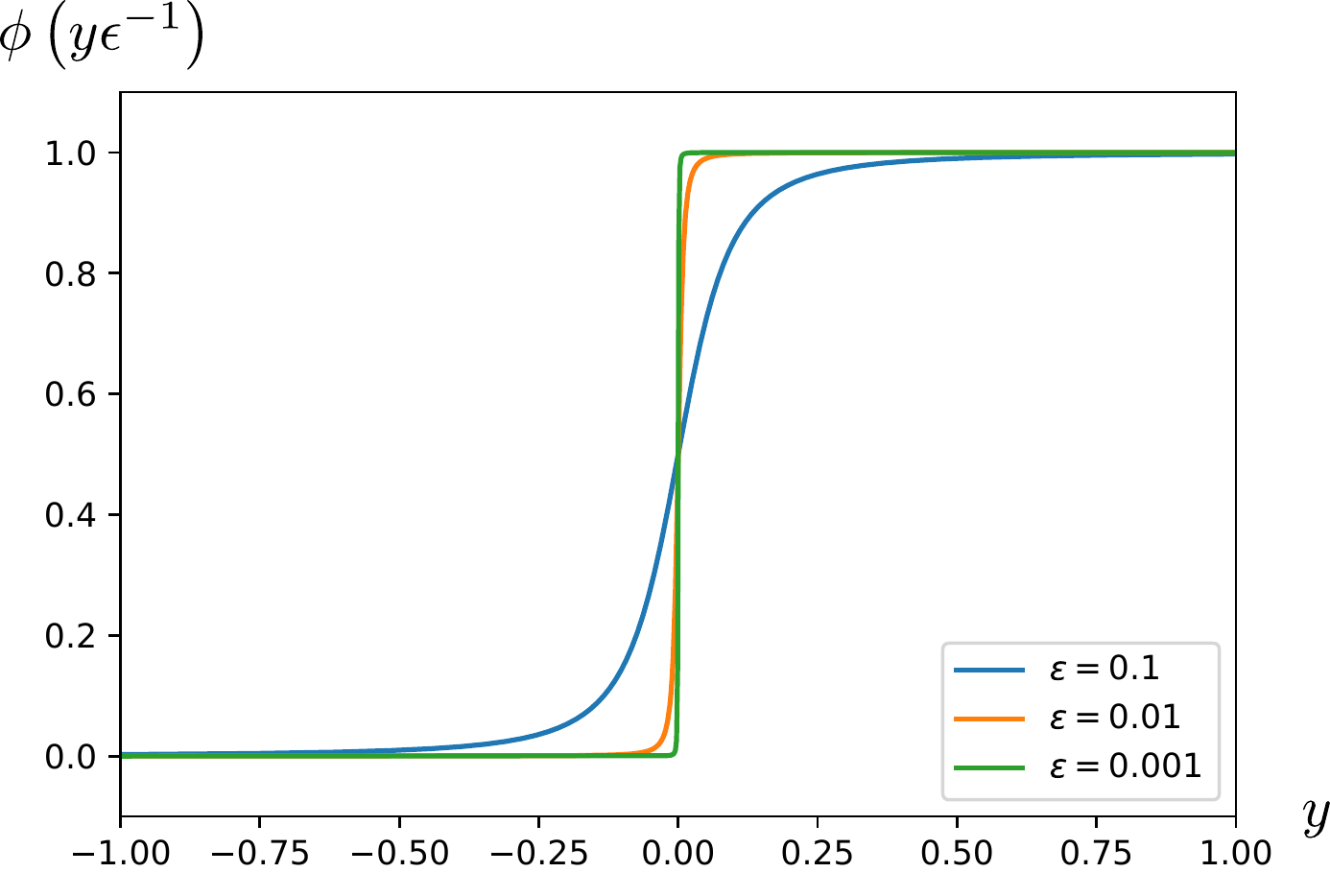}
	\caption{A regularisation function $\phi(y \epsilon^{-1}) = \tfrac{1}{2} \left(1 + \tfrac{y}{\sqrt{y^2+\epsilon^2}}\right)$ satisfying \assumptionref{ass1b} and \assumptionref{ass3} with $k = 2$, $\beta = 1/4$, shown for varying values of $\epsilon$.}
	\figlab{reg_fig}
\end{figure}
Hence, system \eqref{main} under \assumptionref{ass1} and  \assumptionref{ass1b} can be viewed as either
\begin{enumerate}
	\item[(i)] a singularly perturbed system 
	limiting to a PWS system, or 
	\item[(ii)] a regularisation of the PWS system \eqref{PWS}. 
\end{enumerate}
In fact both scenarios are common, and we consider applications of each kind in \secref{bf_bifurcation_in_application}. \assumptionref{ass1b} defines the `regularisation function' $\phi(y \epsilon^{-1})$, which is smooth and switch-like for $\epsilon \in (0,\epsilon_0]$, but singular as $\epsilon \to 0$; see \figref{reg_fig}, which shows one example of a particular regularisation satisfying \assumptionref{ass1b} for varying values of $\epsilon$. \note{We emphasise here an important distinction between systems of type (i) and (ii), however. Namely, for systems of type (i) the regularisation function $\phi(s)$ is \textit{intrinsic}, i.e.~uniquely determined by taking the singular limit in system \eqref{main}. In contrast, there is no intrinsic choice for $\phi(s)$ in the context of systems of type (ii), here one must `choose' a regularisation \cite{Jeffrey2018}.}

It is \note{also} important to distinguish functions satisfying \assumptionref{ass1b} from the well-known class of Sotomayor-Teixeira (ST) regularisations consisting of non-analytic regularisations $\psi(s)$ with
\[
\psi'(s) > 0 , \qquad s \in (-1,1) ,
\]
and
\[
\psi(s) =
\begin{cases}
0, & \text{for } s \leq -1 , \\
1, & \text{for } s \geq 1 ,
\end{cases}
\]
which \textit{do not} satisfy \assumptionref{ass1b}. We emphasise that ST regularisations were initially introduced in \cite{Sotomayor1996} for theoretical purposes including the study of certain singularities and structural stability of PWS systems. In contrast to regularisations satisfying \assumptionref{ass1b}, ST regularisations do not arise naturally in applications due to the artificial cutoffs at the boundary to the \textit{switching layer} \SJnew{$\{(x,y) \in \mathbb R^2 : y = \mathcal O(\epsilon)\}$} at $s= \pm 1$.

\begin{remark}
	\remlab{rem_singular}
	System \eqref{main} under \assumptionref{ass1} and  \assumptionref{ass1b} is a singular perturbation problem in the sense that it loses smoothness in the limit $\epsilon \to 0$. Such problems are related but should be distinguished from the class of slow-fast systems
	\begin{equation}
	z' = H(z;\epsilon) , \qquad 0 < \epsilon \ll 1 ,
	\end{equation}
	for which the set of equilibria $\{z\in\mathbb R^2:H(z;0)=0\}$ contains a regularly embedded submanifold of $\mathbb R^2$. This includes slow-fast problems in the standard form
	\begin{equation}
	\eqlab{stnd}
	\begin{split}
	\dot x &= f(x,y,\epsilon) , \\
	\dot y &= \epsilon g(x,y,\epsilon) .
	\end{split}
	\end{equation}
\end{remark}

For analytical purposes we make one additional (technical) assumption regarding the decay rate of the regularisation function $\phi$\note{, which can be viewed as a nondegeneracy condition as $s \to \infty$}:

\begin{assumption}
	\assumptionlab{ass3}
	The regularisation function $\phi(s)$ \note{has algebraic decay at $s\rightarrow \infty$ in the following sense: There exists a  `transition coefficient' $k \in \mathbb N_+$ as well as a smooth function $\phi_+:[0,\infty] \to [0,\infty)$ such that 
	\begin{equation}
	\eqlab{reg_asymptotics}
	\phi(s)=
	1 - s^{-k} \phi_+\left(s^{-1} \right) , \qquad s >0.\\
	\end{equation}
By \assumptionref{ass1b} we may also write
		\begin{equation}
		\eqlab{beta}
		\beta := \phi_+(0) > 0 .  
		\end{equation}}
\end{assumption}


%

Together with the singular perturbation parameter $\epsilon$, the transition coefficient $k \in \mathbb N_+$ in \eqref{reg_asymptotics} provides a direct quantitative measure of the gradient of the smooth transition between the outer and switching layers in systems \eqref{main}.
Both $k$ and $\epsilon$ will play an important role in the statement of our main results.

\note{
	\begin{remark}
		\assumptionref{ass3} is natural in the context of general systems \eqref{main} where the right-hand-side is sufficiently smooth or analytic for $\epsilon \in (0,\epsilon_0]$. Notice specifically, that the condition ensures that $\phi(u^{-1})$ has a well-defined Taylor-expansion at $u=0$:
		\begin{align*}
		 \phi(u^{-1}) = 1-u^k \beta (1+\mathcal O(u)).
		\end{align*}
		Cases like $\tanh(s)$, where $\phi(s)$ has exponential decay for $s \to \infty$, correspond to $k = \infty$ and are therefore excluded by this assumption. Although $k=\infty$ is common in applications, it is significantly complicated by the presence of an essential singularity. We do not consider this further in the present manuscript, however we refer to \cite{Jelbart2019c,Kristiansen2019b} for details on how to extend the analysis to deal with this case using an adaptation of the blow-up method developed in \cite{Kristiansen2017}. Such an extension is described for singularly perturbed BF in an application in \cite[Ch.6]{Jelbart2020b}.
	\end{remark}
}

\subsection{PWS theory and BF bifurcation} 

We use Filippov theory \cite{Filippov1960} to study the PWS system \eqref{PWS} obtained from system \eqref{main} in the singular limit $\epsilon \to 0$.
%

\begin{figure}[t!]
	\begin{center}
		\includegraphics[width=.6\textwidth]{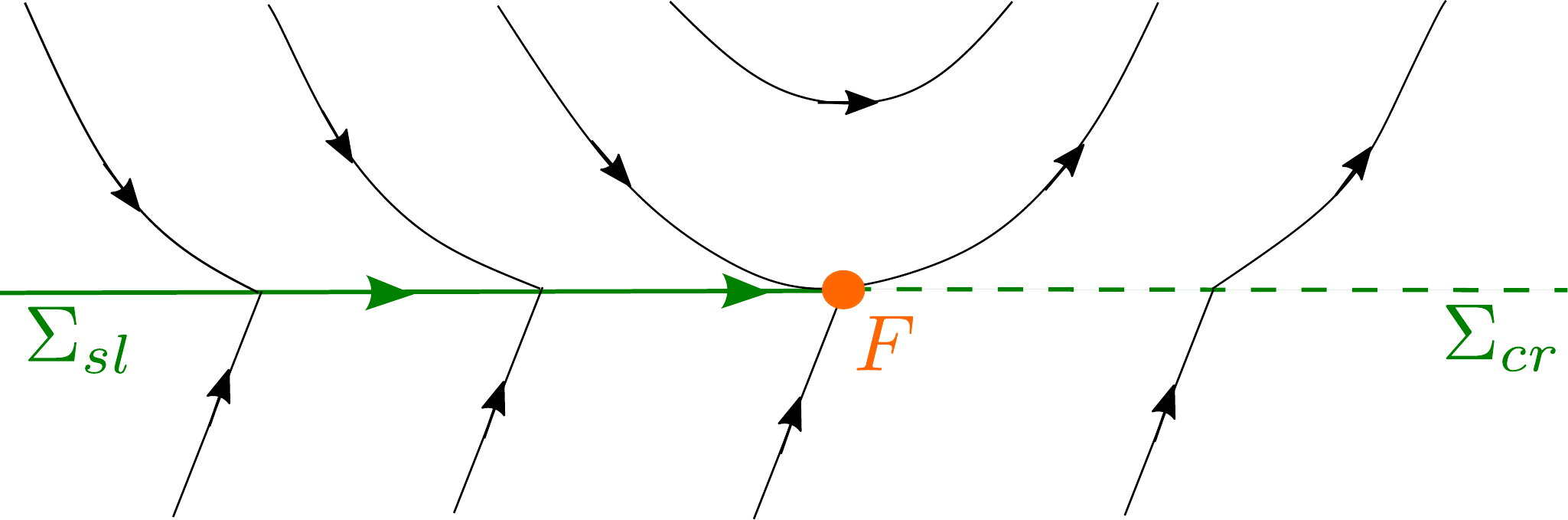}
		\caption{PWS dynamics near a visible fold. The fold point $F$ occurs at a locally quadratic tangency between $\Sigma$ and an orbit of $Z^+$, 
			and separates a sliding region $\Sigma_{sl}$ (bold green), from a crossing region $\Sigma_{cr}$ (dashed green). Trajectories approaching $\Sigma_{sl}$ are `stuck', or `land' on $\Sigma_{sl}$, after which they evolve according to the sliding/Filippov vector field prescribed by \eqref{EqnSlidingVF} until they exit at the fold point. Trajectories colliding with $\Sigma_{cr}$ simply make a non-smooth transition from $\{y \leq 0\}$ to $\{y \geq 0\}$.}
		\figlab{pws_fig}
	\end{center}
\end{figure}

\begin{definition}
	\defnlab{cross_slide}
	Consider system \eqref{PWS}. A point $p \in \Sigma$ is called a crossing point if
	\begin{equation}
	\eqlab{crossing}
	\left(Z^+f(p)\right) \left(Z^-f(p)\right) > 0 , \qquad \qquad \text{(crossing)},
	\end{equation}
	and a sliding point if
	\begin{equation}
	\eqlab{sliding}
	\left(Z^+f(p)\right) \left(Z^-f(p)\right) < 0 , \qquad \qquad \text{(sliding)},
	\end{equation}
	where $Z^\pm f(\cdot) = \langle \nabla f(\cdot), Z^+(\cdot,\alpha) \rangle$ denotes the Lie derivative. We denote the set of crossing points by $\Sigma_{cr}$, and the set of sliding points by $\Sigma_{sl}$.
\end{definition}

Crossing and sliding subsets are separated by tangencies, which occur for points $p \in \Sigma$ for which either one (or both) of $Z^\pm f(p) = 0$. Trajectories are simply assumed to make a non-smooth transition across crossing regions $\Sigma_{cr}$, i.e. to `switch' from $Z^-(p)$ to $Z^+(p)$ (or visa-versa). Trajectories colliding with a sliding region $\Sigma_{sl}$, however, are `trapped' on $\Sigma$ and must evolve in accordance with a vector field defined on $\Sigma_{sl}$. The Filippov convention \cite{Filippov1960} for defining a suitable `sliding vector field' on $\Sigma_{sl}$ is specified in terms of a convex combination of $Z^\pm$ as follows
:
\begin{equation}
Z_{sl}(z,\alpha) := \chi (z,\alpha)Z^+(z,\alpha) + \left(1 - \chi(z,\alpha)\right) Z^-(z,\alpha), \qquad z \in \Sigma_{sl},
\end{equation}
where $z=(x,y)$, and
\[
\chi(z) = \frac{\langle Df(z), Z^-(z,\alpha) \rangle}{\langle Df(z), Z^-(z,\alpha) - Z^+(z,\alpha) \rangle}, \qquad z \in \Sigma_{sl}.
\]
In our case, the simple form of the switching manifold \eqref{switching_manifold} allows us to simplify the above so that the sliding vector field can be expressed straightforwardly in the $x-$coordinate chart as
\begin{equation}
\eqlab{EqnSlidingVF}
\dot x = \frac{\det(Z^+((x,0),\alpha) | Z^-((x,0),\alpha))}{ Z_2^-((x,0),\alpha) - Z_2^+((x,0),\alpha)} = Z_{sl}(x,\alpha) ,  \qquad (x,0) \in \Sigma_{sl},
\end{equation}
where $\det(Z^+((x,0),\alpha) | Z^-((x,0),\alpha))$ denotes the determinant of the $2 \times 2$ matrix with columns $Z^+((x,0),\alpha)$, $Z^-((x,0),\alpha)$. 


Points of tangency between the vector fields $Z^{\pm}$  and the switching manifold $\Sigma$ separate sliding and crossing type regions. The following definition characterises
the least degenerate case: 

\begin{definition}
	\defnlab{fold}
	Consider system \eqref{PWS}. A point $F \in \Sigma$ is called a fold point if either
	\begin{equation}
	\eqlab{fold}
	Z^+f(F) = 0 , \qquad Z^+(Z^+f)(F) \neq 0 , \qquad \text{or} \qquad Z^-f(F) = 0 , \qquad Z^-(Z^-f)(F) \neq 0 .
	\end{equation}
	A fold point $F$ with $Z^+f(F) = 0$ is visible (invisible) if the inequality $Z^+(Z^+f)(F) \neq 0$ is positive (negative). Conversely, a fold point $F$ with $Z^-f(F) = 0$ is visible (invisible) if the inequality $Z^-(Z^-f)(F) \neq 0$ is negative (positive).
\end{definition}

\note{Note that in general, the location of a fold point $F$ also depends on $\alpha$ via the Lie derivatives \eqref{fold}, see again \defnref{cross_slide}.} \figref{pws_fig} shows the dynamics near a visible fold point. 
Given the Filippov convention for defining a vector field on $\Sigma_{sl}$, solutions may be constructed in the PWS approach as concatenations of trajectory segments in $\Sigma^\pm$ and $\Sigma_{sl}$ where solutions can leave $\Sigma_{sl}$ only via its boundary, i.e.~via fold points.
%
\begin{definition}
	\defnlab{pws_cycle}
	A PWS cycle in system \eqref{PWS} is a closed concatenation of orbit segments containing at least one orbit segment from each of $\Sigma_{sl}$ and $\mathbb R^2 \setminus \Sigma_{sl}$.
\end{definition}
It remains to consider bifurcations specific to PWS systems. Our sole focus is on
the \textit{boundary-focus (BF) bifurcation} for PWS systems, which refers to a collision of a hyperbolic focus with the switching manifold under parameter variation. 

\begin{definition}
	\defnlab{bf}
	A PWS system \eqref{PWS} has a BF bifurcation at $z = z_{bf} = (x_{bf},0) \in \Sigma$ for $\alpha = \alpha_{bf}$ if the following conditions hold:
	\begin{equation}
	\eqlab{bf_conds}
	Z^+(z_{bf},\alpha_{bf}) = (0,0)^T , \ \ Z_2^-(z_{bf},\alpha_{bf}) \neq 0 , \ \ \det \left(\frac{\partial Z^+}{\partial \alpha} \big| \frac{\partial Z^+}{\partial x} \right)\bigg|_{(z_{bf},\alpha_{bf})} \neq 0 , \ \ \frac{\partial Z_{sl}}{\partial x}\bigg|_{(x_{bf},\alpha_{bf})} \neq 0 ,
	\end{equation}
	where  $Z^\pm = (Z^\pm_1,Z^\pm_2)^T$, $\det (X | Y)$ denotes the determinant of the matrix with columns $X,Y$, and
	\begin{equation}
	\eqlab{nondegeneracy2}
	\lambda_\pm(\alpha_{bf}) = A \pm B i , \qquad A ,  B \neq 0 ,
	\end{equation}
	where $\lambda_\pm(\alpha_{bf})$ denote the eigenvalues of the Jacobian $(\partial Z^+/\partial z)|_{(z_{bf},\alpha_{bf})}$.
\end{definition}


\defnref{bf} characterises a BF bifurcation in \eqref{PWS} as the transversal collision of a focus-type equilibrium with the switching manifold $\Sigma$ under variation of $\alpha$. Note that the vector field $Z^-(z,\alpha)$ is transverse to $\Sigma$ near the BF point by the genericity condition $Z_2^-(z_{bf},\alpha_{bf}) \neq 0$. In total, there are $2^4=16$ possible cases at the collision value $\alpha = \alpha_{bf}$, depending on 
\begin{enumerate}
	\item[(i)] $\sgn(Z_2^-(z_{bf},\alpha_{bf}))$, 
	\item[(ii)] $\sgn((\partial Z_{sl} / \partial x )|_{(z_{bf},\alpha_{bf})})$, 
	\item[(iii)] orientation of the rotation (clockwise/counter-clockwise), and 
	\item[(iv)] stability/instability of the focal point;
\end{enumerate}
see \cite[Figure 1]{Carvalho2018} for a complete set of normal forms at the collision value $\alpha = \alpha_{bf}$ in the PWS setting. \figref{collision_figs} shows three cases of particular interest. Notice that in general, a BF bifurcation occurs at a (degenerate) tangency point, and is associated with the local transition of a visible fold point (prior to collision), to an invisible fold point (following the collision); see also \cite[Figure~5]{Kuznetsov2003}.
%
Moreover, particular cases of BF bifurcation are correlated with the termination of PWS cycles consisting of orbit segments in $\Sigma_{sl}$ and $\mathbb R^2 \setminus \Sigma$. In \figref{collision_figs}, PWS cycles $\Gamma = \Gamma^1 \cup \Gamma^2$ terminate in the collision limit $\alpha = \alpha_{bf}^+$ in cases BF$_1$ and BF$_3$ (top and bottom panels respectively).

\begin{figure}[h!]
	\begin{center}
		\subfigure[]{\includegraphics[width=.32\textwidth]{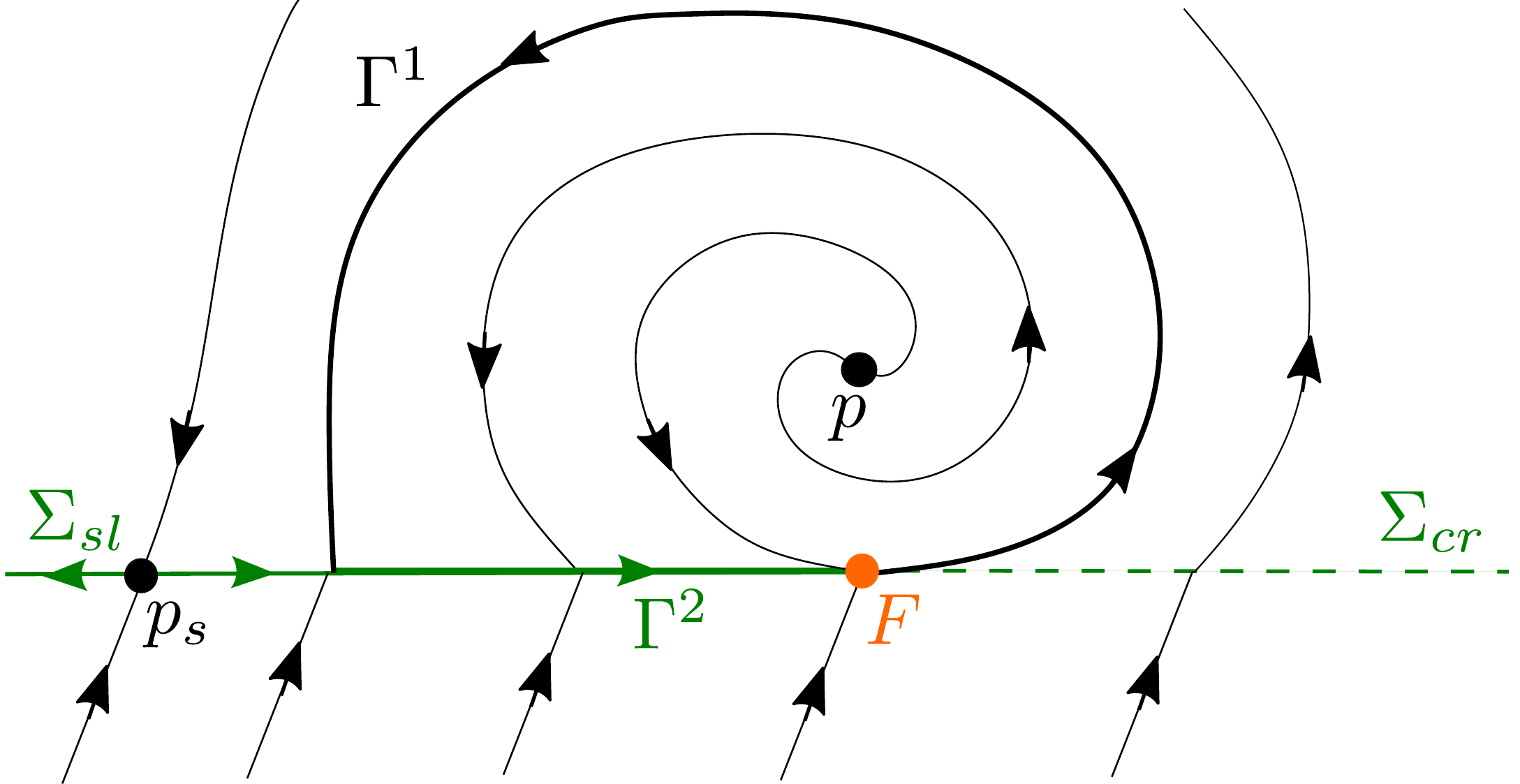}}
		\subfigure[]{\includegraphics[width=.3\textwidth]{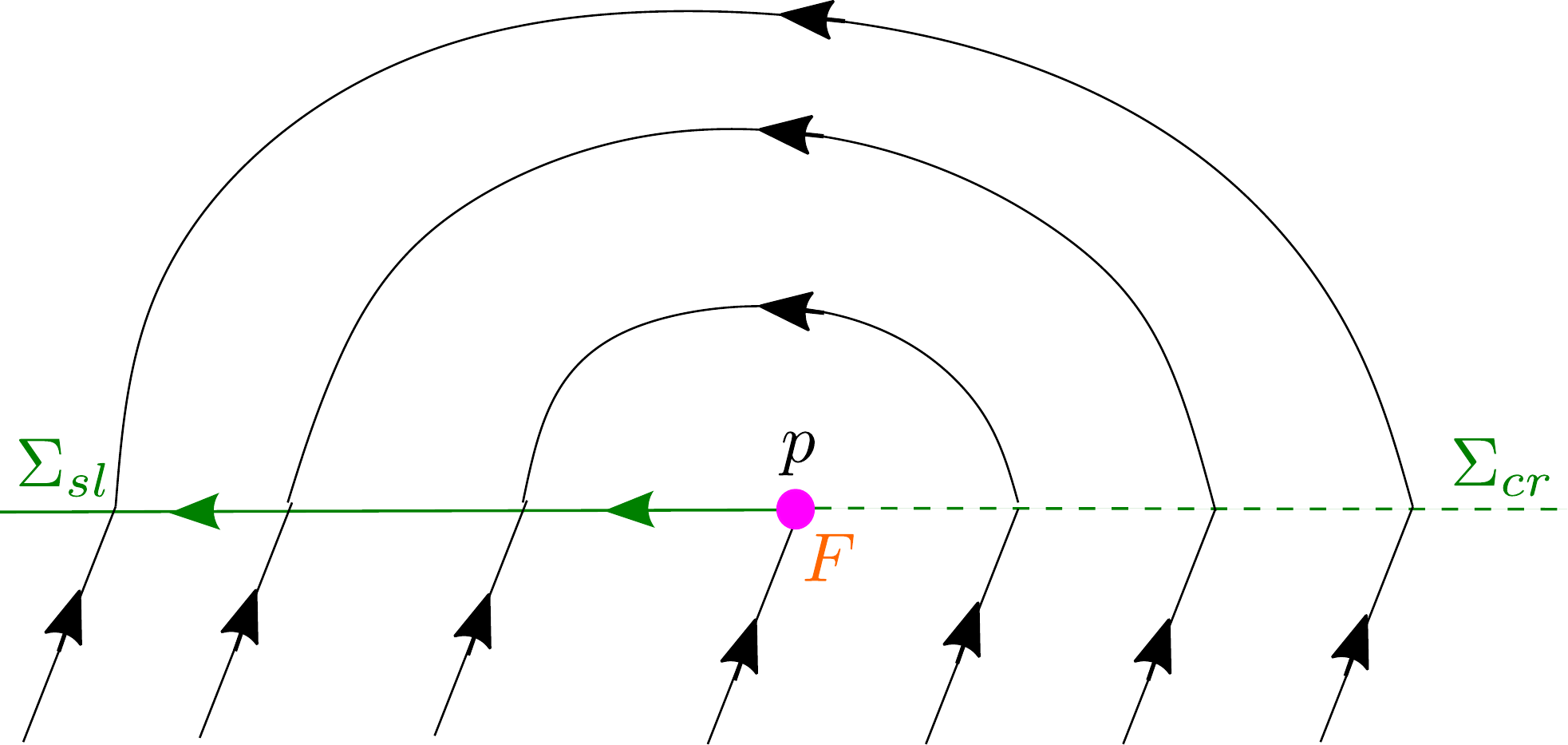}}
		\subfigure[]{\includegraphics[width=.3\textwidth]{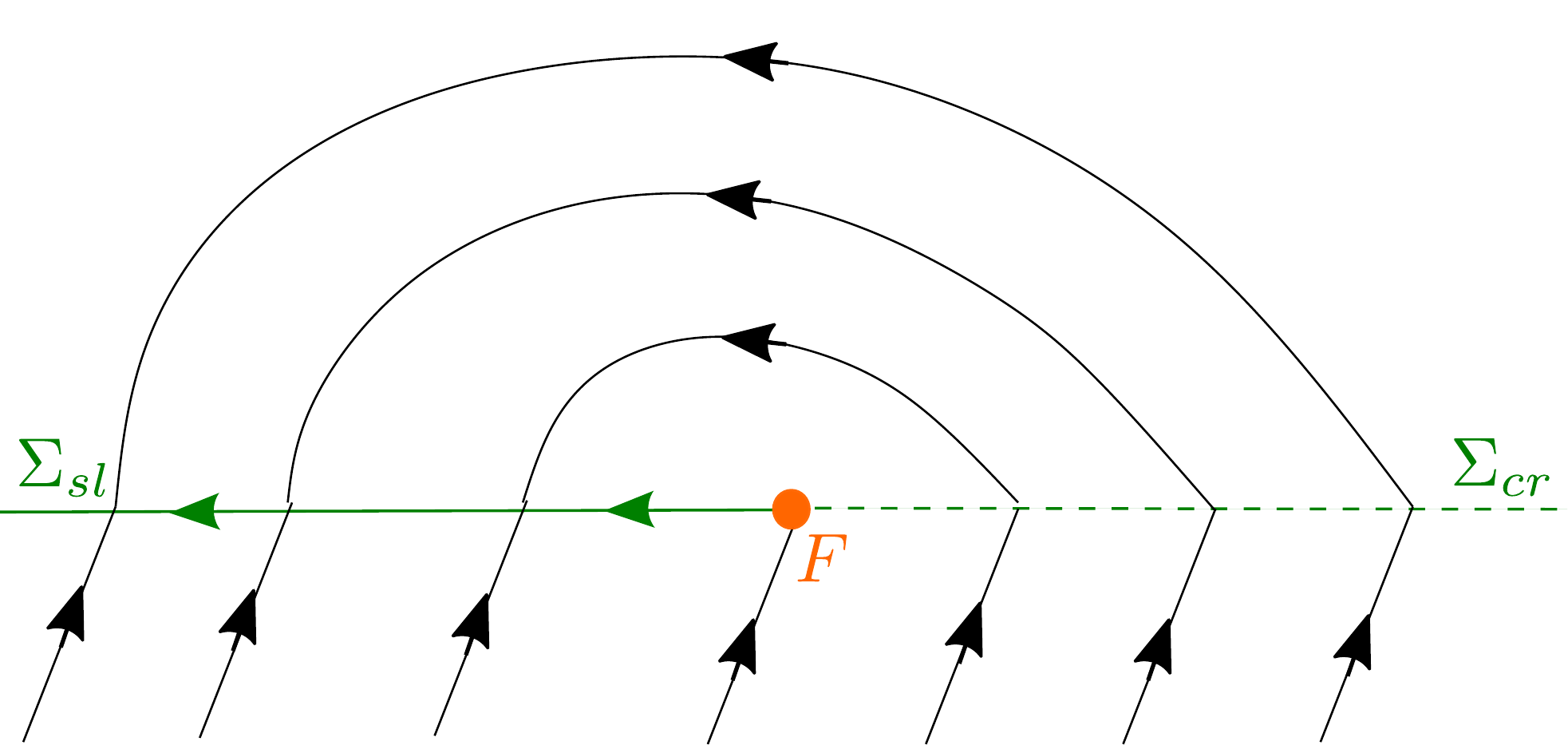}}
		
		\subfigure[]{\includegraphics[width=.305\textwidth]{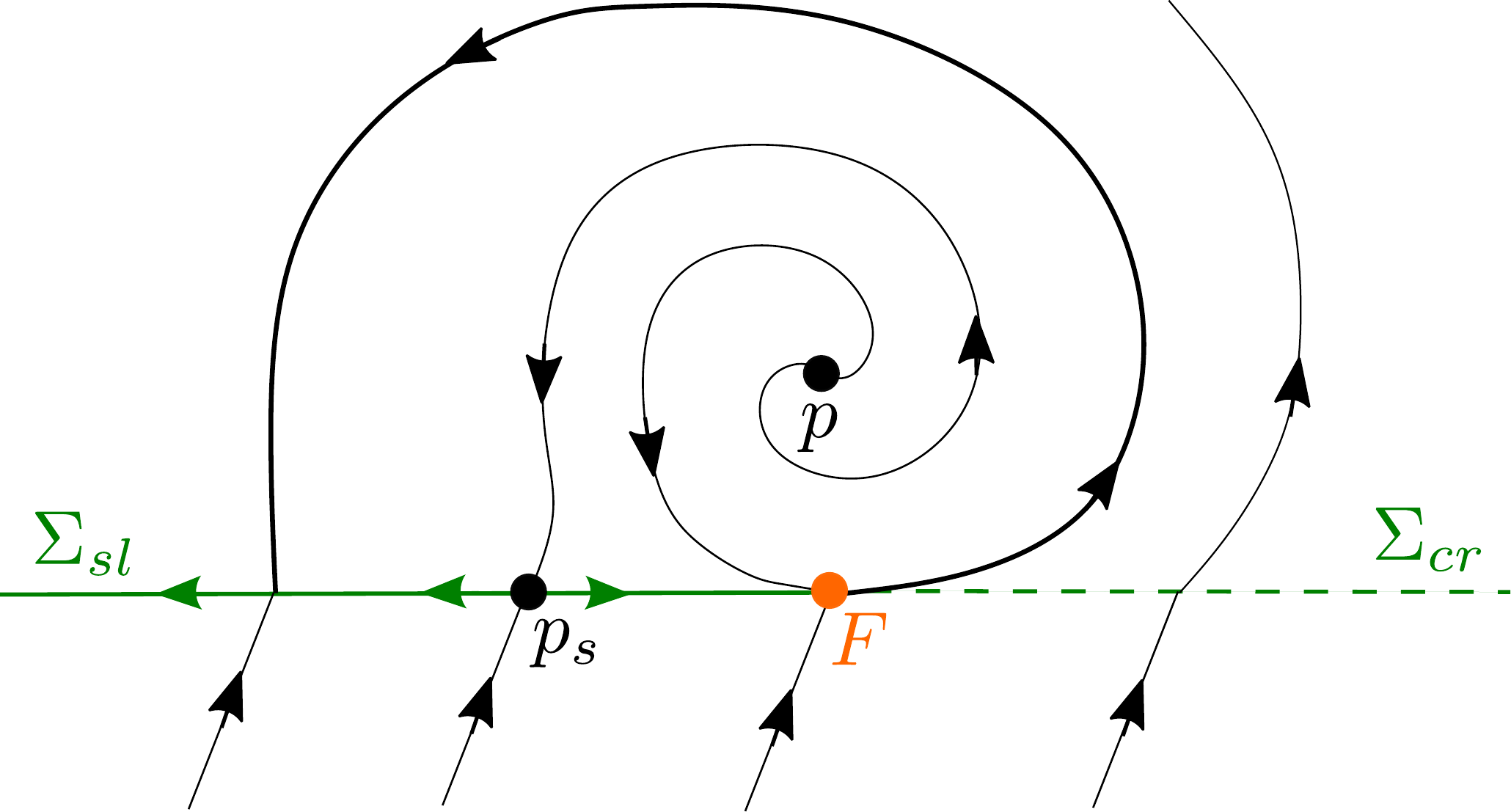}}
		\subfigure[]{\includegraphics[width=.3\textwidth]{./beb_figures/bfb3}}
		\subfigure[]{\includegraphics[width=.3\textwidth]{./beb_figures/pws_neg_12}}
		
		\subfigure[]{\includegraphics[width=.3\textwidth]{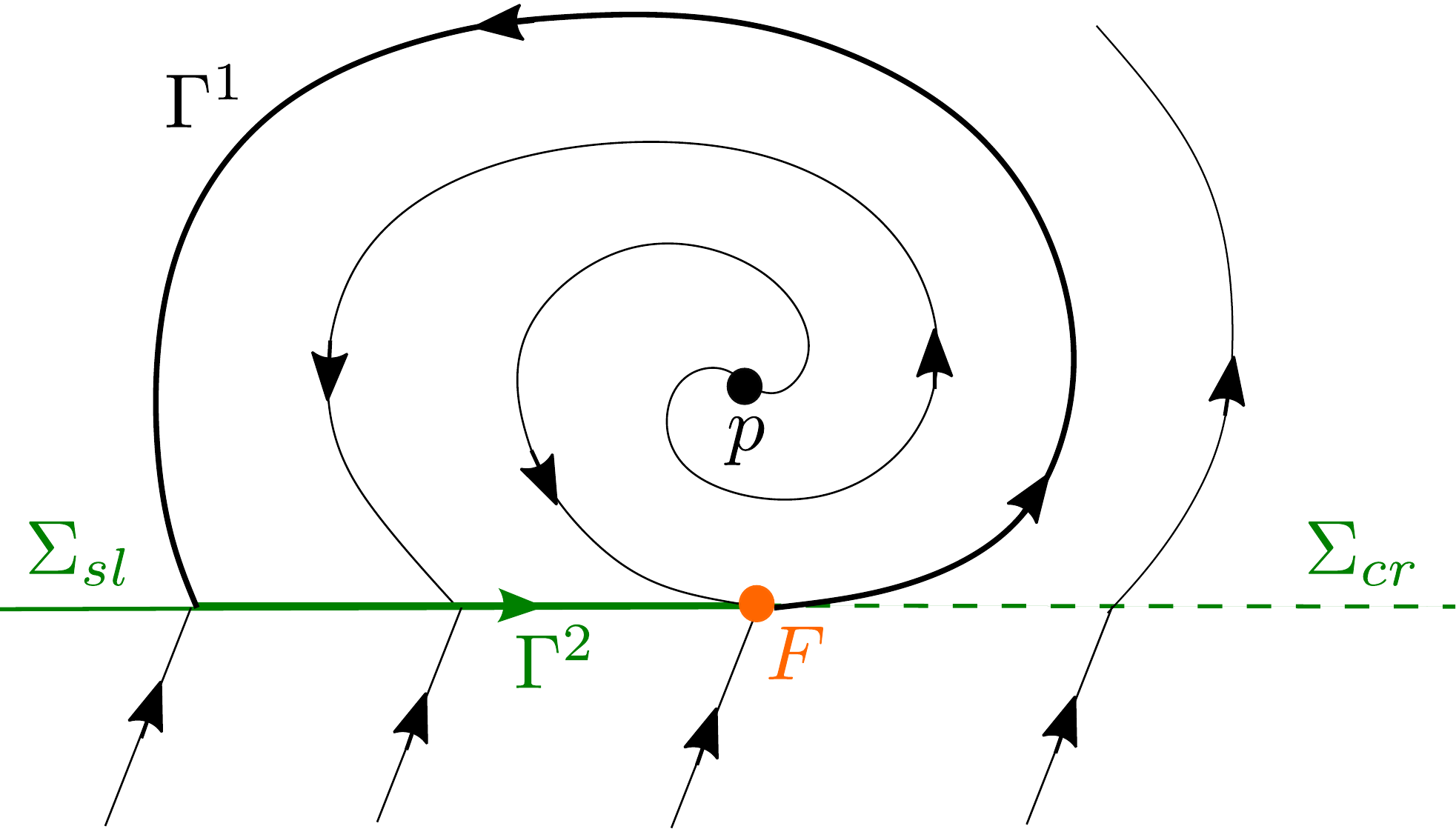}}
		\subfigure[]{\includegraphics[width=.305\textwidth]{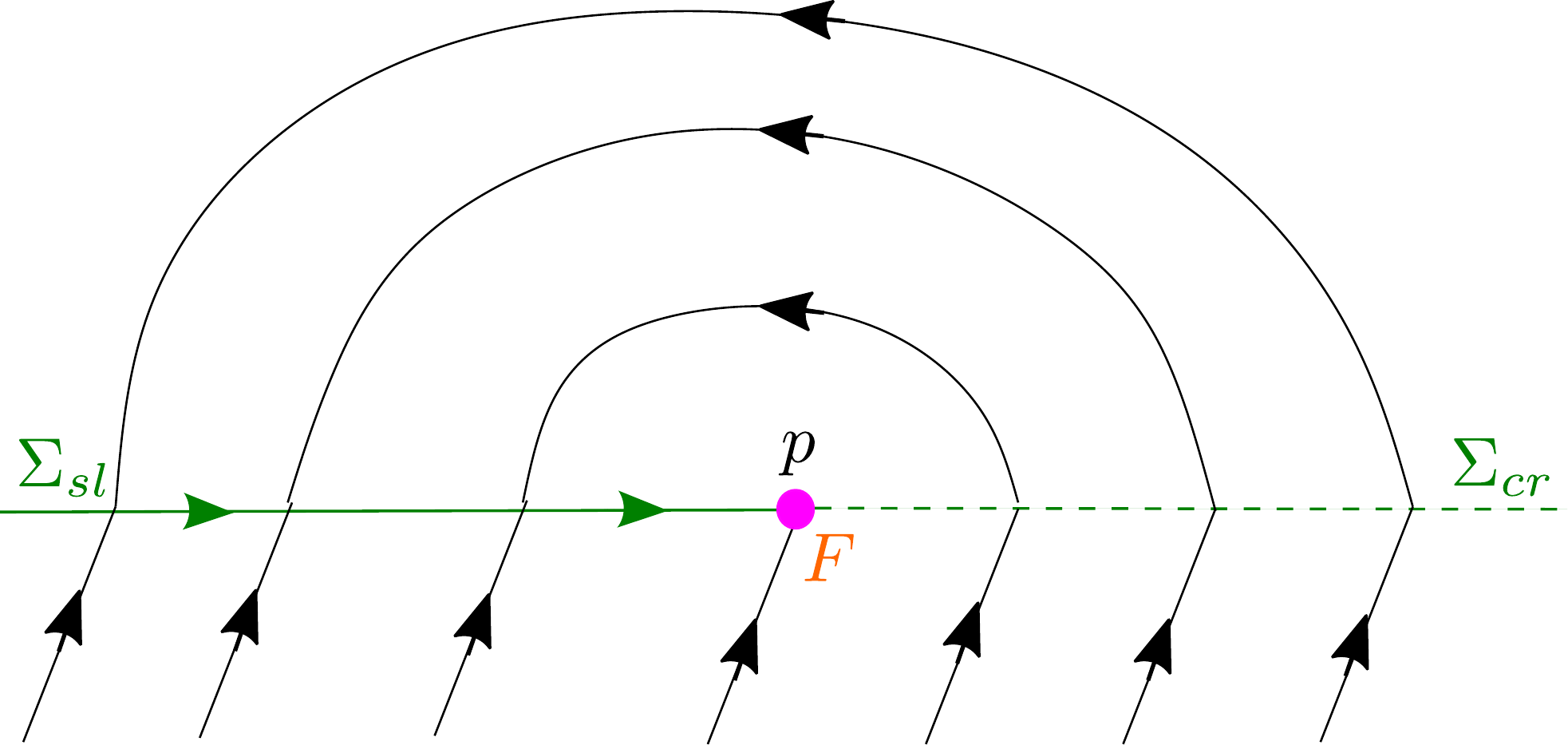}}
		\subfigure[]{\includegraphics[width=.305\textwidth]{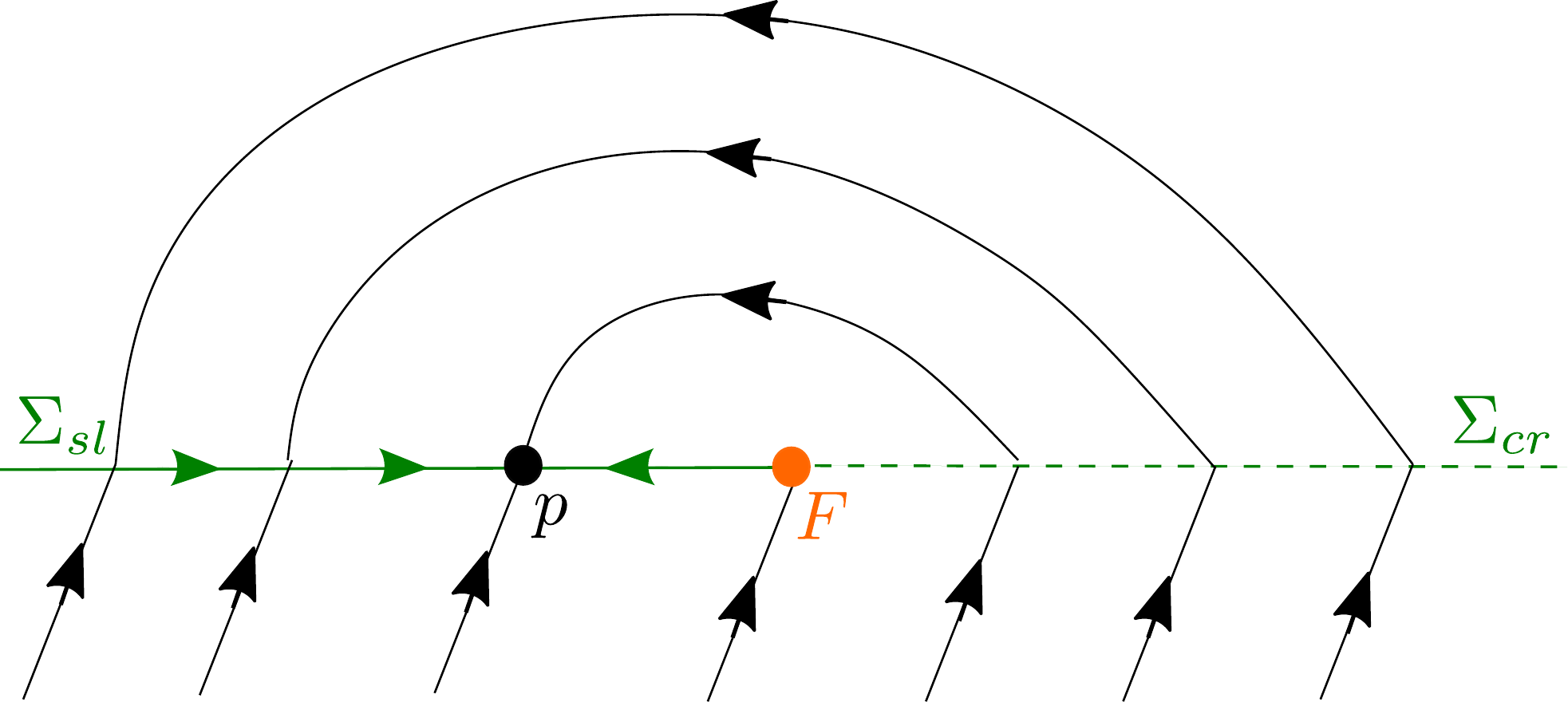}}
		\caption{BF$_i$ bifurcations, $i=1,2,3$, assuming the determinant quantity in \eqref{bf_conds} is positive. The top, middle, and bottom panels show cases $i=1,2,3$ respectively, and $\alpha > \alpha_{bf}$, $\alpha = \alpha_{bf}$, $\alpha_{bf} < \alpha_{bf}$ in the left, middle and right panels respectively. In each case, sliding and crossing regions $\Sigma_{sl}$ and $\Sigma_{cr}$ (bold and dashed green resp.) are separated by a fold singularity $F$ which is visible (invisible) for $\alpha > \alpha_{bf}$ ($\alpha < \alpha_{bf}$), and degenerate for $\alpha = \alpha_{bf}^+$.
			PWS cycles $\Gamma = \Gamma^1 \cup \Gamma^2$ which terminate in the limit $\alpha \to \alpha_{bf}^+$ can be constructed in cases BF$_1$ and BF$_3$ (shown in (a) resp. (g)), but not in case BF$_2$ due to the location of a sliding equilibrium $p_s \in \Sigma_{sl}$.}
		\figlab{collision_figs}
	\end{center}
\end{figure}


\subsection{Local normal form for the {\em smooth} BF bifurcation}
\seclab{normal_form1}

In the PWS literature, see e.g.  \cite{Carvalho2018}, the normal forms obtained are valid only up to `$\Sigma-$equivalence' (roughly speaking, homeomorphism plus preservation of the sliding vector field; see \cite[Definition~4]{Carvalho2018}). For our purposes, however, a \textit{smooth} normal form is required to resolve smooth dynamics limiting on a PWS bifurcation. Unlike \cite{Carvalho2018}, we will also require that our system exhibits a `collision' under parameter variation, as per the determinant condition in \eqref{bf_conds}. 


\begin{proposition}
	\proplab{prop_normal_form}
	Consider system \eqref{main} under \assumptionref{ass1}, \assumptionref{ass1b}, and assume that the PWS system \eqref{PWS} obtained in the limit $\epsilon \to 0$ has a BF bifurcation. Then there exists constants
	\[
	\tau \neq 0, \qquad  \delta > \tau^2 / 4 , \qquad \gamma \in \mathbb R ,
	\]
	such that system \eqref{main} can be smoothly transformed to \SJ{one of the following local normal forms}:
	\begin{equation}
	\eqlab{normal_form}
	\begin{pmatrix}
	\dot{x}\\\dot{y}
	\end{pmatrix}
	=
	\begin{pmatrix}
	\pm (\tau - \gamma) + \phi\left(y \epsilon^{-1} \right) \left( \mp (\tau - \gamma) + \mu + \tau x - \delta y + \theta_1(x,y,\mu) \right) \\
	\pm 1 + \phi\left(y \epsilon^{-1} \right) \left(\mp 1 + x + \theta_2(x,y,\mu) \right)
	\end{pmatrix}
	= :X(x,y,\mu,\epsilon) ,
	\end{equation}
	where
	$\theta_i(x,y,\mu)$, $i=1,2$ are real-valued smooth functions such that
	\note{
	\[
	\theta_1(x,y,\mu) = 
	\mathcal O(x^2,xy,y^2,x\mu,y\mu,\mu^2) , \qquad 
	\theta_2(x,y,\mu) = \mathcal O \left(x^2, xy, y^2, x \mu, y \mu \right) .
	\]
	}
	The PWS system \W{
		\begin{equation}
		\eqlab{PWS_normal_form}
		\begin{pmatrix}
		\dot{x}\\\dot{y}
		\end{pmatrix}
		=
		\left\{
		\begin{aligned}
		&
		\begin{pmatrix}
		\mu + \tau x - \delta y + \theta_1(x,y,\mu)  \\
		x + \theta_2(x,y,\mu) 
		\end{pmatrix}
		=: X^+(x,y,\mu),
		& \quad(y>0) ,
		\\
		&
		\begin{pmatrix}
		\pm (\tau - \gamma)  \\
		\pm 1 
		\end{pmatrix}
		=: X^-(x,y,\mu),
		&\quad(y<0) ,
		\end{aligned}
		\right.
		\end{equation}
	}
	obtained from \eqref{normal_form} in the limit $\epsilon \to 0^+$ \SJ{has a stable (unstable) sliding region in case $X_2^- = +1$ (\SJnew{$X_2^- = -1$}), and a} BF collision at the origin for $\mu = 0$. The Filippov vector field
	is
	\begin{equation}
	\eqlab{Filippov_VF}
	\dot x = \frac{\mu + \gamma x + \theta_1(x,0,\mu) - (\tau - \gamma) \theta_2(x,0,\mu)}{1 - x - \theta_2(x,0,\mu)} =: X_{sl}(x,\mu) , \qquad \forall (x,0) \in \Sigma_{sl} .
	\end{equation}
\end{proposition}

\bpr
	See \appref{proof_of_normal_form}.
\epr

\begin{remark}
	\remlab{normal_forms}
	The normal forms in \eqref{normal_form} and \eqref{PWS_normal_form} assume a locally counter-clockwise rotation near the focus in \propref{prop_normal_form}; the equivalent normal forms with locally clockwise rotation are obtained from \eqref{normal_form} and \eqref{PWS_normal_form} by the transformation $(x,\mu) \mapsto (-x,-\mu)$ (see \appref{proof_of_normal_form}).
\end{remark}



Following the (PWS) topological classification for BF bifurcations given in \cite{Kuznetsov2003}, there are five generic cases (up to orientation), depending on 
\begin{enumerate}
	\item[(i)] the orientation of the `lower' vector field $X^-$ relative to $\Sigma$, 
	\item[(ii)] the orientation of the sliding dynamics on $\Sigma_{sl}$, and 
	\item[(iii)] the existence of PWS cycles depending on the location of equilibria on $\Sigma_{sl}$ relative to the global return. 
\end{enumerate}
We label these cases BF$_i$, $i=1,\ldots,5$ in accordance with the conventions initially put forward in \cite[Figure~5]{Kuznetsov2003}. With respect to our normal form \eqref{normal_form}, we can classify a smooth analogue for each case BF$_i$, $i=1,\ldots,5$.
In particular, $\gamma > 0$ ($\gamma < 0$) is correlated with the existence of an unstable (stable) sliding equilibrium on $\Sigma_{sl}$, and $\sgn(X_2^-)$ determines the local orientation of $X^-$ with respect to $\Sigma$, which is important in determining stability of the sliding manifold $\Sigma_{sl}$.

\begin{classification}
	Bifurcation types BF$_i$, $i = 1 , \ldots , 5$ given in \cite[Figure~5]{Kuznetsov2003} are obtained from the (PWS) normal form \eqref{PWS_normal_form}  as follows:
	\SJ{
		\begin{enumerate}
			\item[(BF$_1$)] $X_2^-=+1, \ \gamma > 0$: stable sliding with an equilibrium $p_s \in \Sigma_{sl}$ and a PWS cycle $\Gamma$ when $\mu > 0$.
			\item[(BF$_2$)] $X_2^-=+1, \ \gamma > 0$: stable sliding with an equilibrium $p_s \in \Sigma_{sl}$ when $\mu > 0$, and no PWS cycles.
			\item[(BF$_3$)] $X_2^-=+1, \ \gamma < 0$: stable sliding with a PWS cycle $\Gamma$ when $\mu > 0$, and an equilibrium $p \in \Sigma_{sl}$ when $\mu < 0$.
			\item[($BF_4$)] \SJnew{$X_2^-=-1, \ \gamma < 0$}: unstable sliding with an  
			\SJnew{equilibrium $p \in \Sigma_{sl}$ when $\mu < 0$, and no PWS cycles}.
			\item[($BF_5$)] \SJnew{$X_2^-=-1, \ \gamma > 0$}: unstable sliding with an equilibrium $p_s \in \Sigma_{sl}$ 
			when $\mu > 0$\SJnew{, and no PWS cycles.}
		\end{enumerate}
	}
\end{classification}

Notice that the `local' conditions for BF$_1$ and BF$_2$ coincide
; these cases are distinguished by global mechanisms, in particular the location of the `drop point' with respect to a sliding equilibrium $p_s$\SJ{; c.f. \figref{collision_figs}a and \figref{collision_figs}d.}

\

In the following, we restrict attention to the cases of interest in our work featuring stable sliding and an unstable focus. Explicitly, we impose the following assumptions: 

\begin{assumption}
	\assumptionlab{ass2}
	The PWS system \eqref{PWS} obtained from \eqref{main} in the singular limit $\epsilon \to 0^+$ has a BF bifurcation at $z_{bf} \in \Sigma$ for $\alpha = \alpha_{bf}$. We also assume the following:
	\begin{enumerate}
		\item[\SJ{(i)}] The lower vector field points towards $\Sigma$, i.e. $Z_2^-(z_{bf},\alpha_{bf}) > 0$. In terms of the normal form \eqref{normal_form},
		\[ 
		\begin{pmatrix}
		\dot{x}\\\dot{y}
		\end{pmatrix}
		=
		\begin{pmatrix}
		\tau - \gamma \\
		1
		\end{pmatrix}
		=X^- ,
		\qquad (y < 0) .
		\]
		\item[\SJ{(ii)}] The colliding equilibrium is of unstable focus type. In terms of the normal form \eqref{normal_form}, $\tau > 0$ so that
		\begin{equation}
		\eqlab{nondegeneracy3}
		\lambda_\pm(0) = A \pm i B , \qquad A ,  B > 0 .
		\end{equation}
	\end{enumerate}
\end{assumption}

Hence in the conventional terminology introduced in \cite{Kuznetsov2003}, we are considering BF bifurcations of type BF$_i$ for $i = 1,2,3$. The cases BF$_i$, $i=4,5$ are similar, only of less interest in our study due to the instability of the sliding region $\Sigma_{sl}$.

\subsection{Singular limit analysis for system \eqref{normal_form}}
\seclab{dynamics_in_the_singular_limit}

A singular limit analysis for the normal form \eqref{normal_form} amounts to a PWS analysis of the PWS system \eqref{PWS_normal_form}.
Our findings 
are summarised in the following result; see again \figref{collision_figs}.

\begin{lemma}
	\lemmalab{PWS_lemma}
	Consider the PWS system \eqref{PWS_normal_form} obtained from system \eqref{normal_form} in the limit $\epsilon \to 0^+$, under \assumptionref{ass2},
	\W{and define the intervals 
		\begin{equation}
		\eqlab{S2_intervals}
		\mu \in \mathcal I_- := (- \mu_+, - \mu_- )  \qquad \text{or} \qquad \mu \in \mathcal I_+ := (\mu_-, \mu_+) ,
		\end{equation}
		where $\mu_+ > \mu_- > 0$ can be chosen arbitrarily small.}
	Then the following assertions hold:
	\begin{enumerate}
		\item[(i)] \SJnew{There exists a smooth function $x_F(\mu)$ such that \eqref{PWS_normal_form} has a visible (invisible) fold for $\mu>0$ ($\mu<0$) at 
		\[
		F := \left(x_F(\mu), 0 \right) \in \Sigma , \qquad x_F(\mu) = \mathcal O\left(\mu^2\right) ,
		\]
		which divides the switching manifold into} 
		the union $\Sigma = \Sigma_{sl} \cup F \cup \Sigma_{cr}$, where
		\[
		\Sigma_{sl} = \left\{(x,0) : x < x_F(\mu) \right\} , \qquad \Sigma_{cr} = \left\{(x,0) : x > x_F(\mu) \right\} .
		\]
		\item[(ii)] For all $\mu \in \mathcal I_+$ 
		with sufficiently small $\mu_+>\mu_->0$, there exists an equilibrium
		\begin{equation}
		\eqlab{focus}
		p(\mu) : \SJ{ \left(\mathcal O(\mu^2) , \frac{\mu}{\delta} + \mathcal O(\mu^2) \right) }
		\end{equation}
		of unstable focus type. \SJ{If additionally $\gamma > 0$, there also exists} an equilibrium
		\begin{equation}
		\eqlab{saddle_eq}
		p_s (\mu): (x_s(\mu),0) = \left(-\frac \mu \gamma + \mathcal O(\mu^2), 0 \right) ,
		\end{equation}
		\SJ{which} is unstable as an equilibrium on $\Sigma_{sl}$. 
		\SJ{Both $p(\mu), p_s(\mu) \to (0,0)$} in the collision limit $\mu \to 0^+$.
		\item[(iii)] For all $\mu \in \mathcal I_-$ 
		with sufficiently small  $\mu_+>\mu_->0$ \SJ{and $\gamma < 0$}, there exists \SJ{an} equilibrium
		\[
		p(\mu) : \left(-\frac \mu \gamma + \mathcal O(\mu^2), 0 \right) ,
		\]
		\SJ{which }is stable as an equilibrium on $\Sigma_{sl}$ \SJ{and satisfies $p(\mu) \to (0,0)$ as $\mu \to 0^-$.}
		
		\item[(iv)] For all $\mu \in \mathcal I_+$ 
		with sufficiently small $\mu_+>\mu_->0$, the $X^+$ trajectory $\Gamma^1$ in $\{y \geq 0\}$ which grazes the visible fold point $F$ in \figref{collision_figs}a,d,g has a first return or `drop point' $(x_d(\mu),0) \in \Sigma_{sl}$, allowing for the construction of the $PWS$ cycles $\Gamma = \Gamma^1 \cup \Gamma^2$ in cases BF$_i$, $i=1,3$ shown in \figref{collision_figs}a and \figref{collision_figs}g respectively.
	\end{enumerate}
\end{lemma}

\bpr
	Assertions (i)-(iii) follow from direct calculations and suitable application of the implicit function theorem.
	
	Assertion (iv): consider the vector field $X^+(x,y,\mu)$ extended to all of $\mathbb R^2$. Since by construction $X^+(x,y,\mu)$ has an hyperbolic unstable focus arbitrarily close to $(0,0)$ for sufficiently small $\mu > 0$, one can always guarantee a return mechanism due to the repulsion and rotation associated with the focal point. Thus the orbit segment $\Gamma^1$ in \figref{collision_figs}a and \figref{collision_figs}g can be identified with the segment of the $X^+$ orbit intersecting $F$ between $F$ and the orbits first return at $(x_d(\mu),0) \in \Sigma_{sl}$. The cases BF$_i$, $i=1,2$ are then distinguished as follows:
	\[
	\text{BF}_1: \ x_d(\mu) > x_s(\mu) , \qquad \text{BF}_2: \ x_d(\mu) < x_s(\mu) ,
	\]
	where the inequalities hold for all $\mu \in (0, \mu_+)$. 
	In cases BF$_1$ and BF$_3$, PWS cycles $\Gamma = \Gamma^1 \cup \Gamma^2$ can be constructed for all $\mu \in \mathcal I_+$, where
	\[
	\Gamma^2 = \left\{(x,0) : x \in [x_d(\mu),x_F(\mu)] \right\} .
	\]
\epr



The preceding analysis also provides us with a useful means of distinguishing between the cases BF$_i$, $i=1,2$ and BF$_3$, which we will use as a shorthand distinction for the remainder of this work: the fact that there exists an equilibrium $p_s \in \Sigma_{sl}$ in the former two cases but not the latter tells us that
\begin{equation}
\eqlab{bf_type}
\text{BF}_i, \ i = 1,2 : \ \gamma > 0 , \qquad  \text{BF}_3 : \ \gamma < 0 .
\end{equation}
The $\gamma = 0$ case `between' BF$_i$, $i=1,2$ and BF$_3$ also shows up in applications, and is considered further in \secref{stick-slip} and \secref{degenerate_cases}.







\section{Main results}
\seclab{main_results}


This section is devoted to presentation of the main results on the smooth normal form
\eqref{normal_form}.
%
The main analytical tool used to obtain these results is the blow-up method developed in \cite{Dum1996,Krupa2001a,Krupa2001b}.
The highly singular nature of this problem requires {\em multiple blow-ups} to unfold the BF bifurcation. Hence, we decided to present the reader with a clear and concise, yet sufficiently detailed geometric picture to illustrate the results while deferring the detailed analysis to the appendix sections. Our aim  is to present the geometric intuition behind this method, which is its major appeal.

There are two important scaling regimes in system \eqref{normal_form} with respect to the parameter $\mu$:
\begin{enumerate}
	\item[(S1):] $\mu  = \hat \mu \epsilon^{k/(k+1)} = \mathcal O \left(\epsilon^{k/(k+1)} \right)$, where $\hat \mu \in \mathcal I := (-\hat \mu_0, \hat \mu_0)$ for arbitrarily large $\hat \mu_0 > 0$;
	\item[(S2):] $\mu = \mathcal O(1)$ or more precisely, $\mu \in \mathcal I_- \cup \mathcal I_+$
	\W{where $\mathcal I_\pm$ are the intervals in \eqref{S2_intervals}.}
\end{enumerate}
\SJ{Notice that (S1) and (S2) are distinct (i.e. non-overlapping) for $0 < \epsilon \ll 1$ sufficiently small.} The unfolding of the BF point $(x,y,\epsilon,\mu) = (0,0,0,0)$ occurs in regime (S1), where no less than three successive blow-ups are required: two blow-ups are necessary to resolve the tangency at $(0,0,0,\mu)$
independently of $\mu$, and  an additional blow-up including $\mu$ is necessary to resolve the singularity persisting for $\mu = 0$. 
%
The dynamics on `either side' of the regularised collision are studied in (S2), including persistence of PWS cycles as relaxation oscillations when $\mu \in \mathcal I_+$. 
%
Finally, \SJ{by including \textit{both $\epsilon$ and $\mu$} in the blow-up, we are able to}  \W{establish} a connection between the dynamics in regimes (S1) and (S2).
This allows for a description of the connection associated with the matching problem\SJnew{s obtained for}
\begin{equation}
\eqlab{matching_limit_neg}
\hat \mu \to - \infty , \qquad \mu \to 0^- , 
\qquad \text{((S1)-(S2) connection, $\mu < 0$)},
\end{equation}
\SJnew{and}
\begin{equation}
\eqlab{matching_limit}
\hat \mu \to + \infty , \qquad \mu \to 0^+ ,
\qquad \text{((S1)-(S2) connection, $\mu > 0$)}\,. 
\end{equation}


\subsection{Blow-up of the switching manifold and fold singularity}
\seclab{Sigma_blowup_outline}

In both scaling regimes (S1) and (S2), we must first resolve 
the degeneracy associated with a loss of smoothness along $\Sigma$ when $\epsilon = 0$, which occurs independently of $\mu$. As is known from e.g. \cite{Jelbart2019c,Kristiansen2019,Kristiansen2019c,Kristiansen2019d}, this can be achieved by means of a cylindrical blow-up in the extended $(x,y,\epsilon)-$space\footnote{See \appref{blow-up_of_Sigma} for details; the blow-up transformation itself is given by \eqref{cyl_bu}.}.
In the blown-up space, the system regains smoothness and the switching manifold $\Sigma \times \{0\} \subset \mathbb R^2 \times \mathbb R_+$ is replaced by a cylinder; see \figref{S1_blowups_i}a and \figref{S1_blowups_i}b. In this approach, the dynamics with rescaled coordinates $(x, Y, \epsilon) = (x, y/\epsilon, \epsilon)$, i.e. within the switching layer, are identified with dynamics on the cylinder itself.
Here, a \SJnew{classical} slow-fast system is found with an attracting critical manifold $S$, which connects to the intersection of the blow-up cylinder with the plane $\{\epsilon = 0\}$ (denoted $l_s$) at a nonhyperbolic point $Q$, see \figref{S1_blowups_i}a-b.

\begin{figure}[t!]
	\begin{center}
		\subfigure[]{\includegraphics[width=.40\textwidth]{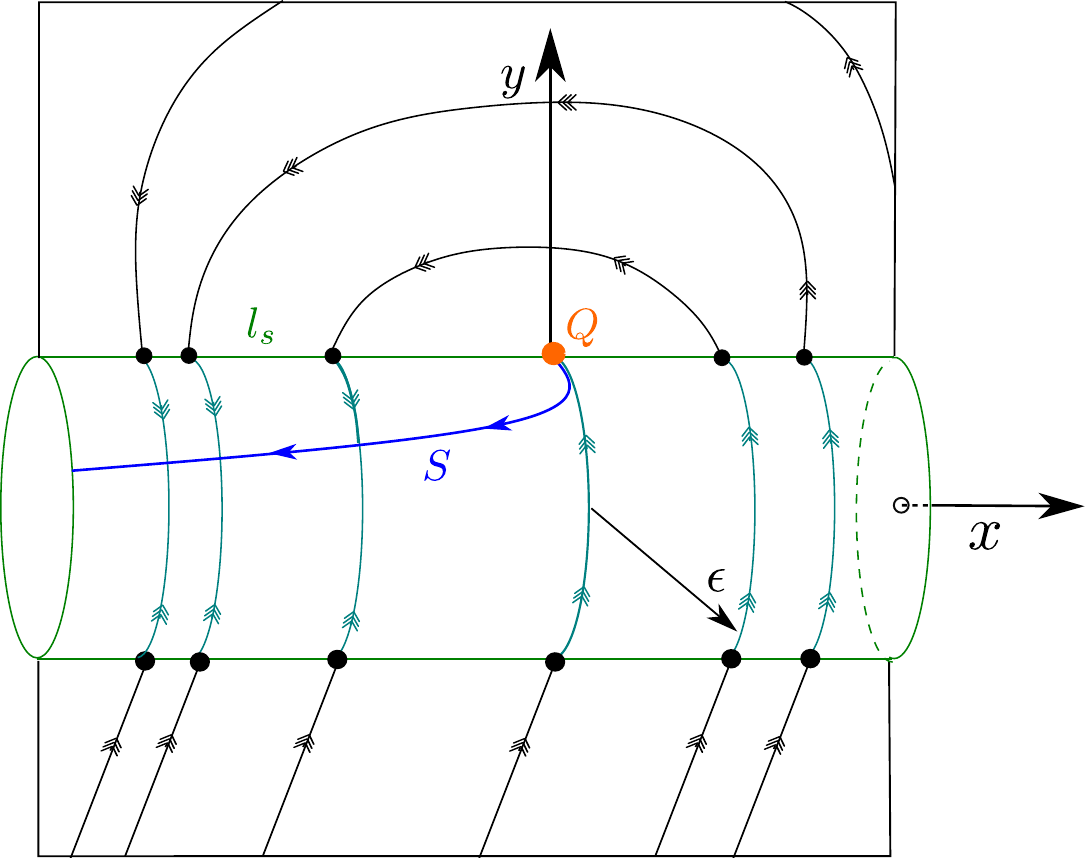}}
		\subfigure[]{\includegraphics[width=.40\textwidth]{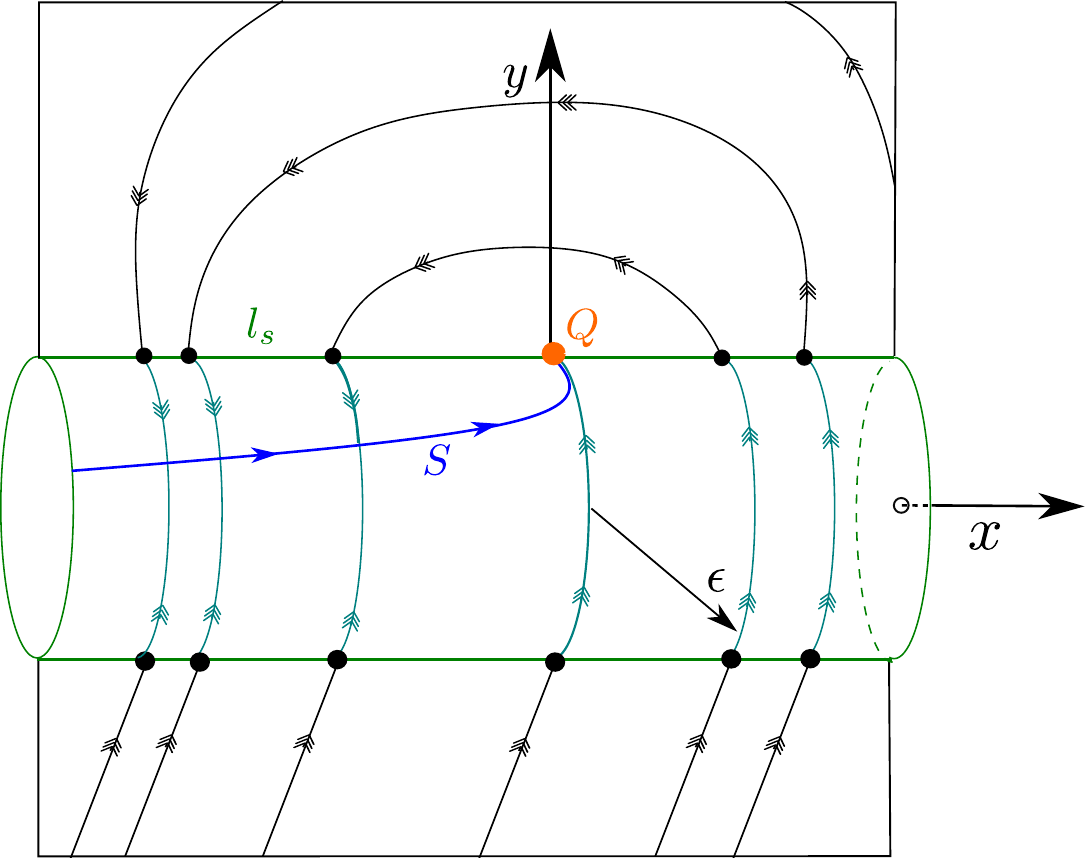}}
		\subfigure[]{\includegraphics[width=.40\textwidth]{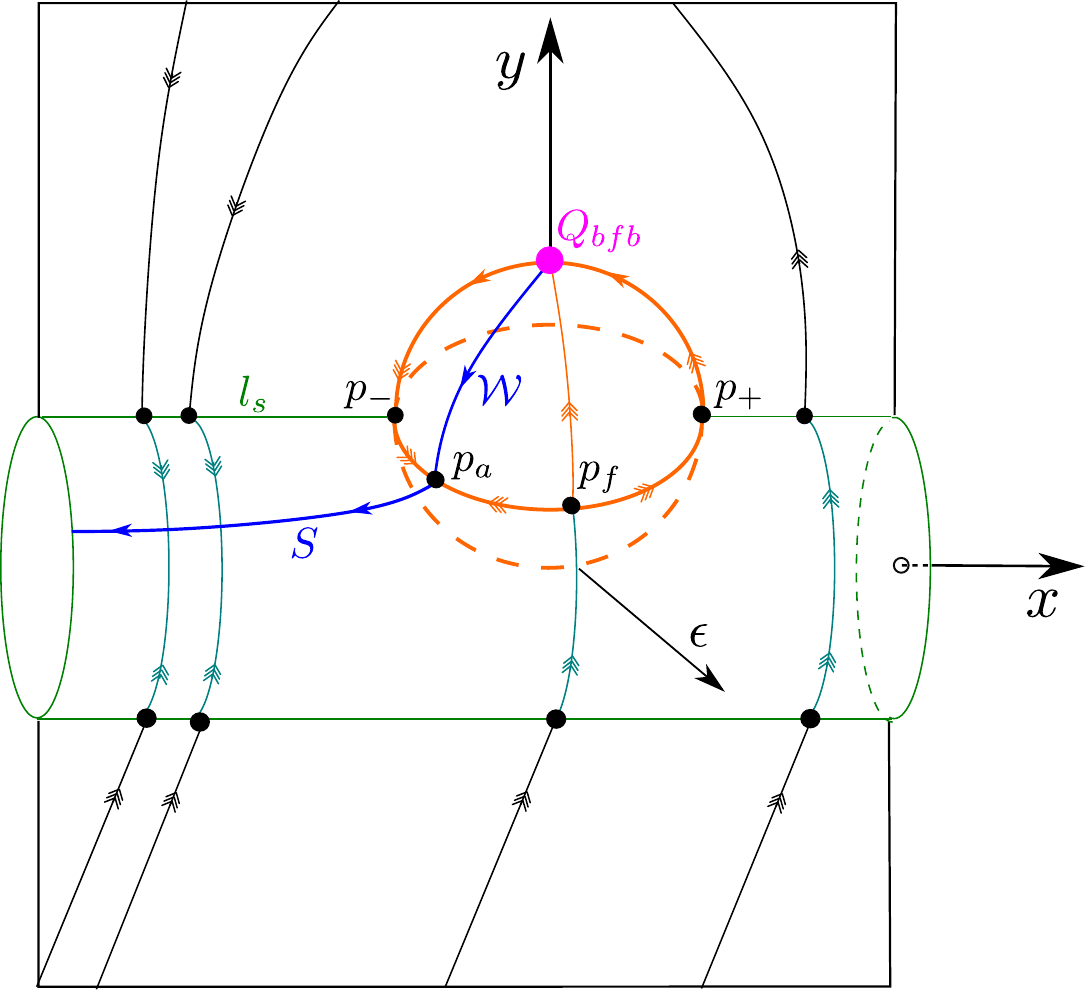}}
		\subfigure[]{\includegraphics[width=.40\textwidth]{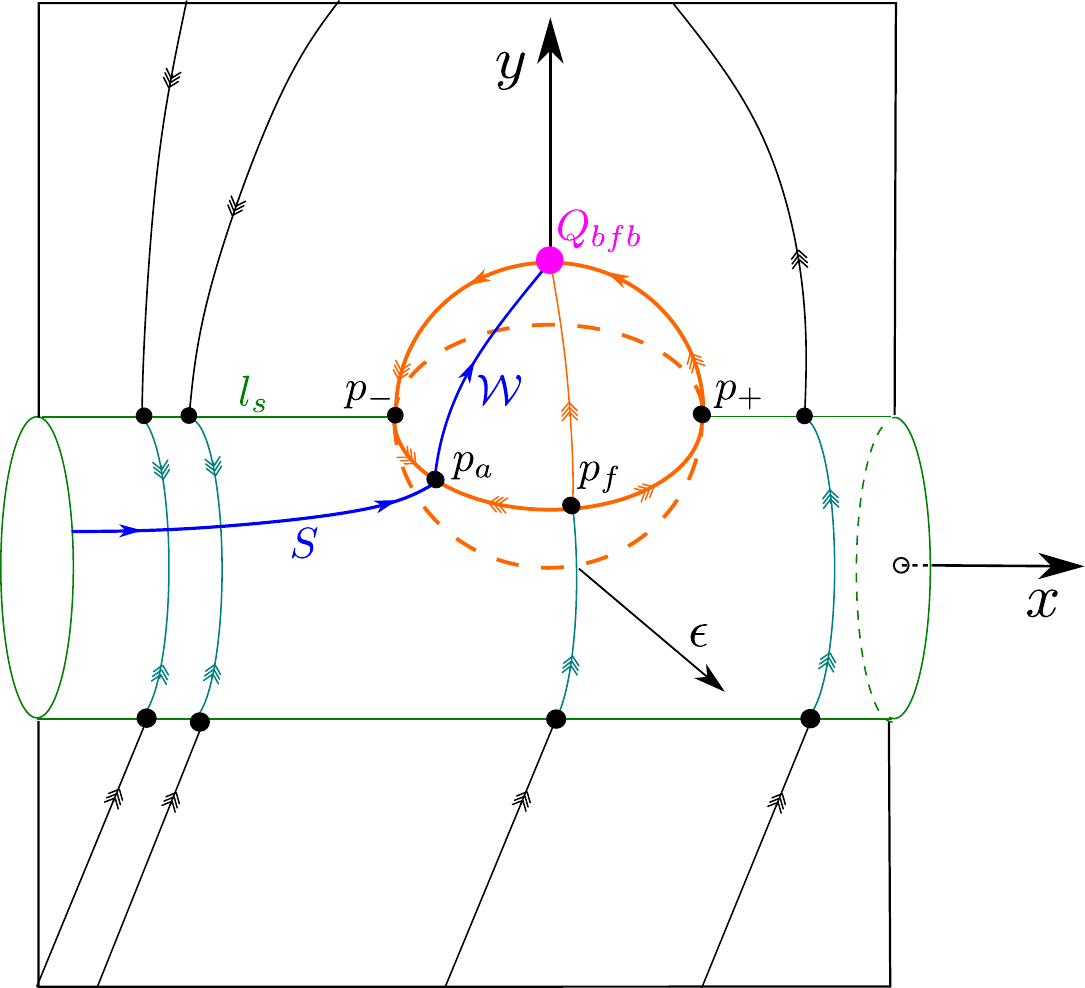}}
		\caption{Blow-up transformations necessary to resolve the switching manifold $\Sigma$ and the fold point $F$ in \eqref{normal_form}, shown with $\mu = 0$.
			The left (right) panel shows the sequence of blow-ups in cases BF$_i$, $i=1,2$ (BF$_3$) for which $\gamma > 0$, ($\gamma < 0$). Smoothness is regained by cylindrical blow-up of $\Sigma$. 
			On the blow-up cylinder we identify a slow-fast system with an attracting critical manifold $S$ which connects to a nonhyperbolic point $Q$ as sketched in (a), (b).
			Degeneracies at $Q$ stemming from the fold singularity $F$ 
			are resolved by a weighted spherical blow-up, leading to the situation in (c), (d). A second attracting critical manifold $\mathcal W$ is identified on the blow-up sphere and connects to a nonhyperbolic point $Q_{bfb}$, which
			requires further blow-up, see \figref{S1_blowups_ii}. \SJ{The} \SJ{orientation of the flow on both $S$ and $\mathcal W$ is determined by the sign of $\gamma$.}}
		\figlab{S1_blowups_i}
	\end{center}
\end{figure}

The degeneracy at $Q$ stems from the fold singularity at $(x,y,\epsilon)=(0,0,0)$ in system \eqref{normal_form}, which also exists for all $\mu$. As is known from, e.g.~\cite{Bonet2016,Kristiansen2019c}, this leads to a loss of hyperbolicity due to 
\begin{enumerate}
	\item[(i)] an alignment between $S$ and the fast fibration on the cylinder, and 
	\item[(ii)] a tangency with the `outer flow' in the upper half-plane $\{(x,y,0) : y \geq 0\}$. 
\end{enumerate}
Applying a successive (weighted) spherical blow-up\footnote{See \appref{blow-up_of_Q} for details; the blow-up transformation itself is given by \eqref{sph_bu}.} which replaces $Q$ with a 2-sphere
resolves the issue (i), with the critical manifold $S$ terminating at a partially hyperbolic and attracting point $p_a$ lying in the intersection between the blow-up sphere and cylinder. We also identify a unique centre manifold $\mathcal W$ emanating from $p_a$ and extending over the blow-up sphere, which limits to an attracting \textit{critical manifold} as $\mu \to 0^+$. Degeneracy due to the issue (ii) is also resolved for $\mu \neq 0$ (see \figref{S2_blowups}), however an additional degeneracy at a nonhyperbolic point $Q_{bfb}$ lying in the intersection between the blow-up sphere and the plane $\{\epsilon = 0\}$ arises due to the collision singularity when $\mu = 0$; this is shown in \figref{S1_blowups_i}c and \figref{S1_blowups_i}d. The manifold $\mathcal W$ connects to $Q_{bfb}$, and the orientation of the reduced flow depends on the sign of $\gamma$ (toward $Q_{bfb}$ for $\gamma < 0$, away from $Q_{bfb}$ for $\gamma > 0$).

\begin{remark}
	The colour-code adopted in \figref{S1_blowups_i} is designed to keep track of blown-up objects and their corresponding counterparts in the original $(x,y,\epsilon)$-coordinates. We represent the blow-up cylinder in green to indicate that it `replaces' the switching manifold $\Sigma$, represented in green in earlier figures. Similarly, the blow-up sphere is shown in orange since it resolves degeneracies at $Q$ stemming from the fold singularity at $F$ (also orange in earlier figures). A similar convention will be adopted in later figures with regard to further blow-up of $Q_{bfb}$, which is shown in magenta.
\end{remark}

\subsection{The regularised BF collision: bifurcations in \SJnew{regime (S1) when $\mu = \mathcal O(\epsilon^{k/(k+1)})$}}
\seclab{results_bifurcation}

The degeneracy at $Q_{bfb}$ is \textit{not} independent of $\mu$, and occurs only for $\mu = 0$. Hence, one must blow-up in the doubly-extended $(x,y,\epsilon,\mu)-$space.
We identify another (weighted) spherical blow-up which replaces $Q_{bfb}$ with a 3-sphere, regaining at least partial hyperbolicity everywhere in the resulting blown-up space (except at bifurcations).\footnote{See \appref{blow-up_of_Qbfb} for details; the blow-up transformation itself is given by \eqref{bu_Qbfb}.} 
\figref{S1_blowups_ii} shows the resulting dynamics in cases $\gamma < 0$ and $\gamma > 0$, restricting to the (invariant) scaling regime (S1) (which allows for a `true' 3-dimensional depiction). In particular, the critical manifold $\mathcal W$ connects to a partially hyperbolic and attracting point $q_a$ lying within the intersection of the first blow-up 2-sphere with the blow-up 3-sphere. A unique centre manifold \SJnew{$\mathcal J'$} emanates from $q_a$, and the orientation of the flow on \SJnew{$\mathcal J'$} is determined by the sign of $\gamma$.
The dynamics on the restricted 3-sphere itself (missing in \figref{S1_blowups_ii}) varies depending on the region of $(\hat \mu, \gamma) -$parameter space.
This is described in the following two lemmas.

\begin{figure}[t!]
	\centering
	\subfigure[]{\includegraphics[width=.49\textwidth]{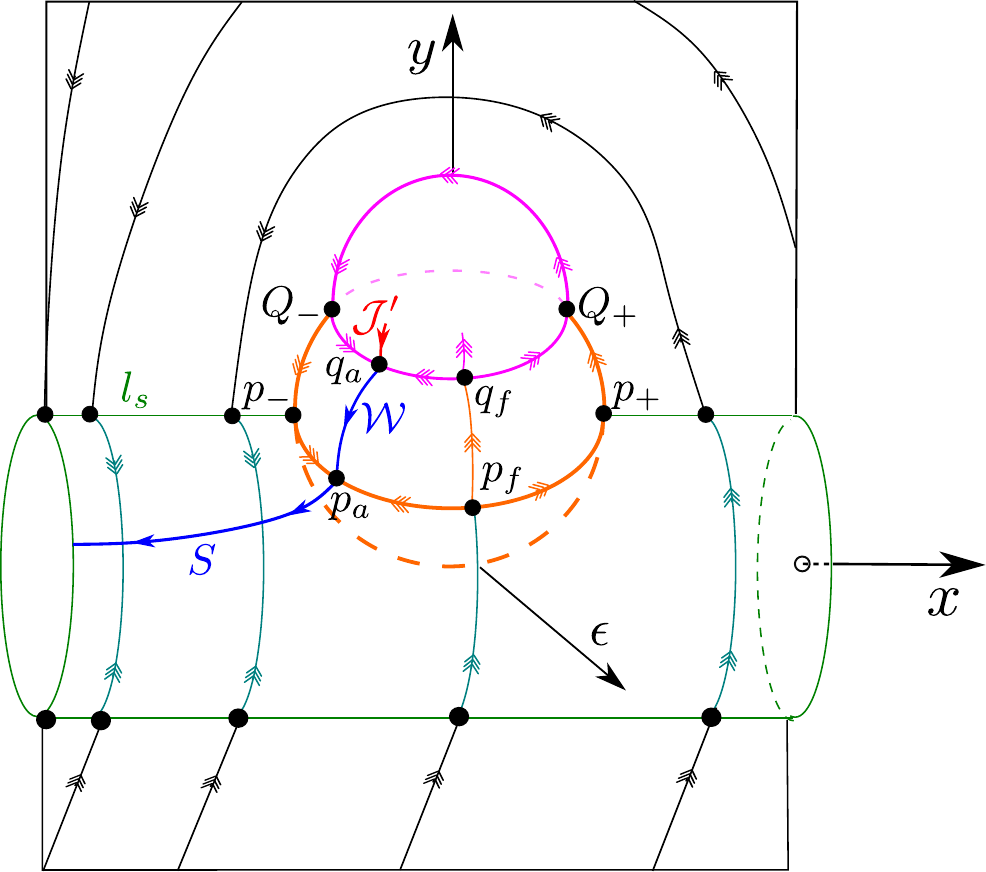}}
	\subfigure[]{\includegraphics[width=.49\textwidth]{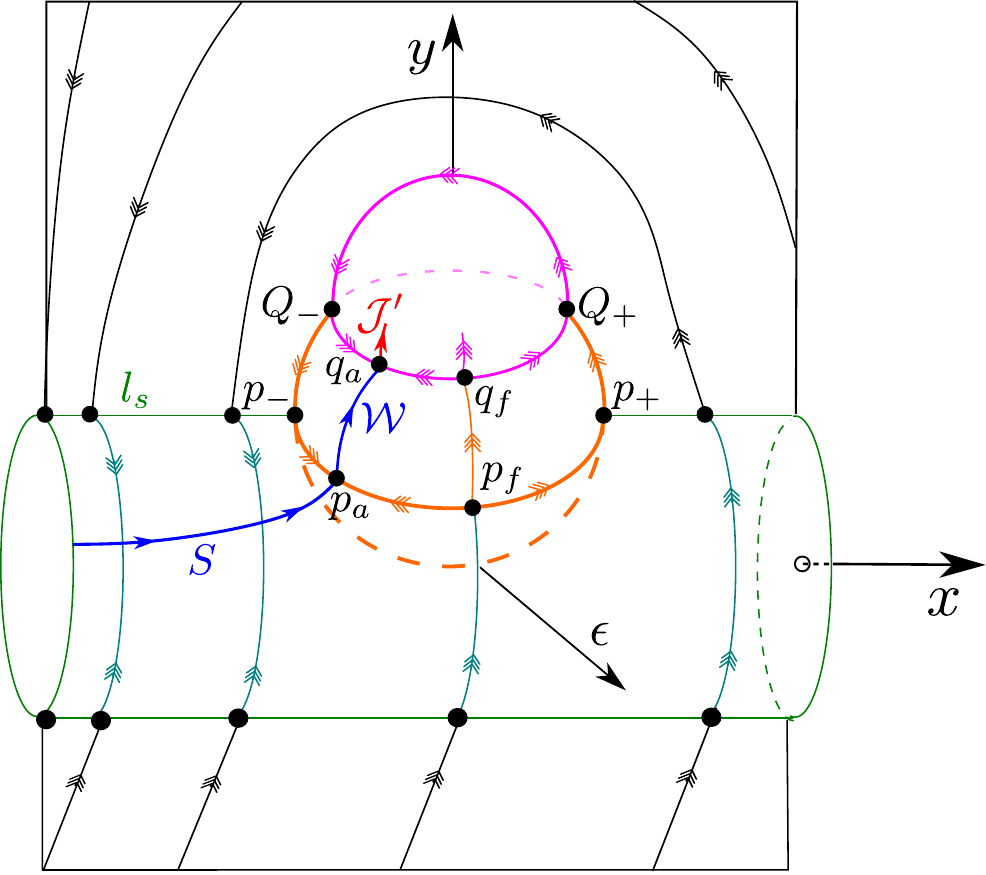}}
	\caption{Key dynamics after blow-up of $Q_{bfb}$ and restriction to regime (S1) in \figref{S1_blowups_i}c and \figref{S1_blowups_i}d. In (a): case $\gamma > 0$ (BF$_i$, $1=1,2$), and in (b): $\gamma < 0$ (BF$_3$). Importantly, partial hyperbolicity has been regained everywhere except at bifurcations. In both (a) and (b), the attracting critical manifold $\mathcal W$ (blue) connects to an attracting partially hyperbolic point $q_a$ at the equator of the restricted blow-up 3-sphere (magenta), and a unique centre manifold \SJnew{$\mathcal J'$} emanates from $q_a$, extending onto the (restricted) 3-sphere. The orientation of the flow on \SJnew{$\mathcal J'$} is determined by the sign of $\gamma$. Dynamics on the restricted 3-sphere itself depend on the region of $(\hat \mu, \gamma) -$parameter space as described in \lemmaref{lem_bifurcations_1}; see also \figref{bfb_bif_diagram}, \figref{bfb_bif_neg} and \figref{bfb_bif_pos}.}
	\figlab{S1_blowups_ii}
\end{figure}

\begin{lemma}
	\lemmalab{lem_desing}
	Applying the parameter-dependent coordinate transformation 
	\begin{equation}
	\eqlab{coord_change}
	\left(x_1, \nu_1, \rho_1, \hat \mu \right) \in \mathbb R \times \mathbb R_+^2 \times \mathbb R \mapsto
	\begin{cases}
	x = \nu_1^{2k(1+k)} \rho_1^{k(1+k)} x_1 , \\
	y = \nu_1^{2k(1+k)} \rho_1^{2k(1+k)} , \\
	\epsilon = \nu_1^{2(1+k)^2} \rho_1^{(1+k)(1+2k)} , \\
	\mu = \hat \mu \nu_1^{2k(1+k)} \rho_1^{k(1+2k)} ,
	\end{cases}
	\end{equation}
	and a desingularisation
	\begin{equation}
	\eqlab{time_change}
	d\tilde t = \rho_1^{k(1+k)} dt 
	\end{equation}
	\W{to the extended system $\{\eqref{normal_form},\varepsilon'=0,\mu'=0\}$},
	one obtains a system of equations for which the subspace $\{\nu_1 = 0\}$ is invariant. The dynamics within $\{\nu_1 = 0\}$ are governed by the planar system
	\begin{equation}
	\eqlab{desing_prob}
	\begin{split}
	x_1' &= \rho_1^{k(1+k)} \left((\tau - \gamma) \beta + \hat \mu \rho_1^{k^2} + \tau x_1 - \delta \rho_1^{k(1+k)} \right) + k x_1 \left(\beta + x_1 \right) , \\
	\rho_1' &= \frac{1}{k} \rho_1 \left(\beta + x_1 \right) ,
	\end{split}
	\end{equation}
	where $\beta > 0$ is given by \eqref{beta} and by a slight abuse of notation $(\cdot)'$ denotes differentiation with respect to the new time $\tilde t$ in \eqref{desing_prob}.
\end{lemma}

\bpr
	The result follows after a direct application of the coordinate and time transformations in \eqref{coord_change} and \eqref{time_change} respectively.
\epr

\begin{remark}
	Determining the coordinate transformation \eqref{coord_change} in \lemmaref{lem_desing} is a nontrivial task involving the determination and composition of the appropriate blow-up transformations. This analysis is
	carried out in the appendix: to obtain the transformation \eqref{coord_change}, we have composed the blow-up transformations \eqref{K1_coord}, \eqref{mathcal_K2_coord}, and \eqref{mathfrak_K1_coord}, and restricted to the set
	\begin{equation}
	\eqlab{mathcal_A_1}
	\mathcal A = \left\{\left(x,y,\epsilon,\hat \mu \epsilon^{k/(k+1)} \right) : \hat \mu \in \mathbb R 
	\right\} ,
	\end{equation}
	which is invariant since both $\epsilon$ and $\mu$ are constants of the motion in system $\{\eqref{normal_form},\varepsilon'=0,\mu'=0\}$; see also \lemmaref{par_sets}.
\end{remark}

\begin{remark}
	\remlab{coord_change}
	\lemmaref{lem_desing} provides an explicit means for obtaining quantitative information in applications, without the need to transform the system into the local normal form \eqref{normal_form}. By applying the coordinate transformation \eqref{coord_change} to a general system \eqref{main} satisfying \assumptionref{ass1}, \assumptionref{ass1b} and \assumptionref{ass2} with a BF bifurcation at $(x,y,\epsilon,\alpha)=(0,0,0,0)$, it is possible to obtain a desingularised system analogous to \eqref{desing_prob} which governs the dynamics in (S1). Precise quantitative information about, e.g. bifurcations, can be obtained for the given application from this desingularised system. We discuss this approach further in the context of a predator-prey model in \secref{gause}, and apply it directly for a mechanical oscillator model in \secref{stick-slip}. 
\end{remark}

\begin{figure}[t!]
	\begin{center}
		\includegraphics[width=.85\textwidth]{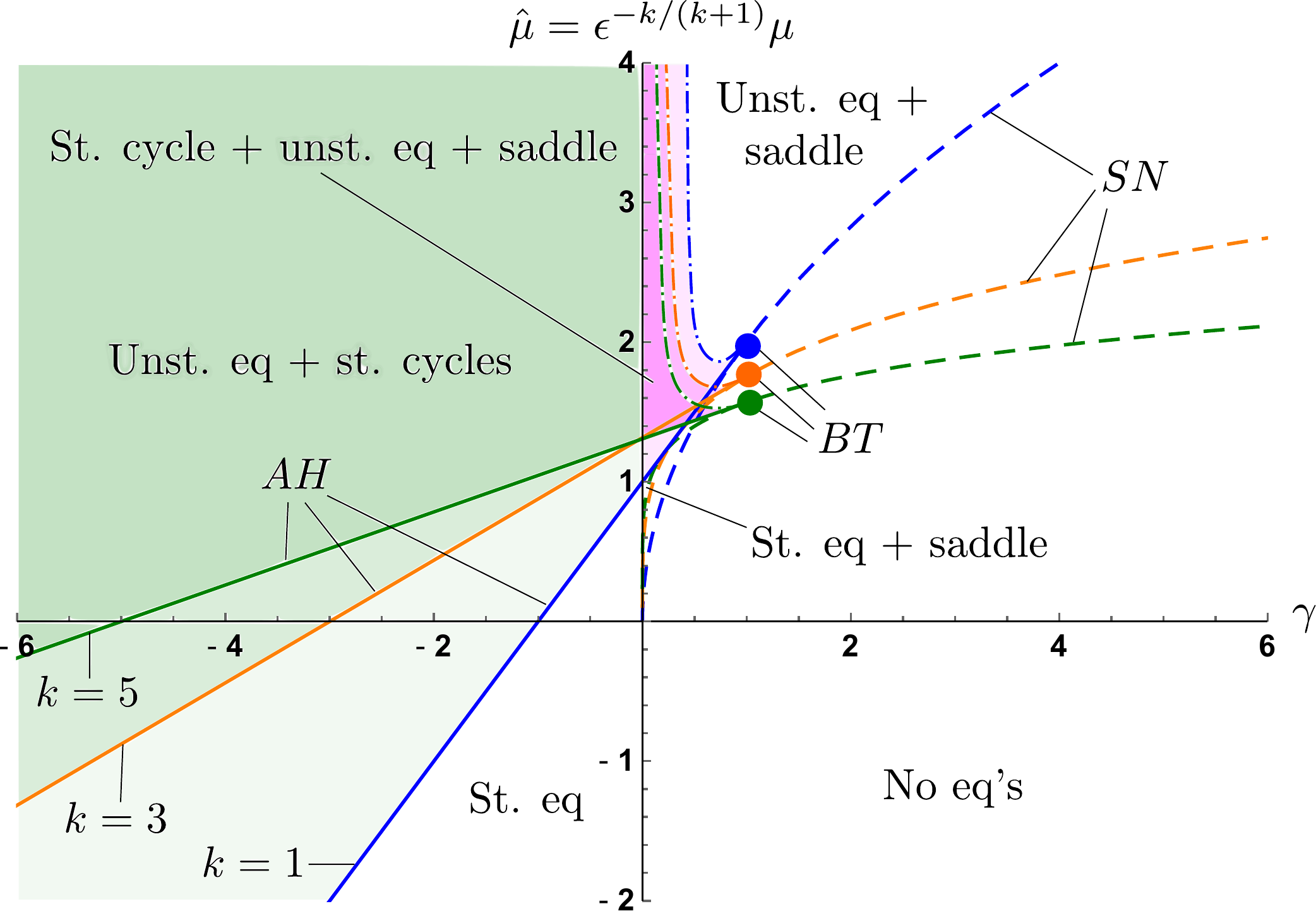}
		\caption{Bifurcation diagram for system \eqref{desing_prob} 
			with 
			$k=1,3,5$ (blue, yellow, green resp.) and parameter values $\beta = \tau = \delta = 1$. The bifurcation diagram for case BF$_3$ is contained within the open half-plane $\gamma < 0$, while BF$_i$, $i=1,2$ cases are contained within the open half-plane $\gamma > 0$ (the degenerate case $\gamma = 0$ is omitted, but discussed in \secref{degenerate_cases}). Supercritical Andronov-Hopf, saddle-node and homoclinic curves (solid, dashed and solid-dashed resp.)~ emanate from a Bogdanov-Takens point. 
			In case BF$_3$, stable oscillations exist within the shaded green region, assuming no terminating global bifurcations.
			Stable oscillations also exist when $\gamma > 0$, in the shaded magenta regions bounded between the Andronov-Hopf and homoclinic curves.
			We emphasise that the indicated homoclinic curves are sketches only, with asymptotic knowledge only in the neighbourhood of Bogdanov-Takens points and for the limit $\gamma \to 0^+$, as described in \thmref{thm_bifurcation} (see also \secref{degenerate_cases}).}
		\figlab{bfb_bif_diagram}
	\end{center}
\end{figure}

System \eqref{desing_prob} governs the dynamics on the restricted 3-sphere in the \W{final} blow-up, i.e.~the `missing dynamics' in \figref{S1_blowups_ii}. Equivalently, one can consider \eqref{desing_prob} as the `singular limit problem' associated with system \eqref{normal_form} in regime (S1). 

\begin{lemma}
	\lemmalab{lem_bifurcations_1}
	System \eqref{desing_prob} has zero, one or two equilibria in the half-plane $\{\rho_1 > 0\}$, depending on the region in $(\hat \mu, \gamma)-$parameter space. 
	\begin{enumerate}
		\item[(i)] 
		There exists an equilibrium $p$ which undergoes a supercritical Andronov-Hopf bifurcation along the parameter-space curve
		\begin{equation}
		\eqlab{hopf_value_1}
		\hat \mu_{ah}(\gamma) = \frac{k \delta + \tau \gamma}{k} \left(\frac {k\beta} \tau \right)^{1/(k+1)} , \qquad  \gamma \in \left(- \infty, \frac{\delta}{\tau} \right) ,
		\end{equation}
		which is stable for $\hat \mu < \hat \mu_{ah}(\gamma)$ and unstable for $\hat \mu > \hat \mu_{ah}(\gamma)$, for all three cases BF$_i$, $i = 1,2,3$. For $\gamma < 0$, i.e. for case BF$_3$, there are no other local bifurcations.
		\item[(ii)] For $\gamma > 0$, i.e. cases BF$_i$, $i=1,2$, saddle-node bifurcations occur along the parameter-space curve
		\begin{equation}
		\eqlab{sn_curve_lem_1}
		\hat \mu_{sn}(\gamma) = \frac{(1+k) \delta} k \left(\frac{k \beta \gamma}{\delta} \right)^{1/(k+1)}, \qquad \gamma > 0.
		\end{equation}
		The focus-node and saddle-type equilibria $p$ and $p_s$ born in the saddle-node bifurcation exist for $\hat \mu > \hat \mu_{sn}(\gamma)$ (and not for $\hat \mu < \hat \mu_{sn}(\gamma)$).
		\item[(iii)] For $\gamma > 0$, i.e. cases BF$_i$, $i=1,2$, there exists a Bogdanov-Takens bifurcation at the intersection of saddle-node and Andronov-Hopf curves described in (i) and (ii) above, for
		\begin{equation}
		\eqlab{bt_point_lem_1}
		(\hat \mu_{bt}, \gamma_{bt} ) = \left(\frac{(1+k) \delta}{k} \left(\frac{k \beta}{\tau} \right)^{1/(1+k)}, \frac \delta \tau \right) .
		\end{equation}
		\item[(iv)] There exists a parameter-space curve $\hat \mu_{hom}(\gamma)$ with a quadratic tangency to $\hat \mu_{ah}(\gamma)$ at $(\hat \mu_{bt}, \gamma_{bt} )$ along which the system has a saddle homoclinic connection. The curve $\hat \mu_{hom}(\gamma)$ is not defined on $\gamma < 0$, and its domain is bounded above (locally) by $\gamma = \delta / \tau$.
	\end{enumerate}
\end{lemma}

\bpr
	See \appref{proof_of_lem_bifurcation}.
\epr

\begin{figure}[t!]
	\centering
	\subfigure[]{\includegraphics[width=.4\textwidth]{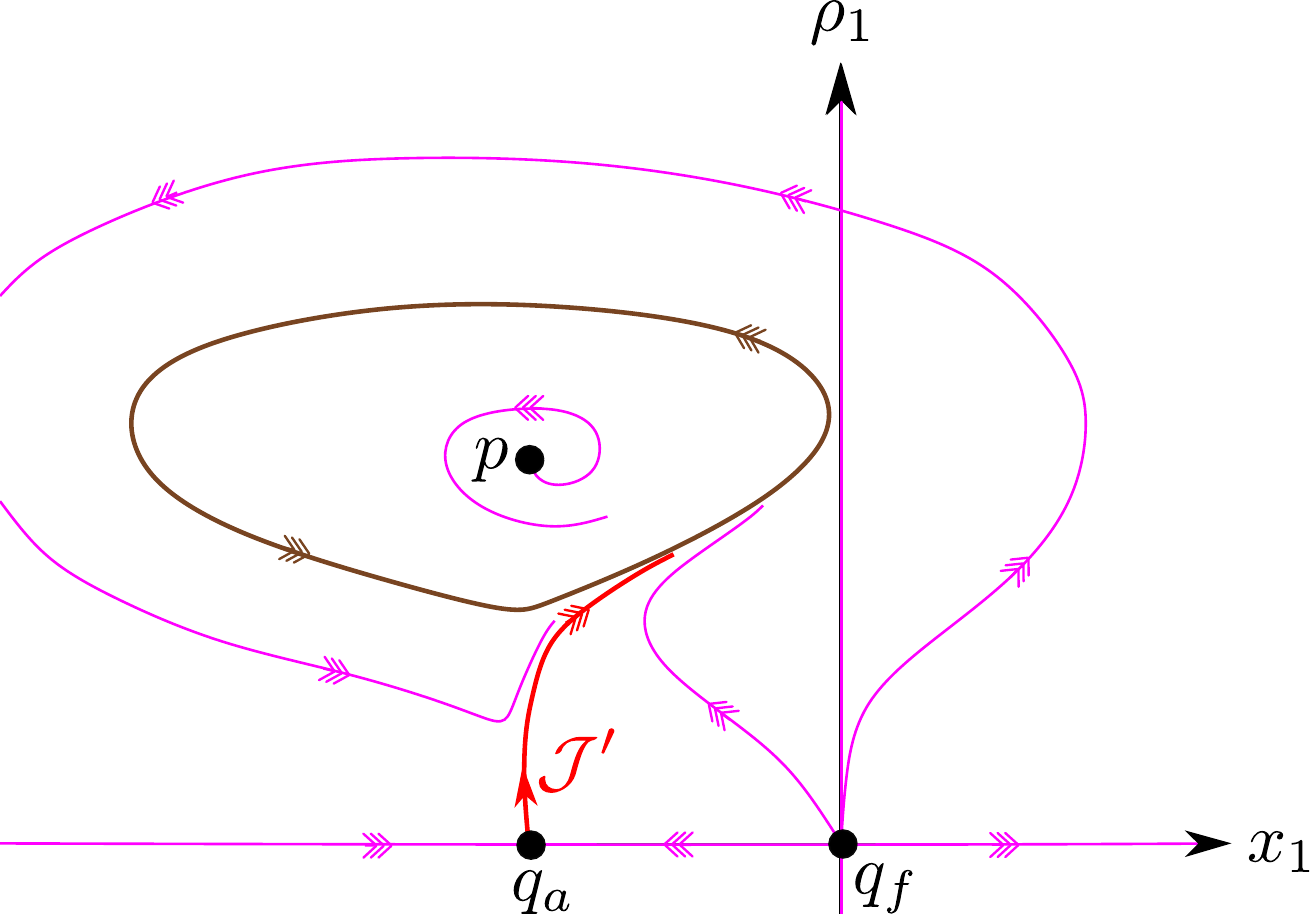}}
	\hspace{0.5cm}
	\subfigure[]{\includegraphics[width=.4\textwidth]{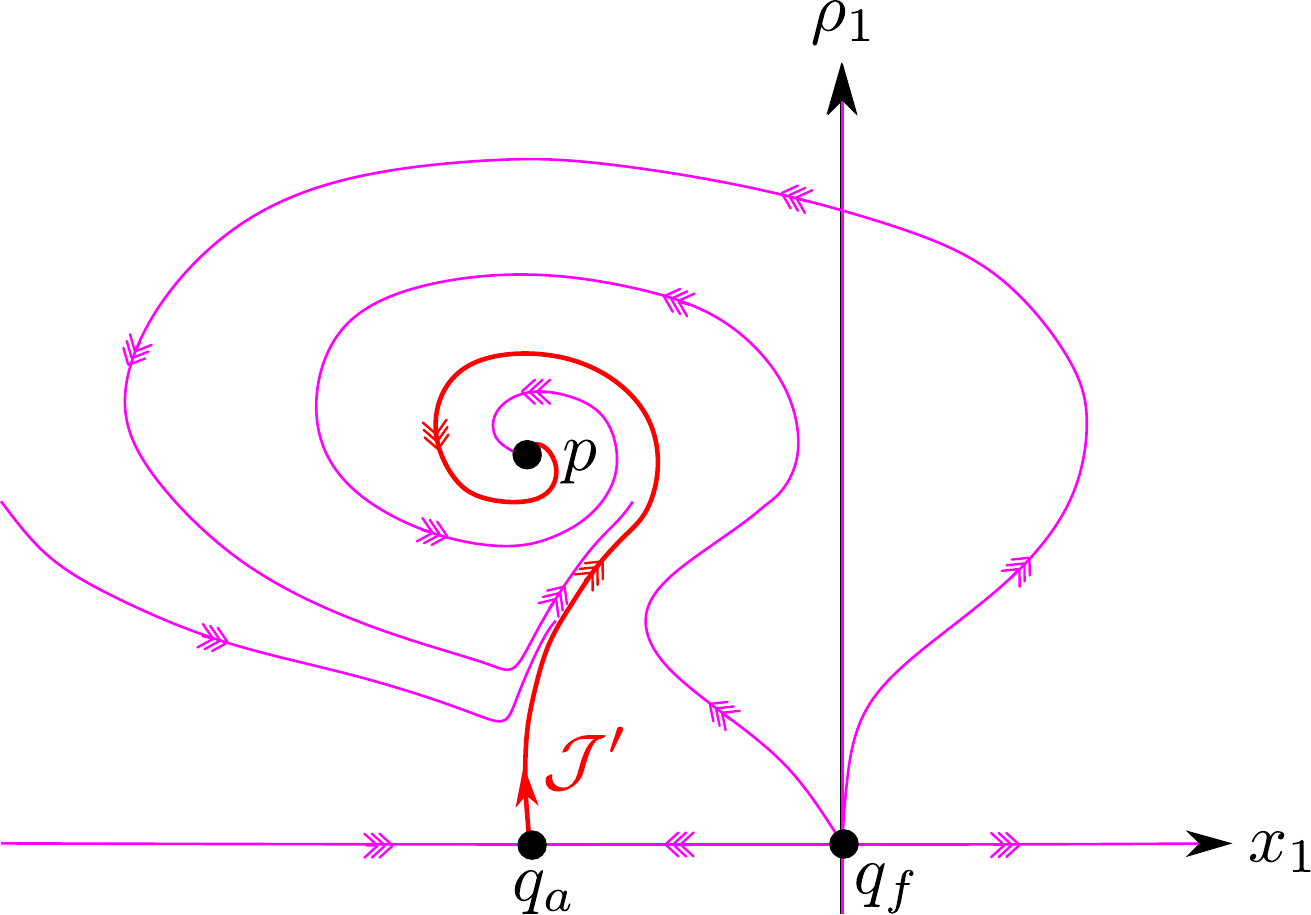}}
	\caption{Phase portraits for system \eqref{desing_prob} 
		with $\gamma < 0$; c.f. \SJ{\figref{S1_blowups_ii} and} \figref{bfb_bif_diagram}. In (a): $\hat \mu > \hat \mu_{ah}(\gamma)$, in which case $p$ is an unstable focus, and the unique centre manifold \SJnew{$\mathcal J'$} (red) emanating from the partially hyperbolic point $q_a$ is forward asymptotic to a stable limit cycle (brown). In (b): $\hat \mu < \hat \mu_{ah}(\gamma)$. Following termination of oscillations in (a) in an Andronov-Hopf bifurcation under $\hat \mu$ variation, $p$ becomes a stable focus.
		}
	\figlab{bfb_bif_neg}
\end{figure}

\begin{figure}[t!]
	\begin{center}
		\includegraphics[width=.9\textwidth]{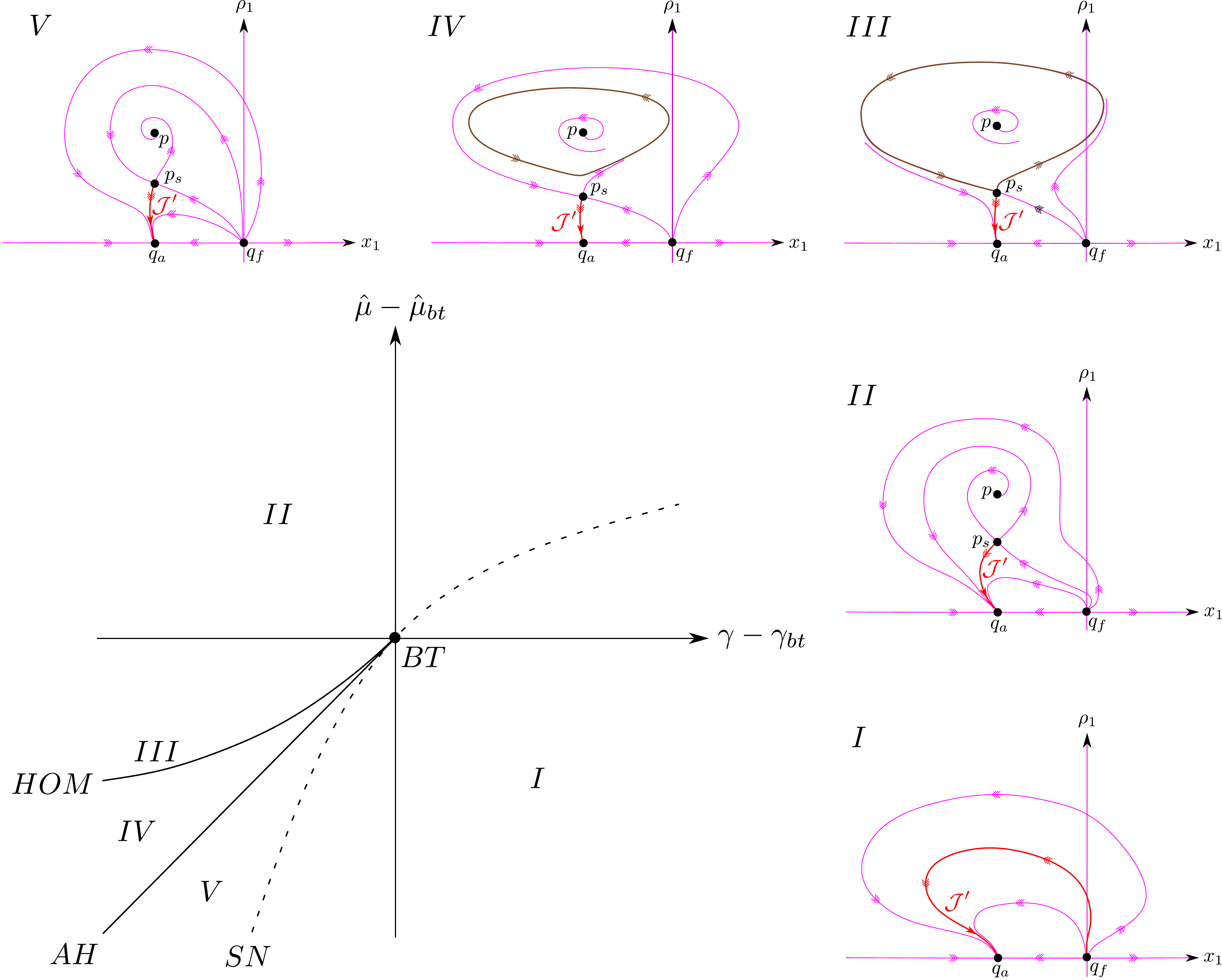}
		\caption{Local bifurcation diagram for system \eqref{desing_prob} 
			in case $\gamma > 0$, near the Bogdanov-Takens value $(\hat \mu_{bt},\gamma_{bt})$; c.f. \SJ{\figref{S1_blowups_ii} and} \figref{bfb_bif_diagram}.}
		\figlab{bfb_bif_pos}
	\end{center}
\end{figure}

A two-parameter bifurcation diagram is shown in \figref{bfb_bif_diagram}, and phase portraits for each case $\gamma < 0$ and $\gamma > 0$ are sketched in \figref{bfb_bif_neg} and \figref{bfb_bif_pos} respectively.  We finally state our main result characterising the dynamics in (S1).





\begin{theorem}
	\thmlab{thm_bifurcation}
	Consider system \eqref{normal_form} under
	\assumptionref{ass3} and \assumptionref{ass2}. Then there exists $\epsilon_0, \mu_0 > 0$ such that for all $\epsilon \in (0,\epsilon_0)$ and $\mu \in (-\mu_0,\mu_0)$, \SJ{the following assertions hold:
		\begin{enumerate}
			\item[(i)] There is a supercritical Andronov-Hopf bifurcation along the parameter-space curve
			\begin{equation}
			\eqlab{hopf_curve}
			\mu_{ah}(\gamma,\epsilon) = \frac{k \delta + \tau \gamma}{k} \left(\frac {k\beta} \tau \right)^{1/(k+1)} \epsilon^{k/(k+1)} + o\left(\epsilon^{k/(k+1)}\right) , \qquad  \gamma \in \left(- \infty, \frac{\delta}{\tau} \right) ,
			\end{equation}
			where the bifurcating equilibrium $p_\epsilon$ is stable for $\mu < \mu_{ah}(\gamma,\epsilon)$, and unstable for $\mu > \mu_{ah}(\gamma,\epsilon)$.
			\item[(ii)] There is a saddle-node bifurcation along the parameter-space curve
			\begin{equation}
			\eqlab{sn_curve}
			\mu_{sn}(\gamma,\epsilon) = \frac{(1+k) \delta} k \left(\frac{k \beta \gamma}{\delta} \right)^{1/(k+1)} \epsilon^{k/(k+1)} + o\left(\epsilon^{k/(k+1)}\right) , \qquad \gamma > 0.
			\end{equation}
			\item[(iii)] Saddle-node and Andronov-Hopf curves intersect in a codimension-2 Bogdanov-Takens point
			\begin{equation}
			\eqlab{bt_point}
			(\mu_{bt}(\epsilon), \gamma_{bt} ) = \left(\frac{(1+k) \delta}{k} \left(\frac{k \beta}{\tau} \right)^{1/(1+k)} \epsilon^{k/(k+1)} + o\left(\epsilon^{k/(k+1)}\right) , \frac \delta \tau \right) .
			\end{equation}
			\item[(iv)] There exists a parameter-space curve $\mu_{hom}(\gamma,\epsilon)$ which is tangent to $\mu_{ah}(\gamma,\epsilon)$ at $(\mu_{bt}, \gamma_{bt} )$, along which the system has a homoclinic-to-saddle connection. The curve $\mu_{hom}(\gamma)$ is not defined on $\gamma < 0$, and its domain is bounded above (locally) by $\gamma = \delta / \tau$.
		\end{enumerate}
		The bifurcations identified above divide $(\epsilon^{-k/(k+1)}\mu,\gamma)-$parameter space into the following regions, depending on the sign of $\gamma$:
		\begin{itemize}
			\item $\gamma < 0$: there exists a single equilibrium $p_\epsilon$, as described by the Andronov-Hopf statement above.
			\item $\gamma > 0$: singularities $p_\epsilon$ and $p_{s,\epsilon}$ of focus-node and saddle-type exists in the region bounded above the saddle-node curve $\mu = \mu_{sn}(\gamma,\epsilon)$.
			There are no equilibria for $\mu < \mu_{sn}(\gamma,\epsilon)$.
	\end{itemize}}
\end{theorem}


\bpr
	Applying the transformation
	\begin{equation}
	\eqlab{par_blow_down_1}
	\mu = \hat \mu \epsilon^{k/(k+1)} 
	\end{equation}
	to \eqref{hopf_value_1}, \eqref{sn_curve_lem_1} and \eqref{bt_point_lem_1} in \lemmaref{lem_bifurcations_1} yields the leading order estimates \eqref{hopf_curve}, \eqref{sn_curve} and \eqref{bt_point} in \thmref{thm_bifurcation} respectively. Existence and the required properties of the homoclinic curve $\mu_{hom}(\gamma,\epsilon)$ follow immediately from \lemmaref{lem_bifurcations_1} (iv) together with another application of the transformation \eqref{par_blow_down_1}, and the fact that no more than two equilibria of system \eqref{desing_prob} coexist on $\{\rho_1 > 0\}$.
\epr

Since $\mu = \mathcal O(\epsilon^{k/{(k+1)}})$ in \thmref{thm_bifurcation} (see the vertical axis in \figref{bfb_bif_diagram}), all observed bifurcations are `singular' in the sense that they occur within an $\epsilon-$dependent neighbourhood which shrinks to zero in the singular limit $\epsilon \to 0^+$. In particular, all local bifurcations have associated non-zero eigenvalues satisfying $\lambda(\epsilon) \to 0$ as $\epsilon \to 0^+$, with a rate of convergence depending on the transition coefficient $k$. This is of particular interest for the Andronov-Hopf bifurcations: the fact that the eigenvalues $\lambda = A(\epsilon) \pm i B(\epsilon) \to 0$ as $\epsilon \to 0^+$ indicates a frequency/amplitude of the resultant oscillations which depends not only on $\epsilon$, but also on $k$. The reader is referred to the references \cite{Dum1996,Krupa2001b,Kuehn2015} and \cite{Maesschalck2011b,Wechselberger2015} for detailed descriptions of singular Andronov-Hopf and Bogdanov-Takens bifurcations in \SJnew{classical} slow-fast systems respectively. \SJnew{We emphasise that our results in \thmref{thm_bifurcation} are qualitatively independent of the choice of regularisation function $\phi$, insofar as the Andronov-Hopf and Bogdanov-Takens bifurcations 
are supercritical \textit{for all} $\phi$ satisfying \assumptionref{ass1b} and \assumptionref{ass3}, and \textit{for all} values of $k \in \mathbb N_+$.}


\subsection{Relaxation oscillation in \SJnew{regime (S2) when $\mu = \mathcal O(1)$}}
\seclab{results_ro}

\begin{figure}[t!]
	\begin{center}
		\subfigure[]{\includegraphics[width=.43\textwidth]{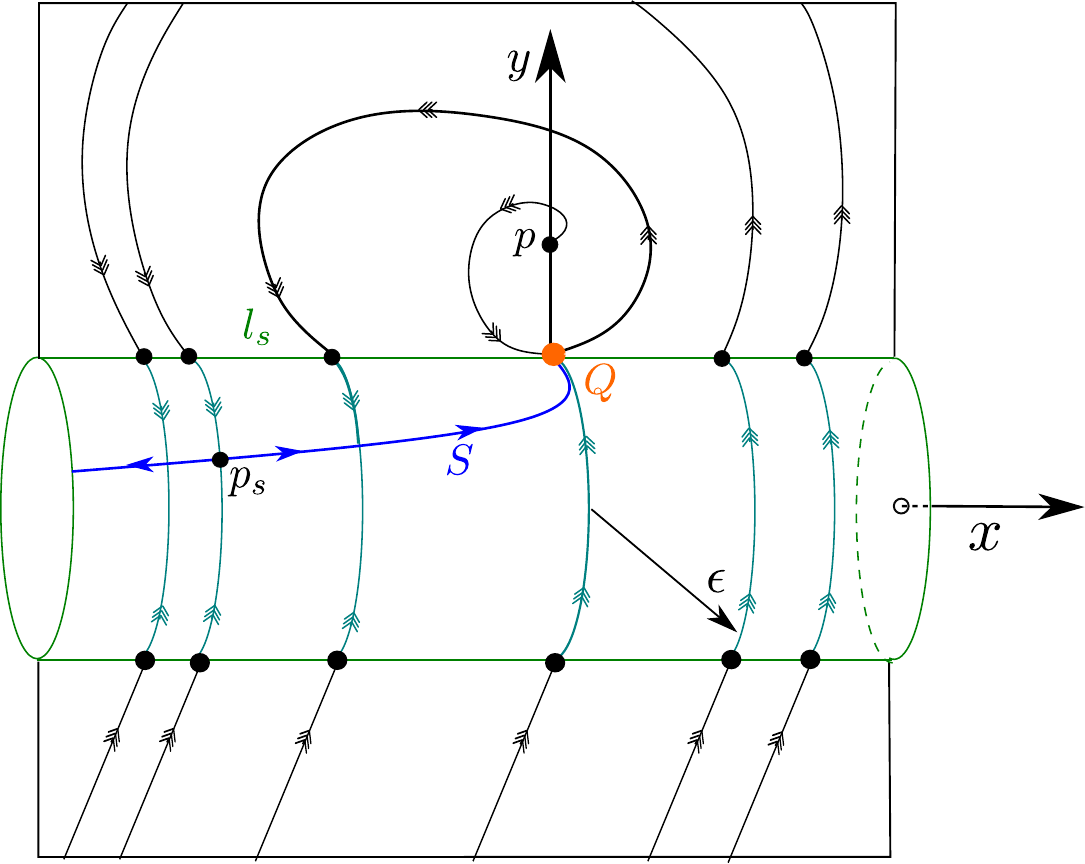}}
		\subfigure[]{\includegraphics[width=.43\textwidth]{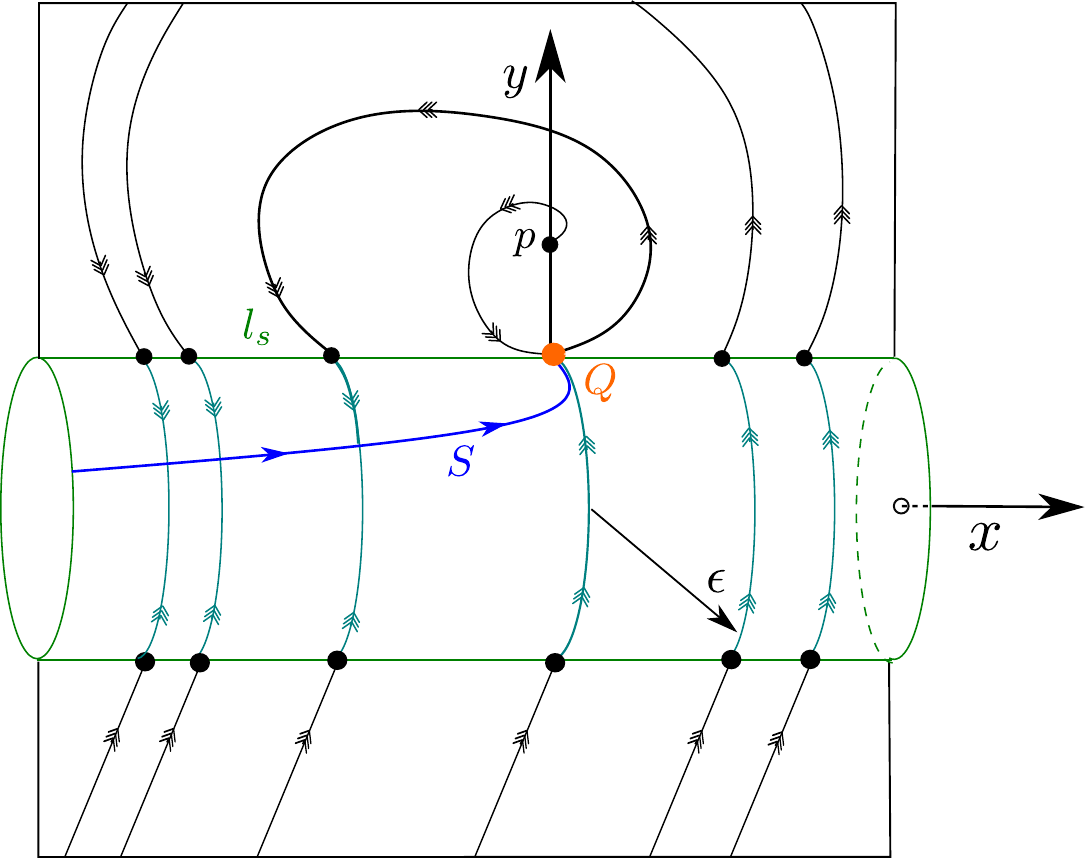}}
		\subfigure[]{\includegraphics[width=.43\textwidth]{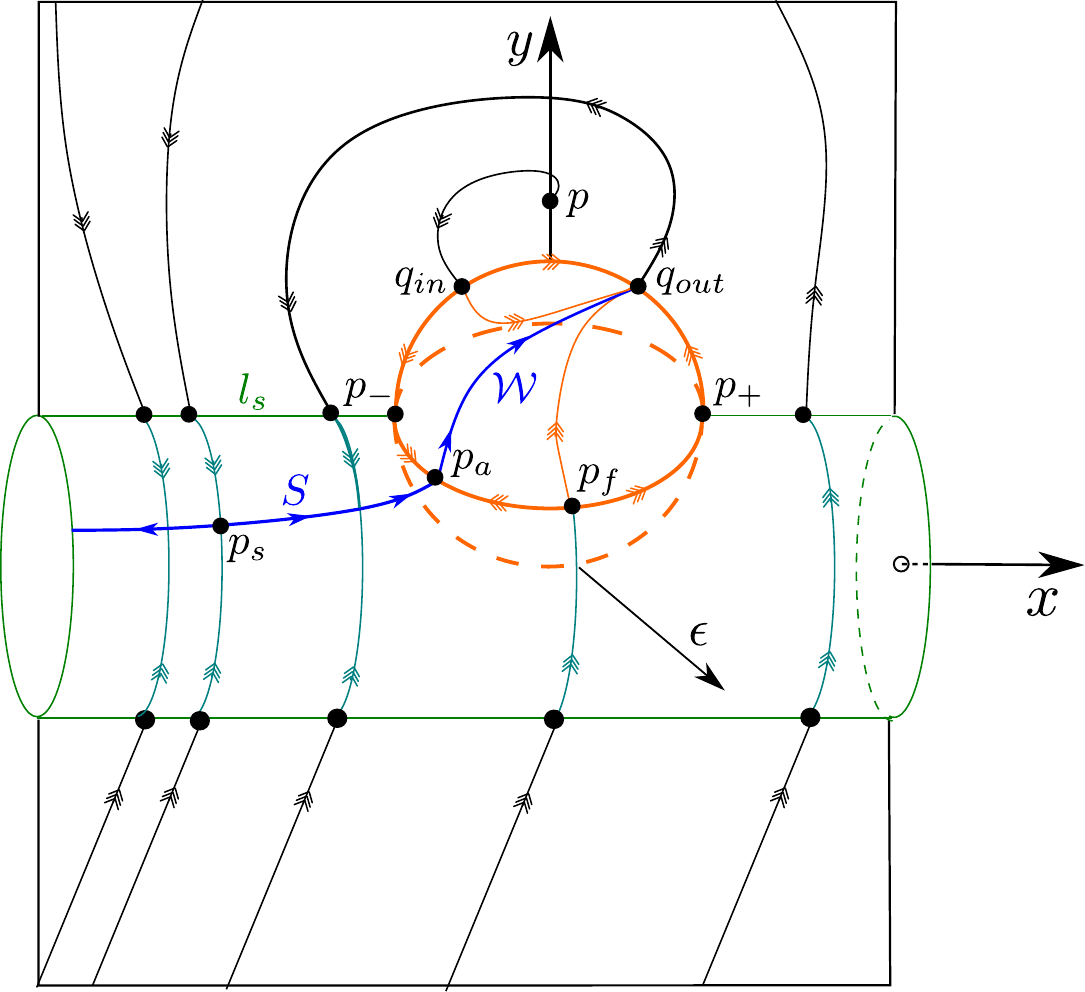}}
		\subfigure[]{\includegraphics[width=.43\textwidth]{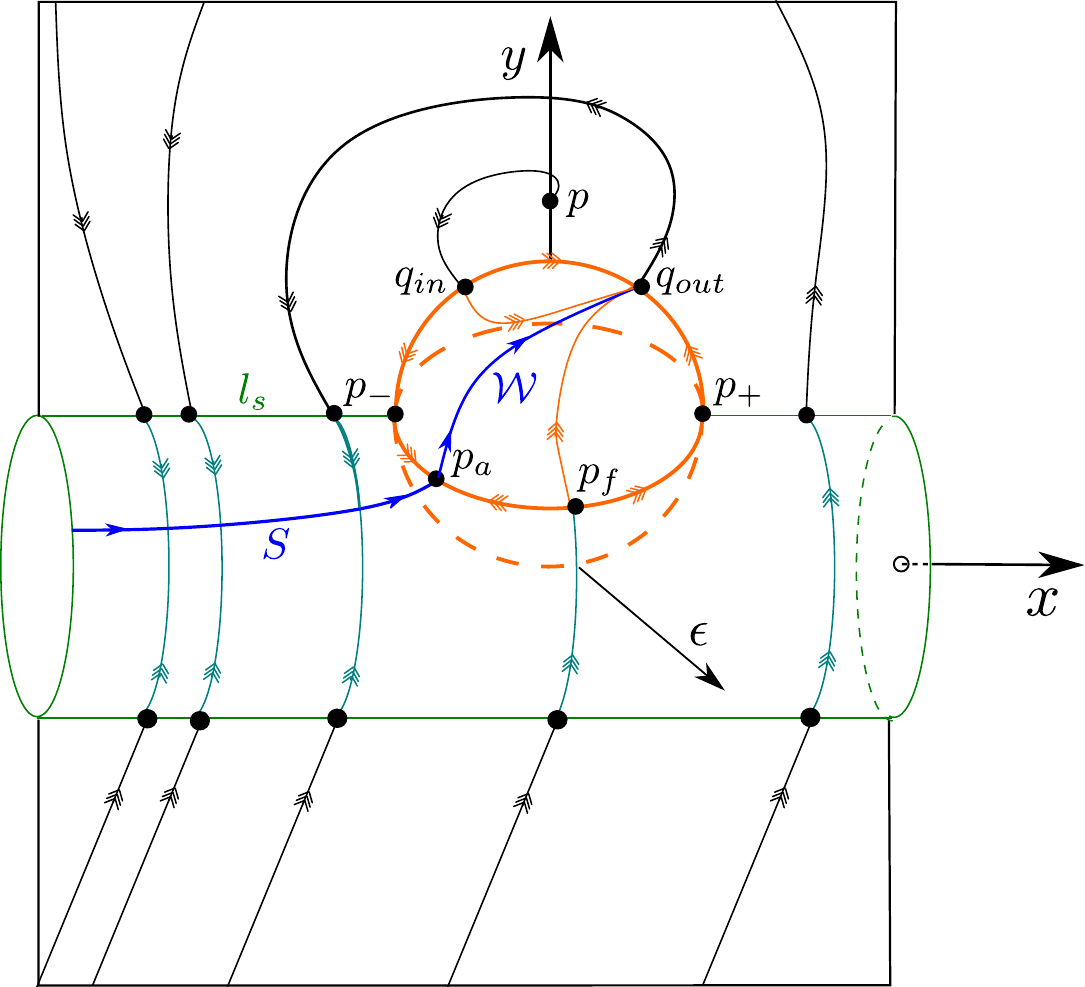}}
		\caption{Blow-up transformations necessary to desingularise system \eqref{normal_form} with fixed $\mu \in \mathcal I_+$. The left (right) panel shows the sequence of blow-ups in case BF$_1$ (BF$_3$). In each case, the loss of smoothness along $\Sigma$ and degeneracy at $Q$ are resolved by the same cylindrical and spherical blow-ups as in the $\mu = 0$ case in \figref{S1_blowups_i}. 
			\SJ{Following spherical blow-up, a unique center manifold $\mathcal W$ (which converges to the critical manifold $\mathcal W$ in \figref{S1_blowups_i}c-d as $\mu \to 0^+$) connects the point $p_a$ to a hyperbolic saddle $q_{out}$ contained within the intersection of the blow-up sphere with $\{\epsilon = 0\}$. Partial hyperbolicity is regained everywhere} \SJ{in (c)-(d), and closed singular cycles can be constructed.}}
		\figlab{S2_blowups}
	\end{center}
\end{figure}

In this case \W{$\mu \in \mathcal I_- \cup \mathcal I_+$,} the point $p \in \{(x,y,0) : y > 0\}$ is an unstable focus \SJ{when $\mu \in \mathcal I_+$,} and $(x,y,\epsilon) = (0,0,0)$ is a visible (invisible) fold point \SJ{when} \W{$\mu \in \mathcal I_+$ ($\mu \in  \mathcal I_-$);} see \lemmaref{PWS_lemma}. 

Applying the same ($\mu-$independent) cylindrical blow-up described in \secref{results_bifurcation} resolves the non-smoothness associated with $\Sigma \times \{0\} \in \mathbb R^2 \times \mathbb R_+$, however this time one identifies
an equilibrium $p_s \in S$ in case $\gamma, \mu > 0$, which is unstable when viewed as an equilibrium on $S$; see \figref{S2_blowups}a, and an equilibrium $p \in S$ for $\gamma, \mu < 0$, which is stable when viewed as an equilibrium on $S$.\footnote{This case is not shown in \figref{S2_blowups}, which only shows the case $\mu \in \mathcal I_+$, but is shown in \figref{complete_BF3}f.} For $\mu > 0$ and either $\gamma > 0$ or $\gamma < 0$, the reduced flow on $S$ is locally increasing in the $x-$coordinate in a neighbourhood of the nonhyperbolic point $Q$ (stemming from the visible fold singularity), where it terminates.

Recall that $Q$ is nonhyperbolic for all $\mu$, and can be resolved by a successive weighted spherical blow-up; see \figref{S2_blowups}c and \figref{S2_blowups}d. As in case $\mu = 0$, the manifold $S$ connects to a partially hyperbolic point $p_a$ at the intersection of blow-up sphere and cylinder. In contrast to the $\mu = 0$ case however, the manifold $\mathcal W$ emanating from $p_a$ is a unique centre manifold (not a critical manifold), which extends over the blow-up sphere and connects to a hyperbolic saddle $q_{out}$ lying within the intersection of the blow-up sphere and the plane $\{\epsilon = 0\}$. Partial hyperbolicity has been regained everywhere in \figref{S2_blowups}c and \figref{S2_blowups}d, and the existence of a connection from $q_{out}$ to a drop point on the edge of the cylinder (denoted $l_s$ in \figref{S2_blowups}) when $\mu \in \mathcal I_+$ allows for the construction of the nondegenerate `singular cycles' (closed loops consisting of distinct orbit and critical manifold segments; c.f. \defnref{pws_cycle}) in \figref{S2_blowups}c and \figref{S2_blowups}d.
\begin{remark}
	The BF collision in the blown-up space correlates with the collision of three equilibria $p, q_{in}, q_{out}$ in the limit $\mu \to 0^+$, with all three points coalescing in a single degenerate point $Q_{bfb}$ for $\mu = 0$, as can be seen by comparing, e.g. \figref{S2_blowups}d and \figref{S1_blowups_i}d.
\end{remark}

We now state the main result for the dynamics in (S2), which (among other things) describes persistence of the PWS cycles associated with (PWS) BF$_i$, $i=1,3$ bifurcations as 
relaxation oscillations in the smooth problem. 

\begin{theorem}
	\thmlab{thm_ro}
	Consider system \eqref{normal_form} under 
	\assumptionref{ass3} and \assumptionref{ass2}. There exists $\epsilon_0 > 0$ and $\mu_\pm > 0$ such that for all $\epsilon \in (0, \epsilon_0)$ the following assertions are true:
	\begin{enumerate}
		\item[(i)] For $\mu \in \mathcal I_+$ there exists an \SJ{unstable focus}
		\[
		p_\epsilon : (x_\epsilon, y_\epsilon) = \left( \mathcal O\left(\mu^2,\epsilon^k \right), \frac{\mu}{\delta} +  \mathcal O\left(\mu^2,\epsilon^k \right) \right) ,
		\]
		and in case $\gamma > 0$ \SJ{there also exists a hyperbolic saddle} 
		\[
		p_{s,\epsilon} : (x_{s,\epsilon}, y_{s,\epsilon}) = \left(-\frac{\mu}{\gamma} + \mathcal O(\mu^2,\epsilon) ,  \mathcal O\left(\mu^2,\epsilon \right) \right)  .
		\]
		\item[(ii)] For and $\mu \in \mathcal I_-$ and $\gamma < 0$, there exists \SJ{a stable node}
		\[
		p_\epsilon : (x_\epsilon, y_\epsilon) = \left(-\frac{\mu}{\gamma} + \mathcal O(\mu^2,\epsilon) ,  \mathcal O\left(\mu^2,\epsilon \right) \right) .
		\]
		\item[(iii)] For $\mu \in \mathcal I_+$, the PWS cycle $\Gamma = \Gamma^1 \cup \Gamma^2$ in both cases BF$_1$ (\figref{collision_figs}a) and BF$_3$ \SJ{(\figref{collision_figs}g)} perturbs to a strongly attracting relaxation oscillation $\Gamma_\epsilon$, which is $\mathcal O(\epsilon^{2k/(2k+1)})-$close to $\Gamma$ in the Hausdorff distance.
		There is no relaxation oscillation corresponding to case BF$_2$ \SJ{(\figref{collision_figs}d)}.
	\end{enumerate}
\end{theorem}

\bpr
	See \appref{proof_of_thm_ro}. Statements (i)-(ii) describe persistence of equilibria identified in \lemmaref{PWS_lemma}, and (iii) describes persistence of the PWS cycles as relaxation oscillations in the perturbed problem. 
\epr
\begin{remark}
	As discussed in \secref{introduction}, the observed `relaxation oscillations' are \textit{not slow-fast}, though they are clearly a consequence of singularly perturbed dynamics. Recall \remref{rem_singular} which emphasises that the systems being considered are `singular' in a different sense (i.e. they lose smoothness); see also \cite{Jelbart2019c,Jelbart2020a,Kosiuk2016,Kristiansen2019b,Kristiansen2019d} for examples of systems exhibiting non-slow-fast but nonetheless `relaxation-type' behaviour.
\end{remark}

\subsection{Connecting the dynamics \W{between regimes} (S1) and (S2)}
\seclab{results_connection}

Having separately described the dynamics in (S1) and (S2), we turn our attention to the `connection' between these regimes, i.e.~we are interested in connections associated with the matching problems obtained in the dual limits \eqref{matching_limit_neg} and \eqref{matching_limit}. The first result describes the motion of equilibria between regimes (S1) and (S2) under $\mu-$variation.

\begin{proposition}
	\proplab{eq_match}
	There exists $\epsilon_0 > 0$ such that for all \SJnew{$\epsilon \in (0,\epsilon_0)$} system \eqref{normal_form} has smooth parameterised families of equilibria
	\begin{equation}
	\eqlab{eq_fams}
	\mu \mapsto p_\epsilon(\mu), \qquad \text{and} \qquad \mu \mapsto p_{s,\epsilon}(\mu) ,
	\end{equation}
	on the interval $\mu \in (-\mu_+,\mu_+)$. For $\mu \in \mathcal I$ the equilibria $p_\epsilon(\mu)$ and $p_{s,\epsilon}(\mu)$ are the same equilibria described in \thmref{thm_bifurcation}, while for $\mu \in \mathcal I_+ \cup \mathcal I_-$ they are the same equilibria described in \thmref{thm_ro}.
\end{proposition}

\bpr
	\SJnew{This requires an understanding of the dynamics associated with both limits \eqref{matching_limit_neg} and \eqref{matching_limit} in the blow-up, where the families of equilibria \eqref{eq_fams} show up as 2-dimensional surfaces in the extended $(x,y,\epsilon,\mu)-$space. Proving \propref{eq_match} amounts to finding suitable parameterisations for these surfaces in multiple overlapping coordinate charts. We omit the details, which can be found in the forthcoming PhD thesis \cite{Jelbart2020b}.}
\epr

Now consider the matching problem associated with \eqref{matching_limit} and restrict attention to the BF$_3$ bifurcations, i.e. we assume $\gamma < 0$ in system \eqref{normal_form}.
\begin{figure}[t!]
	\centering
	\subfigure[]{\includegraphics[width=.49\textwidth]{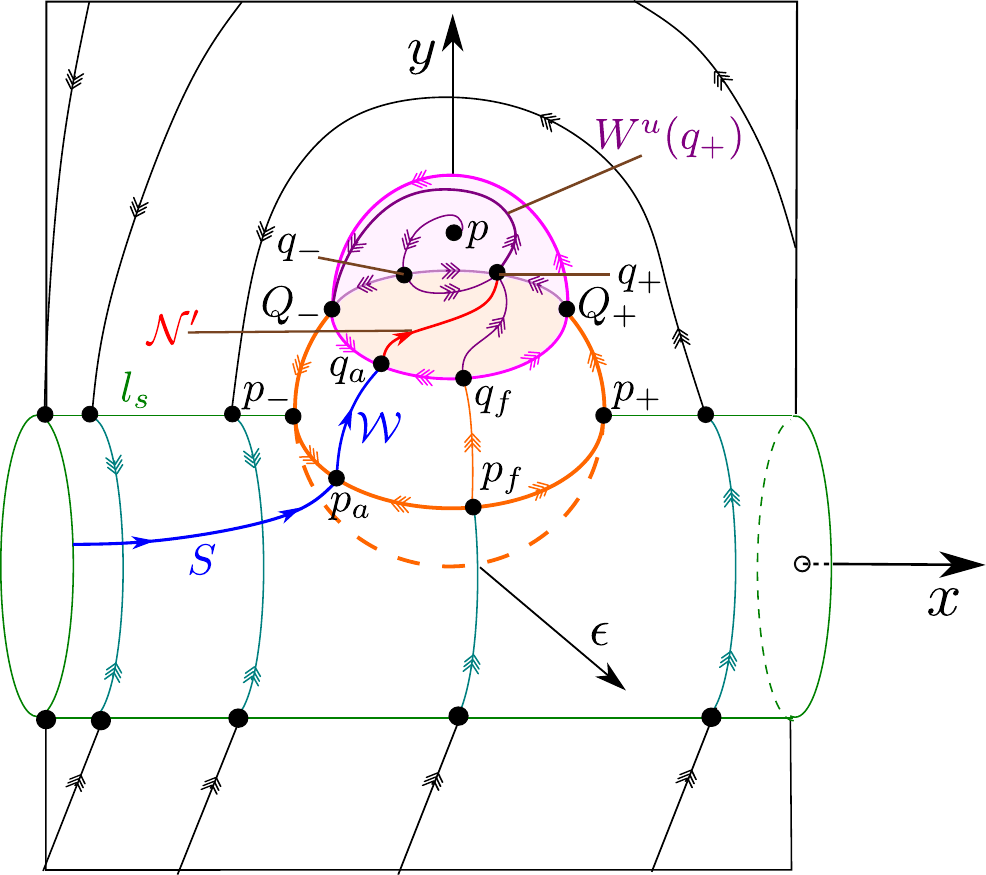}}
	\subfigure[]{\includegraphics[width=.49\textwidth]{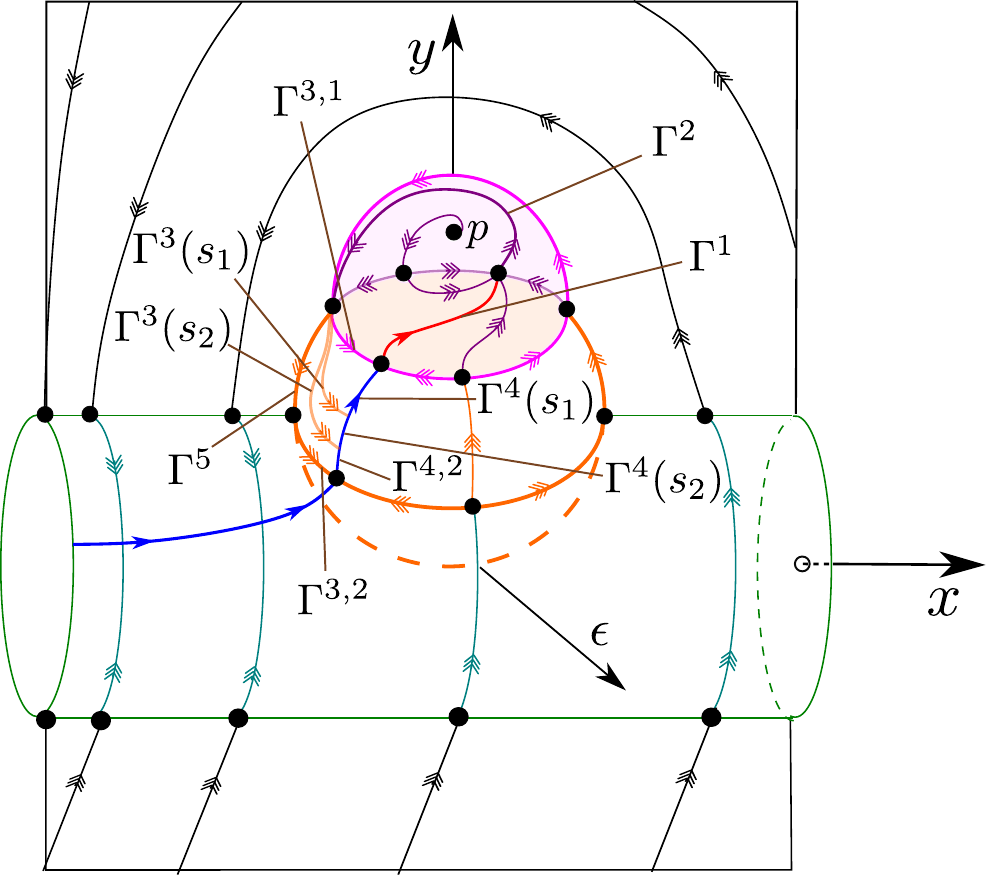}}
	\caption{In (a): Key dynamics following blow-up in the limit \eqref{matching_limit}.
		Invariant subsets corresponding to intersections of the blow-up 3-sphere (which replaces $Q_{bfb}$) with the `first' blow-up sphere (which replaces $Q$) and the plane $\{\epsilon = 0\}$ are shown in shaded orange and magenta respectively. The attracting critical manifold $\mathcal W$ connects to a partially hyperbolic point $q_a$, and a unique centre manifold \SJnew{$\mathcal N'$} contained within the intersection of the blow-up spheres emanates from $q_a$ and terminates at $q_+$. Within the intersection of the 3-sphere and $\{\epsilon = 0\}$ there exists a regular connection $W^u(q_+)$ between $q_+$ and a hyperbolic saddle $Q_-$. In (b): A family of singular cycles $\{\Gamma(s) : s \geq 0\}$ of relaxation type bounded between `small' and `large' singular cycles $\Gamma^s = \Gamma^1 \cup \Gamma^2 \cup \Gamma^{3,1}$ and $\Gamma^l = \Gamma^1 \cup \Gamma^2 \cup \Gamma^{3,2} \cup \Gamma^{4,2} \cup \Gamma^5$ respectively. Since $Q_-$ is an unstable node upon restriction to the invariant region bounded between $\Gamma^s$ and $\Gamma^l$ on the first blow-up sphere, the interior of this region is filled out by trajectories that are backward asymptotic to $Q_-$. This allows for the construction of a family of singular cycles (formally defined in \appref{singular_cycles}). Two candidate cycles $\Gamma(s_i) = \Gamma^1 \cup \Gamma^2 \cup \Gamma^3(s_i) \cup \Gamma^4(s_i)$ for $i=1,2$ are shown here, with $s_2 > s_1 > 0$.}
	\figlab{blowup_connection}
\end{figure}
\figref{blowup_connection}a shows the key dynamics identified in this limit, following the blow-up of $Q_{bfb}$. In order to effectively visualise the relevant dynamics (recalling that the point $Q_{bfb}$ has been replaced by a 3-sphere), only the important dynamics in key invariant subspaces are shown. These regions are in shaded orange and magenta in \figref{blowup_connection}, and correspond to intersections between the blow-up 3-sphere with the blow-up 2-sphere (resulting from blow-up of $Q$), and the plane $\{\epsilon = 0\}$ respectively.
A unique centre manifold \SJnew{$\mathcal N'$} emanating from $q_a$ extends across the intersection of the two blow-up spheres, connecting to $q_+$ which lies within the intersection of both blow-up spheres with $\{\epsilon = 0\}$. Within the intersection of the 3-sphere and $\{\epsilon = 0\}$, one identifies a regular connection $W^u(q_+)$ from $q_+$ to the point $Q_-$, which is a hyperbolic saddle. Considered within the invariant (planar) region bounded by $\Gamma^{3,1}$, $\Gamma^{3,2}$, $\Gamma^{4,2}$ and $\Gamma^5$ (see \figref{blowup_connection}b), however, the equilibrium $Q_-$ is an unstable node, thus providing an entire family of connections back to $\mathcal W$. \figref{blowup_connection}b shows two candidate singular cycles $\Gamma(s_j)$, $j=1,2$, chosen from an entire family 
\begin{equation}
\eqlab{family_cycs}
\{\Gamma(s) : s \geq 0\} , \qquad \Gamma(s) = \Gamma^1 \cup \Gamma^2 \cup \Gamma^3(s) \cup \Gamma^4(s) ,
\end{equation}
bounded between the singular cycles
\begin{equation}
\eqlab{outer_sing_cycs}
\Gamma^s = \Gamma^1 \cup \Gamma^2 \cup \Gamma^{3,1} , \qquad
\Gamma^l = \Gamma^1 \cup \Gamma^2 \cup \Gamma^{3,2} \cup \Gamma^{4,2} \cup \Gamma^5 .
\end{equation}
As part of the proof \thmref{thm_connection} below, it is shown that the singular cycles $\Gamma^s$ and $\Gamma^l$ bounding the family $\{\Gamma(s) : s \geq 0\}$ perturb to regular cycles in regime (S1), and relaxation-type cycles in regime (S2) respectively. Thus by understanding how the family of cycles $\{\Gamma(s) : s \geq 0\}$ perturbs, we can describe a transition from regular oscillation in (S1) to relaxation oscillation in (S2).

\

We now state our final main result, identifying a connection between limit cycles of regular and relaxation type in regimes (S1) and (S2) respectively.

\begin{theorem}
	\thmlab{thm_connection}
	Consider system \eqref{main} under \assumptionref{ass1}, \assumptionref{ass1b}, \assumptionref{ass3} and \assumptionref{ass2}, with fixed $\gamma < 0$. Then there exists $\epsilon_0 > 0$, $K > 0$ and $\mu_+ > 0$ such that \SJ{for all $\epsilon \in (0,\epsilon_0)$} system \eqref{main} has a \SJnew{continuous} parameterised family of stable limit cycles 
	\begin{equation}
	\eqlab{fam_cycs}
	\mu \mapsto \Gamma(\mu,\epsilon), \qquad \mu \in \left(K \epsilon^{k/(k+1)}, \mu_+\right),
	\end{equation}
	containing (but not limited to) the relaxation oscillations identified for $\mu \in \mathcal I_+$ in \thmref{thm_ro}.
\end{theorem}


\begin{figure}[t!]
	\centering
	\subfigure[]{\includegraphics[width=.312\textwidth]{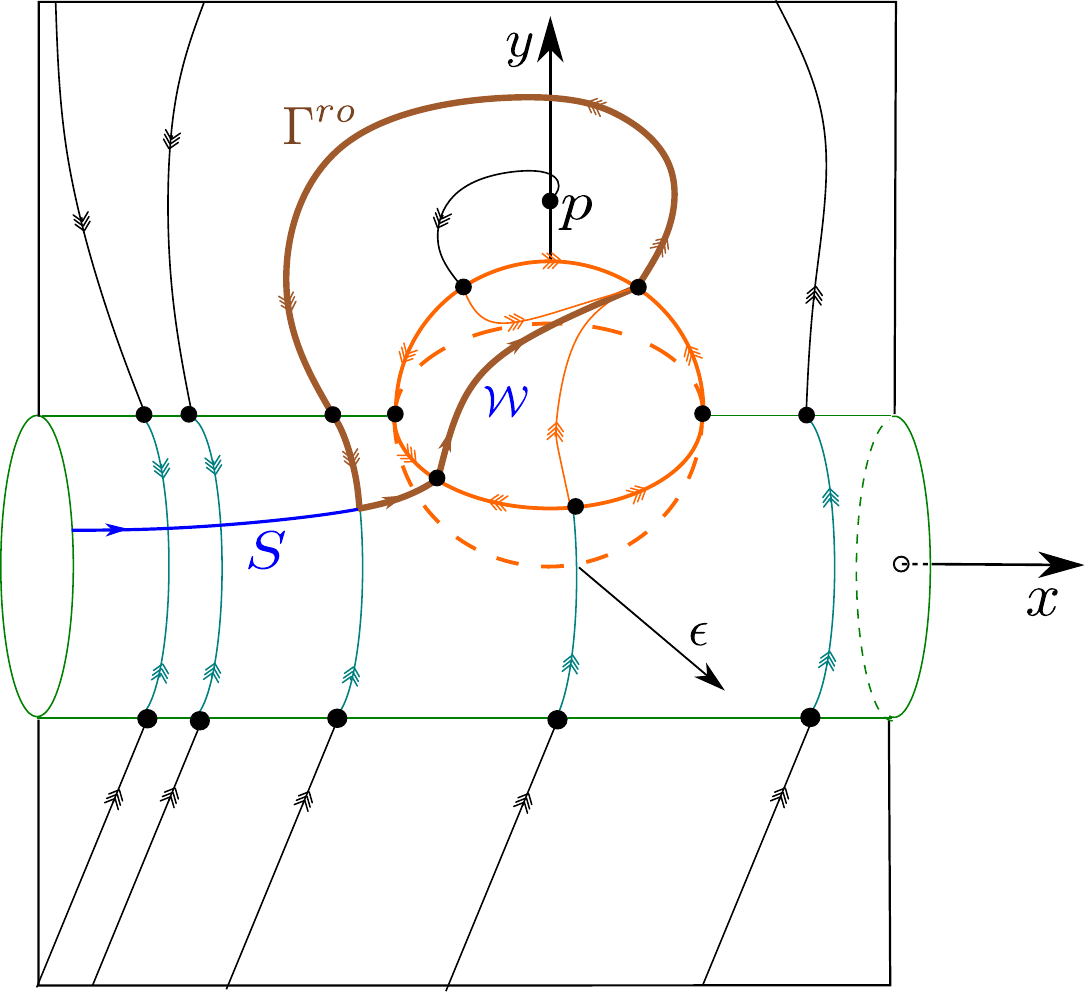}}
	\subfigure[]{\includegraphics[width=.32\textwidth]{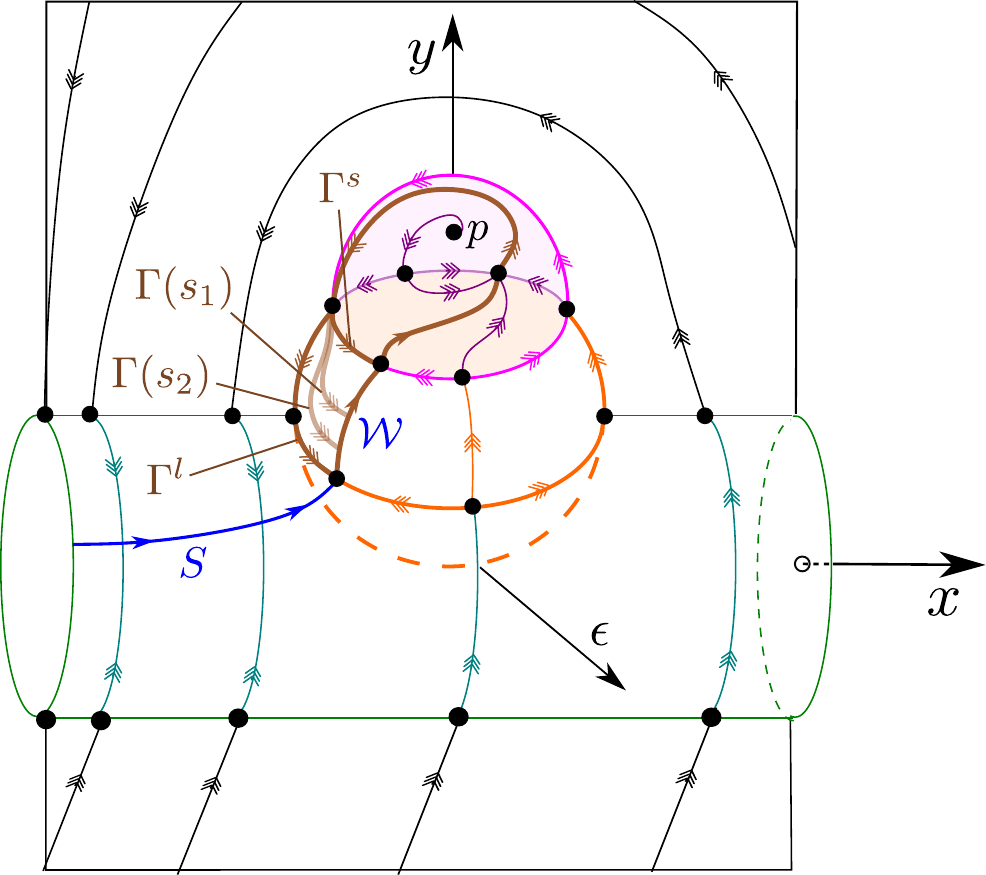}}
	\subfigure[]{\includegraphics[width=.32\textwidth]{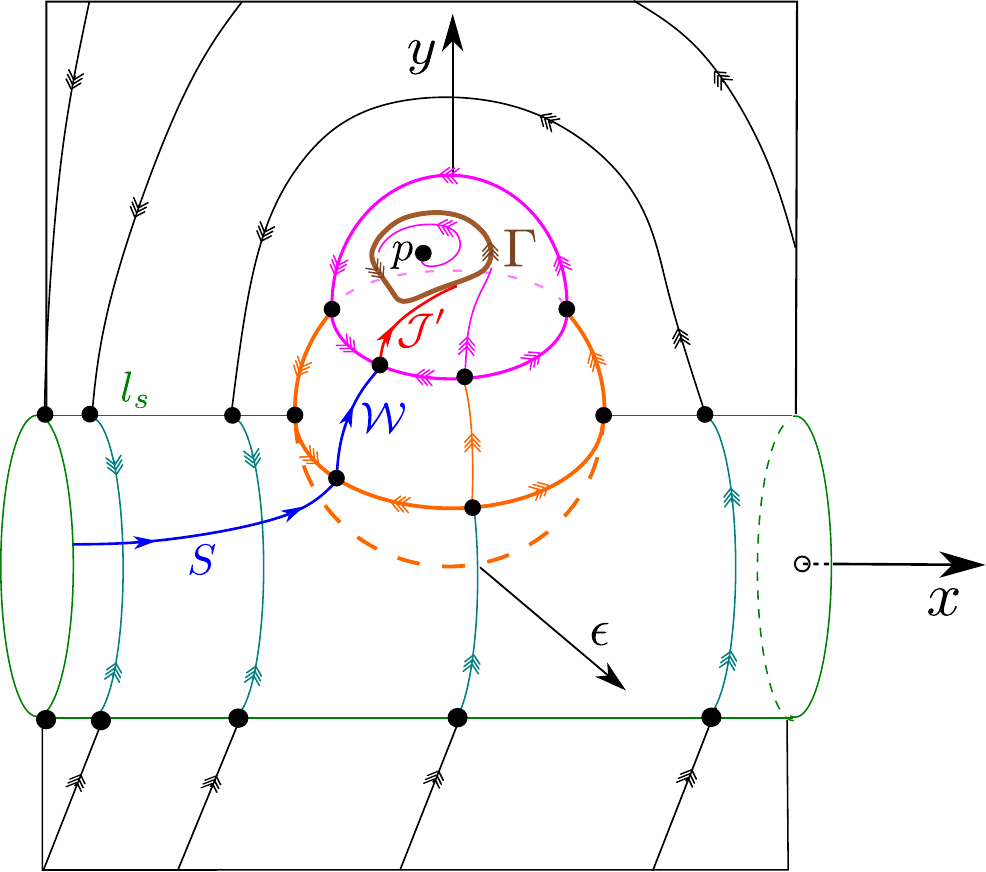}}
	
	\subfigure[]{\includegraphics[width=.32\textwidth]{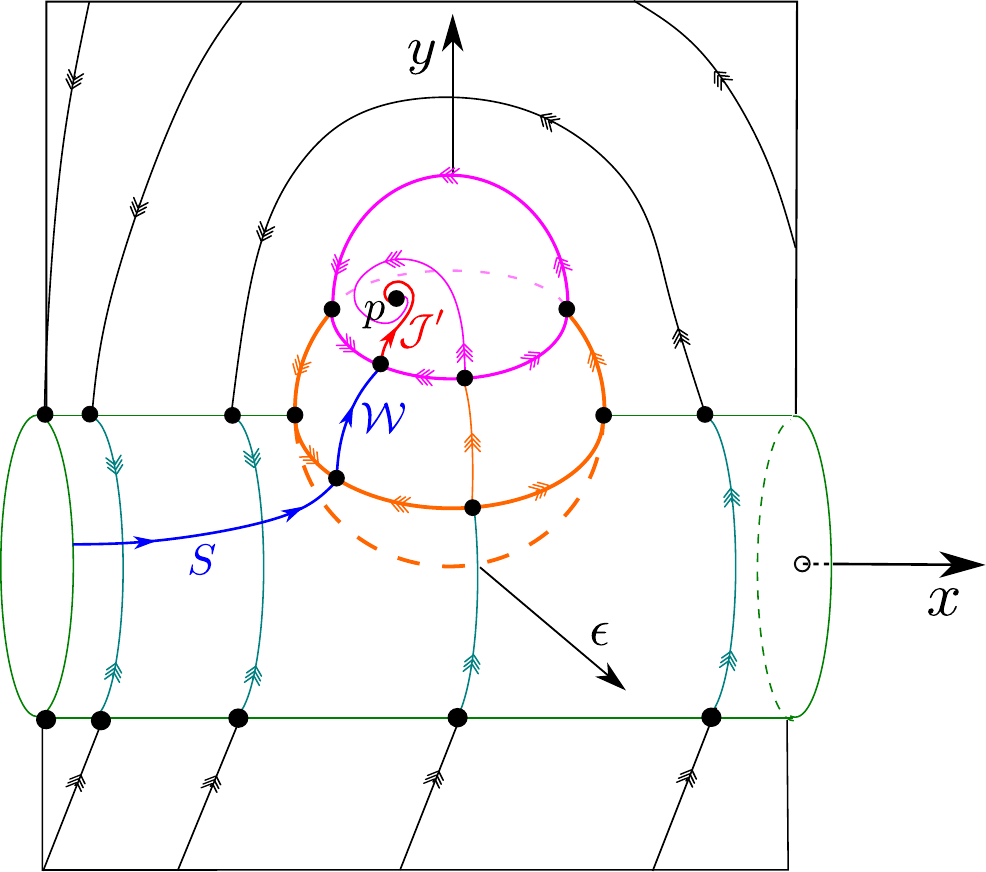}}
	\subfigure[]{\includegraphics[width=.32\textwidth]{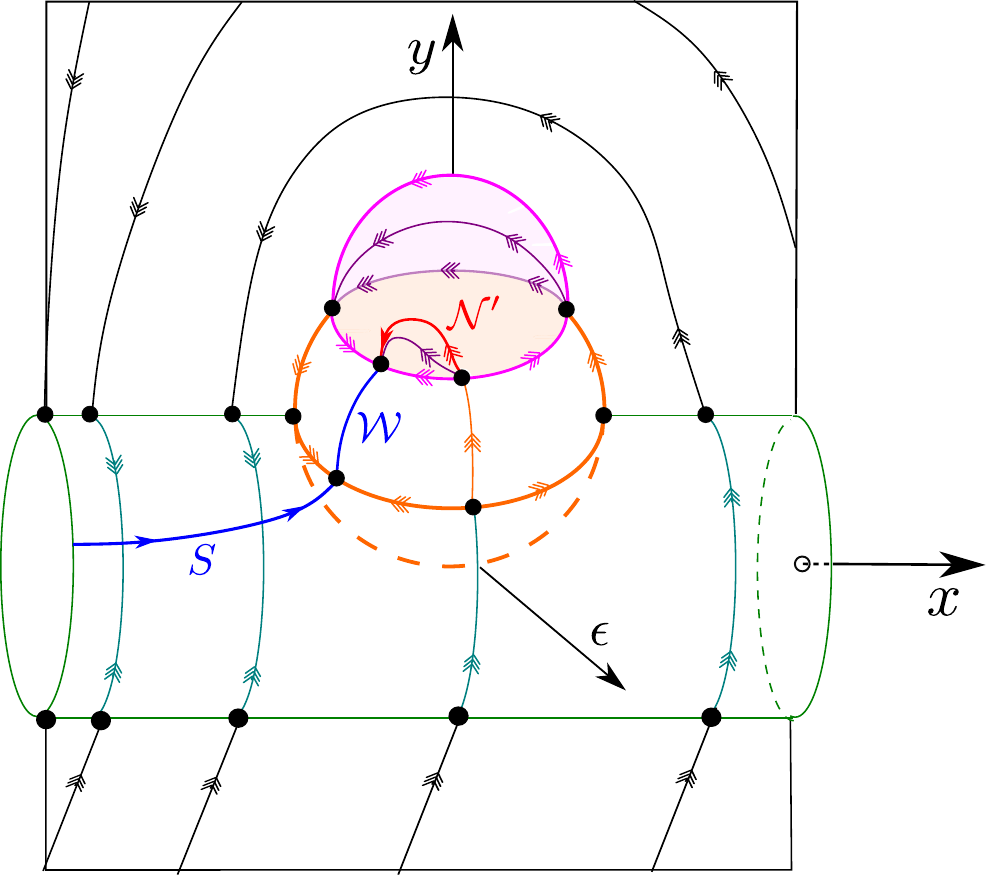}}
	\subfigure[]{\includegraphics[width=.312\textwidth]{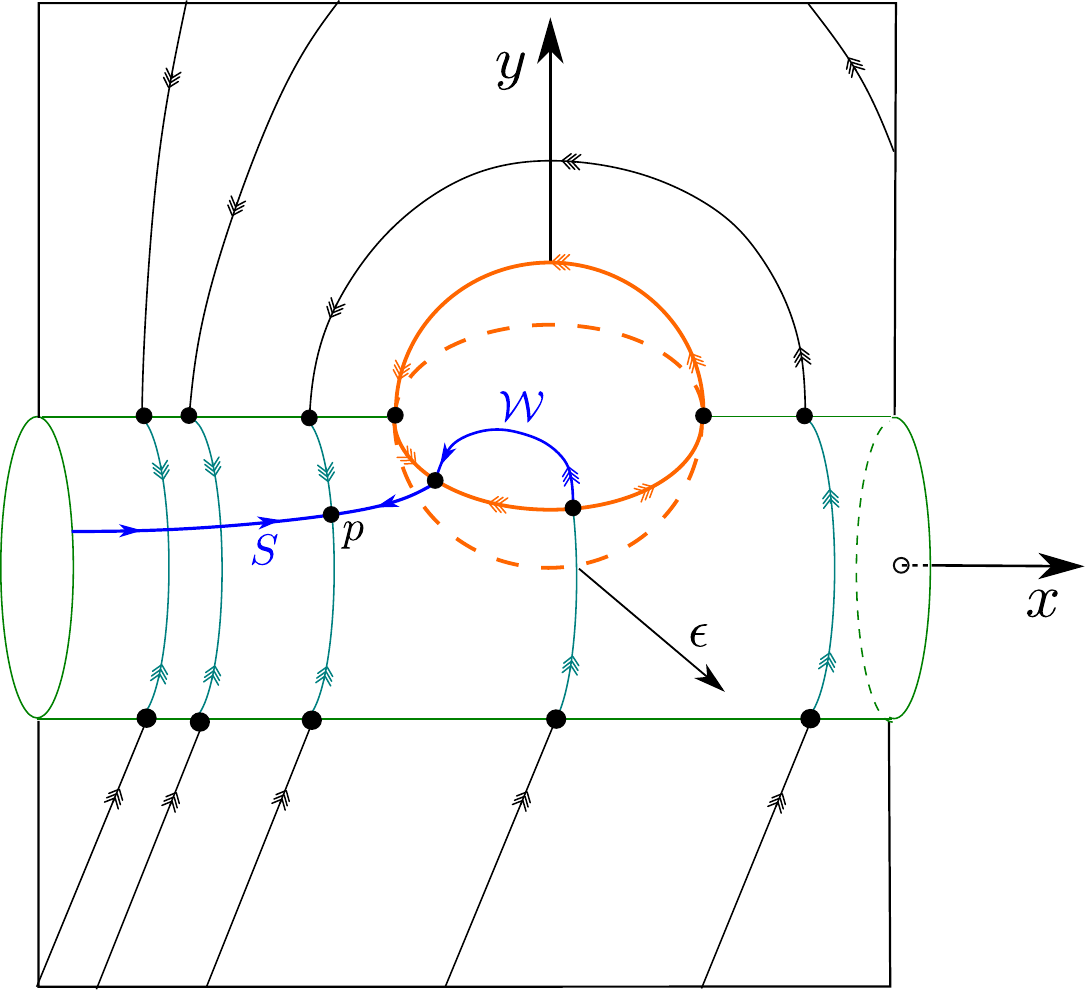}}	
	\caption{The unfolding for systems with type BF$_3$ bifurcation in the nonsmooth limit $\epsilon \to 0^+$. In (a): $\mu \in \mathcal I_+$. There exists a singular relaxation cycle \SJnew{$\Gamma^{ro}$} (brown) which perturbs to a 
		relaxation oscillation by \thmref{thm_ro}. In (b): the limit \eqref{matching_limit}. A family of singular cycles can be constructed as in \figref{blowup_connection}. The largest of these is $\Gamma^l$, which perturbs to relaxation oscillations in (S2) for $\mu > 0$. The smallest is $\Gamma^s$, which perturbs to regular cycles in (S1) for $\hat \mu \gg 1$. (c) 
		$\hat \mu \in (\hat \mu_{ah}, \hat \mu_{ah} + \kappa)$. An attracting regular cycle $\Gamma$ (also in brown) born/terminating in the Andronov-Hopf bifurcation exists in (S1). The cycle is depicted on the surface of the blow-up 3-sphere (magenta) restricted to (S1); see also \figref{S1_blowups_ii} and \figref{bfb_bif_neg}. (d) 
		$\hat \mu \in (-\infty, \hat \mu_{ah})$. The limit cycle in (c) has terminated in the supercritical Andronov-Hopf bifurcation for $\hat \mu = \hat \mu_{ah}$ described in \thmref{thm_bifurcation}, leaving behind a stable focus $p$. (e) The limit \eqref{matching_limit_neg}. Here the equilibrium $p \in \mathcal W$. 
		Flow on the (unique) 1D centre manifold \SJnew{$\mathcal N'$} (red) is now locally toward the intersection of $\mathcal W$ with the blow-up 3-sphere. (f) $\mu \in \mathcal I_-$. The equilibrium $p \in S$ has emerged on the blow-up cylinder, consistent with \propref{eq_match}.}
	\figlab{complete_BF3}
\end{figure}

\bpr
	\SJnew{This follows from \thmref{thm_cycs}, which constitutes a more general description of the family of limit cycles \eqref{fam_cycs}. This is deferred to \appref{proof_of_thm_connection}, and relies on the rigorous blow-up analysis performed in \appref{blow-up_of_Sigma_and_Q} and \appref{blow-up_of_Qbfb}.} 
\epr

\begin{remark}
	\remlab{connection}
	The supercritical Andronov-Hopf bifurcation identified in \thmref{thm_bifurcation} also implies the existence of a family of stable limit cycles $\Gamma_{ah}(\mu,\epsilon)$ for all $\mu \in (\mu_{ah}(\gamma,\epsilon), \mu_{ah}(\gamma,\epsilon) + \kappa \epsilon^{k/(k+1)})$ for some $\kappa > 0$. To prove a smooth, regular connection between the Andronov-Hopf cycles $\Gamma_{ah}(\mu,\epsilon)$ and the cycles $\Gamma(\mu,\epsilon)$ described in \thmref{thm_connection}, one must show that the Andronov-Hopf cycles in the desingularised system \eqref{desing_prob} grow unboundedly for all $k \in \mathbb N_+$. \SJnew{Numerical continuations of system \eqref{desing_prob} with parameter values $\tau = \delta = 1$, $\gamma = -1$, $\beta = 1/4$ and $k = 2$ consistent with the regularisation function in \figref{reg_fig} suggests the existence of such a connection, however} we do not consider this further in the present manuscript.
\end{remark}


Combining \thmref{thm_connection} with \thmref{thm_bifurcation} and \thmref{thm_ro}, we obtain a complete characterisation for general systems \eqref{main} satisfying \assumptionref{ass1}, \assumptionref{ass1b}, \assumptionref{ass3} and \assumptionref{ass2} with BF$_3$ bifurcations in the nonsmooth limit $\epsilon \to 0^+$, barring a connection between regular cycles within (S1) (see \remref{connection} above). This is sketched in the blown-up space in \figref{complete_BF3}, which shows the transition from $\mu \in \mathcal I_+$ through to $\mu \in \mathcal I_-$.

\section{BF bifurcation in applications}
\seclab{bf_bifurcation_in_application}

We consider two applications, 
namely:
\begin{enumerate}
	\item [(M1)]A Gause predator-prey model with a smooth Holling-type response, and a BF$_3$ bifurcation in the nonsmooth limit.
	\item[(M2)] A mechanical oscillator subject to friction, with a degenerate BF bifurcation in the nonsmooth limit which violates the nondegeneracy condition $\gamma = (\partial Z_{sl} / \partial x)|_{(0,0)} \neq 0$.
\end{enumerate}
Our results from \secref{main_results} apply directly to (M1). 
In (M2) we are able to adapt our analysis to describe the birth of oscillations in regime (S1) and the persistence of PWS cycles as relaxation oscillations in regime (S2) respectively, despite the degeneracy due to $\gamma = 0$. 



\subsection{BF$_3$ bifurcation in model (M1)}
\seclab{gause}

We consider the following Gause predator-prey model with refuge originally proposed in \cite{Gause1936}:
\begin{equation}
\eqlab{Gause}
\begin{split}
\dot R &= r R - C f(R) , \\
\dot C &= \left(\tilde e f(R) - m \right) C,
\end{split}
\end{equation}
where the variables $R$ and $C$ denote prey and predator densities respectively, and $r, \tilde e, m > 0$ are fixed parameters; see \cite{Gause1936,Krivan2011} for details on the model and it's interpretation.
The function $f$ is the so-called \textit{Gause functional response}, which is expected to have a \textit{smooth but sharp transition} near a critical prey population threshold $R_c$. In \cite{Krivan2011}, the author considers a PWS approximation for $f(R)$ of the following form:
\begin{equation}
\eqlab{Krivan_response}
f(R) = 
\begin{cases}
0 , & R < R_c , \\
\left[0, \frac{\lambda R_c}{1 + \lambda h R_c} \right] , & R = R_c , \\
\frac {\lambda R}{1 + h \lambda R} , & R > R_c .
\end{cases}
\end{equation}
In order to consider a \textit{smooth} response, we consider the following regularisation:
\begin{equation}
\eqlab{reg_holling}
f(R) = \phi\left((R-R_c)\epsilon^{-1}\right) \frac {\lambda R}{1 + h \lambda R} ,
\end{equation}
where the regularisation function $\phi:\mathbb R \to \mathbb R$ is defined to satisfy \assumptionref{ass1b} and \assumptionref{ass3}. 


\begin{remark}
	\remlab{gause_regularisation}
	In \cite{Krivan2011} the PWS functional response \eqref{Krivan_response} is described as the $\sigma \to \infty$ limiting case of the following Holling type-III response curve:
	\[
	f_{III}(R,\sigma) = \frac{\lambda R^{1+\sigma} / (R^{\sigma} - R_c^{\sigma})}{1 + h \lambda R^{1+\sigma} / (R^{\sigma} - R_c^{\sigma})} .
	\]
	Following an approach put forward in \cite{Szmolyan2017}, a convenient form can be obtained by setting $\sigma = \epsilon^{-1}$, $R = e^u$ and $R_c = e^{u_c}$, which yields
	\[
	f_{III} \left(e^u, \epsilon^{-1} \right) = \phi \left((u - u_c)\epsilon^{-1} , \epsilon \right) \frac{\lambda e^u}{1 + h \lambda e^u} ,
	\]	
	where the regularisation function
	\begin{equation}
	\eqlab{gause_reg}
	\phi(s,\epsilon) = \frac{1 + h \lambda e^{u_c} e^{\epsilon s}}{1 + e^{-s} + h \lambda e^{u_c} e^{\epsilon s}} 
	\end{equation}
	satisfies \assumptionref{ass1b}, although it has a more general form with a second argument $\epsilon$. Such regularisations have been considered in the context of both theory and applications in, e.g. \cite{Kristiansen2019c,Kristiansen2019d}. Since the second argument in $\epsilon$ is not expected to lead to significant dynamical differences in the regularised BF bifurcation, we restrict to the smaller class of regularisations defined by \assumptionref{ass1b} for simplicity.
	
	More significantly, the regularisation function \eqref{gause_reg} is flat for $s \to \pm \infty$, therefore violating \assumptionref{ass3}. This leads to significant (though not insurmountable) complications for the blow-up analysis relating to the loss of hyperbolicity at greater than algebraic rates; see \cite{Kristiansen2017} on how to deal with this using blow-up, and \cite{Bossolini2017b,Jelbart2019c,Kristiansen2019b} for applications.
\end{remark}

We make the following change of coordinates for the sake of consistency with earlier notations:
\[
x = C, \qquad y = R - R_c , \qquad \mu = R_c ,
\]
and consider the regularised problem
\begin{equation}
\eqlab{reg_gause}
\begin{pmatrix}
\dot{x}\\\dot{y}
\end{pmatrix}
=
\begin{pmatrix}
x \left(-m + \phi \left( y \epsilon^{-1} \right) \tilde e w(y, \mu) \right)  \\
r (\mu + y) - \phi \left( y \epsilon^{-1} \right) x w(y, \mu) ,
\end{pmatrix}
=: Z(x,y,\mu,\epsilon) 
\end{equation}
obtained by considering
\[
Z(x,y,\mu,\epsilon) = \phi\left(y \epsilon^{-1}\right) Z^+(x,y,\mu) + \left(1 - \phi\left(y \epsilon^{-1}\right) \right) Z^-(x,y,\mu) ,
\]
where
\begin{equation}
\eqlab{PWS_gause}
\begin{pmatrix}
\dot{x}\\\dot{y}
\end{pmatrix}
= \left\{
\begin{aligned}
&
\begin{pmatrix}
x \left(\tilde e w(y, \mu) - m \right)\\
r \mu + r y - x w(y, \mu) 
\end{pmatrix} 
=:  Z^+(x,y,\mu) &\qquad (y>0) ,\\
&
\begin{pmatrix}
- m x \\
r \mu + r y 
\end{pmatrix}  
=: Z^-(x,y,\mu) & \qquad (y<0) , 
\end{aligned}
\right.
\end{equation}
with
\[
w(y, \mu) = \frac{\lambda (y + \mu)}{1 + h \lambda (y + \mu)} .
\]
Note that the physically meaningful domain in our new coordinates is $x \geq 0$, $y \geq - \mu$. We also assume `low handling times' $h \in (0,h_\ast)$, where
\[
h_\ast = \frac{2e}{r} \left(-1 + \sqrt{1 + \frac{r}{m}}\right) .
\]
\W{Note, the regularised system \eqref{reg_gause} fulfils Assumptions 1-3}.

\W{
	\begin{lemma}
		The PWS system \eqref{PWS_gause} \SJ{with $h \in (0,h_\ast)$} has a BF bifurcation at
		\begin{equation}
		\eqlab{gause_bf_value}
		\mu_{bf} = \frac{m}{\lambda (\tilde e - hm)}\,.
		\end{equation}
	\end{lemma}
	
	\bpr
		The switching manifold $\Sigma = \{(x,y) \in \mathbb R^2 : f_\Sigma(x,y) = y = 0\}$ decomposes into crossing and sliding submanifolds 
		$\Sigma_{cr} = \left\{(x,0) : x \in [0,x_F) \right\}$
		and
		$\Sigma_{sl} = \left\{(x,0) : x > x_F \right\}$
		where
		\[
		x_F(\mu) = \frac{r}{\lambda} \left(1 + h \lambda \mu \right)\,.
		\]
		%
		A calculation at $F : (x_F(\mu),0)$ gives
		\[
		Z^+(Z^+ f_\Sigma)(F) = r \mu \left(m - e w(0,\mu) \right) ,
		\]
		from which it follows that the tangency at $F$ is a visible (invisible) fold for $\mu < \mu_{bf}$ ($\mu > \mu_{bf}$),
		where $\mu_{bf}$ is given by \eqref{gause_bf_value}.
		
		The Filippov vector field on $\Sigma_{sl}$ is given in the $x-$coordinate chart by
		\[
		\dot x = \tilde e r \mu - m x  =: Z_{sl}(x,\mu) , \qquad x \in \Sigma_{sl} .
		\]
		
		The PWS system \eqref{PWS_gause} has an equilibrium at
		\[
		p(\mu) = \left(\frac{\tilde e r}{\lambda (\tilde e - hm)} , -\mu + \frac{m}{\lambda (\tilde e - hm)} \right) ,
		\]
		where $\tilde e - hm > 0$ for all $h \in (0, h_\ast)$. This equilibrium $p(\mu)$ collides with $\Sigma$
		in the collision limit $\mu \to \mu_{bf}^-$\SJ{, where $\mu_{bf}$ is given by \eqref{gause_bf_value}}.
		\SJ{Direct calculations verify the conditions in \eqref{bf_conds}, and the
			eigenvalues of the Jacobian at $p(\mu_{bf})$ are given by
			\[
			\lambda_\pm = \frac{hmr}{2\tilde e} \pm \frac{\sqrt{mr}}{2\tilde e} \sqrt{-4\tilde e^2 + 4m\tilde eh+ mrh^2} ,
			\]
			which take the form $\lambda_\pm=A\pm iB$ with $A,B > 0$ for all $h \in (0,h_\ast)$, as required.}
	\epr
}
\begin{figure}[t!]
	\centering
	\subfigure[]{\includegraphics[width=.45\textwidth]{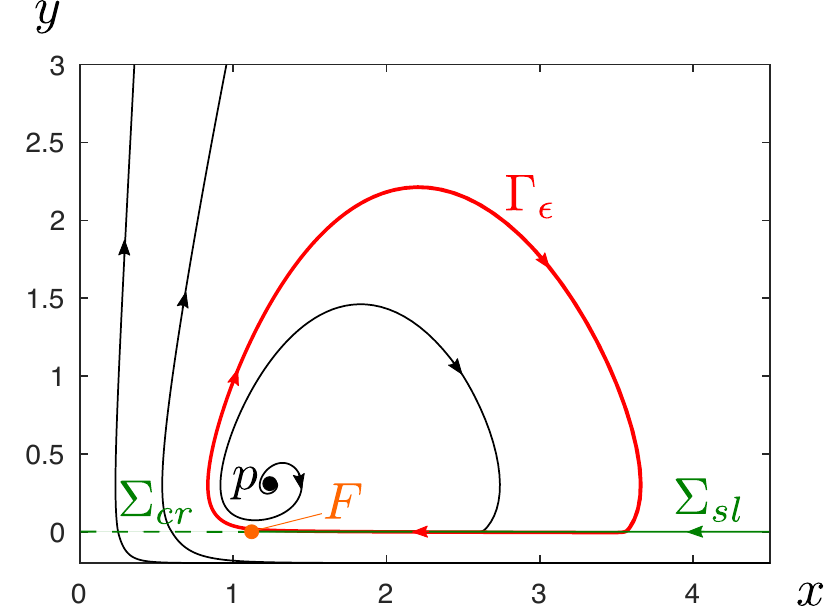}}
	\hspace{0.5cm}
	\subfigure[]{\includegraphics[width=.45\textwidth]{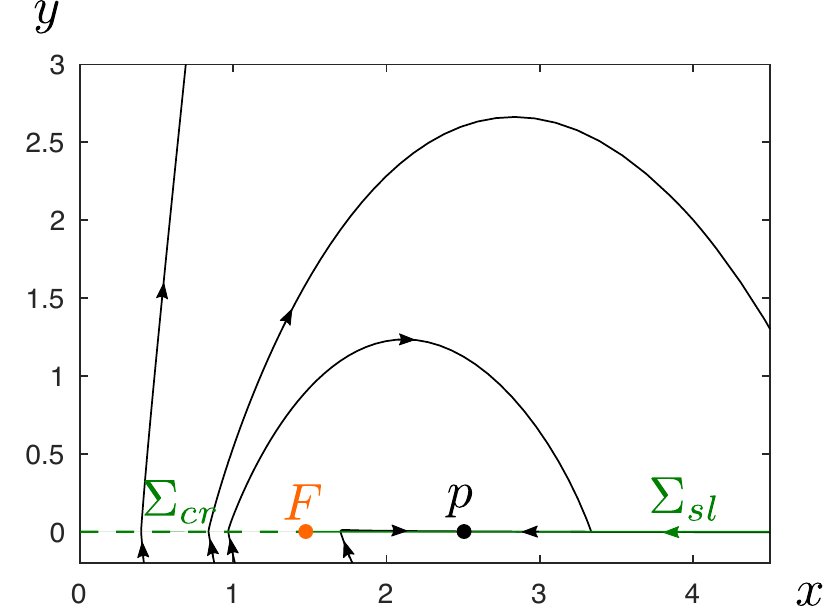}}
	\caption{Phase portraits for system \eqref{reg_gause} with the regularisation function from \figref{reg_fig}, $\epsilon = 10^{-2}$, overlaid with important objects identified in the PWS analysis of system \eqref{PWS_gause}. The parameter values $r=1, \ \lambda =1 , \ m=0.2 , \ \tilde e = 0.5$ are taken from \cite[Figure 2]{Krivan2011}, giving a BF value of $\mu_{bf} = 0.5$. Sliding and crossing submanifolds $\Sigma_{sl}$ and $\Sigma_{cr}$ 
		are separated by a tangency point $F$, 
		which is a visible (invisible) fold for $\mu < \mu_{bf}$ ($\mu > \mu_{bf}$). In (a): $\mu = 0.2 < \mu_{bf}$, for which the equilibrium $p$ (black disk) is an unstable focus and there exists a relaxation oscillation $\Gamma_\epsilon$ (red). In (b): $\mu = 1 > \mu_{bf}$, \SJnew{there are no oscillations and $p$ is asymptotically stable}.}
	\figlab{gause_pws}
\end{figure}

\figref{gause_pws} shows the dynamics for $\epsilon = 10^{-2}$, together with the PWS objects identified in this proof.
We note that for all $\mu_{bf} - \mu > 0$ sufficiently small, the unique $Z^+$ orbit intersecting $F \in \Sigma$ (where we consider the extension of $Z^+$ on $\{y \geq 0\}$) has a first return, or `drop point' $(x_d,0) \in \Sigma_{sl}$. Hence, we can construct PWS cycles $\Gamma = \Gamma^1 \cup \Gamma^2$ where
\[
\Gamma^2 = \left\{(x,0) : x \in [x_F,x_d] \right\} \subset \Sigma ,
\]
and $\Gamma^1 \subset \{y \geq 0\}$ is identified with the segment of the $Z^+$ orbit connecting $F$ and $(x_d,0)$.

\

The following result summarises our findings for system \eqref{reg_gause}.


\begin{theorem}
	\thmlab{gause_thm}
	Consider system \eqref{reg_gause} with $\epsilon \in (0,\epsilon_0)$ and $h \in (0,h_\ast)$. Then the following holds for $\epsilon_0 > 0$ sufficiently small:
	\begin{enumerate}
		\item[(i)] There exists a unique equilibrium $p_\epsilon(\mu)$ such that $p_\epsilon(\mu) \to p(\mu_{bf})$ as $\epsilon \to 0^+$ and $\mu \to \mu_{bf}^+$, which undergoes a supercritical Andronov-Hopf bifurcation for
		\begin{equation}
		\eqlab{gause_hopf}
		\mu_{ah}(\epsilon) = \mu_{bf} + \hat \mu_{ah} \epsilon^{k/(k+1)} + o \left(\epsilon^{k/(k+1)}\right) ,
		\end{equation}
		where $\mu_{bf}$ is given by \eqref{gause_bf_value} and the constant $\hat \mu_{ah}$ can be calculated explicitly in terms of the system parameters.
		\item[(ii)] There exist constants $\mu_\pm \in (0, \mu_{bf})$ such that
		for all $\mu \in (\mu_-, \mu_+)$, there exists a relaxation oscillation $\Gamma_\epsilon$ which is strongly attracting, unique in an arbitrarily large ball $B(0,M)$, and $\mathcal O(\epsilon^{2k/(2k+1)})-$close in Hausdorff distance to a PWS cycle $\Gamma = \Gamma^1 \cup \Gamma^2$.
		\item[(iii)] The branch of relaxation oscillations described in (ii) connects to a branch of regular cycles in (S1), i.e. there exists $K > 0$ such that system \eqref{reg_gause} has a smooth parameterised family of stable limit cycles
		\begin{equation}
		\eqlab{gause_cycles}
		\mu \mapsto \Gamma(\mu,\epsilon) , \qquad \mu \in \left( \mu_-, \mu_{bf} - K \epsilon^{k/(k+1)} \right) ,
		\end{equation}
		which contains (but is exhausted by) the relaxation oscillations in (ii).
	\end{enumerate}
\end{theorem}


\bpr
	\W{Since Assumptions 1-4 are satisfied, the statement (i) follows from \thmref{thm_bifurcation}, noting that the inequality $\gamma = (\partial Z_{sl} / \partial x)|_{(p(\mu_{bf}),\mu_{bf})} < 0$ implies a type BF$_3$ bifurcation.	
		
		The first part of statement (ii) follows from \thmref{thm_ro}}, i.e.~there exists a stable limit cycle of relaxation type $\Gamma_\epsilon$ for $\mu \in (\mu_-,\mu_+)$ with $\mu_\pm > 0$ sufficiently small which converges to a PWS cycle $\Gamma = \Gamma^1 \cup \Gamma^2$ in Hausdorff distance as $\epsilon \to 0^+$; it remains to prove uniqueness on $B(0,M)$. It follows from the contraction mapping argument in the proof for \thmref{thm_ro} \SJnew{(see \appref{proof_of_thm_ro})} that $\Gamma_\epsilon$ is unique in an open tubular neighbourhood $\mathcal T_\epsilon$ of width $\mathcal O(\epsilon^{2k/(2k+1)})$. \SJ{Applying the Dulac criterion \cite{Dulac1937} with the same Dulac function $1/(xy)$ used in \cite{Krivan2011} to prove non-existence of cycles in \eqref{PWS_gause}$|_{\{y>0\}}$,} we can show that no limit cycles can exist on domains $\mathcal D_\epsilon = \{(x,y) \in \mathbb R^2 : y > \epsilon^k \tilde c\}$ with constant $\tilde c > 0$, since
	\[
	\frac{\partial }{\partial x} \left( \frac{1}{xy} Z_1^+(x,y,\mu,\epsilon) \right) + \frac{\partial }{\partial y} \left( \frac{1}{xy} Z_2^+(x,y,\mu,\epsilon) \right) > 0 
	\]
	follows via the relation
	\[
	\frac{df(y+\mu)}{dy} = \frac{\lambda}{(1+\lambda h (y + \mu))^2} + \mathcal O(\epsilon^k) < \frac{f(y+\mu)}{y} = \frac{\lambda}{1+\lambda h (y+\mu)} + \mathcal O(\epsilon^k) .
	\]
	Since the intersection $\mathcal T_\epsilon \cap \mathcal D_\epsilon$ is non-empty and open, the result follows.
	
	\W{Statement (iii)}, i.e. existence of the smooth family of stable cycles \eqref{gause_cycles} containing regular oscillations in (S1) and relaxation oscillations in (S2), follows from \thmref{thm_connection}.
\epr

\begin{remark}
	\remlab{guase_desing}
	In order to obtain an explicit expression for $\hat \mu_{ah}$ in \eqref{gause_hopf}, one must study the desingularised system obtained \SJ{from system \eqref{reg_gause} after} applying the transformation
	\begin{equation}
	\eqlab{gause_coord_change}
	\left(x_1, \nu_1, \rho_1, \hat \mu \right) \in \mathbb R \times \mathbb R_+^2 \times \mathbb R \mapsto
	\begin{cases}
	x = \frac{e r}{\lambda (e - h m)} + \nu_1^{2k(1+k)} \rho_1^{k(1+k)} \left(h r \hat \mu \rho_1^{k^2} + x_1 \right) , \\
	y = \nu_1^{2k(1+k)} \rho_1^{2k(1+k)} , \\
	\epsilon = \nu_1^{2(1+k)^2} \rho_1^{(1+k)(1+2k)} , \\
	\mu = \mu_{bf} + \hat \mu \nu_1^{2k(1+k)} \rho_1^{k(1+2k)} ,
	\end{cases}
	\end{equation}
	\SJ{which is }obtained by composing \eqref{coord_change} with a (parameter-dependent) coordinate translation which moves the BF point so that it occurs at the origin for $\mu = 0$. After desingularising the resulting system by dividing the right-hand-side by a common factor of $\rho_1^{k(k+1)}$ and restricting to the invariant subspace $\{\nu_1 = 0\}$, one can locate the Andronov-Hopf bifurcation by standard methods in the resulting planar system; see \appref{proof_of_lem_bifurcation}, where this is done for system \eqref{desing_prob}.
\end{remark}



\subsection{Degenerate BF bifurcation in model (M2)}
\seclab{stick-slip}

For our second application, we consider the following (regularised) model for a mechanical oscillator subject to friction, in the formalism of \cite{Kristiansen2019c} (see e.g. \cite{Berger2002,Chen2014,Chen2018,Popp1990,Won2016} for PWS formulations):
\begin{equation}
\eqlab{stick-slip}
\begin{pmatrix}
\dot{x}\\\dot{y}
\end{pmatrix}
=
\begin{pmatrix}
y - \alpha, \\
-x - \mu \left(y, \phi\left(y \epsilon^{-1} \right) \right) 
\end{pmatrix}
= Z(x,y,\alpha,\epsilon) ,
\end{equation}
where $x$ and $y$ denote displacement and velocity of a mass relative to a conveyor moving at a constant speed $\alpha$, and $\phi(y \epsilon^{-1})$ is a regularisation function satisfying \assumptionref{ass1b} and \assumptionref{ass3}. The function $\mu(y, \phi(y \epsilon^{-1}))$ is known as the `characteristic of friction', and can take different forms depending on the application.
Following \cite{Kristiansen2019c}, we write
\begin{equation}
\eqlab{char}
\mu\left(y, \phi \left(y \epsilon^{-1}\right) \right) = \mu_+(y) \phi\left(y \epsilon^{-1} \right) - \mu_+(-y) \left(1 - \phi\left(y \epsilon^{-1} \right) \right) ,
\end{equation}
and consider the case where $\mu_+: \mathbb R \to \mathbb R$ is a smooth function satisfying a `Stribeck law', i.e.~there exists some $y_0 > 0$ (possibly infinite) such that
\begin{enumerate}
	\item[(i)] $\mu_+(y) > 0$ for all $y \in [0,y_0)$ and $\mu_+(0) =: \mu_s > 0$;
	\item[(ii)] $\mu_+'(y) \in (-2,0)$ for all $y \in [0,y_0)$.
\end{enumerate}
We also assume that the translated regularisation function $\psi(y \epsilon^{-1}) := \phi(y \epsilon^{-1}) - 1/2$ is odd
, from which it follows that
\begin{equation}
\eqlab{mu_sym}
\mu\left(-y,\phi\left(-y\epsilon^{-1}\right)\right) = \mu\left(-y, 1 - \phi \left(y \epsilon^{-1}\right) \right) =
- \mu\left(y,\phi\left(y \epsilon^{-1}\right)\right) .
\end{equation}
\SJnew{Physically, this is necessary in applications for which the friction force does not depend on the direction of motion.}



System \eqref{stick-slip} satisfies \W{\assumptionref{ass1}, \assumptionref{ass1b} and \assumptionref{ass3}} by construction, and converges in the singular limit $\epsilon \to 0^+$ to the PWS system
\W{
	\begin{equation}
	\eqlab{PWS_stick-slip}
	\begin{pmatrix}
	\dot{x}\\\dot{y}
	\end{pmatrix}
	=\left\{
	\begin{aligned}
	&
	\begin{pmatrix}
	y - \alpha \\
	- x - \mu_+(y) 
	\end{pmatrix}
	= Z^+(x,y,\alpha) 
	\quad& (y>0) ,
	\\
	&
	\begin{pmatrix}
	y - \alpha  \\
	- x + \mu_+(-y) 
	\end{pmatrix}
	= Z^-(x,y,\alpha)  
	\quad& (y<0) ,
	\end{aligned}
	\right.
	\end{equation}}
with switching manifold $\Sigma = \{(x,y) \in \mathbb R^2 : f_\Sigma(x,y) = y = 0\}$. Restricting our focus to belt speeds $\alpha \geq 0$ (see \remref{symmetry} below), we identify a unique equilibrium in system \eqref{PWS_stick-slip} at
\[
p(\alpha) : \left(- \mu_+(\alpha), \alpha \right) ,
\]
which is an unstable focus for all $\alpha \in [0,y_0)$. Notice that $p(\alpha) \to (-\mu_s,0) \in \Sigma$ as $\alpha \to 0^+$, i.e. there is a boundary collision in the limit $\alpha \to 0^+$.

\begin{remark}
	\remlab{symmetry}
	System \eqref{PWS_stick-slip} also has a boundary collision $p(\alpha) \to (0, \mu_s) \in \Sigma$ as $\alpha \to 0^-$. Since both \eqref{stick-slip} and its PWS counterpart \eqref{PWS_stick-slip} exhibit the symmetry
	\begin{equation}
	\eqlab{symmetry}
	(u, y, \alpha) \leftrightarrow (-u,-y,-\alpha) ,
	\end{equation}
	however, the dynamics for $\alpha \leq 0$ are easily inferred from knowledge of the dynamics for $\alpha \geq 0$ after applying the symmetry \eqref{symmetry}. Therefore, we need only consider the case $\alpha \geq 0$.
\end{remark}

\begin{lemma}
	Consider the PWS system \eqref{PWS_stick-slip}. The switching manifold $\Sigma$ is given by the union
	\[
	\Sigma = \Sigma_{cr}^- \cup F^- \cup \Sigma_{sl} \cup F^+ \cup \Sigma_{cr}^+ ,
	\]
	where
	$\Sigma_{cr}^\pm$ ($\Sigma_{sl}$) are crossing (sliding) submanifolds given by
	\begin{equation}
	\eqlab{ss_manifolds}
	\Sigma_{cr}^- = \left\{(x,0) : x < -\mu_s \right\} , \ \ \Sigma_{sl} = \left\{(x,0) : x \in (-\mu_s, \mu_s) \right\} , \ \ \Sigma_{cr}^+ = \left\{(x,0) : x > -\mu_s \right\} ,
	\end{equation}
	and $F^{\pm} : (\pm \mu_s, 0)$. The point $F^-$ is a visible (invisible) fold for $\alpha > 0$ ($\alpha < 0$), and $F^+$ is a visible (invisible) fold for $\alpha < 0$ ($\alpha > 0$).
\end{lemma}

\W{
	\bpr
		The switching manifold $\Sigma$ decomposes into crossing and sliding submanifolds depending on the sign of
		\begin{equation}
		\eqlab{ss_crossing}
		\left(Z^+f_\Sigma(x,0)\right) \left(Z^-f_\Sigma(x,0)\right) = x^2 - \mu_s^2 ,
		\end{equation}
		\SJ{which is strictly positive (negative) for the crossing (sliding) submanifolds $\Sigma_{cr}^\pm$ ($\Sigma_{sl}$) in \eqref{ss_manifolds}. Tangencies occur when the quantity \eqref{ss_crossing} is zero, i.e.~at $F^{\pm}$.}
		Checking the tangency properties, we have that
		\[
		Z^+(Z^+ f_\Sigma)(F^-) = \alpha , \qquad Z^-(Z^- f_\Sigma)(F^+) = \alpha .
		\]
		Hence by \defnref{fold}, $F^-$ is a visible (invisible) fold for $\alpha > 0$ ($\alpha < 0$), while $F^+$ is a visible (invisible) fold for $\alpha < 0$ ($\alpha > 0$); see \figref{ss}. 
	\epr
}
The Filippov vector field on $\Sigma_{sl}$ is given in the $x-$coordinate chart by
\begin{equation}
\eqlab{ss_Filippov}
\dot x = Z_{sl}(x,\alpha) = - \alpha , \qquad (x,0) \in \Sigma_{sl}.
\end{equation}

\begin{figure}[t!]
	\centering
	\subfigure[]{\includegraphics[width=.45\textwidth]{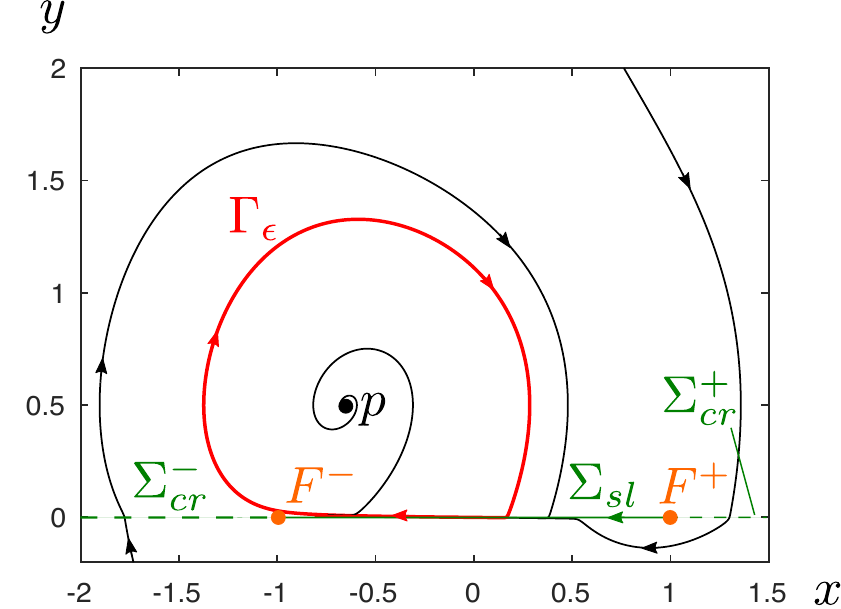}}
	\hspace{0.25cm}
	\subfigure[]{\includegraphics[width=.45\textwidth]{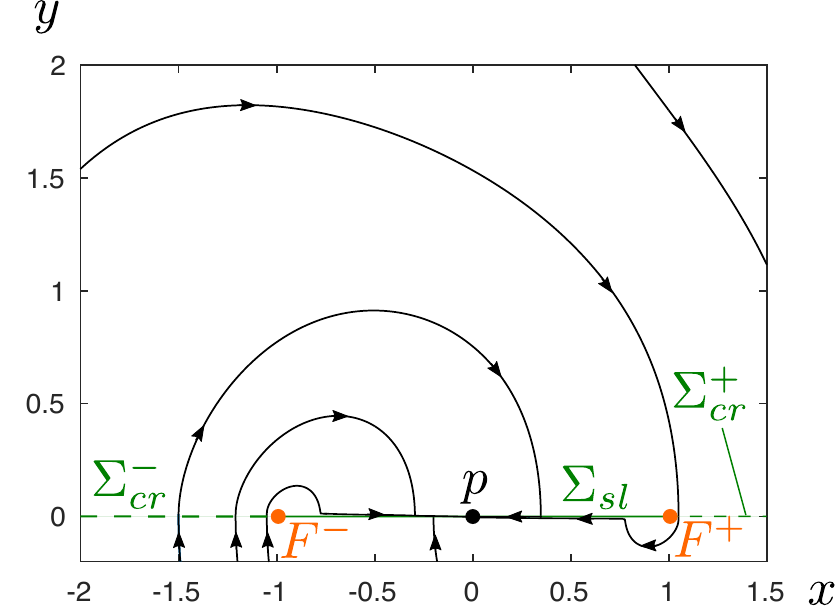}}
	\caption{Phase portraits for system \eqref{stick-slip} with the regularisation function from \figref{reg_fig} with $\epsilon = 10^{-2}$ and \SJ{a polynomial characteristic $\mu_+(y) = \mu_s - \frac{3(\mu_s - \mu_m)}{2 v_m} y + \frac{(\mu_s - \mu_m)}{2 v_m^3} y^3$ appearing in,  e.g.~\cite{Ibrahim1994,Panovko1965,Thomsen2003}, with parameter values $\mu_s = 1, \ \mu_m = 0.5$ and $v_m = 1$.} 
		In (a): $\alpha = 0.5$. In this case $F^-$ ($F^+$) is a visible (invisible) fold, $p$ is an unstable focus, and a relaxation oscillation $\Gamma_\epsilon$ (shown in red) is identified in the regularised problem. In (b): $\alpha = 0$. Here the equilibrium $p:(0,0)$ is a stable node and constitutes the only attractor within the ball $B(0,\mu_s)$. The system exhibits the `creep'  phenomenon: flow along a locally invariant manifold $C_\epsilon$ toward $p$ on a slow timescale $\tilde \tau = \mathcal O(\epsilon)$.}
	\figlab{ss}
\end{figure}

PWS cycles can be constructed for all $\alpha > 0$ sufficiently small, since the unique $Z^+$ orbit intersecting $F^- \in \Sigma$ (where we consider the extension of $Z^+$ on $\{y \geq 0\}$) has a first return, or `drop point' $(x_d,0) \in \Sigma_{sl}$. Hence we can define $\Gamma = \Gamma^1 \cup \Gamma^2$ where
\[
\Gamma^2 = \left\{(x,0) : x \in [-\mu_s,x_d] \right\} \subset \Sigma ,
\]
and $\Gamma^1 \subset \{y \geq 0\}$ is identified with the segment of the $Z^+$ orbit connecting $F^-$ and $(x_d,0)$.

\

Notice that the Filippov vector field \eqref{ss_Filippov} is degenerate for belt speed $\alpha = 0$, i.e.~in the collision limit. This degeneracy indicates a violation of the rightmost nondegeneracy in \eqref{bf_conds}, since
\begin{equation}
\eqlab{degen_1}
\frac{\partial Z_{sl}}{\partial x}\bigg|_{(x,0)} \equiv 0 , \qquad (x,0) \in \Sigma_{sl} .
\end{equation}
In the case of system \eqref{stick-slip}, the degeneracy due to \eqref{degen_1} implies that 
the nullclines
\[
N_1(\alpha) = \left\{(x,\alpha) : x \in \mathbb R \right\} , \qquad N_2(\alpha,\epsilon) = \left\{(- \mu (y, \phi(y \epsilon^{-1})), y) : y\in \mathbb R \right\} ,
\]
which have a unique intersection when $\alpha \neq 0$, $\epsilon \to 0^+$, intersect along the entire line segment $F^- \cup \Sigma_{sl} \cup F^+$ when $\alpha = 0$, $\epsilon \to 0^+$. This leads to the so-called {\em creep phenomenon}, in which flow on a slow timescale arises for belt speeds $\alpha \approx 0$. \SJ{The fact that the nullcline $N_2(\alpha,\epsilon)$ is a graph for all $0 < \epsilon \ll 1$ implies the existence of a \textit{unique} equilibrium in the regularised problem \eqref{stick-slip} even when $\alpha = 0$, i.e. when the conveyor belt is turned off, despite the fact that the Filippov vector field is degenerate\SJnew{; see again \figref{ss}b}. Hence the creep phenomenon is unique to the regularised system \eqref{stick-slip}, and \textit{not} described by the PWS system defined jointly by \eqref{PWS_stick-slip} and \eqref{ss_Filippov} often considered in the literature.}

\begin{lemma}
	\lemmalab{stick-slip}
	\SJ{The PWS system \eqref{PWS_stick-slip} obtained from \eqref{stick-slip} in the singular limit $\epsilon \to 0^+$ has a degenerate BF bifurcation at $F^-$ in the collision limit $\alpha \to 0^+$, in the sense that all conditions in \defnref{bf} are satisfied except that $(\partial {Z_{sl}} / \partial x)|_{(F^-,0)} = 0$.}
\end{lemma}

\bpr
	\SJ{Direct calculations verify the remaining conditions in \defnref{bf}.} 
\epr

We now state our main results for system \eqref{stick-slip}, which describes the dynamics in scaling regimes (S1) and  (S2), as well as an additional regime (S0) relating to the creep dynamics which does not show up in the unfolding of the nondegenerate BF bifurcations.


\begin{theorem}
	\thmlab{ss_thm}
	Consider system \eqref{stick-slip} with $\epsilon \in (0,\epsilon_0)$. Then the following holds for $\epsilon_0 > 0$ sufficiently small:
	\begin{enumerate}
		\item[(i)] There exists a unique equilibrium $p_\epsilon(\alpha)$ such that $p_\epsilon(\alpha) \to F^-$ as $\epsilon \to 0^+$ and $\alpha \to 0^+$, which undergoes a supercritical Andronov-Hopf bifurcation for \note{
		\begin{equation}
		\eqlab{ss_ah_curve}
		\alpha_{ah} = \left(\frac{2 \beta \mu_s k}{-\mu_+'(0)}\right)^{1/(1+k)} \epsilon^{k/(k+1)} + o \left(\epsilon^{k/(k+1)}\right) ,
		\end{equation}
		}where $\beta$ is given by \eqref{beta}. 
		
		\item[(ii)] There exists constants $\alpha_\pm > 0$ such that for all $\alpha \in (\alpha_-, \alpha_+)$, 
		there exists a relaxation oscillation $\Gamma_\epsilon$ which is strongly attracting, unique in a tubular neighbourhood $\mathcal T$ of a PWS cycle $\Gamma = \Gamma^1 \cup \Gamma^2$, and $\mathcal O(\epsilon^{2k/(2k+1)})-$close to $\Gamma$ in Hausdorff distance. If additionally $\mu'(y) < 0$ for all $y \geq 0$ (i.e. $y_0 = \infty$), then there are no other limit cycles in an arbitrarily large ball $B(0,M)$.
		
		\item[(iii)]
		The relaxation oscillations described in (ii) connect to regular cycles in (S1), i.e. there exists $K > 0$ such that system \eqref{stick-slip} has a smooth parameterised family of stable limit cycles
		\begin{equation}
		\eqlab{ss_cycles_pos}
		\alpha \mapsto \Gamma_+(\alpha,\epsilon) , \qquad \alpha \in \left(K \epsilon^{k/(k+1)}, \alpha_+ \right) ,
		\end{equation}
		which contains (but is exhausted by) the relaxation oscillations for $\alpha \in (\alpha_-,\alpha_+)$.
		
		\item[(iv)] There exists an additional scaling regime given by
		\[
		(S0) : \alpha = \mathcal O(\epsilon) ,
		\]
		in which $p_\epsilon$ is the only attractor within the ball $B(0,\mu_s)$ and the system exhibits creep dynamics, i.e. flow along a locally invariant slow manifold $C_\epsilon$ towards $p_\epsilon$ on a slow timescale $\tilde \tau = \mathcal O(\epsilon)$. Here
		\begin{equation}
		\eqlab{ss_slow_manifold}
		C_\epsilon = \{(x,\mathcal O(\epsilon)) : x \in [-\mu_s + \kappa, \mu_s - \kappa] \} ,
		\end{equation}
		where the constant $\kappa > 0$ can be chosen arbitrarily small as $\epsilon \to 0^+$.
	\end{enumerate}
\end{theorem}

\bpr
	In order to prove statement (i), we look in scaling regime (S1). We can gain explicit information about the dynamics in regime (S1) \SJ{by considering system \eqref{stick-slip} and applying} the coordinate transformation
	\note{
		\begin{equation}
		\eqlab{coord_change_ss}
		\left(x_1, \nu_1, \rho_1, \hat \alpha \right) \in \mathbb R \times \mathbb R_+^2 \times \mathbb R \mapsto
		\begin{cases}
		x = -\mu_s + \nu_1^{2k(1+k)} \rho_1^{k(1+k)} x_1 , \\
		y = \nu_1^{2k(1+k)} \rho_1^{2k(1+k)} , \\
		\epsilon = \nu_1^{2(1+k)^2} \rho_1^{(1+k)(1+2k)} , \\
		\alpha = \hat \alpha \nu_1^{2k(1+k)} \rho_1^{k(1+2k)} ,
		\end{cases}
		\end{equation}
	}
	which is obtained by combining \eqref{coord_change} with a coordinate translation which moves the fold $F^-$ to the origin in \eqref{stick-slip}. This yields the system
	\note{
		\begin{equation}
		\eqlab{stick-slip_S1}
		\begin{split}
		x_1' &= \rho_1^{k(1+2k)} \left(\rho_1^k - \hat \alpha \right) - k x_1 \left(x_1 - 2 \mu_s \beta + \mu_+'(0) \rho_1^{k(1+k)} \right) , \\
		\rho_1' &= - \frac{1}{k} \rho_1 \left(x_1 - 2 \mu_s \beta + \mu_+'(0) \rho_1^{k(1+k)} \right) ,
		\end{split}
		\end{equation}
	}
	after a suitable desingularisation (division by $\rho_1^{k(k+1)}$) and restriction to the invariant subspace $\{\nu_1 = 0\}$. System \eqref{stick-slip_S1} has an Andronov-Hopf bifurcation for $\hat \alpha = \hat \alpha_{ah}$, where
	\note{
		\[
		\hat \alpha_{ah} = \left(\frac{2 \beta \mu_s k}{-\mu_+'(0)}\right)^{1/(1+k)} ,
		\]
	}
	which is supercritical with first Lyapunov coefficient
	\note{
		\[
		l_1 =  - \frac{\beta  k^3 (1+k) (2+k)}{16}
		\left(\frac{\beta  k}{\tau }\right)^{-2/k(1+k)}  < 0 .
		\]
	}
	The \SJ{leading order estimate} in \eqref{ss_ah_curve} \SJ{follows from the relation} $\alpha = \hat \alpha \epsilon^{k/(k+1)}$.
	
	\
	
	Existence of a locally unique relaxation oscillation $\Gamma_\epsilon$ converging to the PWS cycle $\Gamma$ in Hausdorff distance as $\epsilon \to 0^+$ follows from 
	the proof for \thmref{thm_ro} \SJnew{(\appref{proof_of_thm_ro})}, which can be adapted for the case $\gamma = 0$; see also  \secref{degenerate_cases}.
	
	Now assume that $\mu_+'(y) < 0$ for all $y \geq 0$. It follows from the contraction mapping argument in the proof for \thmref{thm_ro} that $\Gamma_\epsilon$ is unique in an open tubular neighbourhood $\mathcal T_\epsilon$ of width $\mathcal O(\epsilon^{2k/(2k+1)})$, and we can show that no limit cycles can exist on domains $\mathcal D_\epsilon = \{(x,y) \in \mathbb R^2 : y > \epsilon^k \tilde c\}$ with constant $\tilde c > 0$, since
	\[
	\langle \nabla , Z(x,y,\alpha,\epsilon) \rangle = - \mu_+(y) + \mathcal O \left(\epsilon^k\right) \leq 0 
	\]
	for $\epsilon_0 > 0$ sufficiently small, for all $y \geq 0$. Since the intersection $\mathcal T_\epsilon \cap \mathcal D_\epsilon$ is non-empty and open, the result follows.
	
	\
	
	Existence of the smooth family of stable cycles \eqref{ss_cycles_pos} containing regular oscillations in (S1) and relaxation oscillations in (S2) follows from a generalisation of \thmref{thm_connection} to the case $\gamma = 0$; see again \secref{degenerate_cases}. 
	
	\
	
	It remains to prove the statement (iv). In order to study the dynamics in (S0), we consider the rescaled problem with
	\[
	y = \epsilon Y, \qquad \alpha = \epsilon a .
	\]
	Rewriting system \eqref{stick-slip} and transforming to a fast timescale $T = t / \epsilon$ yields
	\begin{equation}
	\eqlab{creep}
	\begin{split}
	x' &= \epsilon^2 \left(Y - a \right) , \\
	Y' &= - x + \mu_s\left(1 - 2 \phi(Y) \right) - \epsilon \mu_+'(0) Y - \epsilon^2 \varphi(Y,\epsilon) ,
	\end{split}
	\end{equation}
	where $(\cdot)'$ denotes differentiation with respect to $T$, and $\varphi(Y,\epsilon) = \mathcal O((Y\epsilon)^3)$. System \eqref{creep} is a slow-fast system in the standard form \eqref{stnd}, with a normally hyperbolic and attracting critical manifold
	\[
	S = \left\{\left(\mu_s \left(1 - 2 \phi(Y)\right) , Y \right) : Y \in \mathbb R \right\} ,
	\]
	and an associated non-trivial eigenvalue $\lambda(Y) = - 2 \mu_s \phi'(Y) < 0$. 
	Since $x' = \mathcal O(\epsilon^2)$ however, there is no reduced flow on $S$, and one must look to higher orders to observe a (infra)-slow flow. Standard matching arguments lead to the following series expression for the nearby slow manifold $S_\epsilon$ as $\epsilon \to 0^+$:
	\[
	S_\epsilon : x = \mu_s\left(1 - 2 \phi(Y) \right) - \epsilon a Y - \frac{\epsilon^2}{2} \left(\mu_+''(0) - \frac{Y - a}{\phi'(Y)} \right) + \mathcal O(\epsilon^3) .
	\]
	Restricting system \eqref{creep} to $S_\epsilon$ and rewriting it on a new, \textit{infra-slow} timescale $\tilde \tau = \epsilon^2 T = \epsilon t$, we obtain the following for $\epsilon = 0$ in the $Y-$coordinate chart:
	\[
	\frac{d Y}{d \tilde \tau} = - \left(\frac{Y - a}{2 \phi'(Y)}\right) ,
	\]
	which has a unique equilibrium $p : (\mu_s(1-2\phi(Y)), a) \in S$. Asymptotic stability of $p$ follows by a direct calculation, and 
	the expression for $C_\epsilon$ in \eqref{ss_slow_manifold} is obtained by taking a compact submanifold of $S_\epsilon$ and writing it in $(x,y)-$coordinates.
\epr




\section{Discussion and conclusion}
\seclab{summary_and_conclusion}

In this article we considered \textit{smooth} systems \eqref{main} with a BF bifurcation in the nonsmooth limit $\epsilon \to 0$. By considering general regularisations of the kind defined by \assumptionref{ass1} and \assumptionref{ass1b}, rigorous results were obtained for the case $0 < \epsilon \ll 1$ relevant in applications. After deriving a suitable local normal form in \propref{prop_normal_form}, we adopted a blow-up approach in order to identify the scaling regime (S1) relevant for the unfolding of the BF point. Importantly for applications, the scaling regime (S1) scales in an $(\epsilon,k)-$dependent fashion in accordance with $\mu = \mathcal O(\epsilon^{k/(k+1)})$, where \SJ{the transition coefficient} $k \in \mathbb N_+$ is an algebraic decay rate associated with the regularisation function $\phi$.~
The bifurcation structure within (S1) is determined by the type of BF bifurcation in the nonsmooth limit\SJnew{, and is \textit{qualitatively independent} of the choice of regularisation function $\phi$ within the  class of regularisations defined by \assumptionref{ass1b} and \assumptionref{ass3}}; see \thmref{thm_bifurcation}. In case BF$_3$, a supercritical Andronov-Hopf bifurcation is identified as the mechanism for the onset of stable oscillations, while in cases BF$_i$, $i=1,2$, a codimension-2 \SJnew{supercritical} Bogdanov-Takens point organises the \SJnew{local} bifurcation diagram; this is summarised in \figref{bfb_bif_diagram}. In cases BF$_1$ and BF$_3$ we were also able to prove persistence of PWS cycles as relaxation oscillations in the smooth system $0 < \epsilon \ll 1$, within a second scaling regime (S2) with values $\mu = \mathcal O(1)$ bounded away from the bifurcations; see \thmref{thm_ro}. Finally in \thmref{thm_connection}, we showed that in case BF$_3$ there exists a smooth parameterised family of stable limit cycles connecting regular cycles within (S1) to relaxation oscillations within (S2). Our results were finally applied for two applications (M1) and (M2) in \secref{bf_bifurcation_in_application}.

While \thmref{thm_bifurcation}, \thmref{thm_ro} and \thmref{thm_connection} are derived from the specific normal form \eqref{normal_form} and thus qualitative in nature, our blow-up analysis also lead to the identification of the nontrivial coordinate transformation \eqref{coord_change} in \lemmaref{lem_desing}, which can be applied directly in applications in order to obtain a desingularised system which governs the dynamics within the scaling regime (S1). By studying this desingularised problem, one can obtain \textit{quantitative} information regarding the bifurcation structure of the given application. This is outlined in \remref{guase_desing} in the context of the Gause model (M1), and \SJnew{undertaken} explicitly the application (M2) in the proof for \thmref{ss_thm} in \secref{bf_bifurcation_in_application}.

\

In the remainder of this section, we briefly discuss a number of related topics that are not considered in detail in this work.

\subsection{Other BE bifurcations}
\seclab{other_cases}

The applications (M1) and (M2) also feature \textit{boundary-node} BN bifurcations in the nonsmooth limit $\epsilon \to 0$ within physically relevant parameter regimes; for larger handling times $h \in (h_\ast , m / \tilde e)$ in application (M1), and for sufficiently large and negative gradients $\mu_+'(y) \in (-\infty, -2)$ in application (M2). In \cite{Kristiansen2019d}, the authors successfully apply blow-up techniques in order to analyse smooth systems \eqref{main} featuring a BN bifurcation in the nonsmooth limit in the context of a model for substrate-depletion oscillation in \cite{Kristiansen2019d}. \note{In a sequel paper \cite{Jelbart2021} we extend the methods presented herein in order to provide a complete description for the unfolding of hyperbolic boundary-equilibrium bifurcations in the general case where the incident equilibrium a focus, node or saddle.} 

\subsection{Degenerate cases}
\seclab{degenerate_cases}

Our analysis also provides some insight into the \SJnew{degenerate} BF bifurcations which lie `between' the cases BF$_i$, $i=1,2,3$.

\SJ{
	Cases BF$_1$ and BF$_2$ are separated by a degenerate case known as the `homoclinic boundary focus' (HBF) bifurcation; see e.g. \cite{Dercole2011} for an analysis of this case in the PWS setting.
	For regularised systems \eqref{main} featuring a HBF bifurcation in the nonsmooth limit, one can show that PWS homoclinic cycles
	persist in the scaling regime (S2) as singular homoclinic cycles after blow-up, which exist along a curve in $(\mu,\gamma)-$parameter space with $\mu \in \mathcal I_+$, $\gamma > 0$. It seems reasonable that the homoclinic curve $\hat \mu_{hom}(\gamma)$ in $(\hat \mu, \gamma)-$space identified in \lemmaref{lem_bifurcations_1} (see \figref{bfb_bif_diagram}) corresponds to the extension of this curve into the scaling regime (S1). Proving this would, however, require a global understanding of the curve $\hat \mu_{hom}(\gamma)$ in the limit $\hat \mu \to \infty$.
	
	Cases BF$_1$ and BF$_3$ are separated by a degenerate BF bifurcation with $\gamma = 0$, which  
	appeared in the application (M2). Since the normal forms \eqref{normal_form} and \eqref{PWS_normal_form} do not require that $\gamma \neq 0$, 
	many of our results generalise to the case $\gamma = 0$. In particular, one can identify a supercritical Andronov-Hopf bifurcation scaling regime (S1) for $\mu =\mu_{ah}(\gamma,\epsilon) = \mu_{ah}(0,\epsilon)$,
	PWS cycles can be shown to persist as stable relaxation oscillations in the scaling regime (S2), and a connection of the type described in \thmref{thm_connection} can be shown to persist by extending the proof in \appref{proof_of_thm_connection} to all values $\gamma \leq 0$. The degeneracy due to $\gamma = 0$ becomes more pronounced, however, when one considers the matching problem associated with the limit in \eqref{matching_limit_neg}, due to the potential for creep dynamics on an additional slow timescale $\tilde \tau = \epsilon t$.
}

\subsection{The comparison with canards}
\seclab{canards}

The existence of the entire family of singular cycles \SJnew{defined by \eqref{family_cycs} and shown in \figref{blowup_connection}b (see also \figref{complete_BF3}b)} is reminiscent of the well-known \textit{canard explosion} phenomena known to occur in slow-fast systems \cite{Wechselberger2015,Dum1996,Gucwa2009,Krupa2001b}, in which small oscillations 
born in a singular Andronov-Hopf bifurcation grow smoothly to $\mathcal O(1)$-amplitude oscillations of relaxation-type under an exponentially small variation of an additional system parameter. The existence of so-called `canards' -- solutions which lie on the intersection of attracting and repelling slow manifolds -- provides the key underlying mechanism for the canard explosion phenomena (see e.g. \cite{Krupa2001a,Kuehn2015}). Interestingly however, `canard-like' mechanisms leading to `explosions-like' dynamics have also been observed in systems in which one of the intersecting invariant manifolds is not a slow manifold;
see \cite{Kristiansen2019d}, where such an `explosion-like' phenomena occurs in a model for substrate-depletion oscillation with no attracting slow manifold. The situation in system \eqref{normal_form} is different again: here there is neither
a repelling slow manifold, nor an obvious invariant manifold which `replaces' the repelling slow manifold.
Thus, there can be \textit{no canards} and hence \textit{no canard explosion} in the sense described above.

\section*{Acknowledgements}

MW was supported by the Australian Research Council DP180103022 grant. SJ would like to thank the Technical University of Denmark (DTU) and the second author for hospitality during their stay, without which this work may not have been possible.

\newpage

\appendix
\renewcommand{\thesection}{\Alph{section}}

\section*{Appendix}

\SJnew{The appendices are structured as follows: \appref{proof_of_normal_form} contains the proof for the normal form result \propref{prop_normal_form}. \appref{S1S2_proofs} includes the detailed analysis necessary to prove \thmref{thm_bifurcation} and \thmref{thm_ro}. A proof of \thmref{thm_connection} is finally given in \appref{proof_of_thm_connection}, once the detailed blow-up analysis has been presented in \appref{blow-up_of_Sigma_and_Q} and \appref{blow-up_of_Qbfb}.}

\section{Normal form: proof of \propref{prop_normal_form}}
\applab{proof_of_normal_form}

In this section we outline a proof for \propref{prop_normal_form}, referring the reader to \cite{Jelbart2020b} for further details.

\

We assume without loss of generality that system \eqref{PWS} has a BF bifurcation in system \eqref{PWS} at $z_{bf}=(0,0)$ for $\alpha_{bf}=0$. It follows from the second expression in \eqref{bf_conds} that $Z_2^-(z, \alpha) \neq 0$ in sufficiently small neighbourhood $\mathcal U = U \times I_\alpha \ni (0,0,0)$, for suitably small $U \subset \mathbb R^2$ and $I_\alpha = (-\alpha_0,\alpha_0)$. Hence, orbit segments of $Z^-(z,\alpha)$ in $\mathcal U$ can be described as level sets $L(x,y,\alpha) = const.$ where $L(x,y,\alpha)$ is a smooth function satisfying $(\partial L /\partial x)|_{(x,y,\alpha)} \neq 0$ for all $(x,y,\alpha) \in \mathcal U$. Defining a new coordinate $u = L(x,y,\alpha)$ with locally defined inverse $x=M(u,y,\alpha)$, one obtains the following system from \eqref{main} after dividing the right-hand-side of the transformed system by the locally positive factor $|Z^-_2((M(u,y,\alpha), y), \alpha)|$:
\begin{equation}
\eqlab{rectified_system2}
\begin{split}
\dot u &= \phi\left(y \epsilon^{-1} \right) \left\{ \left[\left(\frac{\partial M}{\partial u}\right)^{-1} \frac{\partial L}{\partial u} + \left(\frac{\partial M}{\partial y} \right)^{-1} \frac{\partial L}{\partial y} \right] Z_1^+ + \frac{\partial L}{\partial y} Z_2^+ \right\} / |Z_2^-|,  \\
\dot y &= \pm 1 + \phi\left(y \epsilon^{-1} \right) \left\{ \mp 1 + \frac{Z_2^+}{|Z_2^-|} \right\} ,
\end{split}
\end{equation}
where all quantities are evaluated in terms of $(u,y,\alpha,\epsilon)$. Note that the `$\pm$' sign is chosen to agree with $\sgn(Z_2^-((0,0),0))$, and by a slight abuse of notation the $(\dot{\ })$ notation now denotes differentiation with respect to the new (rescaled) time.

Since $u = L(x,y,\alpha)$ is a diffeomorphism, the nondegeneracy conditions \eqref{bf_conds} imply a BF bifurcation in system \eqref{rectified_system2} at $(u,y)=(0,0)$ when $\alpha = 0$. Hence, we can express \eqref{rectified_system2} as an expansion
\begin{equation}
\eqlab{expanded_sys}
\begin{split}
\dot u &= \phi\left(y \epsilon^{-1} \right) \left[ s_1(\alpha) + a(\alpha) u + b(\alpha) y + \varphi_1\left(u, y, \alpha \right) \right] ,  \\
\dot y &= \pm 1 + \phi\left(y \epsilon^{-1} \right) \left[ \mp1 + s_2(\alpha) + c(\alpha) u + d(\alpha) y + \varphi_2\left(u, y, \alpha  \right) \right] ,
\end{split}
\end{equation}
where $a(\alpha), b(\alpha), c(\alpha), d(\alpha)$ are smooth functions of $\alpha$ such that $a(0) = a, \ b(0) = b, \ c(0) = c, \ d(0) = d$ are constant, and $s_i(\alpha)$, $\varphi_i(y,y,\alpha)$, $i=1,2,$ are smooth functions satisfying $s_i(0) = 0, \ i = 1,2$, $s_2'(0) > 0$, and $||\varphi_i(u,y,\alpha)|| = \mathcal O (||(u,y,\alpha)||^2)$. 
The presence of a hyperbolic focus in system \eqref{expanded_sys} implies trace, determinant and discriminant conditions
\[
\tau := a + d \neq 0 , \qquad \delta := ad-bc > 0, \qquad \tau^2 - 4\delta < 0 ,
\]
with either $b \neq 0 $ or $c \neq 0$ (or both). Without loss of generality, we assume that $c \neq 0$
, and introduce the linear coordinate transformation
\[
v = \sgn(c) \left( c (u + w(\alpha)) + d y \right) , \qquad u = \frac{1}{c} \left( \sgn(c) v - d y\right) - w(\alpha) ,
\]
where $w : I_\alpha \to \mathbb R$ is some smooth yet to be determined function of $\alpha$ satisfying $w(0) = 0$. Applying the coordinate transformation and imposing the requirement that constant terms must vanish in the resulting equation for $\dot y$, we obtain the system
\begin{equation}
\begin{split}
\sgn(c) \dot v &= d + \phi\left(y \epsilon^{-1} \right) \left(-d + \mu + \tau \sgn(c) v - \delta y + \theta_1(\sgn(c) v,y,\alpha) \right) , \\
\dot y &= \pm1 + \phi\left(y \epsilon^{-1} \right) \left(\mp1 + \sgn(c) v + \theta_2(\sgn(c) v,y,\alpha) \right) ,
\end{split}
\end{equation}
after setting $\mu(\alpha) := c s_1(\alpha) + d s_2(\alpha) - c \tau w(\alpha) + \mathcal O(\alpha^2)$, and $w(\alpha) = (s_2'(0) / c) \alpha + \mathcal O(\alpha^2)$. Since $\mu'(0) \neq 0$ follows by the determinant condition in \eqref{bf_conds} applied to \eqref{expanded_sys}, the inverse function theorem implies existence of a unique inverse function $\alpha = \tilde \alpha(\mu)$ such that $\tilde \alpha(0) = 0$ and $\tilde \alpha'(0) \neq 0$. By considering $\mu$ as a new parameter we may write $\theta_i(v,y,\mu) := \tilde \theta_i(v,y, \tilde \alpha(\mu))$, $i = 1,2$, and permitting a slight abuse of notation by writing $x$ in place of $v$, one obtains the desired form \eqref{normal_form} after setting $\gamma = \tau \mp d$ 
and writing down the case $\sgn(c) = + 1$ (locally counter-clockwise rotation; see \remref{normal_forms}). The form of $X_{sl}(x,\mu)$ given in \eqref{Filippov_VF} follows by a direct computation using the expression in \eqref{EqnSlidingVF} with $X^\pm$ in place of $Z^\pm$, where $X^\pm$ are given by \eqref{PWS_normal_form}.
\qed

\section{Proofs for \lemmaref{lem_bifurcations_1} and \thmref{thm_ro}}
\applab{S1S2_proofs}

In this section we prove \lemmaref{lem_bifurcations_1} and \thmref{thm_ro}. We start with \lemmaref{lem_bifurcations_1}, which is sufficient for the proof of \thmref{thm_bifurcation} given in \secref{results_bifurcation}.

\subsection{Proof of \lemmaref{lem_bifurcations_1}}
\applab{proof_of_lem_bifurcation}

We consider system \eqref{desing_prob}. Bifurcating equilibria occur for $\rho_1-$values satisfying the equation
\begin{equation}
\eqlab{eq_eqn_1}
\varphi(\rho_1) = \gamma \beta - \hat \mu \rho_1^{k^2} + \delta \rho_1^{k(1+k)} =0 , \qquad \rho_1 > 0, \ \hat \mu \in \mathbb R ,
\end{equation}
which cannot be solved explicitly in terms of radicals for arbitrary $k \in \mathbb N_+$. We need to show that equation \eqref{eq_eqn_1} has a unique solution when $\gamma \leq 0$, and up to two solutions when $\gamma > 0$. First consider the case $\gamma \leq 0$. Direct calculations show that
\begin{enumerate}
	\item[(i)] $\varphi(\rho_1)$ has a minimum at $\rho_1 = \rho_{1,c} > 0$, and no other extrema on $\rho_1 > 0$;
	\item[(ii)] $\varphi(\rho_1) < 0$ for all $\rho_1 \in (0,\rho_{1,c}]$;
	\item[(iii)] $\lim_{\rho_1 \to \infty}\varphi(\rho_1) = \infty$;
\end{enumerate}
see \figref{number_eqs_1}a. Existence of a unique solution $\rho_{1,\ast} > 0$ satisfying $\varphi(\rho_{1,\ast}) = 0$ follows from (i)-(iii) after an application of the intermediate value theorem.

\begin{figure}[h!]
	\begin{center}
		\subfigure[]{\includegraphics[width=.48\textwidth]{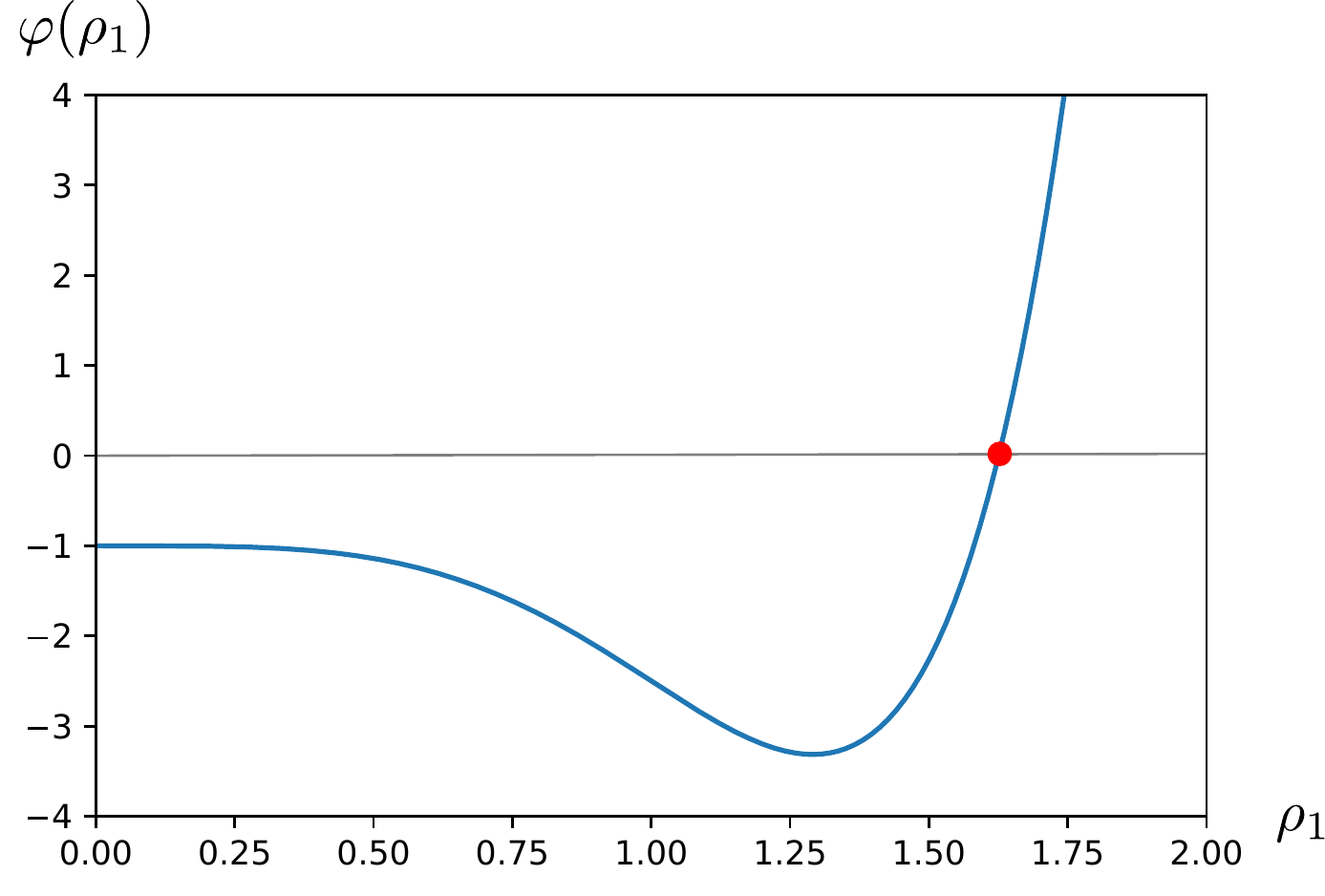}}
		\subfigure[]{\includegraphics[width=.48\textwidth]{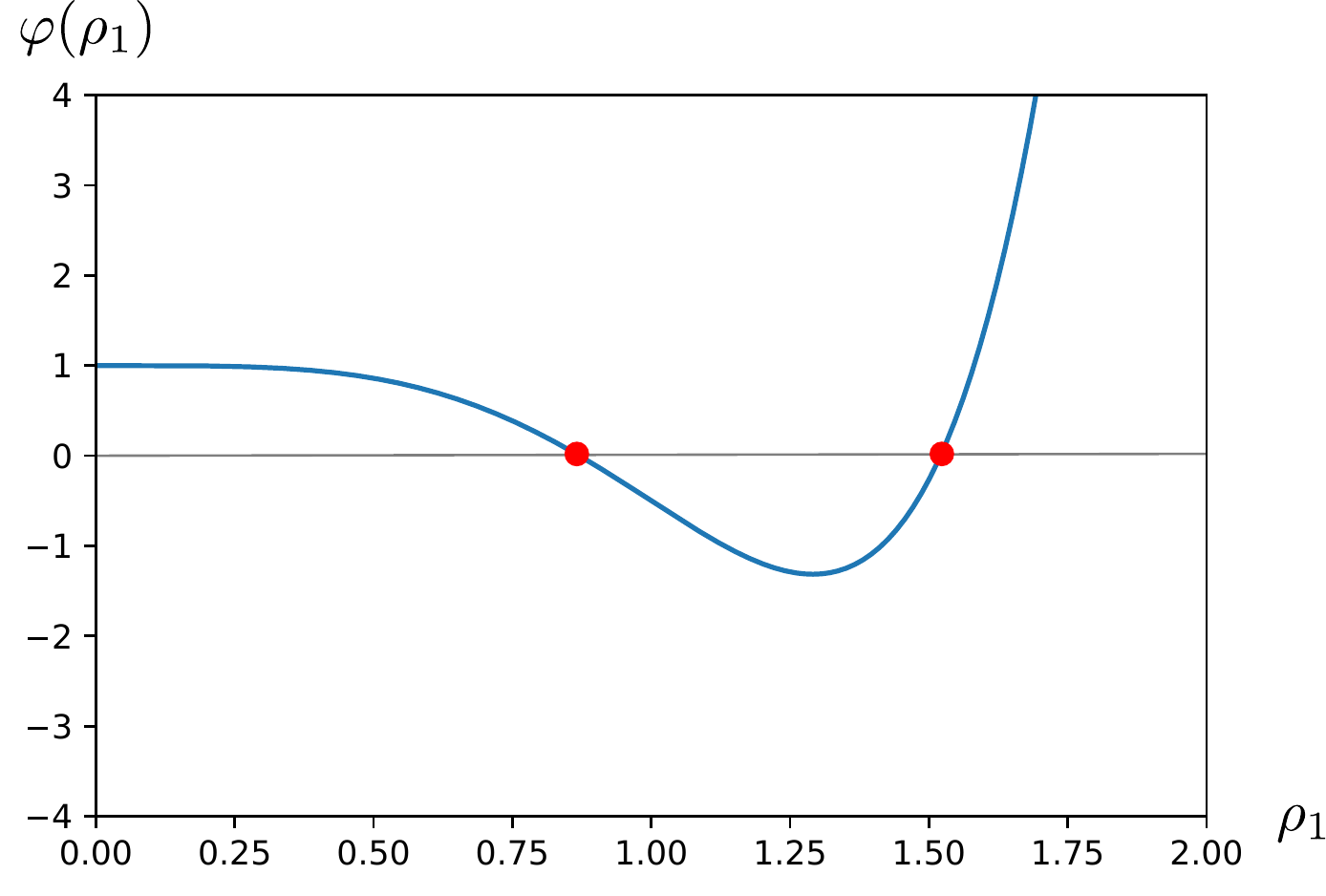}}
		\caption{The curve $\{(\rho_1, \varphi(\rho_1)) : \rho_1 \geq 0\}$ with $\beta = \delta = 1$, $k = 2$ and $\hat \mu = 5/2$. Each $\rho_{1,\ast}$ satisfying $\varphi(\rho_{1,\ast}) = 0$ corresponds to an equilibrium of \eqref{desing_prob}, and the number of possible roots (red dots in the figure) depends on $\sgn(\gamma)$. In (a): $\gamma = -1 < 0$ and there is a unique root. In (b): $\gamma = 1 > 0$ and up to two distinct roots (the case shown here) may exist depending on $\hat \mu$.}
		\figlab{number_eqs_1}
	\end{center}
\end{figure}

Considering the case $\gamma > 0$, one can show that
\begin{enumerate}
	\item[(i)] $\varphi(0) = \beta \gamma > 0$;
	\item[(ii)] $\varphi(\rho_1)$ has a minimum at $\rho_1 = \rho_{1,c}(\hat \mu) > 0$, and no other extrema on $\rho_1 > 0$;
	\item[(iii)] There exists $\hat \mu_{c}$ such that the minimum value satisfies $\varphi(\rho_{1,c}(\hat \mu)) < 0$ for all $\hat \mu > \hat \mu_{c}$;
	\item[(iv)] $\lim_{\rho_1 \to \infty}\varphi(\rho_1) = \infty$;
\end{enumerate}
see \figref{number_eqs_1}b. Hence by choosing $\hat \mu > \hat \mu_{c}$, (i)-(iv) above imply two solutions to $\varphi(\rho_1) = 0$ with $\rho_1 > 0$, as required.

\

We now consider bifurcations. Although \eqref{eq_eqn_1} cannot be solved explicitly, we \textit{can} locate bifurcations by considering suitable intersections in $(\rho_1,\hat \mu, \gamma)-$space.
Calculating the trace associated with the Jacobian $J$ as a function of $\hat \mu, \gamma$ and $\rho_1$, we obtain
\begin{equation}
\eqlab{J_Tr}
\tr J\big|_{x_1 = -\beta} = - k \beta + \tau \rho_1^{k(1+k)} ,
\end{equation}
which after solving for the local Andronov-Hopf condition $\tr J |_{x = -\beta} = 0$ yields a fixed location
\[
p(\hat \mu_{ah}) : \left(x_{1,ah},\rho_{1,ah}\right) = \left(-\beta, \left(\frac{k \beta}{\tau}\right)^{1/{k(1+k})} \right).
\]
Substituting this into the equilibrium condition \eqref{eq_eqn_1} and solving for $\hat \mu$ as a function of $\gamma$ yields the expression for $\hat \mu_{ah}(\gamma)$ in \eqref{hopf_value_1}, and checking the local determinant condition we obtain
\begin{equation}
\eqlab{J_Det}
\det J\big|_{x = - \beta, \hat \mu=\hat \mu_{ah}(\gamma)} = \left(\frac{k \beta}{\tau} \right)^2 \left(\delta - \gamma \tau \right)  > 0 , \qquad \text{for} \qquad \gamma < \frac{\delta}{\tau} ,
\end{equation}
indicating an Andronov-Hopf bifurcation for $\hat \mu = \hat \mu_{ah}(\gamma)$, $\gamma < \delta / \tau$.
We can characterise the asymptotic stability of the bifurcating equilibrium by using the fact that $p(\hat \mu) \in \{x = -\beta\}$ for all $\hat \mu \in \mathbb R$. In particular, the equilibrium occurs for
\[
\hat \mu(\rho_1) = \frac{1}{\rho_1^{k^2}} \left( \gamma \beta + \delta \rho_1^{k(1+k)} \right),
\]
from which it follows that
\[
\frac{\partial \hat \mu}{\partial \rho_1}\bigg|_{\rho_1 = \rho_{1,ah}} = k^2 \beta \left(\frac{\tau}{k \beta} \right)^{(k^2+1)/k(1+k)} \left(\frac{\delta}{\tau} - \gamma \right) > 0 , \qquad \iff \qquad \gamma < \frac{\delta}{\tau}.
\]
Hence, the equilibrium moves up the line $\{x = - \beta\}$ with increasing $\rho_1$ in a neighbourhood of the Andronov-Hopf point. This together with
\[
\frac{\partial }{\partial \rho_1} \tr J \big|_{\rho_1 = \rho_{1,ah}} = (1+k)k^2 \beta \left(\frac{\tau}{k \beta}\right)^{1/{k(1+k)}} > 0 ,
\]
proves the stability claim in assertion (i), i.e. that the equilibrium $p$ transitions from stable to unstable with increasing $\hat \mu$ over $\hat \mu_{ah}(\gamma)$. To determine criticality, a calculation using the formula from \cite[p.211]{Chow1994} (and Mathematica) for the first Lyapunov coefficient $l_1$ yields
\note{
	\[
	l_1 = - \frac{\beta  k^3 (1+k)}{16 (\delta - \gamma  \tau )}
	\left(\frac{\beta  k}{\tau }\right)^{-2/k(1+k)} ((2+k) \delta -\gamma  \tau ) ,
	\]
}implying a supercritical bifurcation for all $k \in \mathbb N_+$, since
\[
\gamma < \frac \delta \tau \qquad \implies \qquad  (2+k) \delta - \gamma \tau > \delta - \gamma \tau > 0 \qquad   \implies \qquad l_1 < 0 .
\]
In order to conclude the statement (i), it remains to show that there are no other local bifurcations when $\gamma < 0$. This fact will follow from the preceding arguments together with assertion (ii); we turn to the proof of (ii) and (iii) now.

\

Restricting to $\{x_1 = - \beta\}$ and calculating the Jacobian determinant as a function of $\hat \mu, \gamma$ and $\rho_1$ yields the saddle-node condition
\begin{equation}
\eqlab{sn_cond_1}
\det J\big|_{x = - \beta} = (1+k) \beta \gamma - (1 + 2 k) \hat \mu \rho_1^{k^2} + 2 (1 + k) \delta \rho^{k(1+k)} = 0.
\end{equation}
This equation can be solved for $\rho_1^{k^2}$, which can in turn be substituted into \eqref{eq_eqn_1} in order to obtain
\begin{equation}
\eqlab{sn_location_1}
p_{sn}(\hat \mu_{sn} (\gamma)) : (x_{1,sn}, \rho_{1,sn}) = \left(-\beta, \left(\frac{k \beta \gamma}{\delta} \right)^{1/k(1+k)} \right) .
\end{equation}
Substituting \eqref{sn_location_1} into \eqref{sn_cond_1} and solving for $\hat \mu$ in terms of $\gamma$ yields the desired expression for $\hat \mu_{sn}(\gamma)$ in \eqref{sn_curve_lem_1}. The fact that the equilibria emerging from the saddle-node bifurcation exist only for $\hat \mu > \hat \mu_{sn}(\gamma)$ follows from the fact that the curve $\hat \mu_{ah}(\gamma)$ is bounded above $\hat \mu_{sn}(\gamma)$ in $(\hat \mu, \gamma)-$parameter space, together with the existence of a Bogdanov-Takens point (which ensures that the equilibrium which undergoes Andronov-Hopf bifurcation is born in the saddle-node bifurcation). The existence and location of the Bogdanov-Takens point given in \eqref{bt_point_lem_1} can be verified by simultaneously setting both trace and determinant expressions in \eqref{J_Tr} and \eqref{J_Det} (or equivalently \eqref{sn_cond_1}) to zero, and solving for $\rho_1, \hat \mu, \gamma$. This proves assertions (ii)-(iii). To conclude the proof of (i), notice that $p_{sn}(\hat \mu(\gamma)) \to (-\beta,0)$ as $\gamma \to 0^+$, implying that saddle-node bifurcation in the invariant regime $\{\nu_1=0\} \cap \mathcal A$ is only possible for $\gamma \geq 0$. Hence there can be no local bifurcations when $\gamma < 0$ in addition to the Andronov-Hopf bifurcation described above, as required.


\


It remains to prove statement (iv). Here we note that existence of the Bogdanov-Takens bifurcation identified in (iii) implies local existence of a curve $\hat \mu_{hom}(\gamma)$ that is quadratically tangent to $\hat \mu_{ah}(\gamma)$ at $(\hat \mu_{bt}, \gamma_{bt} )$, along which a homoclinic-to-saddle connection exists \cite{Kuznetsov2013}. The local upper bound $\gamma_{bt}=\delta / \tau$ on the domain for $\hat \mu_{hom}(\gamma)$ is determined by the upper bound on the domain for $\hat \mu_{ah}(\gamma)$ (i.e. they are the same). Finally, it follows from the fact that there is only one equilibrium $p$ within $\{\rho_1 > 0\}$ that a homoclinic-to-saddle connection cannot exist for $\gamma < 0$, since such a connection requires at least two equilibria within $\{\rho_1 > 0\}$.
\qed

\subsection{Proof of \thmref{thm_ro}}
\applab{proof_of_thm_ro}

In the following we prove \thmref{thm_ro}. We consider system \eqref{normal_form} with $\mu \in \mathcal I_- \cup \mathcal I_+$.

\

Existence of the unstable focus $p_\epsilon$ in statement (i) when $\mu \in \mathcal I_+$ follows from the existence of the unstable focus $p$ in \lemmaref{PWS_lemma} (ii) and regular perturbation theory.

Now consider the dynamics within the switching layer $\{y = \mathcal O(\epsilon)\}$. For general systems \eqref{main} satisfying \assumptionref{ass1} and \assumptionref{ass1b}, we are given as a \textit{general result} that the Filippov vector field associated with the PWS system \eqref{PWS} obtained for $\epsilon \to 0$ 
agrees with the reduced flow on a normally hyperbolic critical manifold contained within the switching layer; see e.g. \cite{Bonet2016,Buzzi2006,Kristiansen2017,Kristiansen2019,Llibre2009,Llibre2007}. In our case, this immediately implies existence of an attracting critical manifold $S$ with a reduced vector field given in the $x-$coordinate chart by \SJ{\eqref{Filippov_VF}.} 
\SJ{The second part of statement (i) and statement (ii) follow from \lemmaref{PWS_lemma} and Fenichel theory \cite{Fenichel1979}.}

\

\begin{figure}[h!]
	\begin{center}
		\includegraphics[width=.5\textwidth]{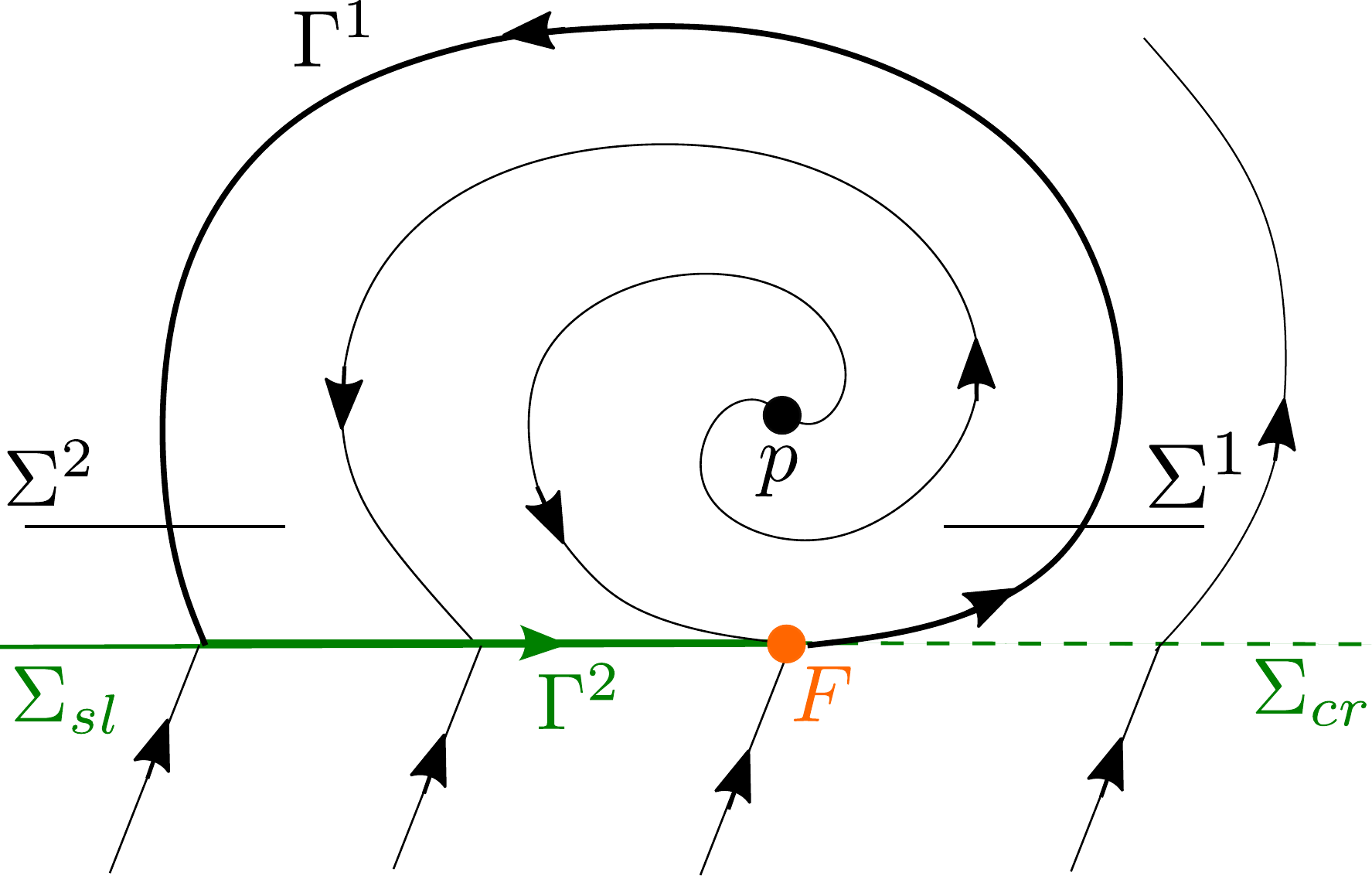}
		\caption{Setup for the proof of \thmref{thm_ro}, with sections $\Sigma^j$, $j=1,2$ transversal to the PWS cycle $\Gamma = \Gamma^1 \cup \Gamma^2$. Note that only the BF$_3$ case is shown; the setup for case BF$_1$ is similar.}
		\figlab{ro_proof}
	\end{center}
\end{figure}

It remains to prove statement (iii), i.e. existence and uniqueness of relaxation oscillations for $\mu \in \mathcal I_+$. The proof follows by a contraction mapping argument: we define sections transversal to $\Gamma$ by
\[
\Sigma^1 = \left\{(x, c): x \in J^1\subset (0, c_+) \right\} , \qquad \Sigma^2 = \left\{(x,c): x \in J^2 \subset (- c_-,0) \right\} ,
\]
where $c, c_\pm > 0$ are \SJ{suitably defined} constants, and show that the Poincar\'e map $\Pi: \Sigma^1 \to \Sigma^1$ defined by the composition $\Pi = \Pi^2 \circ \Pi^1$ \SJ{of flow maps $\Pi^1 : \Sigma^1 \to \Sigma^2$ and $\Pi^2 : \Sigma^2 \to \Sigma^1$} is a strong contraction; the setup is shown in \figref{ro_proof}.

\

First consider the map $\Pi^1$. \lemmaref{PWS_lemma} (iv) ensures the existence an intersection $\Gamma^1 \cap \Sigma^2 = \{(\gamma_L, c)\}$. By choosing $\epsilon_0 > 0$ sufficiently small, we can guarantee that the flow on $\{(x,y) : y \geq \tilde c > 0\}$ is a regular perturbation of the flow induced by \eqref{PWS_normal_form}. In particular, choosing $\tilde c \in (0,c)$ it follows by regular perturbation theory that the map $\Pi^1$ takes the form
\[
\Pi^1 : \Sigma^1 \to \Sigma^2, \qquad (x,c) \mapsto \left(\tilde x(x, \epsilon), c \right) ,
\]
where $\tilde x(x,\epsilon)$ is smooth and, moreover, we have that for any $K>0$ the constants $\epsilon_0$ and $c$ can be chosen sufficiently small so that $|\tilde x(x,\epsilon) - x_d| < K$, where $x_d$ is the `drop point' associated with the $\Gamma^1 \cap \Sigma_{sl}$.

\

Now consider the map $\Pi^2$. It follows from \lemmaref{K1_mu_pos} (v) that in cases BF$_i$, $i=1,3$, the map $\Pi^2$ has the same properties as the local transition map described in \cite[Theorem 1.3]{Kristiansen2019c}.\footnote{\SJ{\lemmaref{K1_mu_pos} (v) and Fenichel theory are required in case BF$_1$ to ensure that the equilibrium $p_{s,\epsilon}$ remains uniformly bounded away from $\Gamma$ for $\epsilon > 0$.}} In particular, if we define $\gamma_R$ via $\Gamma^1 \cap \Sigma^1 = \{(\gamma_R,c)\}$, then the map $\Pi^2$ takes the form
\[
\Pi^2 : \Sigma^2 \to \Sigma^1, \qquad (x,c) \mapsto \left(\gamma_R + \epsilon^{2k/(2k+1)} m_1(\epsilon) , c \right) ,
\]
where $m_1(\epsilon)$ is continuous and, moreover, the restricted map $Q:J^2 \to J^1$ on the $x-$coordinate is a strong contraction, i.e.
\begin{equation}
\eqlab{contraction}
Q(x,\epsilon) = \gamma_R + \epsilon^{2k/(2k+1)} m_1(\epsilon) + \mathcal O(e^{-\kappa/\epsilon}) , \qquad \frac{\partial Q}{\partial x}\bigg|_{(x, \epsilon)} = \mathcal O(e^{-\kappa/\epsilon}) ,
\end{equation}
for some constant $\kappa > 0$ and continuous function $m_1$.

Composing the maps $\Pi^i$, $i=1,2$, we obtain the explicit form
\[
\Pi : \Sigma^1 \to \Sigma^1, \qquad (x,c) \mapsto \left(\gamma_R + \epsilon^{2k/(2k+1)} m_1(\epsilon) , c \right) ,
\]
and using this one can show that the strong contraction property follows from \eqref{contraction}, thus proving the desired result after an application of the contraction mapping principle.
\qed

\section{Blow-up of $\Sigma$ and $Q$}
\applab{blow-up_of_Sigma_and_Q}

The remainder of the manuscript is devoted to the blow-up analysis summarised alongside the statement of our main results in \secref{main_results}, finally culminating with a proof for \thmref{thm_connection} 
in \appref{proof_of_thm_connection}. 
We begin with the blow-up of the switching manifold $\Sigma$ and, subsequently, the degenerate point $Q$ which persists as a consequence of the fold singularity independent of $\mu$.

\

\subsection{Blow-up of $\Sigma$}
\applab{blow-up_of_Sigma}

As discussed in \secref{Sigma_blowup_outline}, system \eqref{normal_form} is singular in the sense that it gives rise to PWS dynamics or `loses smoothness' along $\Sigma$ in the singular limit $\epsilon \to 0$, and blow-up methods can be applied in order to regain smoothness in an extended blow-up space; see e.g. \cite{Jelbart2019c,Kristiansen2019,Kristiansen2019c,Kristiansen2019d}. 
As a first step in the application of such an approach, we consider the extended system obtained by appending system \eqref{normal_form} with the trivial equation $\epsilon' = 0$ and transforming to a fast timescale $\tilde t = t/ \epsilon$:
\begin{equation}
\eqlab{main_extended}
\begin{split}
x' &= \epsilon \left(\tau - \gamma + \phi\left(y \epsilon^{-1} \right) \left(-(\tau - \gamma) + \mu + \tau x - \delta y + \theta_1(x,y,\mu) \right) \right), \\
y' &= \epsilon \left(1 + \phi\left(y \epsilon^{-1} \right) \left(-1 + x + \theta_2(x,y,\mu) \right) \right) , \\
\epsilon' &= 0 ,
\end{split}
\end{equation}
System \eqref{main_extended} has the entire plane $\{\epsilon = 0\}$ as a manifold of equilibria, containing the degenerate subset $\Sigma \times \{0\} = \{(x, 0, 0) : x \in \mathbb R\} \subset \mathbb R^3$ (i.e. the switching manifold) along which the system is non-smooth. 
Smoothness is regained by introducing a blow-up transformation
\begin{equation}
\eqlab{cyl_bu}
r \geq 0 , \ (\bar y, \bar \epsilon ) \in S^1 \mapsto
\begin{cases}
y = r \bar y , \\
\epsilon = r \bar \epsilon ,
\end{cases}
\end{equation}
which blows the switching manifold $\Sigma \times \{0\}$ up to a cylinder $\{r=0\} \times S^1 \times \mathbb R$. Note that since $\epsilon \geq 0$, only the closed half-cylinder corresponding to $\bar \epsilon \geq 0$ is relevant for our analysis. 
For the most part, the analysis can be carried out in a single coordinate chart $K_1$ \SJ{defined by setting $\bar y = 1$.}
In line with the common notational conventions introduced in \cite{Krupa2001a,Krupa2001b}, we define coordinates \SJ{this chart by}
\begin{equation}
\eqlab{K1_coord}
K_1 : \ y = r_1, \qquad \epsilon = r_1 \epsilon_1 .
\end{equation}

We obtain the following equations, after a suitable time desingularisation which amounts to division of the right hand side by a common factor of $\epsilon_1 \geq 0$:
\begin{equation}
\eqlab{K1}
\begin{split}
x' &= r_1 F(x, r_1,\epsilon_1, \mu) , \\
r_1' &= r_1 G(x, r_1,\epsilon_1, \mu) , \\
\epsilon_1' &= - \epsilon_1 G(x, r_1,\epsilon_1, \mu) ,
\end{split}
\end{equation}
where
\begin{equation}
\begin{split}
F(x, r_1,\epsilon_1, \mu) &= (\tau - \gamma) \epsilon_1^k \phi_+(\epsilon_1) + \left(1 - \epsilon_1^k \phi_+(\epsilon_1) \right) \left( \mu + \tau x - \delta r_1 + \theta_1(x,r_1,\mu) \right) , \\
G(x, r_1,\epsilon_1, \mu) &= \epsilon_1^k \phi_+(\epsilon_1) + \left(1 - \epsilon_1^k \phi_+(\epsilon_1) \right) \left( x + \theta_2(x, r_1, \mu) \right) .
\end{split}
\end{equation}
Note that in the above we have used \eqref{reg_asymptotics} in order to write $\phi(\epsilon_1^{-1}) = 1 - \epsilon_1^k \phi_+(\epsilon_1)$.

\begin{remark}
	Formally speaking the `desingularisation' used in the derivation of the equations in \eqref{K1} amounts to a positive transformation of time $d \tilde t_1 = \epsilon_1 d \tilde t$ where by a slight abuse of notation we allow the dash notation in system \eqref{K1} to denote differentiation with respect to the new time $\tilde t_1$. Desingularisations of this kind are an essential feature of the blow-up approach, and commonly understood less abstractly as multiplication or division of the right-hand-side by a common factor following the application of a given blow-up transformation. In the remainder of this work, we simply specify the relevant division or multiplication required to obtain the system of interest.
\end{remark}

System \eqref{K1} has a line of equilibria
\begin{equation}
\eqlab{l_s}
l_s = \left\{(x, 0, 0) : x \in \mathbb R \right\} ,
\end{equation}
along the intersection of the blow-up cylinder with the plane $\{\epsilon_1 = 0\}$; 
see e.g.~\figref{S2_blowups}. Linearising about $l_s$ yields the eigenvalues $\lambda = 0, -x, x$. Hence, the origin $Q:(0,0,0)$ is nonhyperbolic, while $l_s \setminus Q$ is normally hyperbolic and saddle-type.

The planes $\{r_1 = 0\}$ and $\{\epsilon_1 = 0\}$ are invariant under the flow induced by \eqref{K1}. Within $\{r_1 = 0\}$ we identify \SJ{a 1D critical manifold}
\begin{equation}
\eqlab{S1}
S_1 = \left\{\left(x, 0, \epsilon_1 \right) \in \mathbb R \times \mathbb R_+^2 : \epsilon_1^k \phi_+(\epsilon_1) = - \frac{x+\theta_2(x,0,\mu)}{1-x-\theta_2(x,0,\mu)} \right\} ,
\end{equation}
\SJ{which for $\epsilon_1 > 0$ coincides with the attracting critical manifold $S$ identified in the switching layer and described in the proof of \thmref{thm_ro} in \appref{proof_of_thm_ro}.}
Within $\{\epsilon_1 = 0\}$, we identify the image of the equilibrium $p$ in \eqref{focus}, which by an application of the implicit function theorem can be expressed as
\begin{equation}
\eqlab{p1}
p_1(\mu) : \left(\mathcal O(\mu^2), \frac{\mu}{\delta} + \mathcal O(\mu^2) , 0 \right) ,
\end{equation}
in a neighbourhood of $(x,r_1,0,\mu) = (0,0,0,0)$ with $\mu \geq 0$. Note that $p_1(\mu) \to Q$ as $\mu \to 0^+$, i.e. the boundary collision for $\mu \to 0^+$ occurs precisely at the nonhyperbolic point $Q$.

\begin{lemma}
	\lemmalab{K1_mu_pos}
	The following holds for system \eqref{K1}:
	\begin{enumerate}
		\item[(i)] Considered within $\{r_1 = 0\}$, the manifold of equilibria $S_1 \setminus Q$ emanating from the nonhyperbolic point $Q$ is normally hyperbolic and attracting.
		\item[(ii)] There exists a 2D locally invariant, attracting slow manifold $M_1$ with base along $S_1 \setminus Q$.
	\end{enumerate}
	The following holds for system \eqref{K1} with $\mu \in \mathcal I_+$:
	\begin{enumerate}
		\item[(iii)] The slow flow on $M_1$ is increasing (decreasing) in the $x-$ ($\epsilon_1-$) direction locally towards a neighbourhood of $Q$.
		\item[(iv)] Considered within $\{\epsilon_1 = 0\}$, the equilibrium $p_1$ is an unstable focus.
		\item[(v)] In cases BF$_i$, $i=1,2$, there exists an unstable equilibrium $p_{s,1} \in S_1 \setminus Q$, and in cases BF$_i$, $i=1,3$, an improved singular cycle
		\[
		\Gamma = \Gamma^1 \cup \Gamma^2 \cup \Gamma^3
		\]
		can be constructed, where by a slight abuse of notation $\Gamma^2$ in \figref{collision_figs}a,g, has been replaced with $\Gamma^2 \cup \Gamma^3$.
	\end{enumerate}
	The following holds for system \eqref{K1} with $\mu = 0$:
	\begin{enumerate}
		\item[(vi)] In case BF$_3$, the equilibria $p_1$ and $Q$ coalesce in a single non-hyperbolic point at the origin, and the reduced flow on $S_1$ is trivially zero. In cases BF$_i$, $i=1,2,$ the equilibrium $p_{s,1}$ also coalesces at $Q$.
	\end{enumerate}
	Finally, the following holds for system \eqref{K1} with $\mu \in \mathcal I_-$:
	\begin{enumerate}
		\item[(vii)] The slow flow on $M_1$ is decreasing (increasing) in the $x-$ ($\epsilon_1-$) direction locally away from a neighbourhood of $Q$.
		\item[(viii)] In case BF$_3$, there exists a stable equilibrium $p_1 \in S_1 \setminus Q$.
	\end{enumerate}
\end{lemma}

\bpr
	Statements (i) and (iv) are direct calculations, and \SJ{existence of the equilibria $p_{s,1}$ and $p_1$ with stability properties described in (v) and (viii) follows from \thmref{thm_ro} (i)-(ii).} 
	
	\SJ{Statements regarding existence of a 2D center manifold $M_1$ with the properties described in (ii), (iii), and (vii) follow after applying Fenichel theory within switching layer coordinates $(x,Y,\epsilon) = (x,y/\epsilon,\epsilon)$, and extending the results to chart $K_1$ via the transition map
		\begin{equation}
		\eqlab{K_trans}
		r_1 = \epsilon Y, \qquad \epsilon_1 = Y^{-1}, \qquad Y > 0.
		\end{equation}
		
		
		Now consider the second part of statement (v) regarding 
		an improved singular cycle $\Gamma$. The segment $\Gamma^1$ may be identified with $\Gamma^1$ in \lemmaref{PWS_lemma} (iv). In case BF$_3$, the remaining segments} can be explicitly constructed as
	\begin{equation}
	\eqlab{cyc_segs2}
	\begin{split}
	\Gamma^2 &= \left\{(x_d,0,\epsilon_1): \epsilon_1 \in [0,\epsilon_d(x_d) \in S_1] \right\}  , \\
	\Gamma^3 &=  \left\{(x,0,\epsilon_1) \in S_1 : x\in[x_d,0] \right\} ,
	\end{split}
	\end{equation}
	where $x_d$ denotes the $x-$coordinate of the `drop point' associated with the global return mechanism, i.e. the $x-$coordinate of the intersection $\Gamma^1 \cap l_s$; see again \lemmaref{PWS_lemma} on this. Since \SJ{for $\mu \in \mathcal I_+$} there are no equilibria on $S_1$ in case BF$_3$, we need not worry about the location of this drop point. Cases BF$_i$, $i=1,2$ on the other hand, are distinguished by the location of the nonhyperbolic saddle $p_{s,1}$ relative to the global return. In particular it follows from the form and invariance of the fast fibers $\{(c,0, \epsilon_1): c = const\} \in \{r_1 = 0\}$ on the cylinder that if $p_{s,1} : (x_s,0,\epsilon_{1,s})$, we may distinguish between cases BF$_1$ and BF$_2$ as follows:
	\begin{equation}
	\eqlab{BF12_distinction}
	BF_1: \ x_s < x_d , \qquad BF_2: \ x_s > x_d .
	\end{equation}
	In case BF$_1$, a singular cycle $\Gamma$ can be constructed precisely as for case BF$_3$. In case BF$_2$, such a construction is interrupted by the presence of $p_{s,1}$ on the interior of $\Gamma^3$.
	
	\SJ{Now consider statement (vi). The reduced flow within the switching layer is governed by \eqref{Filippov_VF}. Hence $\dot x \equiv 0$ for $\mu = 0$ on $S$, which can be extended to $S_1$ by \eqref{K_trans}. Moreover for $\epsilon = 0$ (again within the switching layer), the $x-$coordinate of the equilibrium (either $p_s$ or $p$, depending on $\mu$ and $\gamma$), is $\tilde x = -\mu/\gamma + \mathcal O(\mu^2)$. 
		The statement (vi) follows, since the $x-$coordinate is the same in chart $K_1$.}
\epr

\SJ{The dynamics are sketched in a representative subset of cases: see \figref{S2_blowups}a-b for cases BF$_i$, $i=1,3$ with $\mu \in \mathcal I_+$ respectively. The case BF$_2$ is similar to BF$_1$, except that the equilibrium $p_{s,1}$ is further to the right on $S_1$, prohibiting the existence of a singular cycle $\Gamma$. \figref{S1_blowups_i}a-b show dynamics in the collision limit $\mu \to 0^+$ in cases BF$_i$, $i=1,2$ and BF$_3$ respectively. Note that in figures like \figref{S1_blowups_i} and \figref{S2_blowups} showing global dynamics we drop the subscript notation on objects like $S_1$ and $p_1$. We adopt this convention throughout the manuscript.}


\begin{remark}
	\SJ{The degenerate case between BF$_1$ and BF$_2$ cases occurs when $x_d = x_s$; see \eqref{BF12_distinction}. This case corresponds to the (PWS) HBF bifurcation discussed in \secref{degenerate_cases}. In the blow-up, an improved singular homoclinic-to-saddle orbit is constructed by setting $x_d=x_s$ in \eqref{cyc_segs2}; see also \cite{Kaleda2011} for results relevant to this.}
\end{remark}

\subsection{Blow-up of $Q$}
\applab{blow-up_of_Q}

In order to resolve degeneracies at the point $Q$ due to both (i) tangency with the flow within $\{\epsilon_1 =0\}$, 
and (ii) loss of normal hyperbolicity 
due to alignment between $S_1$ and the fast fibration, we introduce a successive (spherical) blow-up. We consider system \eqref{K1}, dropping the subscripts for notational convenience, and define a weighted blow-up transformation via the map
\begin{equation}
\eqlab{sph_bu}
\rho \geq 0, \ \left(\bar x, \bar r, \bar \epsilon \right) \in S^2 \mapsto
\begin{cases}
x = \rho^{k(k+1)} \bar x , \\
r = \rho^{2k(k+1)} \bar r, \\
\epsilon = \rho^{k+1} \bar \epsilon ,
\end{cases}
\end{equation}
which blows the nonhyperbolic point $Q$ up to the sphere $\{\rho = 0\} \times S^2$. Note that the blow-up transformation \eqref{sph_bu} does \textit{not} depend on $\mu$, since the degeneracies (i)-(ii) described above exist independently of $\mu$. We will make use of the following coordinate charts:
\[
\mathcal K_1: \bar \epsilon = 1, \qquad \mathcal K_2 : \bar r = 1 , \qquad \mathcal K_3 : \bar x = -1 ,
\]
for which we define chart specific coordinates
\begin{align}
\mathcal K_1: \ x &= \rho_1^{k(k+1)} x_1 , && r = \rho_1^{2k(k+1)} r_1, && \epsilon = \rho_1^{k+1} , \eqlab{mathcal_K1_coord} \\
\mathcal K_2: \ x &= \rho_2^{k(k+1)} x_2 , && r = \rho_2^{2k(k+1)}, && \epsilon = \rho_2^{k+1} \epsilon_2, \eqlab{mathcal_K2_coord} \\
\mathcal K_3: \ x &= - \rho_3^{k(k+1)} , && r = \rho_3^{2k(k+1)} r_3 , && \epsilon = \rho_3^{k+1} \epsilon_3,
\end{align}
related via the transition maps
\begin{equation}
\eqlab{mathcal_transition_map}
\begin{aligned}
\kappa_{12} & : \ x_1 = \epsilon_2^{-k} x_2, && r_1 = \epsilon_2^{-2k} , && \rho_1 = \epsilon_2^{1/(k+1)} \rho_2 , && \qquad \epsilon_2 > 0 , \\
\kappa_{23} & : \ x_2 = -r_3^{-1/2} , && \epsilon_2 = r_3^{-1/2k} \epsilon_3 , && \rho_2 = r_3^{1/2k(k+1)} \rho_3 , && \qquad r_3 > 0 , \\
\kappa_{31} & : \ r_3 = (- x_1)^{-2} r_1 , && \epsilon_3 = (-x_1)^{-1/k} , && \rho_3 =(-x_1)^{1/k(k+1)} \rho_1 , && \qquad x_1 > 0 .
\end{aligned}
\end{equation}

\begin{remark}
	\SJ{In \cite[eqn.~(3.2)]{Kristiansen2019c} the author uses the equivalent weights $(k,2k,1)$ in order to resolve the same degeneracies (i)-(ii). We consider the map \eqref{sph_bu} instead in order to avoid noninteger weights in subsequent blowup transformations.}
\end{remark}

\subsubsection*{Chart $\mathcal K_1$ dynamics}

After a suitable desingularisation (division by $\rho_1^{k(k+1)}$), the equations in chart $\mathcal K_1$ are
\begin{equation}
\eqlab{Ks1}
\begin{split}
x_1' &= r_1 f_1(x_1,r_1,\rho_1,\mu) + k x_1 g_1(x_1,r_1,\rho_1,\mu) ,  \\
r_1' &= r_1 \left(1+2k \right) g_1(x_1,r_1,\rho_1,\mu)  , \\
\rho_1' &= - \frac{\rho_1}{k+1} g_1(x_1,r_1,\rho_1,\mu) ,
\end{split}
\end{equation}
where
\[
\begin{split}
f_1(x_1,r_1,\rho_1,\mu) &= F\left(\rho_1^{k(k+1)} x_1, \rho_1^{2k(k+1)} r_1, \rho_1^{k+1}, \mu \right) , \\
g_1(x_1,r_1,\rho_1,\mu) &= \rho_1^{-k(k+1)} G\left(\rho_1^{k(k+1)} x_1, \rho_1^{2k(k+1)} r_1, \rho_1^{k+1}, \mu \right) ,
\end{split}
\]
and in particular,
\[
f_1(x_1,r_1,0,\mu) = \mu + \mathcal O(\mu^2) , \qquad g_1(x_1,r_1,0,\mu) = \beta + x_1 + \mathcal O(\mu x_1) .
\]

The subspaces $\{\rho_1 = 0\}$, $\{r_1 = 0\}$ and hence $\{\rho_1 = r_1 = 0\}$ are all invariant under the flow induced by \eqref{Ks1}, and correspond to the surface of the blow-up sphere, blow-up cylinder, and their intersection respectively. The dynamics within $\{\rho_1 = r_1 = 0\}$
are governed by the restricted 1-dimensional problem
\begin{equation}
\eqlab{sph_cyl_flow}
x_1' = k x_1 \left(\beta + x_1 + \mathcal O(\mu x_1) \right) ,
\end{equation}
for all $\mu \in (-\mu_0,\mu_0)$ 
with $\mu_0 > 0$ sufficiently small, where we have used $\beta := \phi_{+}(0) > 0$, see \eqref{beta}. 
Hence for each $\mu \in (-\mu_0,\mu_0)$, we identify two equilibria:
\[
p_f : (0,0,0) , \qquad p_a : (-\beta + \mathcal O(\mu), 0, 0) .
\]
The extension of the attracting critical manifold $S_1$ in \eqref{S1} is contained within $\{r_1 = 0\}$, and given implicitly by
\begin{equation}
\eqlab{Ks_S1}
\mathcal S_1 = \left\{\left(x_1,0,\rho_1 \right) : g_1(x_1,0,\rho_1,\mu) = 0  \right\} .
\end{equation}

\

The key dynamics in chart $\mathcal K_1$ are described in the following result, \SJ{and sketched in a representative subset of cases: see \figref{S2_blowups}c-d for cases BF$_i$, $i=1,3$ with $\mu \in \mathcal I_+$ respectively. \figref{S1_blowups_i}c-d show dynamics in the collision limit $\mu \to 0^+$ in cases BF$_i$, $i=1,2$ and BF$_3$ respectively. Case BF$_3$ with $\mu \in \mathcal I_-$ is shown in  \figref{complete_BF3}f.}

\begin{lemma}
	\lemmalab{K1s}
	There exists $\mu_0 > 0$ such that the following holds for system \eqref{Ks1} for all $\mu \in (-\mu_0,\mu_0)$:
	\begin{enumerate}
		\item[(i)] The equilibrium $p_f$ is a hyperbolic saddle, 
		while the equilibrium $p_a$ is partially hyperbolic and attracting with a single non-trivial eigenvalue $\lambda = - k \beta + \mathcal O(\mu)$.
		\item[(ii)] Considered within the invariant plane $\{r_1 = 0\}$, the manifold $\mathcal S_1$ is normally hyperbolic and attracting, terminating at the point $p_a$. Moreover, there exists a local, 1D attracting centre manifold $\mathcal W_1$ emanating from $p_a$ within the invariant plane $\{\rho_1=0\}$.
		\item[(iii)] There exists a 2D attracting slow manifold $\mathcal M_1$ which may be considered as the extension of the manifold $M_1$ described in \lemmaref{K1_mu_pos}, such that
		\[
		\mathcal M_1 \big|_{r_1 = 0} = \mathcal S_1, \qquad \mathcal M_1 \big|_{\rho_1 = 0} = \mathcal W_1 .
		\]
		For $\mu \in \mathcal I_+ \subset (0,\mu_0)$, the slow flow on $\mathcal M_1$ near $p_a$ is increasing in the $r_1-$coordinate, and for $\mu \in \mathcal I_- \subset (-\mu_0, 0)$, the slow flow on $\mathcal M_1$ near $p_a$ is decreasing in the $r_1-$coordinate.
	\end{enumerate}
	The following holds for system \eqref{Ks1} in the limit $\mu \to 0$:
	\begin{enumerate}
		\item[(iv)] The 1D centre manifold $\mathcal W_1$ described above in (iii) limits to
		\begin{equation}
		\eqlab{Ks1_manifold}
		\mathcal W_1 = \left\{(-\beta, r_1, 0) : r_1 \geq 0 \right\} ,
		\end{equation}
		which is a normally hyperbolic attracting critical manifold when considered within the invariant plane $\{\rho_1 = 0\}$, with a single non-trivial eigenvalue given by $\lambda = - k \beta < 0$.
		\item[(v)] The manifold $\mathcal M_1$ described above in (iii) persists, and the slow-flow on $\mathcal M_1 \setminus \mathcal W_1$ is increasing (decreasing) in the $r_1-$direction in case $\mu \to 0^+$ ($\mu \to 0^-$).
	\end{enumerate}
\end{lemma}

\bpr
	%
	Statement (i) and the first half of statement (ii) are direct calculations.
	
	Existence of the centre manifolds $\mathcal M_1$ and $\mathcal W_1$ with the properties described in statement (iii) follow from centre manifold theory. 
	In particular, standard matching arguments lead to following the local graph form for $\mathcal W_1$ for $\mu \in (0,\mu_+)$ sufficiently small:
	\[
	\mathcal W_1 : x_1 = - \frac{\beta}{1 + A_1 \mu} + \frac{1}{\beta k} \mu r_1 + \mathcal O(r_1 \mu^2, r_1^2) ,
	\]
	where $A_1 = ( \partial \omega_1 / \partial \mu )|_{(x_1,r_1,0,0)}$. The flow on $\mathcal W_1'$ is given in the $r_1-$coordinate chart by
	\[
	r_1' \big|_{\mathcal W_1} = (1 + 2k) \mu r_1^2 + \mathcal O(r_1^2 \mu^2, r_1^3) ,
	\]
	which is positive (negative) when $\mu \in \mathcal I_+$ ($\mu \in \mathcal I_-$) and $r_1 > 0$, as required.
	
	In order to prove statements (iv)-(v) we consider system \eqref{Ks1} as a singularly perturbed problem with $|\mu| \ll 1$. First consider the system \eqref{Ks1}$|_{\{\rho_1 = 0\}}$, which can be written in the general form of a perturbation problem
	\begin{equation}
	\eqlab{Ks1_sing}
	\begin{pmatrix}
	x_1' \\ 
	r_1' \\ 
	\end{pmatrix}
	=
	\begin{pmatrix}
	k x_1 \\ 
	r_1 (1 + 2k) \\ 
	\end{pmatrix}
	g_1(x_1,r_1,0,0) + \mu
	\begin{pmatrix}
	r_1 + A_1 k x_1^2 + \mathcal O(\mu) \\ 
	A_1 (1 + 2k) x_1 r_1 + \mathcal O(\mu) \\ 
	\end{pmatrix} .
	\end{equation}
	System \eqref{Ks1_sing} has a critical manifold $\mathcal W_1$ when $\mu = 0$, which by solving the simple equation $g_1(x_1,r_1,0,0) = \beta + x_1 = 0$ can be seen to take the explicit form in \eqref{Ks1_manifold}. Direct calculation shows that $\mathcal W_1$ is normally hyperbolic and attracting with non-trivial eigenvalue $\lambda = - k \beta$, as required. 
	In order to prove statement (v) we consider the reduced problem on $\mathcal W_1$, which is given in the $r_1-$coordinate chart by
	\[
	\dot r_1 = \pm \left( \frac{1+2k}{k \beta} \right) r_1^2 \geq 0 ,
	\]
	where the right-hand-side is positive (negative) in the limiting case $|\mu| = \mu \to 0^+$ ($|\mu| = - \mu \to 0^+$). This combined with centre manifold theory applied along $\mathcal W_1$ is sufficient to prove the desired result.
\epr

\subsubsection*{Chart $\mathcal K_2$ dynamics}

After a suitable desingularisation (division by $\rho_2^{k(k+1)}$), the equations in chart $\mathcal K_2$ are
\begin{equation}
\eqlab{Ks2}
\begin{split}
x_2' &= f_2(x_2,\epsilon_2,\rho_2,\mu) - \frac{x_2}{2} g_2(x_2,\epsilon_2,\rho_2,\mu) , \\
\epsilon_2' &= -\epsilon_2 \left(\frac{1+2k}{2k}\right) g_2(x_2,\epsilon_2,\rho_2,\mu) , \\
\rho_2' &= \frac{\rho_2}{2k(1+k)} g_2(x_2,\epsilon_2,\rho_2,\mu) ,
\end{split}
\end{equation}
where 
\[
\begin{split}
f_2(x_2,\epsilon_2,\rho_2,\mu) &= F\left(\rho_2^{k(k+1)} x_2, \rho_2^{2k(k+1)}, \rho_2^{k+1} \epsilon_2, \mu \right) , \\
g_2(x_2,\epsilon_2,\rho_2,\mu) &= \rho_2^{-k(k+1)} G\left(\rho_2^{k(k+2)} x_2, \rho_2^{2k(k+1)}, \rho_2^{k+1} \epsilon_2, \mu \right) ,
\end{split}
\]
and in particular,
\[
f_2(x_1,\epsilon_2,0,\mu) = \mu + \mathcal O(\mu^2) , \qquad g_2(x_2,\epsilon_2,0,\mu) = \beta \epsilon_2^k + x_2 + \mathcal O(\mu x_2) .
\]

The subspaces $\{\rho_2 = 0\}$, $\{\epsilon_2 = 0\}$ and $\{\rho_2 = \epsilon_2 = 0\}$ are all invariant under the flow induced by \eqref{Ks1}, and correspond to the blow-up sphere, the plane $\{(x,y,\epsilon) : y \geq 0 , \epsilon = 0\}$, and their intersection respectively. The dynamics within $\{\rho_2 = \epsilon_2 = 0\}$
are governed by the restricted 1-dimensional problem
\begin{equation}
\eqlab{up_equator}
x_2' = \mu - \frac{x_2}{2} \left(x_2 + \mathcal O(\mu x_2) \right) + \mathcal O(\mu^2) ,
\end{equation}
for all $\mu \in (-\mu_0,\mu_0)$ with $\mu_0 > 0$ sufficiently small. We identify two equilibria when $\mu > 0$:
\begin{equation}
\eqlab{Ks2_eqs}
q_{in} : \left(- \sqrt{2 \mu} + \mathcal O(\mu), 0, 0 \right), \qquad q_{out} : \left(\sqrt{2 \mu } + \mathcal O(\mu), 0, 0 \right) .
\end{equation}
Notice in particular that $q_{in}$ and $q_{out}$ collide in a single point $Q_{bfb}:(0,0,0)$ in the collision limit $\mu \to 0^+$.

System \eqref{Ks2} also has an equilibrium $P_2$ within $\{\epsilon_2 = 0\}$ (the image of $p$ in chart $\mathcal K_2$), which can be located by applying the blow-up transformation \eqref{mathcal_K2_coord} to the expression \eqref{p1} to obtain
\begin{equation}
\eqlab{true_eq_Ks2}
P_2(\mu) : \left( \mathcal O(\mu^{3/2}) , 0, \mathcal O\left(\mu^{1/(2k(1+k))} \right) \right) ,
\end{equation}
in a neighbourhood of $(x_2,0,\rho_2,\mu)=(0,0,0,0)$ with $\mu \geq 0$.

\begin{lemma}
	\lemmalab{K2s}
	There exists $\mu_0 > 0$ such that the following holds for system \eqref{Ks1} for all $\mu \in (-\mu_0,\mu_0)$:
	\begin{enumerate}
		\item[(i)] For all $\mu \in \mathcal I_+ \subset (0,\mu_0)$, the equilibria $q_{in}$ and $q_{out}$ are hyperbolic saddles. 
		\item[(ii)] Considered within the invariant plane $\{\epsilon_2 = 0\}$, the equilibrium $P_2(\mu)$ is an unstable focus with a regular heteroclinic connection to $q_{in}$ for all $\mu \in \mathcal I_+ \subset (0,\mu_0)$.
		\item[(iii)] The extension of the 1D centre manifold $\mathcal W_1$ identified in \lemmaref{K1s}, call it $\mathcal W_2$, is contained within $\{\rho_2 = 0\}$ and connects to $q_{out}$ in case $\mu \in \mathcal I_+ \subset (0,\mu_0)$ (\figref{S2_blowups}c-d). In case $\mu \in \mathcal I_- \subset (-\mu_0,0)$ the manifold $\mathcal W_1$ connects to $p_f$ (\figref{complete_BF3}f). The extension of the manifold $\mathcal M_1$, call it $\mathcal M_2$, satisfies
		\[
		\mathcal M_2 \big|_{\rho_2 = 0} = \mathcal W_2 .
		\]
	\end{enumerate}
	The following holds for system \eqref{Ks2} in the limit $\mu \to 0$:
	\begin{enumerate}
		\item[(iv)] The equilibria $q_{in}, q_{out}, P_2$ all coalesce in a single point $Q_{bfb} : (0,0,0)$, and the orbit $\mathcal W_2$ described above in (iii) limits to a 1D critical manifold
		\begin{equation}
		\eqlab{Ks2_manifold}
		\mathcal W_2 = \left\{(-\beta \epsilon_2^k, \epsilon_2, 0) : \epsilon_2 \geq 0 \right\} ,
		\end{equation}
		with a single non-trivial eigenvalue given by $\lambda = - \beta k \epsilon_2^k \leq 0$. Hence $Q_{bfb}$ is nonhyperbolic, while the manifold $\mathcal W_2 \setminus Q_{bfb}$ is a normally hyperbolic and attracting when considered within the invariant plane $\{\rho_2 = 0\}$.
		\item[(v)] The manifold $\mathcal M_2$ described above in (iii) persists, and the slow flow on $\mathcal M_2 \setminus \mathcal W_2$ is decreasing (increasing) in the $\epsilon_2-$direction in case $\mu \to 0^+$ ($\mu \to 0^-$).
	\end{enumerate}
\end{lemma}

\bpr
	Statement (i) follows after linearisation  of system \eqref{Ks2} and restriction to $q_{in/out}$.
	
	The fact that $P_2(\mu)$ is an unstable focus when considered within $\{\epsilon_2 = 0\}$ follows directly from \lemmaref{K1_mu_pos} (v) and the form of the coordinates given by \eqref{mathcal_K2_coord}. The existence of a regular heteroclinic connection between $P_2$ and $q_{in}$ follows from the Poincar\'e-Bendixon theorem applied within the compact region in $\{\bar \epsilon = 0\}$ bounded between $\Gamma^1$, $l_s$, and the equator $\{\rho = 0\} \times \{(\bar x, \bar r, 0) \in S^2\}$ of the blow-up sphere, see \figref{S2_blowups}c-d or \figref{complete_BF3}a. In particular, $P_2$ is the only repeller with an unstable manifold intersecting the interior of this region, and there are no limit cycles for $\mu \in (0,\mu_0)$ sufficiently small since $\langle \nabla, X^+(x,y,\mu) \rangle \neq 0$ in a neighbourhood of $(0,0,0)$. Hence the (unique) stable manifold $W^s(q_{in})$ emanating from $q_{in}$ connects to $P_2$, proving the statement (ii).
	
	In order to prove statement (iii), consider the (unique) 1D centre manifold $\mathcal W_1$ identified in \lemmaref{K1s}. Since $\mathcal W_1$ enters the visible region in chart $\mathcal K_2$ transversally, we can define its extension $\mathcal W_2$ via transition map \eqref{mathcal_transition_map}, i.e. $\mathcal W_2 = \kappa_{12}(\mathcal W_1)$. Another application of the map $\kappa_{12}$ shows that $\mathcal W_1 \subset \{\rho_1 = 0\}$ implies $\mathcal W_2 \subset \{\rho_2 = 0\}$. In case $\mu \in \mathcal I_+ \subset (0,\mu_0)$, $q_{out}$ is the only attractor for trajectories intersecting the interior of the compact sphere-segment $\{\rho = 0\} \times \{(\bar x, \bar r, \bar \epsilon) : \bar r \geq 0, \bar \epsilon \geq 0 \}$,
	it follows by the Poincar\'e-Bendixon theorem that $\mathcal W_2$ connects to $q_{out}$. In case $\mu \in \mathcal I_- \subset (-\mu_0,0)$, the equilibria $q_{in}$, $q_{out}$ do not exist and the flow along the equator given in \eqref{up_equator} satisfies $x_2' < 0$ uniformly with respect to $\mu$. In this case $p_a$ ($p_f$) is the only attractor (repeller) for trajectories intersecting the interior of the compact sphere-segment $\{\rho = 0\} \times \{(\bar x, \bar r, \bar \epsilon) : \bar r \geq 0, \bar \epsilon \geq 0 \}$, implying a connection by the Poincar\'e-Bendixon theorem; see again \figref{complete_BF3}f.
	
	In order to prove statements (iv)-(v), we consider system \eqref{Ks2} in the limit $\mu \to 0$, treating $|\mu| \ll 1$ as a singular perturbation parameter. In particular, we have the following singularly perturbed problem within $\{\rho_2 = 0\}$:
	\begin{equation}
	\eqlab{Ks2_sing}
	\begin{pmatrix}
	x_2' \\ 
	\epsilon_2' \\ 
	\end{pmatrix}
	=
	\begin{pmatrix}
	- \frac {x_2} 2 \\ 
	- \left(\frac{1+2k}{2k}\right) \epsilon_2  \\ 
	\end{pmatrix}
	g_2(x_2,\epsilon_2,0,0) + \mu
	\begin{pmatrix}
	1 - A_2 \frac{x_2^2}{2}  + \mathcal O(\mu) \\ 
	-A_2 \left(\frac{1+2k}{2k}\right) \epsilon_2 x_2 + \mathcal O(\mu) \\ 
	\end{pmatrix} ,
	\end{equation}
	where $A_2 = ( \partial \omega_2 / \partial \mu )|_{(x_1,0,0)}$.
	The fact that all three points $q_{in}, q_{out}, P_2$ coalesce at a single point $Q_{bfb} : (0,0,0)$ follows directly from expressions \eqref{Ks2_eqs} and \eqref{true_eq_Ks2}, and the form for $\mathcal W_2$ given in \eqref{Ks2_manifold} follows after solving $g_2(x_2,\epsilon_2,0,0) = \beta \epsilon_2^k + x_2 = 0$. The second half of statement (iv) follows from direct calculations. 
	
	Finally, consider statement (v). The reduced problem associated with \eqref{Ks2_sing} is given in the $\epsilon_2-$coordinate chart by
	\begin{equation}
	\eqlab{mathcal_K2_reduced}
	\dot \epsilon_2 = \mp \left( \frac{1+2k}{2k^2 \beta} \right) \epsilon_2^{1-k} < 0 ,
	\end{equation}
	where the right-hand-side is negative (positive) in the limiting case $|\mu| = \mu \to 0^+$ ($|\mu| = - \mu \to 0^+$). The manifold $\mathcal M_2$ (which has its base along $\mathcal W_2$) therefore has a slow flow which is also decreasing in the $\epsilon_2-$direction.
\epr

\begin{remark}
	The fact that the reduced problem \eqref{mathcal_K2_reduced} has a finite time blow-up for $k \geq 2$ reflects the fact that there is an order $k$ tangency between the critical manifold $\mathcal W_2$ and the invariant fast `fiber' contained within $\{(0,\epsilon_2,0) : \epsilon_2 \geq 0\}$. Generically, tangencies of this kind are associate with finite time blow-up due to a factor $\epsilon_2^{-k}$, however in this case the equilibrium $Q_{bfb}$ when $\epsilon_2 = 0$ leads to the factor $\epsilon_2^{1-k}$ in \eqref{mathcal_K2_reduced} instead.
\end{remark}

\subsubsection*{Chart $\mathcal K_3$ dynamics}

After a suitable desingularisation (division by $\rho_3^{k(k+1)}$), the equations in chart $\mathcal K_3$ are
\begin{equation}
\eqlab{Ks3}
\begin{split}
r_3' &= r_3 \left(2 r_3 f_3(r_3,\epsilon_3,\rho_3,\mu) + g_3(r_3,\epsilon_3,\rho_3,\mu) \right) ,  \\
\epsilon_3' &= \epsilon_3 \left(\frac{r_3}{k} f_3(r_3,\epsilon_3,\rho_3,\mu) - g_3(r_3,\epsilon_3,\rho_3,\mu) \right) ,  \\
\rho_3' &= - \frac{\rho_3 r_3}{k(k+1)} f_3(r_3,\epsilon_3,\rho_3,\mu) ,
\end{split}
\end{equation}
where
\[
\begin{split}
f_3(r_3,\epsilon_3,\rho_3,\mu) &= F\left(-\rho_3^{k(k+1)}, \rho_3^{2k(k+1)} r_3, \rho_3^{k+1} \epsilon_3, \mu \right) , \\
g_3(r_3,\epsilon_3,\rho_3,\mu) &= \rho_1^{-k(k+1)} G\left(-\rho_3^{k(k+1)}, \rho_3^{2k(k+1)} r_3, \rho_3^{k+1} \epsilon_3, \mu \right) ,
\end{split}
\]
and in particular,
\[
f_3(r_3,\epsilon_3,0,\mu) = \mu + \mathcal O(\mu^2) , \qquad g_3(r_3,\epsilon_3,0,\mu) = \beta \epsilon_3^k + x_1 + \mathcal O(\mu ) .
\]

One can identify the image of $q_{in}$, $p_a$, $\mathcal W$ and $\mathcal S$ here in the relevant invariant subspaces, however our main interest here concerns the existence of a line of fixed points $\mathcal L_s = \{(0,0,\rho_3) : \rho_3 \geq 0\}$ which coincides with $l_s \cap \{x < 0\}$ (with $l_s$ given by \eqref{l_s}) where domains overlap.
We have the following result:

\begin{lemma}
	\lemmalab{Ks3}
	There exists $\mu_0 > 0$ such that for all $\mu \in (-\mu_0,\mu_0)$, the line $\mathcal L_s$ (including $p_-$) is normally hyperbolic and saddle-type. 
\end{lemma}

\bpr
	The result follows after linearisation along $\mathcal L_s$.
\epr

\subsection{$\mu>0$ singular cycles}
\applab{singular_cycles_pos}


The analysis thus far is sufficient to construct a family of singular relaxation cycles by a concatenation $\Gamma^{ro} = \Gamma^1 \cup \Gamma^2 \cup \Gamma^3 \cup \Gamma^4$ when $\mu > 0$; see \figref{S2_blowups}c-d and \figref{complete_BF3}a. The cycle segments $\Gamma^j$, $j=1,2,3$ can be obtained by extending previous definitions from \lemmaref{K1_mu_pos} and \eqref{cyc_segs2} into the blow-up space obtained after applying \eqref{sph_bu}, where in particular $\Gamma^1$ can be identified with the global unstable manifold $W^u(q_{out}) \subset \{\epsilon = 0\}$ and $\Gamma^3$ connects to $p_a$. The final segment $\Gamma^4$ can be identified with $\mathcal W$, see again \figref{S2_blowups}c-d. We obtain a singular relaxation cycle for each fixed $\mu \in (0,\mu_+)$ (and hence for each $\mu \in \mathcal I_+$), which is \textit{nondegenerate} in the sense that partial hyperbolicity has been regained everywhere. For analytical purposes we parameterise the family of singular cycles obtained by this construction by the chart $\mathcal K_3$ coordinate $\rho_3$:
\begin{equation}
\eqlab{large_ro}
\Gamma^{ro}(\rho_3) = \left(\Gamma^{1}\cup \Gamma^{2} \cup \Gamma^{3} \cup \Gamma^{4}\right)(\rho_3) , \qquad \rho_3 \in (0,R) , 
\end{equation}
which is possible since the map $\mu \mapsto \rho_{3,d}(\mu)$, where $\rho_{3,d}(\mu)$ is determined by the drop point $\Gamma^1(\mu) \cap l_s$, is a diffeomorphism satisfying $\rho_{3,d}'(\mu) > 0$ for all $\mu \in (0,\mu_0)$ with $\mu_0 > 0$ sufficiently small. In particular, the constants $R$ and $\mu_\pm$ can be chosen so that $\rho_{3,d}(\mathcal I_+) \in (0,R)$, i.e. so that the family \eqref{large_ro} contains all singular cycles which  
perturb to the relaxation oscillations described by \thmref{thm_ro}. 

\section{Blow-up of $Q_{bfb}$}
\applab{blow-up_of_Qbfb}

By \lemmaref{K2s} (iv), the equilibria $q_{in}, \ q_{out}, \ P_2$ identified in \eqref{Ks2_eqs} and \eqref{true_eq_Ks2}, chart $\mathcal K_2$, all coalesce in a single non-hyperbolic point $Q_{bfb}$ when $\mu = 0$. Since the system is only degenerate at $Q_{bfb}$ when $\mu = 0$, we must consider an \textit{extended} phase space in which $\mu$ is also considered as a variable. Hence we consider system \eqref{Ks2}, dropping the subscripts for notational convenience and appending the equations with the trivial equation $\mu' = 0$:
\begin{equation}
\eqlab{main_extended2}
\begin{split}
x' &= f(x,\epsilon,\rho,\mu) - \frac{x}{2} g(x,\epsilon,\rho,\mu) , \\
\epsilon' &= -\epsilon \left(\frac{1+2k}{2k}\right) g(x,\epsilon,\rho,\mu) , \\
\rho' &= \frac{\rho}{2k(1+k)} g(x,\epsilon,\rho,\mu) , \\
\mu' &= 0 .
\end{split}
\end{equation}

\SJ{\begin{remark}
		Alternatively, one can blow-up the point $(\epsilon,\mu) = (0,0)$ in parameter space prior to the coordinate blow-ups in \eqref{cyl_bu} and \eqref{sph_bu}, thus allowing for a representation of the problem in a three-dimensional space; see e.g. \cite{Maesschalck2011b,Dum1996,Dumortier2009} for applications of arguments of this kind. On this approach, however, additional blow-up transformations in both parameter and coordinate spaces are required in order to desingularise the system \eqref{normal_form}, which significantly complicates the geometry and analysis. For these reasons, we opt to consider the doubly extended system \eqref{main_extended2}.
\end{remark}}

The point $Q_{bfb}:(0,0,0,0)$ is nonhyperbolic in the extended system \eqref{main_extended2}. To resolve this degeneracy, we define a weighted blow-up transformation
\begin{equation}
\eqlab{bu_Qbfb}
\nu \geq 0 , \ \left(\bar x, \bar \epsilon, \bar \rho, \bar \mu \right) \in S^3 \mapsto
\begin{cases}
x = \nu^{k(k+1)} \bar x , \\
\epsilon = \nu^{k+1} \bar \epsilon , \\
\rho = \nu \bar \rho , \\
\mu = \nu^{2k(k+1)} \bar \mu ,
\end{cases}
\end{equation}
which transforms the point $Q_{bfb}$ into the 3-sphere $\{\nu = 0\} \times S^3$. We define three coordinate charts
\[
\mathcal K_1' : \bar \epsilon = 1 , \qquad \mathcal K_2' : \bar \mu = 1, \qquad \mathcal K_3' : \bar x = - 1 ,
\]
with chart-specific coordinates
\begin{align}
\mathcal K_1' : x &= \nu_1^{k(k+1)} x_1, \ && \epsilon = \nu_1^{k+1} , \ && \rho = \nu_1 \rho_1, \ && \mu = \nu_1^{2k(k+1)} \mu_1 , \eqlab{mathfrak_K1_coord} \\
\mathcal K_2' : x &= \nu_2^{k(k+1)} x_2, \ && \epsilon = \nu_2^{k+1} \epsilon_2 , \ && \rho = \nu_2 \rho_2, \ && \mu = \nu_2^{2k(k+1)} ,   \\
\mathcal K_3' : x &= - \nu_3^{k(k+1)} , \ && \epsilon = \nu_3^{k+1} \epsilon_3 , \ && \rho = \nu_3 \rho_3, \ && \mu = \nu_3^{2k(k+1)} \mu_3 , \eqlab{mathfrak_K3_coord}
\end{align}
and transition maps
\begin{equation}
\eqlab{dash_transition_maps}
\begin{aligned}
\kappa_{12}': x_1 &= \epsilon_2^{-k} x_2, \ && \rho_1 = \epsilon_2^{-1/(1+k)} \rho_2 , \ && \nu_1 = \epsilon_2^{1/(1+k)} \nu_2 , \ && \mu_1 = \epsilon_2^{-2k} , \\ 
\kappa_{23}': x_2 &= - \mu_3^{-1/2}, \ && \rho_2 = \mu_3^{-1/2k(1+k)} \rho_3 , \ && \epsilon_2 = \mu_3^{-1/2k} \epsilon_3 , \ && \nu_2 = \mu_3^{1/2k(1+k)} \nu_3 , \\ 
\kappa_{31}': \epsilon_3 &= (-x_1)^{-1/k} , \ && \rho_3 = (- x_1)^{-1/k(1+k)} \rho_1 , \ && \nu_3 = (-x_1)^{1/k(1+k)} \nu_1 , \ && \mu_3 = (-x_1)^{-2} \mu_1 , \\ 
\end{aligned}
\end{equation}
where $\kappa_{12}'$, $\kappa_{23}'$ and $\kappa_{31}'$ are defined for $\epsilon_2 > 0$, $\mu_3 > 0$ and $x_1 < 0$ respectively. 


\

In the following we consider the dynamics within each coordinate chart $\mathcal K'_i$, $i=1,2,3$, with a particular emphasis on the dynamics relevant for understanding the matching problem associated with the limit in \eqref{matching_limit}. This will prove to be sufficient for a proof of \thmref{thm_connection} in \appref{proof_of_thm_connection}. 

\subsubsection*{$\mathcal K_1'$ chart}

Following a suitable desingularisation (division by $\nu_1^{k(k+1)}$), we obtain the following equations in chart $\mathcal K_1'$:
\begin{equation}
\eqlab{mathcal_Ks1}
\begin{split}
x_1' &= \tilde f_1(x_1,\rho_1,\nu_1,\mu_1) + k x_1 \tilde g_1(x_1,\rho_1,\nu_1,\mu_1), \\
\rho_1' &= \frac{\rho_1}{k} \tilde g_1(x_1,\rho_1,\nu_1,\mu_1), \\
\nu_1' &= - \nu_1 \left(\frac{1+2k}{2k(1+k)} \right) \tilde g_1(x_1,\rho_1,\nu_1,\mu_1), \\
\mu_1' &= \mu_1 \left(1+2k\right) \tilde g_1(x_1,\rho_1,\nu_1,\mu_1), \\
\end{split}
\end{equation}
where
\[
\begin{split}
\tilde f_1(x_1,\rho_1,\nu_1,\mu_1) &= \nu_1^{-2k(k+1)} f\left(\nu_1^{k(k+1)} x_1, \nu_1^{k+1} , \nu_1 \rho_1, \nu_1^{2k(k+1)} \mu_1 \right) , \\
\tilde g_1(x_1,\rho_1,\nu_1,\mu_1) &= \nu_1^{-k(k+1)} g \left(\nu_1^{k(k+1)} x_1, \nu_1^{k+1} , \nu_1 \rho_1, \nu_1^{2k(k+1)} \mu_1 \right) ,
\end{split}
\]
and in particular,
\[
\begin{split}
\tilde f_1 & (x_1,\rho_1,0,\mu_1) = \mu_1 + \rho_1^{k(k+1)} \left( \beta (\tau - \gamma) + \tau x_1 - \delta \rho_1^{k(k+1)} \right) 
, \qquad
 \tilde g_1(x_1,\rho_1,0,\mu_1) = \beta + x_1 .
\end{split}
\]

System \eqref{mathcal_Ks1} has a number of invariant subspaces; the hyperplanes $\{\rho_1=0\}$, $\{\nu_1 =0\}$, $\{\mu_1 = 0\}$ and any intersection of these are all invariant. In particular, the dynamics within the invariant line $\{\rho_1 = \nu_1 = \mu_1 =0\}$ are governed by
\[
x_1' = k x_1 \left(\beta + x_1 \right) ,
\]
with two equilibria
\begin{equation}
\eqlab{mathfrak_K1_eqs}
q_f: \left(0, 0, 0, 0 \right) , \qquad q_a: \left(-\beta, 0, 0, 0 \right) .
\end{equation}
We also have non-trivial invariant sets which define the scaling regime (S1), as described by the following result:

\begin{lemma}
	\lemmalab{par_sets}
	The one-parameter family of sets
	\begin{equation}
	\eqlab{mathcal_A}
	\mathcal A_{1} = \left\{\left(x_1, \rho_1, \nu_1, \hat \mu \rho_1^{k(1+2k)} \right) : \hat \mu \in \mathbb R \right\} ,
	\end{equation}
	is invariant under the flow induced by system \eqref{mathcal_Ks1}. Hence, $\hat \mu$ can be considered as a bifurcation parameter for the restricted system \eqref{mathcal_Ks1}$|_{\mathcal A_1}$.
\end{lemma}

\bpr
	Since $\epsilon$ and $\mu$ are constants of the motion for \eqref{main_extended2}, so are the right-hand-sides in the expressions
	\begin{equation}
	\begin{split}
	\epsilon &= r_1 \epsilon_1 = \rho_2^{(k+1)(1+2k)} \epsilon_2 = \nu_1^{2(1+k)^2} \rho_1^{(k+1)(1+2k)}, \\
	\mu &= \nu_1^{2k(k+1)} \mu_1.
	\end{split}
	\end{equation}
	Restricting $\mu_1$ to $\mathcal A_{1}$, we obtain
	\[
	\mu = \hat \mu \nu_1^{2k(k+1)} \rho_1^{k(1+2k)} = \hat \mu \epsilon^{k/(k+1)} ,
	\]
	so that in particular,
	\begin{equation}
	\eqlab{par_scale}
	\hat \mu = \mu \epsilon^{-k/(k+1)}
	\end{equation}
	is constant.
\epr

\begin{remark}
	Equation \eqref{par_scale} 
	defines the scaling regime (S1), i.e.~
	the invariant set $\mathcal A_1$ identified in \eqref{mathcal_A} is the image of the set $\mathcal A$ in \eqref{mathcal_A_1} in chart $\mathcal K_1'$. Typically one identifies invariant sets like $\mathcal A_1$ first, using them to infer the associated and dynamically important scaling regime.
\end{remark}

It follows from \lemmaref{par_sets} and invariance of the hyperplane $\{\rho_1 = 0\}$ that the intersection $\mathcal A_{1} \cap \{\rho_1 = 0\}$ is also invariant. Within $\mathcal A_{1} \cap \{\rho_1 = 0\}$, we identify the extension of the critical manifold $\mathcal W_2$ described in \lemmaref{K2s}, 
\begin{equation}
\eqlab{mathcal_W1_dash}
\mathcal W_1' = \left\{(-\beta, 0, \nu_1, 0) : \nu_1 \geq 0 \right\} ,
\end{equation}
which connects to the point $q_a$; see \figref{S1_blowups_ii}.

Similarly, $\{\nu_1 = 0\} \cap \mathcal A_1$ is invariant, with dynamics in this set governed by
the system \eqref{desing_prob} identified in \lemmaref{lem_desing} and described in \lemmaref{lem_bifurcations_1}. Recall that by \lemmaref{lem_bifurcations_1} system \eqref{desing_prob} has either zero, one, or two equilibria on $\{\rho_1 > 0\}$, depending on the region in $(\hat \mu, \gamma)-$parameter space; see again \figref{bfb_bif_diagram}, \figref{bfb_bif_neg} and \figref{bfb_bif_pos}.


\begin{lemma}
	\lemmalab{mathfrak_K1}
	The following holds for system \eqref{mathcal_Ks1}:
	\begin{enumerate}
		\item[(i)] The equilibrium $q_f$ is a hyperbolic saddle, 
		and the equilibrium $q_a$ is partially hyperbolic with a single non-zero eigenvalue $\lambda = - k b < 0$.
		\item[(ii)] Considered within $\{\rho_1 =0 \} \cap \mathcal A_1$, the critical manifold $\mathcal W_1'$ is normally hyperbolic and attracting with non-trivial eigenvalue $\lambda = - k \beta$.
		\item[(iii)] There exists a 3D attracting centre manifold $\mathcal M_1'$ such that
		\[
		\mathcal M_1'\big|_{\{\rho_1=0\} \cap \mathcal A_1} = \mathcal W_1' , \qquad \mathcal M_1' \big|_{\{\nu_1=0\} \cap \mathcal A_1} = \mathcal J_1' , \qquad \mathcal M_1' \big|_{\{\rho_1 = \nu_1 = 0\}} = \mathcal N_1' ,
		\]
		where $\mathcal J_1' \subset \{\nu_1 = 0\} \cap \mathcal A_1$ and $\mathcal N_1' \subset \{\nu_1 = \nu_1 = 0\}$ denote local attracting 1D centre manifolds emanating from $q_a$.
		In case $\gamma < 0$ ($\gamma > 0$) the slow flow on $\mathcal M_1' \setminus \mathcal W_1'$ is (decreasing) increasing in the $\rho_1-$direction.
		\item[(iv)] \SJ{In case $\gamma < 0$, system \eqref{desing_prob}$|_{\rho_1 > 0}$ has \SJ{an} equilibrium $P_1'(\hat \mu)$. Moreover,
			\begin{equation}
			\eqlab{P_lims_pos}
			\lim_{\hat\mu \to - \infty} P_1'(\hat \mu) = q_a, \qquad \lim_{\hat\mu \to \infty} P_1'(\hat \mu) = (-\beta, \infty, 0, 0 ) .
			\end{equation}
			In case $\gamma > 0$ with $\hat \mu > \hat \mu_{sn}(\gamma)$, where $\mu_{sn}(\gamma)$ is the saddle-node value given by \eqref{sn_curve_lem_1}, system \eqref{desing_prob}$|_{\rho_1 > 0}$ has two equilibria $P_1'(\hat \mu)$ and $P_{s,1}'(\hat \mu)$. Moreover,
			\begin{equation}
			\eqlab{P_lims_neg}
			\lim_{\hat\mu \to \infty} P_1'(\hat \mu) = (-\beta, \infty, 0, 0 ) , \qquad \lim_{\hat\mu \to \infty} P_{s,1}'(\hat \mu) = q_a .
			\end{equation}}
	\end{enumerate}
\end{lemma}

\bpr
	Statements (i)-(ii) follow after linearisation at $q_f, \ q_a$, and $\mathcal W_1'$, and statement (iii) follows from centre manifold theory. To show that the flow on $\mathcal M_1'$ is increasing in the $\rho_1-$direction, it suffices to prove the statement for the flow on $\mathcal J_1'$. Making a power series ansatz for the form of the local centre manifold $\mathcal J_1'$ yields the following graph form near $q_a$:
	\[
	\mathcal J_1' : x_1 = - \beta - \frac{\gamma \beta}{k} \rho_1^{k(k+1)} + \mathcal O\left(\rho_1^{k(k+1)+1}\right) .
	\]
	Hence, the flow on $\mathcal J_1'$ is described in the $\rho_1-$coordinate chart by
	\[
	\rho_1'\big|_{\mathcal J_1'} = - \frac{\gamma \beta}{k} \rho_1^{k(k+1) + 1} + \mathcal O\left(\rho_1^{k(k+1)+2}\right) ,
	\]
	which is positive (negative) when $\gamma < 0$ ($\gamma > 0$) in a neighbourhood of $q_a$.
	
	Now consider statement (iv). It follows from the properties of $\varphi(\rho_1)$ described in \appref{proof_of_lem_bifurcation} that for $\gamma \neq 0$ the function $\varphi(\rho_1)$ has a minimum $\rho_{1,c} > 0$ and a solution to $\varphi(\rho_{1,\ast}) = 0$ with $\rho_{1,\ast} \in (\rho_{1,c},\infty)$ corresponding to the equilibrium $P_1'$; see \figref{number_eqs_1}a. A direct calculation shows that the location of the minima $\rho_{1,c}$ is given by the following function of $\hat \mu$:
	\[
	\rho_{1,c}(\hat \mu) = \frac{k \hat \mu}{\delta (1 + k)} .
	\]
	Since $ P_1'(\hat \mu) : (-\beta , \rho_{1,\ast}(\hat \mu), 0, 0)$ where $\rho_{1,\ast}(\hat \mu) > \rho_{1,c}(\hat \mu)$, 
	the expressions for $\lim_{\hat \mu \to \infty}P_1'(\hat \mu)$ in both \eqref{P_lims_pos} and \eqref{P_lims_neg} (i.e. in  both cases $\gamma < 0$ and $\gamma > 0$) follow.
	
	Now restrict to case $\gamma < 0$. Then uniqueness of $P_1'(\hat \mu)$ on $\{\rho_1 > 0\}$ follows from the properties (i)-(iii) of the function $\varphi(\rho_1)$ with $\gamma < 0$ given in the proof of \lemmaref{lem_bifurcations_1}, \appref{proof_of_lem_bifurcation}; see \figref{number_eqs_1}a. To see that $P_1'(\hat \mu) \to q_a$ as $\hat \mu \to -\infty$, notice that $\rho_{1,\ast}(\hat \mu)$ satisfies the equilibrium condition
	\[
	-\gamma \beta + \hat \mu \rho_{1,\ast}(\hat \mu)^{k^2} - \delta \rho_{1,\ast}(\hat \mu)^{2k(k+1)} = 0 \qquad \implies \qquad
	\hat \mu = \frac{\gamma \beta + \delta \rho_{1,\ast}(\hat \mu)^{2k(k+1)}}{\rho_{1,\ast}(\hat \mu)^{k^2}} .
	\]
	By the equality on the right, $\rho_{1,\ast}(\hat \mu) \to 0^+$ as $\hat \mu \to -\infty$, implying that  $P_1'(\hat \mu) \to q_a$ as required.
	
	Finally, we consider case $\gamma > 0$. The existence of a second equilibrium $P_{s,1}'(\hat \mu)$ for $\hat \mu > \hat \mu_{sn}(\gamma)$ follows from the properties (i)-(iv) of the function $\varphi(\rho_1)$ with $\gamma > 0$ given in the proof of \lemmaref{lem_bifurcations_1}, noting that $\hat \mu_{sn}$ is the solution to $\varphi(\rho_{1,c}(\hat \mu_{sn})) = 0$; see \figref{number_eqs_1}b. To see that $P_{s,1}'(\hat \mu) \to q_a$ as $\hat \mu \to \infty$, note that when $0 < \hat \mu^{-1} \ll 1$ equation \eqref{eq_eqn_1} has a solution
	\[
	\rho_{s,1}(\hat\mu) = \gamma \beta \hat \mu^{-1/k^2} + \mathcal O\left(\hat \mu^{-2/k^2}\right) 
	\]
	corresponding to $P_{s,1}(\hat \mu)$. Since for $\gamma > 0$ we have $\rho_{s,1}(\hat\mu) \to 0^+$ as $\hat \mu \to \infty$, the expression in \eqref{P_lims_pos} follows.
\epr




\subsubsection*{$\mathcal K_2'$ rescaling chart}

Following a suitable desingularisation (division by $\nu_2^{k(k+1)}$), we obtain the following equations in the rescaling chart $\mathcal K_2'$:
\begin{equation}
\eqlab{mathcal_Ks2}
\begin{split}
x_2' &= \tilde f_2(x_2,\epsilon_2,\rho_2,\nu_2) - \frac{x_2}{2} \tilde g_2(x_2,\epsilon_2,\rho_2,\nu_2) , \\
\epsilon_2' &= -\epsilon_2 \left(\frac{1+2k}{2k}\right) \tilde g_2(x_2,\epsilon_2,\rho_2,\nu_2) , \\
\rho_2' &= \frac{\rho_2}{2k(1+k)} \tilde g_2(x_2,\epsilon_2,\rho_2,\nu_2) , 
\end{split}
\qquad \nu_2 \ll 1 ,
\end{equation}
where we have opted to use the fact that $\nu_2' = 0$ in order to consider the problem as a perturbed system in three variables. Here we have defined
\[
\begin{split}
\tilde f_2(x_2,\epsilon_2,\rho_2,\nu_2) &= \nu_2^{-2k(k+1)} f \left(\nu_2^{k(k+1)} x_2, \nu_2^{k+1} \epsilon_2 , \nu_2 \rho_2, \nu_2^{2k(k+1)} \right) , \\
\tilde g_2(x_2,\epsilon_2,\rho_2,\nu_2) &= \nu_2^{-k(k+1)} g \left(\nu_2^{k(k+1)} x_2, \nu_2^{k+1} \epsilon_2 , \nu_2 \rho_2, \nu_2^{2k(k+1)} \right) ,
\end{split}
\]
and in particular
\[
\tilde f_2 (x_2,\epsilon_2,\rho_2,0) = 1 + \rho_2^{k(k+1)} \left( \beta (\tau - \gamma) \epsilon_2^k + \tau x_2 - \delta \rho_2^{k(k+1)} \right) 
, \qquad
\tilde g_2 (x_2,\epsilon_2,\rho_2,0) = \beta \epsilon_2^k + x_2 .
\]

System \eqref{mathcal_Ks2} has a number of invariant subspaces, including the planes $\{\rho_2 = 0\}$, $\{\epsilon_2 = 0 \}$, and the 1-parameter family of sets
\begin{equation}
\eqlab{mathcal_A2}
\mathcal A_2 = \left\{(x_2, \epsilon_2, \rho_2) : \rho_2^{k(1+2k)} \epsilon_2^{k/(1+k)} = \hat \mu^{-1} \right\}
\end{equation}
obtained from the sets $\mathcal A_1$ in \eqref{mathcal_A} after an application of the transition map $\kappa_{12}'$ in \eqref{dash_transition_maps}. Notice that $\mathcal A_2$ limits to the union $\{\epsilon_2 = 0\} \cup \{ \rho_2 = 0\}$ as $\hat \mu \to \infty$, which is relevant insofar as we are interested in understanding the dynamics associated with the matching problem defined by the limits \eqref{matching_limit}. Hence in the following, we are motivated to consider the dynamics within each invariant plane $\{\epsilon_2 = 0\}$ and $\{ \rho_2 = 0\}$. Within $\{\rho_2 = 0\}$, the limiting system \eqref{mathcal_Ks2}$|_{\nu_2 = 0}$ is
\begin{equation}
\eqlab{mathcal_Ks2_rho}
\begin{split}
x_2' &= 1 - \frac{x_2}{2} \left( x_2 + \beta \epsilon_2^k \right) , \\
\epsilon_2' &= -\epsilon_2 \left(\frac{1+2k}{2k} \right) \left(x_2 + \beta \epsilon_2^k \right) .
\end{split}
\end{equation}
Within $\{\epsilon_2 = 0\}$, the limiting system \eqref{mathcal_Ks2}$|_{\nu_2 = 0}$ is
\begin{equation}
\eqlab{mathcal_Ks2_eps}
\begin{split}
x_2' &= 1 - \frac{x_2^2}{2} + \rho_2^{k(k+1)} \left(\tau x_2 - \delta \rho_2^{k(1+k)} \right) , \\
\rho_2' &= \frac{ \rho_2 x_2}{2 k (1 + k)} ,
\end{split}
\end{equation}
where we identify the equilibrium
\begin{equation}
\eqlab{P2_fam3}
P_2' : \left(0, 0, \delta^{-1/(2k(1+k))} \right) .
\end{equation}
Within the intersection $\{\rho_2 = \epsilon_2 = 0\}$, we identify two equilibria in \eqref{mathcal_Ks2}$|_{\nu_2 = 0}$, namely
\begin{equation}
\eqlab{q_pm}
q_+ : \left(\sqrt 2, 0, 0 \right) , \qquad q_- : \left(-\sqrt 2, 0, 0 \right) .
\end{equation}

\begin{remark}
	\remlab{reg_pert}
	Since we only identify isolated singularities in the limiting systems \eqref{mathcal_Ks2_rho} and \eqref{mathcal_Ks2_eps}, the perturbed systems \eqref{mathcal_Ks2}$|_{\{\rho_2 = 0\}}$ and \eqref{mathcal_Ks2}$|_{\{\epsilon_2 = 0\}}$ with $\nu_2 \ll 1$ are regular perturbations of systems  \eqref{mathcal_Ks2_rho} and \eqref{mathcal_Ks2_eps} respectively.
\end{remark}

Key dynamics of \eqref{mathcal_Ks2} are sketched in \figref{mathcal_K2}, and summarised in the following result. See also \figref{blowup_connection}a for a global representation.

\begin{figure}[h!]
	\centering
	\includegraphics[scale=0.5]{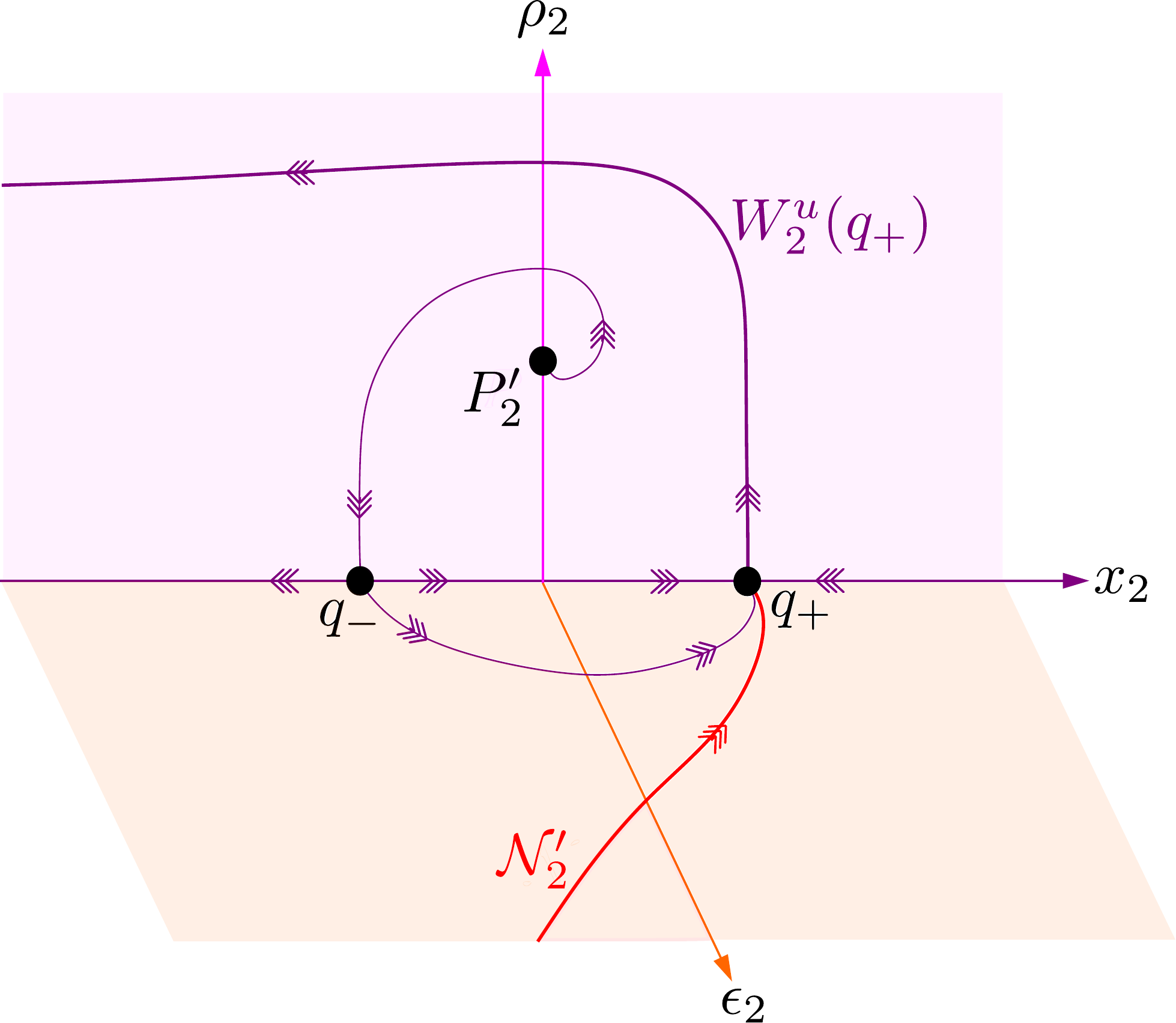}
	\caption{Dynamics in the rescaling chart $\mathcal K_2'$ in the limit $\nu_2 = 0$; see \lemmaref{mathfrak_K2}. Here the extension of the 1D manifold $\mathcal N_1'$ (denoted $\mathcal N_2'$ and contained in $\{\rho_2 = 0\}$) identified in \lemmaref{mathfrak_K1} connects to the hyperbolic saddle $q_+$. The point $P_2'$ is an unstable focus with a regular heteroclinic connection to $q_-$ and, consequently, the (global) unstable manifold $W_2^u(q_+) \subset \{\epsilon_2 = 0\}$ associated with $q_+$ is unbounded.
	}
	\figlab{mathcal_K2}
\end{figure}

\begin{lemma}
	\lemmalab{mathfrak_K2}
	The following holds for system \eqref{mathcal_Ks2}:
	\begin{enumerate}
		\item[(i)] For all $\nu_2 \geq 0$ sufficiently small the equilibria $q_\pm$ are hyperbolic saddles.
		\item[(ii)] Considered within the invariant plane $\{\epsilon_2 = 0\}$, the equilibrium $P_2'$ is an unstable focus with a regular heteroclinic connection to $q_-$.
		\item[(iii)] The extension of the centre manifold $\mathcal N_1'$ identified in \lemmaref{mathfrak_K1}, call it $\mathcal N_2'$, connects to $q_+$, and the flow on $\mathcal N_2'$ is forward asymptotic to $q_+$.
	\end{enumerate}
\end{lemma}

\bpr
	Statement (i) follows from direct calculations.
	
	Statement (ii): Linearising system \eqref{mathcal_Ks2_eps} and evaluating the eigenvalues at $P_2'$ we obtain $\lambda_{1,2} = \left(\tau \pm \sqrt{\tau^2 - 4 \delta} \right) / 2$, where $\tau, \delta > 0$, and $\tau^2 - 4 \delta < 0$ by \assumptionref{ass2}. Hence $P_2'$ is an unstable focus. To see that a regular heteroclinic connection between $P_2'$ and $q_-$ exists, first note that $q_-$ perturbs to a nearby equilibrium in system \eqref{mathcal_Ks2}$|_{\{\epsilon_2 = 0\}}$
	given by
	\[
	\tilde q_-(\nu_2) : \left(-\sqrt 2 + \mathcal O \left(\nu_2^{k(k+1)} \right), 0 ,0 \right)
	\]
	for all $\nu_2 > 0$ sufficiently small. Notice that $\tilde q_-(0) = q_-$, and that $\tilde q_-(\nu_2)$ is the image of the equilibrium $q_{in}$ identified in \eqref{Ks2_eqs} under the blow-up transformation $\mu = \nu_2^{2k(k+1)}$. Similarly, one can show that $P_2'$ perturbs to a nearby equilibrium $\tilde P_2'(\nu_2)$ corresponding the the image of the equilibrium $P_2(\mu)$ in \eqref{true_eq_Ks2} such that $\tilde P_2'(0) = P_2'$. Since by \lemmaref{K2s} (ii) there exists a regular heteroclinic connection between $P_2(\mu)$ and $q_{in}(\mu)$ for all $\mu \in (0,\mu_0)$ with $\mu_0 > 0$ sufficiently small, it follows that there is a regular heteroclinic connection between $\tilde P_2'(\nu_2)$ and $\tilde q_-(\nu_2)$ for all $\nu_2 > 0$ sufficiently small. Recalling that system \eqref{mathcal_Ks2}$|_{\{\epsilon_2 = 0\}}$, is a regular perturbation of system \eqref{mathcal_Ks2_eps} (see \remref{reg_pert}), a connection for $\nu_2 = 0$ (i.e. between $P_2'$ and $q_-$) follows by the requirement that both systems are topologically equivalent for sufficiently small $\nu_2 > 0$.

	
	To prove the statement (iii), notice that centre manifold $\mathcal N_1'$ described in \lemmaref{mathfrak_K1} lies within $\{\rho_1 = \nu_1 = 0\}$. In fact, making a power series ansatz leads to the following local graph form near $q_a$ in chart $\mathcal K_1'$:
	\begin{equation}
	\eqlab{mathcal_N1}
	\mathcal N_1' : x_1 = -\beta + \frac{1}{\beta k} \mu_1 + \mathcal O(\mu_1^2) .
	\end{equation}
	Hence its extension $\mathcal N_2' = \kappa'_{12}(\mathcal N_2') \subset \{\rho_2 = \nu_2 = 0\}$, with its dynamics therefore governed by \eqref{mathcal_Ks2_rho}. Simple calculations (some of which must be undertaken in additional charts $\bar x = \pm 1$) together with invariance for the flow on the sphere $B_\rho := \{\nu = 0\} \times \{(\bar x, \bar \epsilon, 0, \bar \mu) \in S^3\}$ and its boundary $\partial B_\rho$ show that $q_+$ is the only attractor for trajectories intersecting $\text{int} (B_{\rho})$. Hence the Poincar\'e-Bendixon theorem applies, and $\mathcal N_2'$ connects to $q_+$. Given that there are no equilibria on $\text{int}(B_\rho)$, the the flow on $\mathcal N_2'$ is regular and thus forward asymptotic to $q_+$ since
	\begin{equation}
	\eqlab{mu1_N1}
	\mu_1' \big|_{\mathcal N_1'} = \frac{1+2k}{\beta k} \mu_1^2 + \mathcal O(\mu_1^3) \geq 0 
	\end{equation}
	near $q_a$, and the orientation of the flow is preserved.	
\epr


It follows from \lemmaref{mathfrak_K2} (statement (ii) in particular) that the (global) unstable manifold $W^u_2(q_+)$ associated with the saddle $q_+$ is unbounded within $\{\epsilon_2 = 0\}$; see \figref{mathcal_K2}. In order to understand its extension, we must look in chart $\mathcal K_3'$.

\subsubsection*{$\mathcal K_3'$ chart}

In chart $\mathcal K_3'$ we obtain the following equations, after division by a common factor of $\nu_3^{k(k+1)}$:
\begin{equation}
\eqlab{mathcal_Ks3}
\begin{split}
\epsilon_3' &= \frac{\epsilon_3}{k} \left(\tilde f_3(\epsilon_3, \rho_3, \nu_3, \mu_3) - k \tilde g_3(\epsilon_3, \rho_3, \nu_3, \mu_3) \right) , \\
\rho_3' &= \frac{\rho_3}{k(1+k)} \left(\tilde f_3(\epsilon_3, \rho_3, \nu_3, \mu_3) + \tilde g_3(\epsilon_3, \rho_3, \nu_3, \mu_3) \right) , \\
\nu_3' &= - \frac{\nu_3}{2k(k+1)} \left(2 \tilde f_3(\epsilon_3, \rho_3, \nu_3, \mu_3) + \tilde g_3(\epsilon_3, \rho_3, \nu_3, \mu_3)  \right) , \\
\mu_3' &= \mu_3 \left(2 \tilde f_3(\epsilon_3, \rho_3, \nu_3, \mu_3) + \tilde g_3(\epsilon_3, \rho_3, \nu_3, \mu_3)  \right) ,
\end{split}
\end{equation}
where
\[
\begin{split}
\tilde f_3(\epsilon_3, \rho_3, \nu_3, \mu_3) &= \nu_3^{-2k(k+1)} f \left(-\nu_3^{k(k+1)} , \nu_3^{k+1} \epsilon_3 , \nu_3 \rho_3, \nu_3^{2k(k+1)} \mu_3 \right) , \\
\tilde g_3(\epsilon_3, \rho_3, \nu_3, \mu_3) &= \nu_3^{-k(k+1)} g \left(- \nu_3^{k(k+1)}, \nu_3^{k+1} \epsilon_3 , \nu_3 \rho_3, \nu_3^{2k(k+1)} \mu_3 \right) ,
\end{split}
\]
and in particular
\[
\tilde f_3  (\epsilon_3, \rho_3, 0, \mu_3) = \mu_3 + \rho_3^{k(k+1)} \left( \beta (\tau - \gamma) \epsilon_3^k - \tau - \delta \rho_3^{k(k+1)} \right) ,
\qquad \tilde g_3(\epsilon_3, \rho_3, 0, \mu_3) = \beta \epsilon_3^k - 1 .
\]

System \eqref{mathcal_Ks3} has a number of invariant subspaces, including the planes $\{\nu_3 =0 \}$, $\{\rho_3 = 0\}$, $\{\epsilon_3 = 0 \}$, and the 1-parameter family of sets
\begin{equation}
\eqlab{mathcal_A3}
\mathcal A_3 = \left\{\left(\epsilon_3, \rho_3, \nu_3, \hat \mu \epsilon_3^{k/(k+1)} \rho_3^{k(1+2k)} \right) :  \epsilon_3 \geq 0, \rho_3 \geq 0, \nu_3 \geq 0  \right\}
\end{equation}
obtained from the sets $\mathcal A_2$ in \eqref{mathcal_A2} after an application of the transition map $\kappa_{23}'$ in \eqref{dash_transition_maps}.

\

Since we are interested in the extension of the unstable manifold $W_2^u(q_+)$ contained within $\{\nu_2 = \epsilon_2 = 0\}$ in chart $\mathcal K_2'$ (see again \figref{mathcal_K2}), we are motivated to consider \eqref{mathcal_Ks3} within $\{\nu_3 = \epsilon_3 = 0\}$. Here we obtain the planar system
\begin{equation}
\eqlab{mathcal_Ks3_eps}
\begin{split}
\rho_3' &= \frac{\rho_3}{k(1+k)} \left(-1 + \mu_3 \left(1 - \tau \rho_3^{k(k+1)} - \delta \rho_3^{2k(k+1)} \right) \right) , \\
\mu_3' &= \mu_3 \left(-1 + 2 \mu_3 \left(1 - \tau \rho_3^{k(k+1)} - \delta \rho_3^{2k(k+1)} \right) \right) ,
\end{split}
\end{equation}
for which we identify two equilibria; the same equilibrium $q_-:(0,0,0,1/2)$ identified in \eqref{q_pm}, and $Q_-:(0,0,0,0)$.

Dynamics within $\{\nu_3 = \rho_3 = 0\}$ are important for understanding the limit $\hat \mu \to \infty$ associated with the (S1)-(S2) matching. Here we obtain the system
\begin{equation}
\eqlab{mathcal_Ks3_rho}
\begin{split}
\epsilon_3' &= \frac{\epsilon_3}{k} \left(\mu_3 + k \left(1- \beta \epsilon_3^k\right) \right)  , \\
\mu_3' &= \mu_3 \left(2 \mu_3 - 1 +  \beta \epsilon_3^k \right) ,
\end{split}
\end{equation}
which has also has two equilibria: $Q_-$ and $q_a:(\beta^{-1/k},0,0,0)$, where the latter coincides with $q_a$ in \eqref{mathfrak_K1_eqs}.

\begin{figure}[t!]
	\centering
	\includegraphics[scale=0.5]{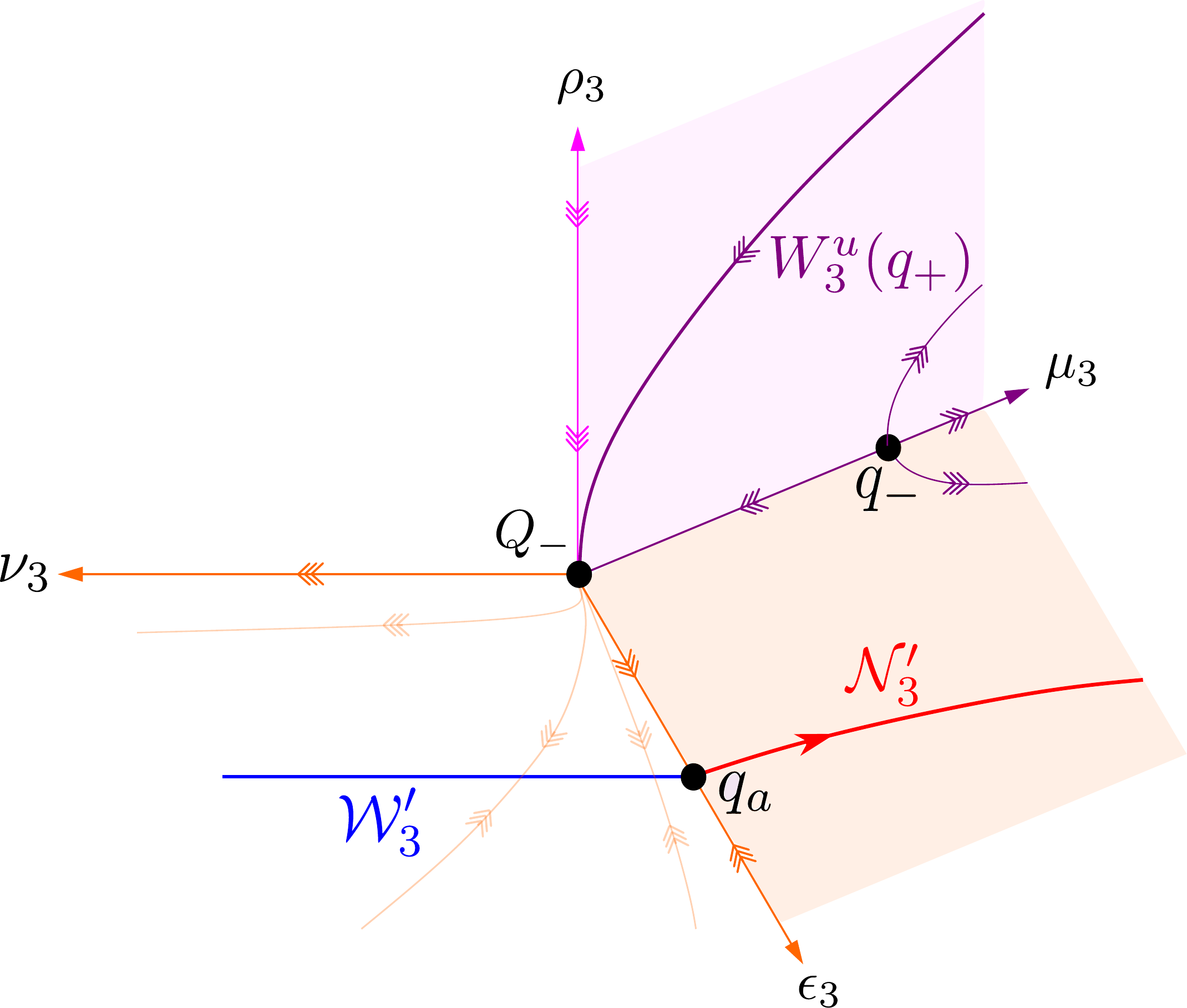}
	\caption{Dynamics in chart $\mathcal K_3'$. Here we sketch the case $\mu_3 \geq 0$ only, allowing us to represent part of the 4-dimensional dynamics of system \eqref{mathcal_Ks3} visible for $\epsilon, \rho_3, \nu_3, \mu_3 \geq 0$. The unstable manifold $W_3^u(q_+) \subset \{\epsilon_3 = \nu_3 = 0\}$ (purple) connects to a hyperbolic saddle $Q_-$ at the origin Also shown: the critical manifold $\mathcal W_3' \subset \{\rho_3 = \mu_3 = 0\}$ (blue), and the unique centre manifold $\mathcal N_3' \subset \{\rho_3 = \nu_3 = 0\}$ (red).}
	\figlab{mathcal_K3}
\end{figure}

Finally, we note that system \eqref{mathcal_Ks3} has a 1D critical manifold
\[
\mathcal W_3' = \left\{\left( \beta^{-1/k}, 0, \nu_3, 0 \right) : \nu_3 \geq 0  \right\}
\] 
within $\{\rho_3 = 0\} \cap \mathcal A_3 = \{\rho_3 = \mu_3 = 0\}$,
which coincides with the critical manifold $\mathcal W_1'$ in \eqref{mathcal_W1_dash}.

\

Key dynamics for system \eqref{mathcal_Ks2} are sketched in \figref{mathcal_K3}, and summarised in the following result. See also \figref{blowup_connection}a for a global representation.

\begin{lemma}
	\lemmalab{mathfrak_K3}
	The following holds for system \eqref{mathcal_Ks3}:
	\begin{enumerate}
		\item[(i)] The point $Q_-$ is a hyperbolic saddle. 
		\item[(ii)] The extension of the unstable manifold $W_2^u(q_+)$ in \figref{mathcal_K2}, call it $W_3^u(q_+)$, connects to $Q_-$ tangent to the $\rho_3-$axis, and the trajectory contained within $W_3^u(q_+)$ is forward asymptotic to $Q_-$.
		\item[(iii)] There exists a unique 3D attracting centre manifold $\mathcal M_3'$ at $q_a$ which coincides with the centre manifold $\mathcal M_1'$ identified in \lemmaref{mathfrak_K1} where domains overlap. In particular, we have
		\begin{equation}
		\eqlab{mathcal_M3}
		\mathcal M_3' \big|_{\{\rho_3 = 0\} \cap \mathcal A_3} = \mathcal W_3' , \qquad \mathcal M_3' \big|_{\{\nu_3 = \rho_3 = 0\}} = \mathcal N_3' ,
		\end{equation}
		where $\mathcal W_3'$ is a 1D critical manifold which is normally hyperbolic and attracting when considered within $\{\rho_3 = \mu_3 = 0\}$, and $\mathcal N_3'$ is the 1D (unique) centre manifold emanating from $q_a$ within $\{\nu_3 = \rho_3 = 0\}$. The manifolds $\mathcal W_3'$ and $\mathcal N_3'$ coincide with the manifolds $\mathcal W_1'$ and $\mathcal N_1'$ described in \lemmaref{mathfrak_K1}, and the slow-flow on $\mathcal M_3' \setminus \mathcal W_3'$ is increasing in the $\mu_3-$direction.
		\item[(iv)] Let $W^s(p)$ denote the stable manifold associated with a point $p \in \mathcal W_3'$, considered within the invariant plane $\{\rho_3 = 0\} \cap \mathcal A_3 = \{\rho_3 = \mu_3 = 0\}$. Then each one-sided stable manifold $W^s(p) \cap \{\epsilon_3 \in [0, \beta^{-1/k}]\}$ contains a regular heteroclinic connection between $p \in \mathcal W_3'$ and $Q_-$.
	\end{enumerate}
\end{lemma}

\bpr
	Statement (i) follows after linearisation of system \eqref{mathcal_Ks3} at $Q_-$.
	
	In order to prove the statement (ii), consider the (invariant) compact subset $\{\nu = 0\} \times \{(\bar x, 0, \bar \rho, \bar \mu) \in S^3 \}$ of the blow-up sphere $\{\nu = 0 \} \times S^3$. It follows from \lemmaref{mathfrak_K2} (ii) that there are no (regular) limit cycles within this region, since such a cycle would need to enclose the point $P_2'$ which is known by  \lemmaref{mathfrak_K2} (ii) to have a regular heteroclinic connection to $q_-$. Hence the point $Q_-$ constitutes the only attractor for (non-equilibrium) trajectories intersecting the interior of this region, and it follows by an application of the Poincar\'e-Bendixon theorem that the unstable manifold $W^u(q_-)$ connects to $Q_-$.


	Statement (iii): existence of a centre manifold $\mathcal M_3'$ satisfying the relevant properties
	follows from \lemmaref{mathfrak_K1} after an application of the inverse transition map associated with $\kappa'_{31}$ in \eqref{dash_transition_maps}.

	Finally, statement (iv) follows by an application of the Poincar\'e-Bendixon theorem: consider the compact region bounded by $\mathcal W$ and the three invariant subspaces
	\[
	\left\{(0,\epsilon_3,0,0) : \epsilon_3 \in [0,\beta^{-1/k}] \right\}, \ \ 
	\left\{(0,0,0,\nu_3) : \nu_3 \geq 0 \right\} , \ \ 
	\left\{(0,\epsilon_3,0) : \epsilon_3 \in \left[0, \beta^{-1/k}\right] \right\},
	\]
	the last of which is written in chart $\mathcal K_3$ coordinates $(r_3,\epsilon_3,\rho_3)$. Since $W^u(p_-) \cap \{\bar \rho = 0, \bar \epsilon \geq 0\}$ is unique, connects to $p_a$ and the flow is regular on the interior of the domain (see e.g. \figref{blowup_connection}), each $W^s(p) \cap \{\epsilon_3 \leq \beta^{-1/k}\}$ with $p \in \mathcal W_3'$ must connect with $Q_-$ ($Q_-$ being the only remaining repeller).
\epr


\subsection{$\mu=0$ singular cycles}
\applab{singular_cycles}

We turn our attention now to the family of singular cycles `connecting' the regular and relaxation-type cycles occurring in regimes (S1) and (S2) respectively. To do so, we consider the dynamics for the matching problem \eqref{matching_limit} shown in \figref{blowup_connection}a. Singular cycles can be constructed with either three, four, or five distinct orbit segments as shown in \figref{blowup_connection}b. We first identify the smallest and largest cycles
\begin{equation}
\eqlab{bound_cycles}
\Gamma^s = \Gamma^1 \cup \Gamma^2 \cup \Gamma^{3,1} , \qquad
\Gamma^l = \Gamma^1 \cup \Gamma^2 \cup \Gamma^{3,2} \cup \Gamma^{4,2} \cup \Gamma^5 ,
\end{equation}
where the subscripts `$s$' and `$l$' stand for `small' and `large' respectively, we have permitted a slight abuse of notation in \eqref{int_cycs_def} by reusing our earlier notation $\Gamma^j$, and the superscript notation $\Gamma^{i,j}$ is explained below. Explicitly, we define
\begin{equation}
\eqlab{int_cycs_def}
\begin{split}
\Gamma^1 &= \mathcal N' , \\
\Gamma^2 &= W^u(q_+) , \\
\Gamma^{3,1} &= \left\{(0,\epsilon_3,0,0) : \epsilon_3 \in [0,\beta^{-1/k}] \right\} , \qquad \text{chart } \mathcal K_1' , \\
\Gamma^{3,2} &= \left\{(0,\epsilon_3,0) : \epsilon_3 \in \left[0, \beta^{-1/k}\right] \right\} , \qquad \ \ \text{chart } \mathcal K_3 , \ \mu = 0, \\
\Gamma^{4,2} &= \mathcal W , \\
\Gamma^5 &=  \left\{(0,0,0,\nu_3) : \nu_3 \geq 0 \right\} , \qquad \qquad  \ \ \ \  \text{chart } \mathcal K_3' , \\
\end{split}
\end{equation}
see again \figref{blowup_connection}b. The singular cycles $\Gamma^s$ and $\Gamma^l$ bound a smooth family of singular orbits $\{\Gamma(s)\}_{s\in\mathbb R_+}$ which can be parameterised as follows: given $s \in \mathbb R_+$, define a point $T(s) = (-\beta, 0, s, 0) \in \mathcal W_1'$ in chart $\mathcal K_1'$ coordinates. It follows from \lemmaref{mathfrak_K3} (iv) that
\[
\Gamma^3(s) = W^s(T(s)) \cap \left\{x_1 \leq -\beta \right\} , \qquad \text{chart } \mathcal K_1' ,
\]
is a regular connection from $T(s)$ to $Q_-$, where $W^s(T(s)) \subset \{\rho_1 = \mu_1 = 0\} = \{\rho_1 = 0\} \cap \mathcal A_1$ denotes the stable manifold of $T(s)$ considered as an equilibrium in $\mathcal W'_1$. We define $\Gamma^4(s)$ at the same $s-$value by
\[
\Gamma^4(s) = \left\{(-\beta, 0, \nu_1,0) : \nu_1 \in [0,s] \right\} , \qquad \text{chart } \mathcal K_1' .
\]
In particular, we have
\[
\lim_{s \to 0^+} \Gamma^3(s) = \Gamma^{3,1} , \ \ \ \  \lim_{s \to \infty} \Gamma^3(s) = \Gamma^{3,2} , \ \ \ \ 
\lim_{s \to 0^+} \Gamma^4(s) = q_a , \ \ \ \  \lim_{s \to \infty} \Gamma^4(s) = \Gamma^{4,2},
\]
allowing us to define the family $\{\Gamma(s)\}_{s\in\mathbb R_+}$ by
\begin{equation}
\eqlab{int_cycles}
\Gamma(s) := \left(\Gamma^1 \cup \Gamma^2 \cup \Gamma^3 \cup \Gamma^4\right)(s) , \qquad s \in \mathbb R_+ ,
\end{equation}
where
\[
\lim_{s \to 0^+} \Gamma(s) = \Gamma^s , \qquad \lim_{s \to \infty} \Gamma(s) \cup \Gamma^5 = \Gamma^l ,
\]
follows by construction.

\begin{remark}
	The family \eqref{int_cycles} is distinct from the family of singular cycles $\{\Gamma^{ro}(\rho_3) : \rho_3 \in (0,L)\}$ defined in \eqref{large_ro}, with $\Gamma^l$ at the boundary between the two families.
\end{remark}

\section{Proof of \thmref{thm_connection}}
\applab{proof_of_thm_connection}


In this section we prove \thmref{thm_connection}. To do so we define two transverse sections $\Sigma^L$ and $\Sigma^R$ as in \figref{bfb_sections} and study the maps 
\begin{equation}
\eqlab{flow_maps}
\Pi_f : \Sigma^R \to \Sigma^L \qquad \text{and} \qquad \Pi_b : \Sigma^R \to \Sigma^L ,
\end{equation}
induced by forward and backward flow respectively. We must consider the flow induced by either system $\{$\eqref{main_extended}$, \ \mu' = 0\}$ or system \eqref{main_extended2}, changing charts whenever necessary. 

\begin{figure}[h!]
	\begin{center}
		\includegraphics[scale=1.75]{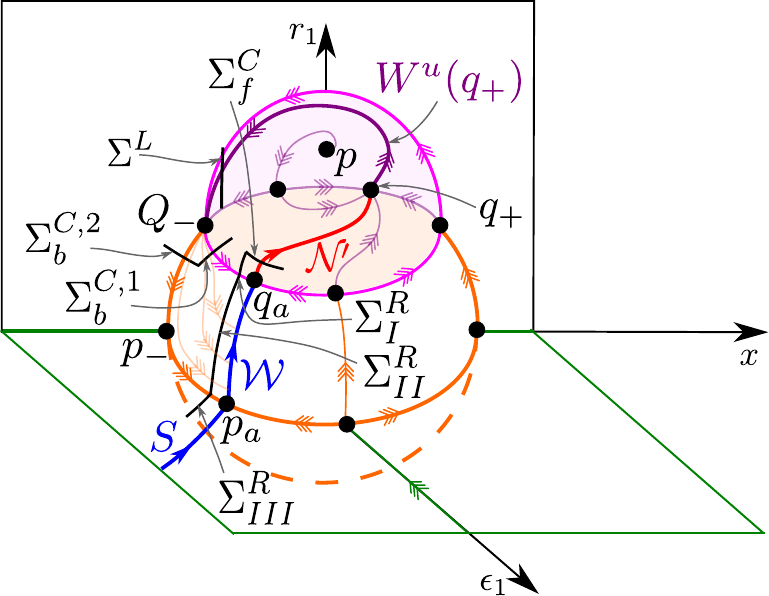}
		\caption{Setup for the proof of \thmref{thm_cycs}. The blow-up spheres are shown in chart $K_1$ coordinates, where the plane $\{r_1=0\}$ (green) corresponds to the face of the blow-up cylinder, c.f. \figref{blowup_connection}a-b. In order to prove \thmref{thm_cycs}, we evolve trajectories with initial conditions on $\Sigma^R = \Sigma^R_I \cup \Sigma^R_{II} \cup \Sigma_{III}^R$ in forward and backward time and match them on $\Sigma^L$. Also shown: additional sections $\Sigma^C_f$, $\Sigma^{C,1}_{b}$ and $\Sigma_b^{C,2}$ used to break the flow maps up into smaller component maps. See \figref{K3_sections} for local figures centred on the points $Q_-$ and $p_-$ in chart $\mathcal K_3'$ and chart $\mathcal K_3$ respectively.}
		\figlab{bfb_sections}
	\end{center}
\end{figure}


\

We begin with formal definitions for $\Sigma^L$ and $\Sigma^R$. 
We will work primarily in charts $\mathcal K_3$ and $\mathcal K_3'$, in which we write coordinates as follows:
\[
\text{Extended chart } \mathcal K_3 : (r_3, \epsilon_3, \rho_3, \mu) ; \qquad \text{Chart }\mathcal K_3' : (\epsilon_3, \rho_3, \nu_3, \mu_3) . 
\]
Note that chart $\mathcal K_3$ is now considered within an \textit{extended} space, where we consider $\mu$ as a variable. Coordinates in $\mathcal K_3$ and $\mathcal K_3'$ are smoothly related where their domains overlap by the following transformation, obtained by composing the transition map $\kappa_{23}^{-1}$ (see \eqref{mathcal_transition_map}) with the blow-up transformation \eqref{mathfrak_K3_coord}:
\begin{equation}
\eqlab{33_trans}
r_3 = \nu_3^{-2k(k+1)} , \qquad \epsilon_3 = \epsilon_3 , \qquad \rho_3 = \rho_3, \qquad \mu = \nu_3^{2k(k+1)} \mu_3 , \qquad \nu_3 > 0 .
\end{equation}
Note that we permit a slight abuse of notation in the above by allowing $\epsilon_3$, $\rho_3$ denote coordinates in both charts $\mathcal K_3$, $\mathcal K_3'$; this is not problematic since these coordinates coincide on their overlapping domain. Now define
\begin{equation}
\eqlab{Sigma_R}
\Sigma^L := \left\{(\epsilon_3, \rho_3 , \nu_3, \sigma_L) : \epsilon_3 \in [0,c_L], \rho_3 \in [0,R_L] ,\nu_3 \in [0, \Lambda_L] 
\right\} , \qquad \mathcal K_3',
\end{equation}
and $\Sigma^R := \Sigma^{R}_{I} \cup \Sigma^{R}_{II} \cup \Sigma^{R}_{III}$, where
\begin{equation}
\eqlab{Sigma_L}
\begin{aligned}
\Sigma^{R}_I &= \left\{(c_R, \rho_3, \nu_3, \mu_3) : \rho_3 \in [0,R_I], \nu_3 \in [0,\Lambda_I] , \mu_3 \in [0,\sigma_I] \right\} , && \text{ chart }  \mathcal K_3' , \\
\Sigma^{R}_{II} &= \left\{(c_R, \rho_3, \nu_3, \mu_3) : \rho_3 \in [0,R_{II}], \nu_3 \in [\Lambda_I,\Lambda_{II}], \mu_3 \in [0,\sigma_{II}] \right\} , && \text{ chart }  \mathcal K_3' , \\
\Sigma^{R}_{III} &= \left\{(c_R, r_3, \rho_3, \mu) : r_3 \in [0,L_{III}], \rho_3 \in [0,R_{III}], \mu \in [0,\mu_0] \right\} , && \text{ chart }  \mathcal K_3 ,
\end{aligned}
\end{equation}
where the constants $R_I, R_{II}, R_{III}, \Lambda_I, \sigma_I, \sigma_{II}, L_{III} > 0$ are fixed and, in the following, often assumed fixed sufficiently small for the validity of local results. The constant $\Lambda_{II} > 0$ can be fixed arbitrarily large and in particular large enough so that $\Sigma^R$ is connected as in \figref{bfb_sections}. Note also that $\Sigma_I^R$ ($\Sigma_{III}^R$) is properly visible in chart $\mathcal K_3'$ ($\mathcal K_3'$) but not in $\mathcal K_3$ ($\mathcal K_3'$); see also \figref{K3_sections}. The section $\Sigma_{II}^R$ is properly visible in both $\mathcal K_3'$ and $\mathcal K_3$.

\

In order to state our main result for this section, we combine families of singular cycles for $\mu > 0$ and $\mu = 0$ defined in \eqref{large_ro} and \eqref{int_cycles} respectively by introducing a new parameter $\theta \in [0,1]$ such that
\begin{itemize}
	\item \SJnew{For $\theta \in [0,1/3)$, the map $\theta \mapsto \mu_3 \in (0,\sigma_I]$ is a bijection.} 
	\item For $\theta \in [1/3,2/3]$ the map \SJnew{$\theta \mapsto s \in [0,\infty]$} 
	is a bijection \SJnew{such that $\theta = 1/3 \mapsto 0$ and $\theta = 2/3 \mapsto \infty$}, so that there exists a family of singular cycles
	\begin{equation}
	\eqlab{cycs1}
	\{\bar \Gamma(\theta) : \theta \in [1/3,2/3]\}
	\end{equation}
	in 1-1 correspondence with the $\mu = 0$ singular cycles defined in \eqref{int_cycles}, \appref{singular_cycles}. \SJnew{In particular, $\bar \Gamma(1/3) = \Gamma^s$ and $\bar \Gamma(2/3) = \Gamma^l$.}
	\item For $\theta \in (2/3,1)$ the map $\theta \mapsto \rho_3 \in (0,R)$ is a bijection, so that there exists a family of singular \SJnew{relaxation} cycles
	\begin{equation}
	\eqlab{cycs2}
	\{\bar \Gamma(\theta) : \theta \in (2/3,1)\}
	\end{equation}
	in 1-1 correspondence with the $\mu > 0$ singular \SJnew{relaxation} cycles defined in \eqref{large_ro}, \appref{singular_cycles_pos}.
\end{itemize}
\SJnew{Note that there are no singular cycles corresponding to $\theta \in [0,1/3)$, since there are no singular cycles intersecting $\Sigma^R \cap \{\nu_3=0, \mu_3 \in (0,\sigma_I]\}$ (the part of $\Sigma^R$ contained within the shaded orange plane in \figref{K3_sections}a).}

\begin{figure}[t!]
	\centering
	\subfigure[]{\includegraphics[width=.49\textwidth]{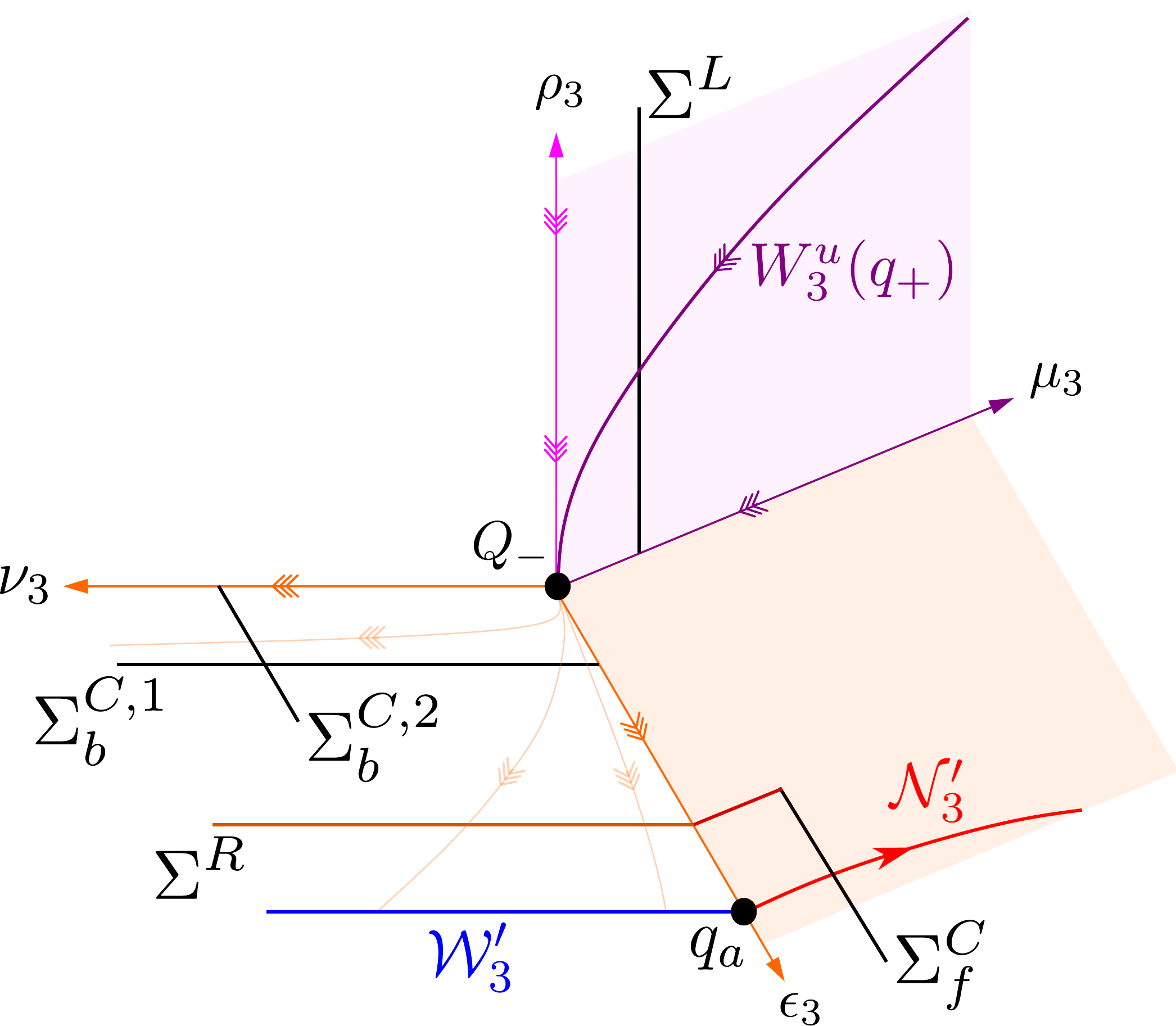}}
	\subfigure[]{\includegraphics[width=.49\textwidth]{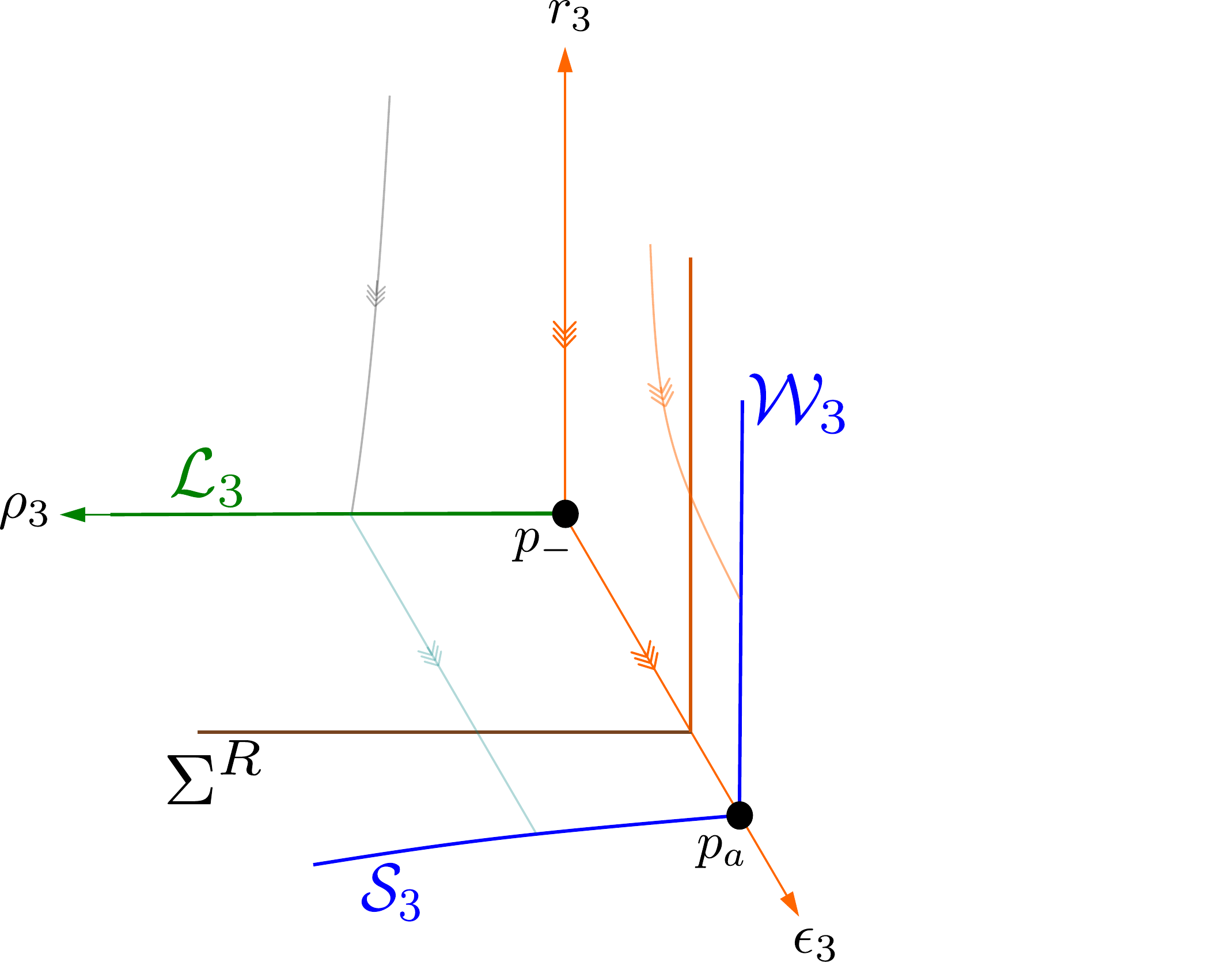}}
	\caption{Local sketch of the setup for the proof of \thmref{thm_cycs}. In (a): local setup near $Q_-$ in chart $\mathcal K_3'$. In (b): local setup near $p_-$ in chart $\mathcal K_3$. Singular cycles $\bar \Gamma(\theta)$ with $\theta \in [1/3,2/3]$ intersect $\Sigma^R \cap \{\rho_3=\mu_3= 0\}$ or $\Sigma^R \cap \{\rho_3 = \mu = 0\}$, which are shown in orange in (a) and (b) respectively. Singular cycles $\bar \Gamma(\theta)$ with $\theta \in (2/3,1)$ intersect $\Sigma^R \cap \{r_3 = \mu = 0\}$, which is shown in brown in (b). \SJnew{There are no singular cycles corresponding to orbits $\bar \Gamma(\theta)$ with $\theta \in [0,1/3)$ which intersect $\Sigma^R \cap \{\rho_3 = \nu_3 = 0\}$ (the red part of $\Sigma^R$ within the shaded orange plane in (a)).} Note that the orange segments $\Sigma^R \cap \{\rho_3 = \mu_3 = 0\}$ in (a) and $\Sigma^R \cap \{\rho_3 = \mu = 0\}$ in (b) are the same by \eqref{33_trans}. See \figref{bfb_sections} for a global comparison.
	}
	\figlab{K3_sections}
\end{figure}

\

We are now in a position to state our main result for this section, which may be considered a more comprehensive analogue of \thmref{thm_connection}.

\begin{theorem}
	\thmlab{thm_cycs}
	Fix $\kappa \in (0,1/3)$ arbitrarily small. Then there exists $\epsilon_0 > 0$ such that for all $\epsilon \in (0, \epsilon_0)$ there exists a parameterised family of stable limit cycles
\begin{equation}
\eqlab{fam2}
\theta \mapsto \left(\mu(\theta, \epsilon ), \bar \Gamma(\theta,\epsilon )\right) , \qquad \theta \in (0,1) ,
\end{equation}
which is continuous in $(\theta,\epsilon)$ and such that:
	\begin{enumerate}
		\item[(i)] The map $\theta\mapsto \mu(\theta,\epsilon)$ is $C^1-$smooth for each $\epsilon\in (0,\epsilon_0)$, satisfying 
		\begin{equation}
		\eqlab{monotonicity}
		\frac{\partial \mu}{\partial \theta} > 0 , \qquad \theta \in (0,1) .
		\end{equation}
		\item[(ii)] Fix $\theta \in (0,1/3)$. Then $\mu(\theta,\epsilon) \to 0^+$ as $\epsilon \to 0$. If $\theta \in (0,1/3-\kappa]$ 
		then $\bar \Gamma(\theta,\epsilon)$ is a regular cycle (i.e.~not a relaxation oscillation), and
		\[
		\mu(\theta,\epsilon) = \mathcal O\left(\epsilon^{k/(k+1)} \right) .
		\]
		\item[(iii)] Fix $\theta \in [1/3, 2/3]$. Then $\mu(\theta,\epsilon) \to 0^+$ as $\epsilon \to 0$, and $\bar \Gamma(\theta,\epsilon)$ is a relaxation oscillation such that $\bar \Gamma(\theta,\epsilon) \to \bar \Gamma(\theta)$ in Hausdorff distance as $\epsilon \to 0$, where $\bar \Gamma(\theta)$ is a member of the family \eqref{cycs1}.
		If $\theta \in [1/3+\kappa, 2/3-\kappa]$, 
		then
		\[
		\mu(\theta,\epsilon) = \mathcal O\left(\epsilon^{k/(2k+1)} \right) .
		\]
		\item[(iv)] Fix $\theta \in (2/3, 1)$. Then $\mu(\theta,\epsilon) \to \mu(\theta,0) > 0$ 
		as $\epsilon \to 0$, and $\bar \Gamma(\theta,\epsilon)$ is a relaxation oscillation 
		such that $\bar \Gamma(\theta,\epsilon) \to \bar \Gamma(\theta)$ in Hausdorff distance as $\epsilon \to 0$, 
		where $\bar \Gamma(\theta)$ is a member of the family \eqref{cycs2}. If $\theta \in [2/3+\kappa, 1-\kappa]$, 
		then
		\[
		\mu(\theta,\epsilon) = \mathcal O(1) ,
		\]
		\SJnew{and $\bar \Gamma(\theta, \epsilon)$ is contained in the family of relaxation oscillations described in \thmref{thm_ro}.}
		\item[(v)] Any limit cycle passing sufficiently close to the union of critical manifolds $S \cup \mathcal W$ (see again \figref{bfb_sections}) is a member of $\bar \Gamma(\theta,\epsilon)$.
	\end{enumerate}
\end{theorem}

\thmref{thm_connection} follows directly from \thmref{thm_cycs}. 
The limit cycles $\bar \Gamma(\theta,\epsilon)$ with $\theta \in [1/3,2/3]$ described in statement (iii) constitute the `connection' between regular oscillations in scaling regime (S1) (described in statement (ii)) and relaxation oscillations in the scaling regime (S2) (described in statement (iv)). It follows from statement (iv) that the relaxation oscillations described in \thmref{thm_ro} (iii) are contained within the family defined by \eqref{fam2}.

\


In order to prove \thmref{thm_cycs}, we need to understand the flow maps $\Pi_{f/b}$ in \eqref{flow_maps}. For a given input $q \in \Sigma^R$, we write
\begin{equation}
\eqlab{maps}
\Pi_{f/b}(q) = 
\left(E_{f/b}(q) , P_{f/b}(q) , N_{f/b}(q) , \sigma_L \right) , \\
\end{equation}
where the right-hand-side is expressed in chart $\mathcal K_3'$ coordinates $(\epsilon_3,\rho_3,\nu_3,\mu_3)$. Since both $\mathcal K_3'$ and $\mathcal K_3'$ are required to cover $\Sigma^R$, the input $q \in \Sigma^R$ will be considered in either $\mathcal K_3'$ or $\mathcal K_3$. We look for limit cycles corresponding to solutions of the fixed point equation
\begin{equation}
\eqlab{fixed_point1}
\left(\Pi_{f} - \Pi_{b}\right)(q) = \left(E_{f} - E_{b} , P_{f} - P_{b} , N_{f} -  N_{b} , 0 \right)(q) = (0,0,0,0) .
\end{equation}
Since $\epsilon$ and $\mu$ are constants of the motion for system \eqref{main_extended2}, the following define invariant sets in chart $\mathcal K_3'$: 
\begin{equation}
\eqlab{consts_motion_m3}
\begin{split}
\mathcal V_\epsilon' &:= \left\{(\epsilon_3, \rho_3, \nu_3, \mu_3) \in \mathbb R_+^3 \times \mathbb R: \nu_3^{2 (k+1)^2} \rho_3^{(k+1) (2k+1)} \epsilon_3 = \epsilon \right\} , \\
\mathcal V_\mu' &:=\left\{(\epsilon_3, \rho_3, \nu_3, \mu_3) \in \mathbb R_+^3 \times \mathbb R :  \nu_3^{2k(k+1)} \mu_3 = \mu \right\} , \\
\mathcal V' &:= \mathcal V_\epsilon' \cap \mathcal V_\mu .
\end{split}
\end{equation}
By restricting \eqref{fixed_point1} to $\mathcal V'$ in particular, it is sufficient to solve \textit{any one} of the equations
\[
\left(E_{f} - E_{b}\right)(q) = 0 , \qquad \left(P_{f} - P_{b}\right)(q) = 0 , \qquad \left(N_{f} - N_{b}\right)(q) = 0 ,
\]
and in the following we will focus on solutions for the second equation $(P_f - P_b)(q) = 0$. 
Inputs $q\in \Sigma^R$ expressed in chart $\mathcal K_3'$ can also be restricted to the invariant sets 
\eqref{consts_motion_m3}, while inputs expressed $\mathcal K_3$ can be restricted to $\epsilon = const.$ sets defined by
\begin{equation}
\eqlab{consts_motion_3}
\mathcal V_\epsilon := \left\{(r_3, \epsilon_3, \rho_3,\mu) \in \mathbb R_+^3 \times \mathbb R :  \rho_3^{(1+k)(1+2k)} r_3 \epsilon_3 = \epsilon \right\} .
\end{equation}
By restricting initial conditions in $\Sigma^R_I$, $\Sigma_{II}^R$ to $\mathcal V'$ and initial conditions in $\Sigma_{III}^R$ to $\mathcal V_\epsilon$, we can define the \textit{restricted maps}
\begin{equation}
\eqlab{restricted_maps}
\mathcal P_{f/b} : \Sigma^R \times [0,\epsilon_0] \times [0,\mu_0] \to \Sigma^L , \ \ \
\mathcal P_{f/b}(q,\epsilon, \mu) = 
\begin{cases}
P_{f/b}(q)\big|_{\mathcal V'} , & q \in (\Sigma_I^R \cup \Sigma_{II}^R) \cap \mathcal V' , \\
P_{f/b}(q) \big|_{\mathcal V_\epsilon}, & q \in \Sigma_{III}^R \cap \mathcal V_\epsilon .
\end{cases}
\end{equation}
Since it will frequently be useful to consider sets of inputs on $\Sigma^R_I$, $\Sigma_{II}^R$ and $\Sigma^R_{III}$ independently, we also introduce the following (further restricted) maps:
\begin{equation}
\eqlab{restricted_maps2}
\mathcal P_{f/b}^I := \mathcal P_{f/b}\big|_{\Sigma_I^R} , \qquad 
\mathcal P_{f/b}^{II} := \mathcal P_{f/b}\big|_{\Sigma_{II}^R} , \qquad 
\mathcal P_{f/b}^{III} := \mathcal P_{f/b}\big|_{\Sigma_{III}^R} ,
\end{equation}
noting that a complete description of all three amounts to a complete description of $\mathcal P_{f/b}$ in \eqref{restricted_maps}. Following the restrictions above, our approach in the following is to solve \eqref{fixed_point1} by solving 
the corresponding restricted fixed point equation
\begin{equation}
\eqlab{fixed_point}
\mathcal D(q,\epsilon,\mu) := \mathcal P_f(q, \epsilon, \mu) - \mathcal P_b(q, \epsilon, \mu) = 0 
\end{equation}
for $\mu$. 

\

The structure of the proof is as follows: (i) characterise the maps $\Pi_{f/b}$; (ii) solve the fixed point equation $(P_f-P_b)(q)=0$ for separate sets of initial conditions $q \in \Sigma^R_I$, $q\in \Sigma_{II}^R$ and $q \in \Sigma_{III}^R$; (iii) restrict solutions obtained for separate sets of inputs in step (ii) to the invariant sets $\mathcal V'$ and $\mathcal V_\epsilon$ defined in \eqref{consts_motion_m3} and \eqref{consts_motion_3}. This yields solutions to \eqref{fixed_point} for three separate cases depending on whether $q \in \Sigma_I^R$, $q \in \Sigma_{II}^R$ or $q \in \Sigma_{III}^R$. Step (iv) is to combine the solutions obtained in (iii) using a global parameterisation in $\theta \in [0,1]$, thus obtaining a  \eqref{fixed_point1} function $\mu(\theta,\epsilon)$ which solves the global fixed point equation \eqref{fixed_point}. The final step (v) (\appref{proof_of_thm_cycs}) is to combine findings in steps (i)-(iv) and prove \thmref{thm_cycs}. We begin with the map $\Pi_f$.

	\begin{lemma}
		\lemmalab{Pf}
		The map $\Pi_f$ is well-defined and at least $C^1-$smooth. In particular we have
		\begin{equation}
		\eqlab{Pf}
		P_f(q) = \gamma_L + J(q) , 
		\end{equation}
		where $\gamma_L \in (0,R_L)$ is defined by the intersection $W^u(q+) \cap \Sigma^L = \{(0, \gamma_L, 0, \sigma_L )\}$, and 
		\[
		J(q) =
		\begin{cases}
		J_{I,II}(\rho_3,\nu_3,\mu_3) , & q = (c_R,\rho_3, \nu_3, \mu_3) \in \Sigma_I\cup \Sigma_{II}^R , \\
		J_{III}(r_3,\rho_3,\mu) , & q = (r_3, c_R,\rho_3, \mu) \in \Sigma_{III}^R ,
		\end{cases}
		\] 
		is continuous such that $J_{I,II}(0,\nu_3,0) = J_{III}(r_3,0,0) = 0$. Moreover, $\Pi_f$ is a strong contraction when restricted to the invariant sets defined by \eqref{consts_motion_m3} and \eqref{consts_motion_3}. More precisely, there exists $\epsilon_0, \mu_0 >0$ such that for each fixed $(\epsilon,\mu) \in (0,\epsilon_0) \times (0,\mu_0)$ \SJnew{the restricted maps defined in \eqref{restricted_maps} and \eqref{restricted_maps2} satisfy} 
		\begin{equation}
		\eqlab{rho1_contraction}
		\frac{\partial \mathcal P^I_f}{\partial \mu_3} = \mathcal O\left(e^{-K_I/\mu_3}\right) , \qquad 
		\frac{\partial \mathcal P^{II}_f}{\partial \nu_3} = \mathcal O\left(e^{-K_{II}/\mu}\right) , \qquad
		\frac{\partial \mathcal P^{III}_f}{\partial \rho_3} = \mathcal O\left(e^{-K_{III}\rho_3^{k(k+1)}/\epsilon \mu}\right) ,
		\end{equation}
		where $K_I, K_{II}, K_{III} > 0$ are constant and the values of $\mu_3$ and $\rho_3$ in the first and third estimates respectively are fixed within the intervals
		\begin{equation}
		\eqlab{mu_rho_bounds}
		\mu_3 \in \left[\mu \Lambda_I^{-2k(k+1)} , \sigma_I \right] , \qquad 
		\rho_3 \in \left[\left(\epsilon \left(c_R L_{III}\right)^{-1}\right)^{1/(k+1)(2k+1)} , L_{III} \right] .
		\end{equation}
	%
	%
	\end{lemma}

\bpr
	Consider $\Pi_f:\Sigma^R \to \Sigma^L$ 
	as a composition $\Pi_f = \Pi_f^2 \circ \Pi_f^1$ of `smaller' component maps $\Pi_f^1 : \Sigma^R \to \Sigma^C_f$ and $\Pi_f^2:\Sigma_f^C \to \Sigma^L$, where
	\[
	\Sigma^C_f = \left\{(x_1,\rho_1,\nu_1,\sigma_f) : |x_1 + \beta| \leq \alpha_f, \rho_1 \in [0,R_f], \nu_1 \in [0,\Lambda_f] \right\} , \qquad \text{chart } \mathcal K_1' ,
	\]
	see \figref{K3_sections} and also \figref{bfb_sections} for a global comparison. Note that $\Sigma_f^C$ is also visible in $\mathcal K_3'$; we define it here in chart $\mathcal K_1'$ only in order to simplify some of the calculations. 
	We consider each component map in turn.
	
	\
	
	\textit{The map $\Pi_f^2$}. 
	We work in chart $\mathcal K_2'$ and consider the map $\Pi_f^2 \circ \kappa_{12}'$, where the transition map $\kappa_{12}'$ is given in \eqref{dash_transition_maps}. Consider the flow within a tubular neighbourhood $\mathcal B_f$ of the connected component of $\mathcal N' \cup W^u(q_+)$ between $\mathcal N' \cap \Sigma^C_f$ and $W^u(q_+) \cap \Sigma^L$; see \figref{bfb_sections}. Assuming $\mathcal B_f$ is sufficiently narrow, 
	the flow in $\mathcal B_f$ is regular everywhere except within a neighbourhood of the point $q_+$, which is a 
	hyperbolic saddle for the limiting system \eqref{mathcal_Ks2}$|_{\{\nu_2=0\}}$. 
	Standard techniques exist which allow one to 
	derive sufficiently precise estimates for the local transition (Dulac) map, and we refer the reader to e.g. \cite{Gucwa2009,Kosiuk2011,Kristiansen2019,Kristiansen2019c,Krupa2001a} where similar transitions occur. We omit the \SJnew{details} here, noting that similar arguments applied to system \eqref{mathcal_Ks2}$|_{\{\nu_2=0\}}$ are sufficient to show that the restricted map $(\Pi^2_f \circ \kappa_{12}')|_{\{\nu_2 = 0\}}$ is (at least) $C^1-$smooth such that the image $(\Pi^2_f \circ \kappa_{12}')(\Sigma_f^C) \subset \Sigma^L$ is a neighbourhood of $W^u(q_+) \cap \Sigma^L$, given $\Sigma_f^C$ suitably small. By regular perturbation theory, similar properties follow for the `perturbed' map $\Pi^2_f \circ \kappa_{12}'$ with sufficiently small $\nu_2 \in [0,\Lambda_f]$; recall \remref{reg_pert}, in turn implying the same properties for the map $\Pi_f^2$ since $\kappa_{12}'$ is a diffeomorphism.
	
	\
	
	\textit{The map $\Pi_f^1$}. 
	In order to describe the map $\Pi_f^1$ we consider 
	the restricted maps $\Pi_f^1|_{\Sigma^R_{I}}$, $\Pi_f^1|_{\Sigma^R_{II}}$ and $\Pi_f^1|_{\Sigma^R_{III}}$ separately. Fortunately the details are similar in each case, and so in the following we provide details for the map $\Pi_f^2|_{\Sigma^R_{I}}$, subsequently describing only those aspects of the analysis which differ in the remaining cases. 
	
	Consider the map $\Pi_f^1|_{\Sigma^R_{I}}\circ \kappa_{31}'$, which has range and domain in chart $\mathcal K_1'$. It follows from \lemmaref{mathfrak_K1} and centre manifold theory that the 3D centre manifold $\mathcal M_1'$ is base for a stable foliation $\mathcal F_1'$ such that contraction along fibers in $\mathcal F_1'$ is stronger that $e^{-K_{I} T}$ for all $K_{I} < k \beta$ during the time interval $[0,T]$, where $T$ denotes the time taken for solutions intersecting $\Sigma^R_{I}$ at time $t=0$ to reach $\Sigma^C_f$. An estimate for $T$ can be obtained by considering the flow on $\mathcal M_1'$ and in particular, the following equation for the evolution of base points $p = \{(\tilde x_1, \tilde \rho_1, \tilde \nu_1, \tilde \nu_1)\} \in \mathcal M_1'$:
	\begin{equation}
	\eqlab{mu1_M1}
	\tilde \mu_1'\big|_{\mathcal M_1'} = \frac{1+2k}{\beta k} \tilde \mu_1^2 + \mathcal O\left(\tilde \mu_1 \tilde \rho_1^{k(k+1)} \right)  , 
	\end{equation}
	which is obtained by lifting the expression in \eqref{mathcal_N1} and hence \eqref{mu1_N1} from $\{\rho_1=0\}$ to $\rho_1 \in [0,R_f]$ for $R_f > 0$ sufficiently small. Let $\tilde T$ denote the transition time defined by $\mu_1(\tilde T) = \sigma_f$, and let $\tilde \sigma_f$ be defined by $\tilde \mu_1(T) = \tilde \sigma_f$. By choosing the constants defining the sections $\Sigma^R_I$ and $\Sigma^C_f$ sufficiently small, the estimates
	\begin{equation}
	\eqlab{est1}
	a_- \tilde \sigma_f < \sigma_f < a_+ \tilde \sigma_f  , \qquad b_- \tilde T < T < b_+ \tilde T  ,
	\end{equation}
	for some constants $a_+ > a_- > 0$ and $b_+ > b_- > 0$ follow from the fact that the flow on $\mathcal M_1'$ is regular between transversal sections contained within $\{\mu_1 = \sigma_f\}$ and $\{\mu_1 = \tilde \sigma_f\}$. Integrating \eqref{mu1_M1} we obtain
	\begin{equation}
	\eqlab{est2}
	\tilde T = \mathcal O \left(\frac{1}{\tilde \mu_1} - \frac{1}{\sigma_f} \right) ,
	\end{equation}
	and by combining \eqref{est1} and \eqref{est2}, we obtain an estimate $T = \mathcal O(\tilde \mu_1^{-1})$. Since the fiber projection map $\pi(q): \mathcal F'_1(p) \to p \in \mathcal M_1'$ is smooth and near-identity with respect to $(\rho_1,\nu_1,\mu_1)$, we obtain the estimate
	\begin{equation}
	\eqlab{final_T}
	T = \mathcal O\left(\frac{1}{\mu_1}\right) = \mathcal O\left(\frac{1}{\mu_3}\right) ,
	\end{equation}
	where we have used $\kappa_{31}'$ from \eqref{dash_transition_maps} in the rightmost equality.
	It follows that contraction along fibers is $\mathcal O(e^{-K_I /\mu_3})$. 
	Hence the image $\Pi^2_f(\Sigma^R_I) \subset \Sigma^C_f$ is wedge-shaped with width $\mathcal O(e^{-K_I /\mu_3})$ about $\mathcal M_1' \cap \Sigma^C_f$ in the $x_1-$coordinate.
	
	
	
	

	
	\
	
	Composing $\Pi^1_f$ and $\Pi_f^2$ yields \eqref{Pf}. It follows from the properties of the maps $\Pi^1_f$ and $\Pi_f^2$ (in particular the $\mathcal O(e^{-K_I /\mu_3})$ contraction along $\mathcal F_1'$ near $q_a$) that the $\epsilon_3-$entry of the restricted map $\Pi_f|_{\Sigma_I^R \cap \mathcal V'}(q)$, call it $\mathcal E_F(q,\epsilon,\mu)$, satisfies
	\[
	\mathcal E_f\big|_{\Sigma_I^R}(q,\epsilon,\mu) = \mathcal E_I(\epsilon,\mu) + \mathcal O(e^{-K_I /\mu_3}) , \qquad 
	\frac{\partial \mathcal E_f}{\partial \mu_3}\bigg|_{\Sigma_{I}^R} = \mathcal O\left(e^{-K_I/\mu_3}\right) ,
	\]
	where $\mathcal E_I$ is continuous and satisfies $\mathcal E_I(0,0)= 0$. 
	Restricting to $\mathcal V'$ in \eqref{consts_motion_m3} with $\mu_3 = \sigma_L$ in order to write $\rho_3$ in terms of $\epsilon$, $\mu$ and $\epsilon_3$, simple algebraic manipulations show that
	\begin{equation}
	\eqlab{restrict_rel}
	\frac{\partial \mathcal P^I_f}{\partial \mu_3} = 
	v(\epsilon,\mu) \frac{\partial \mathcal E_f}{\partial \mu_3}\bigg|_{\Sigma_{I}^R} ,
	\end{equation}
	where the function $v(\epsilon,\mu)$ is at most algebraically expanding. The estimate in \eqref{rho1_contraction} follows. To obtain bounds for $\mu_3$
	, notice that for all inputs $q \in \Sigma^R_I$ we have
	\[
	\nu_3 \in [0,\Lambda_I] , \ \mu_3 \in [0, \sigma_I] \ \ \ \implies \ \ \
	\mu = \mu_3 \nu_3^{2k(k+1)} \in \left[0, \mu_3 \Lambda_I^{2k(k+1)} \right] . 
	\]
	The $\mu_3-$interval in \eqref{mu_rho_bounds} follows after some algebraic manipulation.
	
	\
	
	The desired results for inputs $q\in \Sigma_{II}^R$ and $q \in \Sigma_{III}^R$ follow via similar arguments applied to the restricted transition maps $\Pi^1_f|_{\Sigma^R_{II}}$ and $\Pi^1_f|_{\Sigma^R_{III}}$, respectively. The contraction rates in \eqref{rho1_contraction} differ since the transition time estimates
	\begin{equation}
	\eqlab{trans_times}
	T_{II} = \mathcal O\left(\frac{1}{\mu}\right), \qquad T_{III} = \mathcal O\left(\frac{1}{\mu r_1}\right) = \mathcal O\left(\frac{\rho_3^{k(k+1)}}{\epsilon \mu}\right) ,
	\end{equation}
	(c.f. equation \eqref{final_T}) associated with the maps $\Pi^2_f|_{\Sigma^R_{II}}$ and $\Pi^2_f|_{\Sigma^R_{III}}$ differ. Note that in the above we have restricted to \eqref{consts_motion_3} in the rightmost equality. This leads to
	\[
	\mathcal E_f\big|_{\Sigma_{II}^R}(q) = \mathcal E_{II}(\epsilon,\mu) + \mathcal O\left(e^{-K_{II}/\mu}\right) , \qquad \mathcal E_f\big|_{\Sigma_{III}^R}(q) = \mathcal E_{III}(\epsilon,\mu) + \mathcal O\left(e^{-K_{III}\rho_3^{k(k+1)}/\epsilon\mu}\right) ,
	\]
	for continuous functions $\mathcal E_{II}, \mathcal E_{III}$ satisfying $\mathcal E_{II}(0,0)= \mathcal E_{III}(0,0)=0$. The contraction rates in \eqref{rho1_contraction} follow since the relation \eqref{restrict_rel} still holds for inputs $q \in \Sigma^R_{II}$, and a similar relation obtained by restricting to $\mathcal V_\epsilon$ in \eqref{consts_motion_3} holds for inputs $q \in \Sigma_{III}^R$. Bounds for $\rho_3$ are obtained by noting that for each $q \in \Sigma_{III}^R$ we have
	\[
	r_3 \in [0,L_{III}] , \ \rho_3 \in [0,R_{III}] \ \ \ \implies \ \ \ \epsilon = \rho_3^{(1+k)(1+2k)} r_3 c_R \in \left[0 , c_R L_{III} \rho_3^{(1+k)(1+2k)} \right] , 
	\]
	where we have used \eqref{consts_motion_3}. The $\rho_3-$interval in \eqref{mu_rho_bounds} follows after some algebraic manipulation.
\epr



We turn our attention now to the  
backward flow map $\Pi_b:\Sigma^R \to \Sigma^L$ in \eqref{maps}. 

\begin{lemma}
	\lemmalab{Pb}
	Fix all constants defining the sections $\Sigma^L$ and $\Sigma^R$ in \eqref{Sigma_L} and \eqref{Sigma_R} sufficiently small except $\Lambda_{II}$, which is fixed sufficiently large. Then the map $\Pi_b$ is well-defined, at least $C^1-$smooth, and given by
	\begin{equation}
	\eqlab{Pb}
	\Pi_b(q) = 
	\begin{cases}
	\left(\mathcal O\left(\mu_3\right), \mu_3^{-1/k(k+1)} \rho_3 Q_{I,II}(\rho_3,\nu_3,\mu_3), \mathcal O\left(\nu_3 \mu_3^{1/2k(k+1)}\right), \sigma_L \right) , & q \in \Sigma_I^R \cup \Sigma_{II}^R , \\
	\left(\mathcal O\left(\mu r_3\right) , \mu^{-1/k(k+1)} \rho_3 Q_{III}(r_3,\rho_3,\mu), \mathcal O\left(\mu^{1/{2k(k+1)}} \right) , \sigma_L \right) ,   & q \in \Sigma_{III}^R ,
	\end{cases}
	\end{equation}
	where $Q_{I,II}$ and $Q_{III}$ are continuous functions such that $Q_{I,II}(0,0,0), Q_{III}(0,0,0) > 0$.
\end{lemma}

\bpr
	Consider $\Pi_b$ 
	as a composition $\Pi_b = \Pi_b^2 \circ \Pi_b^1$ of component maps $\Pi_b^1 : \Sigma^R \to \Sigma_b^{C,1} \cup \Sigma_2^{C,2}$ and $\Pi_b^2:\Sigma_b^{C,1} \cup \Sigma_b^{C,2} \to \Sigma^L$, where
	\begin{equation}
	\begin{aligned}
	\Sigma^{C,1}_b &= \left\{(c_{b,1},\rho_3,\nu_3,\mu_3) : \rho_3 \in [0,R_{b,1}], \nu_3 \in [0,\Lambda_{II}], \mu_3 \in [0,\sigma_{b,1}] \right\} , && \text{chart } \mathcal K_3' , \\
	\Sigma^{C,2}_b &= \left\{(\epsilon_3,\rho_3,\Lambda_{I},\mu_3) : \epsilon_3 \in [0,c_{b,2}],  \rho_3 \in [0,R_{b,2}], \mu_3 \in [0,\sigma_{b,2}] \right\} , && \text{chart } \mathcal K_3' ,
	\end{aligned}
	\end{equation}
	see \figref{bfb_sections} and \figref{K3_sections}. 
	We consider each map in turn.
	

	
	\
	
	\textit{The map $\Pi_b^1$}. 
	Direct calculations show that the restricted map $\Pi_b^1|_{\Sigma^R_I \cup \Sigma^R_{II}}$ 
	is a diffeomorphism of the form
	\begin{equation}
	\eqlab{Pib11}
	\Pi_b^1|_{\Sigma^R_I \cup \Sigma^R_{II}} (c_R, \rho_3, \nu_3, \mu_3) = \left(c_{b,1}, \rho_3 A^1(\rho_3, \nu_3, \mu_3), \nu_3 A^2(\rho_3, \nu_3, \mu_3), \mu_3 A^3(\rho_3, \nu_3, \mu_3) \right),
	\end{equation}
	where the functions $A^i$ are continuous with $A^i(0,0,0) > 0$ 
	for each $i=1,2,3$. The map $\Pi_b^1|_{\Sigma^R_{III}}$ is regular outside of a neighbourhood of $p_-$ 
	(see \figref{bfb_sections} and \figref{K3_sections}), and the local transition near $p_-$ can be analysed using standard techniques known from e.g. \cite{Gucwa2009,Kristiansen2019c,Krupa2001a}. Here we simply outline the steps:
	\begin{itemize}
		\item In chart $\mathcal K_3$, consider the system obtained after dividing system \eqref{Ks3} by the locally positive factor $-2r_3f_3(r_3,\epsilon_3,\rho_3) - g_3(r_3,\epsilon_3,\rho_3)$. This system has $r_3'=-r_3$, allowing for a direct calculation of the transition time $T$ taken for trajectories to travel from $\Sigma^{in} \subset \{\epsilon_3 = c_{in}\}$ to $\Sigma^{in} \subset \{r_3 = L_{out}\}$, where $c_{in} , L_{out} > 0$ are assumed sufficiently small in the following.
		\item $C^1$-linearise the restricted planar system in $\{\rho_3=0\}$ for which $p_-$ is a hyperbolic saddle by a near-identity map $\tilde \epsilon_3 \mapsto \epsilon_3 = \tilde \epsilon_3(1+\mathcal O(\tilde \epsilon_3,r_3))$, and subsequently apply this transformation to the full 3-dimensional system. This is sometimes referred to as \textit{partial linearisation} \cite{Kristiansen2019c}.
		\item Define $\tilde \epsilon_3(t) = e^t(c_{in} + z)$ and substitute the solution $r_3(t)=r_{in}e^{-t}$ with $r_{in}=r_3(0)$ into the resulting equations for $\tilde \epsilon_3'$. This leads to
		\begin{equation}
		\eqlab{p_minus_eq}
		z' = -r_{in} z \mathcal T\left(\rho_3^{k(k+1)} \right) ,
		\end{equation}
		where the function $\mathcal T(\rho_3^{k(k+1)})$ is smooth and (one can show) uniformly bounded in a neighbourhood of $p_-$.
		\item Integrate \eqref{p_minus_eq} and apply Gronwall's inequality to obtain $z(T) = \mathcal O(r_{in} \ln r_{in})$. This leads to
		\begin{equation}
		\eqlab{p_map}
		\epsilon_3(T) = r_{in} \left(\frac{c_{in}}{L_{out}}\right) \left(1 + \mathcal O\left(r_{in} \ln r_{in}\right)\right) ,
		\end{equation}
		which can subsequently be used to derive $\rho_3(T) = \rho_{in} \left(1 + \mathcal O\left(r_{in} \ln r_{in}\right)\right)$.
	\end{itemize}
	The maps from $\Sigma^R$ to $\Sigma^{in}$ and from $\Sigma^{out}$ to $\Sigma^L$ are diffeomorphisms which take a similar form to \eqref{Pib11}, leading to the following $C^1$ estimate for $\Pi^1_b|_{\Sigma_{III}^R}$:
	\begin{equation}
	\eqlab{Pib12}
	\Pi^1_b|_{\Sigma_{III}^R}(r_3,c_R,\rho_3,\mu) = \left(r_3 B^1(r_3,\rho_3,\mu), \rho_3 B^2(r_3,\rho_3,\mu) , \Lambda_I, \Lambda_I^{-2k(1+k)} \right) ,
	\end{equation}
	where the functions $B^i$ are continuous and satisfy $B^i(0,0,0) > 0$ for each 
	$i=1,2$.
	
	
	
	\
	
	\textit{The map $\Pi_b^2$}. 
	In order to describe $\Pi_b^2$ we describe the restricted maps $\Pi_b^2|_{\Sigma_b^{C,i}}$ for $i=1,2$ separately, starting with $\Pi_b^2|_{\Sigma_b^{C,1}}$. We are interested in the local transition near $Q_-$ and, hence, we work in chart $\mathcal K_3'$. Dividing system 
	\eqref{mathcal_Ks3} 
	by the locally positive factor $-2\tilde f_3(\epsilon_3, \rho_3, \nu_3,\mu_3) - \tilde g_3(\epsilon_3, \rho_3, \nu_3,\mu_3)$, we obtain the locally equivalent system
	\begin{equation}
	\eqlab{Q_sys}
	\begin{split}
	\epsilon_3' &= \epsilon_3 \psi(\epsilon_3,\rho_3,\nu_3, \mu_3), \\
	\rho_3' &= - \frac{\rho_3}{k(k+1)(2k+1)} \left(k+1 + k \psi(\epsilon_3,\rho_3,\nu_3, \mu_3) \right) , \\
	\nu_3' &= \frac{\nu_3}{2k(k+1)} , \\
	\mu_3' &= - \mu_3 ,
	\end{split}
	\end{equation}
	where $\psi(\epsilon_3,\rho_3,\nu_3, \mu_3) = 1 + \mathcal O(\mu_3,\rho_3^{k+1})$ is smooth. We consider solutions $(\epsilon_3,\rho_3,\nu_3,\mu_3)(t)$ of system \eqref{Q_sys} satisfying
	\begin{equation}
	\eqlab{ICs}
	(\epsilon_3, \rho_3, \nu_3,\mu_3)(0) = (c_{b,1}, \rho_{in}, \nu_{in},\mu_{in}) , \qquad 
	(\epsilon_3, \rho_3, \nu_3,\mu_3)(T) = (\epsilon_{out}, R_L, \nu_{out},\mu_{out})  ,
	\end{equation}
	for some transition time $T > 0$. The equations for $\nu_3'$ and $\mu_3'$ decouple; direct integration gives $\nu_3(t) = \nu_{in} e^{t/2k(1+k)}$ and $\mu_3(t) = \mu_{in} e^{-t}$. Using $\mu_3(T) = \sigma_L$ to solve for $T$, we obtain
	\begin{equation}
	\eqlab{trans_T}
	T = \ln \left(\frac{\mu_{in}}{\sigma_L}\right) ,
	\end{equation}
	and hence $\nu_3(T) = \nu_{in} (\mu_{in} / \sigma_L)^{1/2k(1+k)}$. It remains to determine estimates for $\rho_3(T)$ and $\nu_3(T)$.
	
	Since $Q_-$ is a hyperbolic saddle for system \eqref{Q_sys}, with corresponding eigenvalues
	\begin{equation}
	\eqlab{Q_evs}
	\lambda_1 = 1 , \ \ \lambda_2 = - \frac{1}{k(k+1)} , \ \ \lambda_3 = \frac{1}{2k(k+1)} , \ \ \lambda_4 = -1 ,
	\end{equation}
	satisfying $\lambda_i \neq \lambda_{j_1} + \lambda_{j_2}$ for all $i \in \{1,2,3,4\}$, $j_1 \in \{1,3\}$ and $j_2 \in \{2,4\}$, a result due to Belitskii \cite{Belitskii1973} (see also \cite[Theorem 3.1]{Homburg2010}) guarantees existence of near-identity transformation
	\begin{equation}
	\eqlab{bel_trans}
	\left(\tilde \epsilon_3, \tilde \rho_3, \tilde \nu_3, \tilde \mu_3 \right) \mapsto 
	\begin{cases}
	\epsilon_3 = \tilde \epsilon_3 \tilde {\mathcal R}^1(\tilde \epsilon_3 , \tilde \rho_3, \tilde \nu_3, \tilde \mu_3) , \\
	\rho_3 = \tilde \rho_3 \tilde{\mathcal R}^2(\tilde \epsilon_3 , \tilde \rho_3, \tilde \nu_3, \tilde \mu_3) , \\
	\nu_3 = \tilde \nu_3 , \\
	\mu_3 = \tilde \mu_3 , \\
	\end{cases}
	\end{equation}
	where $\tilde {\mathcal R}^i$ are continuous functions such that $\tilde{\mathcal R}^i(0 , 0, \tilde \nu_3, \tilde \mu_3) = 1$ for each $i=1,2$, which $C^1-$linearises system \eqref{Q_sys} near $Q_-$.
	Applying Belitskii's theorem, integrating the resulting (linear) equations for $\tilde \epsilon_3'$ and $\tilde \rho_3'$ and subsequently inverting \eqref{bel_trans}, we obtain 
	the following estimate for the restricted map $\Pi_b^2|_{\Sigma_b^{C,1}}$:
	\begin{equation}
	\eqlab{Pi_b2}
	\begin{split}
	& \Pi_b^2|_{\Sigma_b^{C,1}} (c_{b,1}, \rho_3, \nu_3, \mu_3 ) = \\
	& \left( \left(\frac{c_{b,1}\mu_3 }{\sigma_L}\right)\mathcal R^{1}(\rho_3,\nu_3, \mu_3) , \rho_3 \left(\frac{\sigma_L}{ \mu_3}\right)^{1/k(1+k)} \mathcal R^{2}(\rho_3,\nu_3, \mu_3), \nu_3 \left( \frac{\mu_3}{\sigma_L}\right)^{1/2k(1+k)}, \sigma_L \right) ,
	\end{split}
	\end{equation}
	where $\mathcal R^{i}$ are continuous with $\mathcal R^i(0,\nu_3, \mu_3) = 1$ for $i=1,2$.
	
	\
	
	Details of the map $\Pi_b^2|_{\Sigma_b^{C,2}}$ are similar, except with $(\epsilon_3, \rho_3, \nu_3,\mu_3)(0) = (\epsilon_{in}, \rho_{in}, \Lambda_{b,2},\mu_{in})$ 
	in place of the initial conditions in \eqref{ICs}. In this case one obtains
	\begin{equation}
	\eqlab{Pi_b22}
	\begin{split}
	& \Pi_b^2|_{\Sigma_b^{C,2}} (\epsilon_3, \rho_3, \Lambda_{b,2}, \mu_3 ) = \\
	& \left(\left(\frac{\epsilon_3 \mu_3}{\sigma_L}\right)\mathcal R^{3}(\epsilon_3,\rho_3,\mu_3) , \rho_3 \left(\frac{\sigma_L}{ \mu_3}\right)^{1/k(1+k)} \mathcal R^{4}(\epsilon_3,\rho_3,\mu_3), \Lambda_{b,2} \left( \frac{\mu_3}{\sigma_L}\right)^{1/2k(1+k)}, \sigma_L \right) ,
	\end{split}
	\end{equation}
	where $\mathcal R^{i}$ are continuous with $\mathcal R^i(0,0, \mu_3) = 1$ for $i=3,4$.
	
	\
	
	The estimates in \eqref{Pb} are obtained by composing \eqref{Pib11} with \eqref{Pi_b2} for inputs $q \in \Sigma^R_I \cup \Sigma^R_{II}$, and \eqref{Pib12} with \eqref{Pi_b22} for inputs $q \in \Sigma_{III}^R$. The $C^1$ property follows since the component maps \eqref{Pib11}, \eqref{Pib12}, \eqref{Pi_b2} and \eqref{Pi_b22} are all $C^1$.
\epr

Having described the maps $\Pi_f$ and $\Pi_b$ in \lemmaref{Pf} and \lemmaref{Pb}, we look for solutions for the fixed point equation
\begin{equation}
\eqlab{fixed_point3}
D(q) = \left(P_f - P_b\right) (q) = 0
\end{equation}
corresponding to limit cycles in either system $\{$\eqref{main_extended}$, \ \mu' = 0\}$ or system \eqref{main_extended2}.

\begin{lemma}
		\lemmalab{fixed_point1}
		The following assertions are true:
			\begin{enumerate}
			\item[(i)] Let $q = (c_R,\rho_3, \nu_3,\mu_3) \in\Sigma_I^R$. Then
			\begin{equation}
			\eqlab{fp_sol_I}
			\rho_3 = \mu_3^{1/k(k+1)} H_{I}(\nu_3,\mu_3) 
			\end{equation}
			is a solution of \eqref{fixed_point3}, where $H_{I}$ is a continuous function such that $H_{I}(0,0) > 0$.
			\item[(ii)] Let $q = (c_R,\rho_3, \nu_3,\mu_3) \in \Sigma_{II}^R$. Then
			\begin{equation}
			\eqlab{fp_sol_II}
			\mu_3 = \rho_3^{k(k+1)} H_{II}(\rho_3,\nu_3) 
			\end{equation}
			is a solution of \eqref{fixed_point3}, where $H_{II}$ is a continuous function satisfying $H_{II}(\nu_3,0) = b_{II}$ for a constant $b_{II} > 0$.
			\item[(iii)] Let $q = (r_3,c_R,\rho_3, \mu) \in \Sigma_{III}^R$. Then
			\begin{equation}
			\eqlab{fp_sol_III}
			\mu = \rho_3^{k(k+1)} H_{III}(r_3,\rho_3) 
			\end{equation}
			is a solution of \eqref{fixed_point3}, where $H_{III}$ is a continuous function such that $H_{III}(0,0) > 0$.
			\item[(iv)] Solutions given by \eqref{fp_sol_I}, \eqref{fp_sol_II} and \eqref{fp_sol_III} coincide where their domains overlap.
		\end{enumerate}
\end{lemma}

\bpr
It follows from \lemmaref{Pf} and \lemmaref{Pb} that
\begin{equation}
\eqlab{fp_again}
D (q) = 
\begin{cases}
\gamma_L + J_{I,II}(\rho_3,\nu_3,\mu_3) - \mu_3^{-1/k(k+1)} \rho_3 Q_{I,II}(\rho_3,\nu_3,\mu_3) , & q \in \Sigma^R_I \cup \Sigma^R_{II}, \\
\gamma_L + J_{III}(r_3,\rho_3,\mu) - \mu^{-1/k(k+1)} \rho_3 Q_{III}(r_3,\rho_3,\mu) , & q \in \Sigma^R_{III} .
\end{cases}
\end{equation}
Setting the first expression equal to zero 
and multiplying by $\mu_3^{1/k(k+1)}$ 
gives
\begin{equation}
\eqlab{aux_fp}
\bar D\big|_{\Sigma_I^R\cup \Sigma_{II}^R}  \left(\rho_3,\nu_3,\mu_3 \right) =  \mu_3^{1/k(k+1)}\left(\gamma_L + J_{I,II}\left(\rho_3,\nu_3,\mu_3\right)\right) - \rho_3 Q_{I,II}\left(\rho_3,\nu_3,\mu_3\right) = 0 ,
\end{equation}
where $\bar D|_{\Sigma_I^R \cup \Sigma_{II}^R} := \mu_3^{1/k(k+1)} D|_{\Sigma_I^R \cup \Sigma_{II}^R}$. Notice that
\[
\bar D\big|_{\Sigma_I^R}\left(0,0,0\right) = 0 , \qquad 
\frac{\partial \bar D|_{\Sigma_I^R}}{\partial \rho_3} \left(0,0,0\right) = - Q_{I,II}(0,0,0) < 0 ,
\]
where the last inequality follows by \lemmaref{Pb}. Hence \eqref{aux_fp} can be solved for $\rho_3$ 
by the implicit function theorem, and the expression \eqref{fp_sol_I} in statement (i) follows. 

The expression \eqref{fp_sol_II} in statement (ii) is similarly obtained by solving \eqref{aux_fp} via the implicit function theorem, except in this case with respect to $\mu_3$. In particular, setting $m_{II} = \mu_3^{1/k(k+1)}$ in \eqref{aux_fp}, one obtains
\[
\bar D\big|_{\Sigma_{II}^R}\left(0,\nu_3,0\right) = 0 , \qquad 
\frac{\partial \bar D|_{\Sigma_{II}^R}}{\partial m_{II}} \left(0,\nu_3,0\right) = \gamma_L > 0 ,
\]
and the result follows.
Statement (iii) is proved by setting $\mu^{1/k(k+1)} = m_{III}$ in the second expression in \eqref{fp_again} and applying similar (implicit function theorem) arguments to obtain \eqref{fp_sol_III}.

\SJnew{Finally}, statement (iv) follows since the solutions \eqref{fp_sol_I}, \eqref{fp_sol_II} and \eqref{fp_sol_III} are locally unique (being derived by the implicit function theorem) on overlapping domains.
\epr

Intersecting solutions \eqref{fp_sol_I}, \eqref{fp_sol_II} with $\mathcal V'$ and \eqref{fp_sol_III} with $\mathcal V_\epsilon$ appropriate invariant sets will allow us to solve three separate fixed point equations
\begin{equation}
\eqlab{fp_3}
\mathcal D_I(q, \epsilon, \mu) = 0, \qquad \mathcal D_{II}(q, \epsilon, \mu) = 0, \qquad \mathcal D_{III}(q, \epsilon, \mu) = 0,
\end{equation}
for $\mu$, where $\mathcal D_I := \mathcal D|_{\Sigma_I^R}$, $\mathcal D_{II} := \mathcal D|_{\Sigma_{II}^R}$ and $\mathcal D_{III} := \mathcal D|_{\Sigma_{III}^R}$ and $\mathcal D$ is defined by \eqref{fixed_point}.

\begin{lemma}
	\lemmalab{fixed_point2}
	Fix $\epsilon_0 > 0$ sufficiently small. Then there exist unique, continuous functions
	$\mu_I(\mu_3,\epsilon)$, $\mu_{II}(\nu_3,\epsilon)$ and $\mu_{III}(\rho_3,\epsilon)$ such that for all $\epsilon \in (0,\epsilon_0)$ we have
	\begin{equation}
	\eqlab{D_eqn}
	\mathcal D_I(q, \epsilon, \mu_I(\mu_3,\epsilon)) = 0, \qquad
	\mathcal D_{II}(q, \epsilon, \mu_{II}(\nu_3,\epsilon)) = 0, \qquad
	\mathcal D_{III}(q, \epsilon, \mu_{III}(\rho_3,\epsilon)) = 0 ,
	\end{equation}
	\SJnew{for all
	\begin{equation}
	\eqlab{mu_ints}
	\mu_3 \in \left[\bar a_I \epsilon^{k/(2k+1)} , \sigma \right] , \ \ \ 
	\nu_3 \in [\Lambda_I, \Lambda_{II}] , \ \ \ 
	\rho_3 \in \left[\left(\epsilon \left(c_R L_{III}\right)^{-1}\right)^{1/(k+1)(2k+1)} , L_{III} \right] ,
	\end{equation}
	respectively, where $\sigma := \sigma_I^{(k+1)/(2k+1)}$ and $\bar a_I > 0$ is a sufficiently large constant.} In particular,
	\begin{equation}
	\eqlab{mu}
	\begin{split}
	\mu_I(\mu_3,\epsilon) &= 
	\left(\epsilon \mu_3^{-1}\right)^{k/(k+1)} \bar {\mathcal H}_I(\mu_3,\epsilon) ,  \\
	\mu_{II}(\nu_3,\epsilon) &= 
	\left(\epsilon \nu_3^{2k(k+1)}\right)^{k/(2k+1)} \bar {\mathcal H}_{II}(\nu_3,\epsilon)  ,  \\
	\mu_{III}(\rho_3,\epsilon) &= 
	\rho_3^{k(k+1)} \bar {\mathcal H}_{III}(\rho_3,\epsilon)  , 
	\end{split}
	\end{equation}
	where $\bar {\mathcal H}_I(\mu_3,\epsilon)$, $\bar {\mathcal H}_{II}(\nu_3,\epsilon)$ and $\bar {\mathcal H}_{III}(\rho_3,\epsilon)$ are continuous and satisfy 
	\[
	\bar {\mathcal H}_I(\mu_3,0) = a_I , \qquad \bar {\mathcal H}_{II}(\nu_3,0) = a_{II}, 
	\qquad \bar {\mathcal H}_{III}(\rho_3,0) = a_{III},
	\]
	where $a_I$, $a_{II}$ and $a_{III}$ are positive constants. 
\end{lemma}

\bpr
First consider inputs $q\in\Sigma_I^R$. In this case, we need to solve \eqref{fp_sol_I} 
for $\mu$. 
First, we restrict to $\mathcal V_\epsilon'$ in order to write $\rho_3$ as
\begin{equation}
\eqlab{rho3_res}
\rho_3 = \left(\epsilon c_R^{-1} \nu_3^{-2(k+1)^2} \right) ^{1/(k+1)(2k+1)} .
\end{equation}
After inserting \eqref{rho3_res} into \eqref{fp_sol_I}, we obtain the following after some algebraic manipulations:
\begin{equation}
\eqlab{fp_I_again}
\left(\epsilon c_R^{-1}\right)^{1/(k+1)(2k+1)} - \nu_3^{2(k+1)/(2k+1)} \mu_3^{1/k(k+1)} H_I(\nu_3,\mu_3) = 0 .
\end{equation}
\SJnew{We solve \eqref{fp_I_again} for $l := \nu_3^{2(k+1)/(2k+1)}$. By \lemmaref{fixed_point1}, the function $H_I$ is continuous with $H_I(0,0) > 0$. Hence choosing $\Sigma_I^R$ sufficiently small guarantees uniform bounds $c_\pm > 0$ such that $0 < c_- \leq H(\nu_3, \mu_3) \leq c_+$. Using these bounds and choosing the constant $\bar a_I>0$ defining the lower bound for the ($\epsilon-$dependent) $\mu_3-$interval in \eqref{mu_ints} sufficiently large, \eqref{fp_I_again} can be solved by the implicit function theorem locally in a neighbourhood of $(l,\mu_3,\epsilon) = (0, \mu_3, 0)$, for all $\mu_3$ in the first interval in \eqref{mu_ints}. We obtain}
%
\[
l = \nu_3^{2(k+1)/(2k+1)} = \left(\epsilon c_R^{-1} \right)^{1/(k+1)(2k+1)} \mu_3^{-1/k(k+1)} \tilde H_I(\mu_3,\epsilon) ,
\]
where $\tilde H_I$ is a continuous function which is strictly positive in a neighbourhood about \SJnew{$(0, \mu_3, 0)$.} 
The expression for $\mu_I(\mu_3,\epsilon)$ follows \SJnew{after} using $\mu = \nu_3^{2k(k+1)} \mu_3$. 

\

Now consider inputs $q \in \Sigma^R_{II}$. We solve \eqref{fp_sol_II} for $\mu$ by restricting to $\mathcal V'$. In particular we have 
\begin{equation}
\eqlab{restricty}
\rho_3 = \left(\epsilon c_R^{-1}\right)^{1/(1+k)(1+2k)} \nu_3^{-2(k+1)/(2k+1)} , \qquad \mu_3 = \mu \nu_3^{-2k(k+1)} ,
\end{equation}
where we note that $\nu_3 \in [\Lambda_I, \Lambda_{II}]$ for inputs $q \in \Sigma_{II}^R$. Using \eqref{restricty} to write \eqref{fp_sol_II} in terms of $\nu_3$, $\epsilon$ and $\mu$ and solving the resulting expression for $\mu$ for all $\epsilon, \mu \ll 1$ and $\nu_3 \in [\Lambda_I,\Lambda_{II}]$ via the implicit function theorem, one obtains the desired expression $\mu = \mu_{II}(\nu_3,\epsilon)$ in \eqref{mu}.

\

Finally, for inputs $q \in \Sigma^R_{III}$ we solve \eqref{fp_sol_III} for $\mu$ by restricting to $\mathcal V_\epsilon$ in \eqref{consts_motion_3}. 
Using
\[
r_3 = \epsilon c_R^{-1} \rho_3^{-(k+1)(2k+1)} 
\]
to write the right-hand-side of \eqref{fp_sol_III} in terms of $\rho_3$ and $\epsilon$ yields the desired expression for $\mu = \mu_{III}(\rho_3,\epsilon)$ in \eqref{mu}. \SJnew{The form of the interval for $\rho_3$ in \eqref{mu_ints} follows from \lemmaref{Pf}; see again \eqref{mu_rho_bounds}.}
\epr

\lemmaref{fixed_point2} implies existence of \SJnew{locally unique} parameterised families of limit cycles
\[
\mu_3 \mapsto \left(\mu_I(\mu_3,\epsilon),\Gamma_I(\mu_3,\epsilon)\right), \ \ \
\nu_3 \mapsto  \left(\mu_{II}(\nu_3,\epsilon),\Gamma_{II}(\nu_3,\epsilon)\right), \ \ \
\rho_3 \mapsto \left(\mu_{III}(\rho_3,\epsilon),\Gamma_{III}(\rho_3,\epsilon)\right) ,
\]
which intersect either $\Sigma_I^R$, $\Sigma_{II}^R$ or $\Sigma_{III}^R$ respectively. \SJnew{It follows from \lemmaref{fixed_point1} and \lemmaref{fixed_point2} that these families coincide on overlapping domains and, hence, that each family constitutes part of a \textit{single family} of limit cycles which are related across charts $\mathcal K_3'$ and $\mathcal K_3$ by the smooth change of coordinates \eqref{33_trans}. In the following, we introduce a parameterisation for this family that is consistent with the convergence and monotonicity properties given in \thmref{thm_cycs}.} 
Let
\begin{equation}
\eqlab{q_par}
q(\theta) =
\begin{cases}
(c_R, \rho_{R,\epsilon}(\theta), \nu_{R,\epsilon}(\theta), \mu_{R,\epsilon}(\theta)) \in \Sigma_I^R\cup\Sigma_{II}^R  & \text{if } \theta \in [0,2/3) ,  \\
(r_{R,\epsilon}(\theta), c_R, \rho_{R,\epsilon}(\theta), \mu) \in \Sigma_{II}^R\cup\Sigma_{III}^R  & \text{if } \theta \in (1/3,1] ,  \\
\end{cases}
\end{equation}
where in the first (second) expression $q(\theta)$ is expressed in chart $\mathcal K_3'$ ($\mathcal K_3$) coordinates. The subscript $\epsilon$ indicates $\epsilon-$dependence which results from restriction to the invariant sets $\mathcal V'$ and $\mathcal V_\epsilon$ in \eqref{consts_motion_m3} and \eqref{consts_motion_3}. By \eqref{33_trans}, the two expressions in \eqref{q_par} coincide on their overlapping domain $\theta \in (1/3,2/3)$ corresponding to $\{\nu_3 > 0\}$ ($\{r_3 > 0\}$) in chart $\mathcal K_3'$ ($\mathcal K_3$). We further impose the following defining properties:
\begin{enumerate}
	\item[(P1)] Smoothness/monotonicity: For each fixed $\epsilon \in (0,\epsilon_0)$ the map $\theta \mapsto q(\theta)$ is smooth and satisfies
	\[
	\mu_{R,\epsilon}'(\theta) <0, \qquad \nu_{R,\epsilon}'(\theta) > 0 \qquad  \text{for all } \theta \in (0,2/3),
	\]
	\[
	r_{R,\epsilon}'(\theta)<0, \qquad \rho_{R,\epsilon}'(\theta)>0 \qquad \text{for all } \theta \in (1/3,1).
	\]
	\item[(P2)] Fixed endpoints corresponding to $\theta \in\{0,1\}$ as $\epsilon \to 0$:
	\[
	\lim_{\epsilon \to 0} \mu_{R,\epsilon}(0) = \sigma , 
	\qquad 
	\lim_{\epsilon \to 0} \rho_{R,\epsilon}(1) = L . 
	\]
	\item[(P3)] Agreement at the intersections $\Sigma_I^R \cap \Sigma_{II}^R$ and $\Sigma_{II}^R \cap \Sigma_{III}^R$: If $\theta = \theta_- \in (1/3,2/3)$ corresponds to the intersection $\Sigma_I^R \cap \Sigma_{II}^R$, then 
	\[
	\nu_{R,\epsilon}(\theta_-) = \Lambda_I , \qquad
	\mu_{R,\epsilon}(\theta_-) = 
	\left(\epsilon \Lambda_I^{-2(1+k)^2} \right)^{k/(2k+1)} \bar {\mathcal H}_{II}(\Lambda_I,\epsilon) ,
	\]
	where the latter can be derived by inserting the expression for $\mu_{II}(\nu_3,\epsilon)$ in \eqref{mu} into the lower bound for the $\mu_3-$interval in \eqref{mu_rho_bounds} and using $\mu = \nu_3^{2k(k+1)} \mu_3$. Similarly, if $\theta = \theta_+ \in (1/3,2/3)$ corresponds to the intersection $\Sigma_{II}^R \cap \Sigma_{III}^R$, then 
	\[
	\nu_{R,\epsilon}(\theta_+) = \Lambda_{II} , \qquad
	\rho_{R,\epsilon}(\theta_+) = \left(\epsilon \left(c_R L_{III}\right)^{-1}\right)^{1/(k+1)(2k+1)} ,  
	\]
	where the latter coincides with the lower bound of the $\rho_3-$interval in \eqref{mu_rho_bounds} and \eqref{mu_ints}.
\end{enumerate}
It follows from (P1)-(P3) and the intermediate value theorem that there exists some $\theta_c \in (0, \theta_-)$ such that $\lim_{\epsilon\to 0}\mu_{R,\epsilon}(\theta)$ is strictly positive for all $\theta \in [0,\theta_c)$, and zero for all $\theta \in [\theta_c,2/3)$. Without loss of generality, we may define the parameterisation so that $\theta_c = 1/3$. Similar arguments regarding the behaviour of $\lim_{\epsilon\to 0}\rho_{R,\epsilon}(\theta)$ allow us to impose on final defining property (P4):
\begin{enumerate}
	\item[(P4)] 
	Limit properties: $\mu_{R,\epsilon}(\theta)$ and $\nu_{R,\epsilon}(\theta)$ satisfy
	\[
	\lim_{\epsilon \to 0 }\mu_{R,\epsilon}(\theta) =
	\begin{cases}
	\mu_R(\theta) > 0 , & \theta \in [0,1/3) , \\
	=0 , & \theta \in [1/3,2/3) , \\
	\end{cases}
	\ \ \ \
	\lim_{\epsilon \to 0 }\nu_{R,\epsilon}(\theta) =
	\begin{cases}
	=0 , & \theta \in [0,1/3) , \\
	\nu_R(\theta)>0 , & \theta \in [1/3,2/3) , \\
	\end{cases}
	\]
	where $\mu_R(\theta)$ and $\nu_R(\theta)$ are diffeomorphisms. Similarly, $r_{R,\epsilon}(\theta)$ and $\rho_{R,\epsilon}(\theta)$ satisfy
	\[
	\lim_{\epsilon \to 0 }r_{R,\epsilon}(\theta)  =
	\begin{cases}
	r_R(\theta)>0 , & \theta \in (1/3,2/3) , \\
	=0 , & \theta \in [2/3,1] , \\
	\end{cases}
	\ \ \ \
	\lim_{\epsilon \to 0 }\rho_{R,\epsilon}(\theta) 
	\begin{cases}
	=0 , & \theta \in (1/3,2/3] , \\
	\rho_R(\theta)>0 , & \theta \in (2/3,1] , \\
	\end{cases}
	\]
	where $r_R(\theta)$ and $\rho_R(\theta)$ are diffeomorphisms.
\end{enumerate}

Using properties (P1)-(P4), we can combine the expressions \eqref{mu} to obtain \textit{single} function
\begin{equation}
\eqlab{global_mu}
\mu : [0,1] \times [0,\epsilon_0] \to \mathbb R , \qquad 
\mu(\theta,\epsilon) = 
\begin{cases}
\mu_I(\mu_{R,\epsilon}(\theta),\epsilon) , & \theta \in [0,1/3], \\
\mu_{II}(\nu_{R,\epsilon}(\theta),\epsilon) , & \theta \in [1/3,2/3], \\
\mu_{III}(\rho_{R,\epsilon}(\theta),\epsilon) , & \theta \in [2/3,1] ,
\end{cases}
\end{equation}
where the component functions are given by \eqref{mu}. We obtain the following result.

\begin{lemma}
	\lemmalab{fixed_point}
	Fix $\kappa \in (0,1/3)$ arbitrarily small. Then there exists $\epsilon_0 > 0$ such that the following assertions hold:
	\begin{enumerate}
		\item[(i)] The function $\mu(\theta,\epsilon)$ in \eqref{global_mu} satisfies the restricted fixed point equation \eqref{fixed_point}.
		\item[(ii)] The function $\mu(\theta,\epsilon)$ satisfies the following estimates:
	\begin{equation}
	\eqlab{mu_estimates}
	\mu(\theta,\epsilon) = 
	\begin{cases}
	\mathcal O\left(\epsilon^{k/(k+1)} \right) , & \text{if } \theta \in \left[0, \frac{1}{3} - \kappa \right] , \\
	\mathcal O\left(\epsilon^{k/(2k+1)}\right) , & \text{if } \theta \in \left[\frac{1}{3} + \kappa, \frac{2}{3} - \kappa \right] ,  \\
	\mathcal O\left(1\right) , & \text{if } \theta \in \left[\frac{2}{3} + \kappa, 1\right] .
	\end{cases}
	\end{equation}
	\end{enumerate}
\end{lemma}

\bpr
Statement (i) is immediate from \lemmaref{fixed_point2}. To see that statement (ii) is true, notice that by property (P4) the functions $\mu_{R,\epsilon}(\theta)$, $\nu_{R,\epsilon}(\theta)$ and $\rho_{R,\epsilon}(\theta)$ are $\mathcal O(1)$ with respect to $\epsilon$ on the respective $\theta-$intervals $[0,1/3-\kappa]$, $[1/3+\kappa,2/3-\kappa]$ and $[2/3+\kappa, 1]$. Hence the estimates \eqref{mu_estimates} follow from \eqref{global_mu} and \eqref{mu}.
\epr

\subsection{Proof of \thmref{thm_cycs}}
\applab{proof_of_thm_cycs}

We now combine \lemmaref{Pf}--\lemmaref{fixed_point} in order to 
to prove \thmref{thm_cycs}, and thus \thmref{thm_connection}.


Existence of the family of limit cycles \eqref{fam2} follows from \lemmaref{fixed_point} (i). 
Stability of the cycles follows by considering the Poincar{\'e} map defined by the composition $\Pi = \Pi^{-1}_b \circ \Pi_f : \Sigma^R \to \Sigma^R$, where $\Pi_b^{-1} : \Sigma^L \to \Sigma^R$ denotes the inverse flow map associated with the map $\Pi_b$. By applying arguments similar to those used to characterise the map $\Pi_b$ in the proof of \lemmaref{Pb}, $\Pi_b^{-1}$ can be shown to be at least $C^1-$smooth and at most algebraically expanding. Since by \lemmaref{Pf} the map $\Pi_f$ is a strong contraction upon restriction to invariant subsets $\mathcal V'$ and $\mathcal V_\epsilon$ defined in \eqref{consts_motion_m3} and \eqref{consts_motion_3}, the restricted map 
\[
\tilde \Pi : \Sigma^R \times [0,\epsilon_0] \times [0,\mu_0] \to \Sigma^L , \ \ \
\tilde \Pi(q,\epsilon, \mu) = 
\begin{cases}
\Pi(q)\big|_{\mathcal V'} , & q \in (\Sigma_I^R \cup \Sigma_{II}^R) \cap \mathcal V' , \\
\Pi(q) \big|_{\mathcal V_\epsilon}, & q \in \Sigma_{III}^R \cap \mathcal V_\epsilon ,
\end{cases}
\]
is also a strong contraction. Stability of the cycles follows after an application of the contraction mapping principle.

Smoothness of $\theta \mapsto \mu(\theta,\epsilon)$ follows from the fact that by the component functions in \eqref{global_mu} are smooth with respect to $\theta$ (see property (P1)), 
and the monotonicity property \eqref{monotonicity} follows from \eqref{global_mu} after differentiation, using $\mu_{R,\epsilon}'(\theta) < 0$, $\nu_{R,\epsilon}'(\theta) > 0$ and $\rho_{R,\epsilon}'(\theta) > 0$ from property (P1). This proves (i).

Now consider assertions (ii)-(iv). The estimates for $\mu(\theta,\epsilon)$ in statements (ii), (iii), and (iv) follow directly from \lemmaref{fixed_point} (ii), so it remains to consider the behaviour the cycles $\bar \Gamma(\theta, \epsilon)$ in the limit $\epsilon \to 0$. \lemmaref{Pf} and \lemmaref{fixed_point} imply that the restricted maps $\mathcal P_{f/b}$ in \eqref{restricted_maps} satisfy 
\[
\lim_{\epsilon \to 0} \mathcal P_{f/b}\left(\theta,\epsilon,\mu(\theta,\epsilon)\right) = \gamma_L + 
\begin{cases}
\mathcal O \left(e^{-K / \mu_R(\theta)}\right) , & \theta \in \left(0,1/3 \right) , \\
0 , & \theta \in \left[2/3, 2/3 \right) ,
\end{cases}
\]
where $\mu_R(\theta) > 0$ is defined in property (P4). It follows that the cycles for $\theta \in [1/3, 1)$ converge in Hausdorff distance to singular cycles in $\{\bar \Gamma(\theta) : \theta \in [1/3,1)\}$ as $\epsilon \to 0$, proving statements (iii) and (iv). In case $\theta \in (0,1/3 - \kappa]$, the expression for $\mu(\theta,\epsilon)$ can be recast in terms of $\hat \mu = \mu \epsilon^{-k/(k+1)}$ using the expression for $\mu_I(\mu_3,\epsilon)$ in \eqref{mu} as follows:
\begin{equation}
\eqlab{mu_I2}
\hat \mu(\theta,\epsilon) = \mu_{R,\epsilon}(\theta)^{-k/(k+1)} \bar {\mathcal H}_I(\theta,\epsilon) , \qquad \theta \in \left(0,1/3-\kappa\right] .
\end{equation}
Hence
\[
\lim_{\epsilon \to 0}\hat \mu(\theta,\epsilon) = \hat \mu(\theta,0) = a_I \mu_R(\theta)^{-k/(k+1)} , \qquad \theta \in \left(0,1/3-\kappa\right] ,
\]
which is fixed and positive since $\mu_R(\theta) > 0$ for all $\theta \in (0,1/3 - \kappa]$ and $a_I>0$ by \lemmaref{fixed_point2}, 
indicating convergence to regular (i.e. not relaxation-type) oscillations within the scaling regime (S1).
This proves statement (ii).

Finally, the uniqueness property (v) follows from uniqueness of the functions $\mu_I(\mu_3,\epsilon)$, $\mu_{II}(\nu_3,\epsilon)$ and $\mu_{III}(\rho_3,\epsilon)$ in \eqref{mu}, which by \lemmaref{fixed_point1} (iv) coincide on overlapping domains.
\qed

\begin{remark}
	The estimates \eqref{rho1_contraction} and \eqref{mu_estimates} can be used to estimate Floquet exponents $\mathcal F l(\theta,\epsilon)$ for limit cycles in \eqref{fam2}. We also obtain the following estimates:
	\begin{equation}
	\eqlab{floq}
	\mathcal F l(\theta,\epsilon) =
	\begin{cases}
	\mathcal O\left(\mu_R(\theta)^{-1} \right) , & \text{if } \theta \in \left(0,\frac{1}{3}-\kappa\right] ,  \\
	\mathcal O\left(\epsilon^{-k/(2k+1)} \right) , & \text{if } \theta \in \left[\frac{1}{3}+\kappa,\frac{2}{3}-\kappa\right] ,  \\
	\mathcal O\left(\epsilon^{-1} \right) , & \text{if } \theta \in \left[\frac{2}{3} + \kappa , 1\right) ,
	\end{cases}
	\end{equation}
	which are consistent with a transition from regular oscillations 
	in (S1) with finite contraction as $\epsilon \to 0$, to \SJnew{`classical'} relaxation-oscillations 
	in (S2) characterised by $\mathcal O(e^{-K/\epsilon})$ contraction.
\end{remark}

\newpage




\bibliography{bibnew}
\bibliographystyle{plain}





\end{document}